\documentclass[a4paper]{amsart}
\usepackage[margin=3cm, marginpar=3cm]{geometry}

\title[The Heisenberg algebra of a vector space and Hochschild
homology]{The Heisenberg algebra of a vector space \\ and Hochschild homology}

\author{\'Ad\'am Gyenge}
\address{Budapest University of Technology and Economics, Department of Algebra and Geometry, M\H{u}egyetem rakpart 3-9., 1111, Budapest, Hungary}
\email{Gyenge.Adam@ttk.bme.hu}

\author{Timothy Logvinenko}
\address{School of Mathematics, Cardiff University, Senghennydd Road, CF24 4AG, Cardiff UK}
\email{LogvinenkoT@cardiff.ac.uk}

\usepackage[T1]{fontenc}
\usepackage[utf8]{inputenc}
\usepackage{amsmath,amsfonts,amssymb,amsthm,stackrel,leftidx}
\usepackage{mathtools,bbold}
\usepackage{lmodern,microtype,csquotes}
\usepackage{graphicx}
\usepackage{mleftright}

\usepackage{tikz, tikz-cd}
\usetikzlibrary{decorations.markings}
\usetikzlibrary{shapes.geometric}
\usetikzlibrary{decorations.pathmorphing}
\usetikzlibrary{calc}

\usepackage{comment}

\usepackage{enumitem}
\setlist[enumerate,1]{label=(\arabic*)}

\usepackage[colorlinks=true, pdfpagemode=none, pdfmenubar=false, pdfstartview=FitH, linkcolor=blue, citecolor=blue, urlcolor=blue, bookmarksdepth =2]{hyperref}

\theoremstyle{plain}

\newtheorem*{Theorem*}{Theorem}
\newtheorem*{Conjecture*}{Conjecture}

\newtheorem{Theorem}{Theorem}[section]

\newtheorem{Proposition}[Theorem]{Proposition}
\newtheorem{Corollary}[Theorem]{Corollary}
\newtheorem{Lemma}[Theorem]{Lemma}
\newtheorem{Conjecture}[Theorem]{Conjecture}

\theoremstyle{definition}
\newtheorem{Definition}[Theorem]{Definition}

\newtheorem*{Definition*}{Definition}

\theoremstyle{remark}

\newtheorem{Example}[Theorem]{Example}

\numberwithin{equation}{section}
\usepackage{suffix}  % To define starred commands
\usepackage{ifthen}
\usepackage{xspace}

\newif\ifdgcal
\dgcaltrue    % Use caligraphic uppercase and italic lowercase.
\newcommand\basecat{{\cat V}}       % The base category for our constructions.

\let\simeq\cong
\newcommand\dg{DG\xspace}

\mathchardef\mhyphen="2D    % hyphen in math environments (used in some macros below)

\newcommand\kk{{\mathbb{k}}}            % base field

\newcommand\ZZ{{\mathbb{Z}}}

\newcommand\cat\mathcal             % 1-categories
\newcommand\bicat\mathbf            % 2-categories

\DeclareMathOperator{\obj}{Ob}
\newcommand\Hom{\operatorname{Hom}}
\newcommand\End{\operatorname{End}}

\newcommand\id{\operatorname{id}}

\newcommand\catgrVect{\mathrm{gr\mhyphen Vect_{\kk}}} % cat of graded vector spaces

\newcommand\Ind{\operatorname{Ind}}     % induction functor
\newcommand\Res{\operatorname{Res}}     % restriction functor
\DeclareMathOperator\Sym{Sym}           % Symmetric powers of a vector space
\newcommand\hochhom{\mathrm{HH}}
\newcommand\hochcx{\mathrm{HC}}

\newcommand{\fix}{{\mathrm{Fix}}}       % fixed point set
\newcommand{\cts}{{\mathrm{cts}}}    % Reduction to continuous chains map on Hochschild complex
\newcommand{\redcts}{{\mathrm{\underline{Rdc}}}}    % Reduction to continuous chains map on Hochschild complex

\newcommand\Algonek{{\bicat{Alg}^1_\kk}}        % bicategory of dg categories
\newcommand\grCatonek{{\bicat{grCat}^1_\kk}}    % bicategory of dg categories
\newcommand\DGCat{{\bicat{dgCat}}}                  % bicategory of dg categories
\newcommand\DGCatone{{\bicat{dgCat}^1}}             % bicategory of dg categories
\newcommand\DGCatdg{{\bicat{dgCat}^{dg}}}           % bicategory of dg categories
\newcommand\HoDGCat{{\bicat{Ho}(\DGCat)}}           % homotopy cat of dg categories
\newcommand\HoDGCatone{{\bicat{Ho}(\DGCatone)}}     % homotopy cat of dg categories
\newcommand\MoDGCatone{{\bicat{Mor}(\DGCatone)}}    % dg cats up to Morita equivalence
\def\EnhCatKC{{\bicat{EnhCat}_{kc}}}                % karoubi complete enhanced triangulated categories
\def\EnhCatKCdg{{\bicat{EnhCat}_{kc}^{dg}}}         % karoubi complete enhanced triangulated categories
\ifdgcal
  \DeclareMathOperator{\modd}{\mathcal{M}\mkern-1.5mu\mathit{o\mkern-1mud}}                    % dg module category
  \DeclareMathOperator{\hperf}{\mathcal{H}\mathit{perf}}                 % h-projective perfect modules (1-category version)
  \DeclareMathOperator{\perf}{{\mathcal{P}\mkern-1.5mu\mathit{erf}}}
  
  \DeclareMathOperator{\pretriag}{{\mathcal{P}\mkern-1.5mu\mathit{re\mhyphen}\mkern-2mu\mathcal{T}\mkern-3mu\mathit{r}}}
  \DeclareMathOperator{\DGFun}{{\mathcal{DGF}\mkern-1.5mu\mathit{un}}}
  
\else
  \DeclareMathOperator{\hperf}{Hperf}                 % h-projective perfect modules (1-category version)
  \DeclareMathOperator{\modd}{Mod}                    % dg module category
  \DeclareMathOperator{\perf}{Perf}
  
  \DeclareMathOperator{\pretriag}{Pre\mhyphen Tr}
  \DeclareMathOperator{\DGFun}{DGFun}
  
\fi

\DeclareMathOperator{\bihperf}{\bicat{Hperf}}               % h-projective perfect modules (2-category version)
\newcommand\sym{{\cat{S}}}                              % symmetric power of a dg category

\newcommand\numGgp[1]{K_0^{\mathrm{num}}\ifthenelse{\equal{#1}{}}{}{(#1)}}      % numerical Grothendieck group

\newcommand\modbar{{\overline{\modd}}}              % dg module bar category
\newcommand\rightmod[1]{{\modd\mkern-2mu\mhyphen #1}}   % right modules
\newcommand\leftmod[1]{{#1\mhyphen\mkern-2.5mu\modd}}   % left modules
\newcommand\bimod[2]{{#1\mhyphen\mkern-2.5mu\modd\mkern-2mu\mhyphen #2}}   % bimodules
\newcommand\bimodbar[2]{{#1\mhyphen\modbar\mhyphen #2}} % bimodules

\newcommand\halg[1]{H_{#1}}              % Heisenberg algebra
\newcommand\chalg[1]{\underline{H}_{#1}} % classical (not idempotent-modified) Heisenberg algebra
\newcommand\chalga[1]{\underline{H}^{A}_{#1}} % classical (not idempotent-modified) Heisenberg algebra, A definition
\newcommand\falg[1]{F_{#1}}              % Fock space
\newcommand\hcat[1]{\bicat{H}_{#1}}                     % The Heisenberg 2-category
\WithSuffix\newcommand\hcat*[1]{\bicat{H}'_{#1}}        % The non-completed Heisenberg 2-category
\newcommand\hcatadd[1]{\bicat{H}^{\mathrm{add}}_{#1}}                       % The additive Heisenberg 2-category
\WithSuffix\newcommand\hcatadd*[1]{\bicat{H}^{\mathrm{add}\prime}_{#1}}     % The non-completed additiveHeisenberg 2-category
\newcommand\monhcat[1]{\underline{\bicat{H}}_{#1}}                          % The monoidal Heisenberg category
\WithSuffix\newcommand\monhcat*[1]{\underline{\bicat{H}}'_{#1}}             % The non-completed monoidal Heisenberg category
\newcommand\PP{{\mathsf{P}}}            % 1-morphism P
\newcommand\QQ{{\mathsf{Q}}}            % 1-morphism Q
\newcommand\hAA{{\mathsf{A}}}           % Hocschild chain A
\newcommand\hunit{{\mathbb{1}}}         % The empty string 1-morphism
\newcommand\ho{N}                       % an object in \hcat
\WithSuffix\newcommand\ho*{N'}          % another object in \hcat

\newcommand\fcat[1]{\bicat{F}_{#1}}                   % The Fock space 2-category
\WithSuffix\newcommand\fcat*[1]{\bicat{F}'_{#1}}      % The non-completed Fock space 2-category
\newcommand\fcatadd[1]{\bicat{F}^{\mathrm{add}}_{#1}}                         % The additive Fock space 2-category
\WithSuffix\newcommand\fcatadd*[1]{\bicat{F}^{\mathrm{add}\prime}_{#1}}       % The non-completed additive Fock space 2-category
\tikzset{dot/.style={circle, fill, inner sep=0, minimum size=4pt}}
\tikzset{serre/.style={star, fill, star points=10, star point ratio = 2, scale=0.4}}
\tikzset{up/.style={color=blue}}
\tikzset{down/.style={color=red}}
\tikzset{cc/.style={color=green!50!black}}
\tikzset{desc/.style={fill=white}}
\tikzset{sym/.style={rectangle, draw, minimum width=1cm}}
\tikzset{many/.style={very thick}}
\newcommand\Spec{\operatorname{Spec}}   % spectrum of a ring
\newcommand\sheaf\mathcal               % sheaves
\DeclareMathOperator{\supp}{Supp}

\DeclareMathOperator{\img}{Im}
\DeclareMathOperator{\homm}{Hom}
\DeclareMathOperator{\tor}{Tor}
\DeclareMathOperator{\eend}{End}
\DeclareMathOperator{\cone}{Cone}
\DeclareMathOperator{\opp}{{opp}}

\DeclareMathOperator{\lder}{\bf L}

\DeclareMathOperator{\ldertimes}{\overset{\lder}{\otimes}}

\DeclareMathOperator{\eval}{eval}

\DeclareMathOperator{\alghh}{{AHH}}
\DeclareMathOperator{\algcat}{{{\mathcal{A}}lg}}
\DeclareMathOperator{\catalg}{{{\mathcal{C}}at}}
\DeclareMathOperator{\forget}{{Frg}}

\def\barperf{{\it\mathcal{P}\overline{er}f}}
\def\barperfA{{\barperf(\A)}}
\def\barperfB{{\barperf(\B)}}

\def\partn{{Part}}
\def\ordpartn{{OrdPart}}

\def\A{{\cat{A}}}
\def\B{{\cat{B}}}
\def\C{{\cat{C}}}
\def\D{{\cat{D}}}

\def\Aopp{{\cat{A}^{\opp}}}

\def\hperfA{{\hperf(\A)}}

\def\modk{{\rightmod{\kk}}}
\def\modA{{\rightmod \A}}

\def\Amod{{\leftmod \A}}

\def\AmodA{{\bimod{\A}{\A}}}

\def\AbimA{{\A\mhyphen\A}}

\def\AmodbarB{{\bimodbar\A\B}}

\def\Bperf{{\B\mhyphen\perf}}

\def\perfA{{\perf(\A)}}

\def\barA{{\bar{\A}}}

\def\bartimes{\mathrel{\overline{\otimes}}}

\def\degzero{{\mathrm{deg.\,0}}}

\newcommand\symbc[1]{{\sym^{#1}\basecat}}           % symmetric powers of \basecat (dg setting)
\def\symbcn{{\symbc{\ho}}}

\def\euler{{\mathrm{eu}}}
\def\strace{{\mathrm{strace}}}

\definecolor{purple}{rgb}{0.5, 0.0, 0.5}

\usepackage{mparhack}
\begin{document}
\begin{abstract}
For any noncommutative smooth and proper variety we construct
three actions of the Heisenberg algebra on the total Hochschild 
homology of its symmetric quotient stacks. 
One is defined algebraically using the orbifold decomposition of
Hochschild homology. One decategorifies the 
$2$-categorical Heisenberg action of Gyenge--Koppensteiner--Logvinenko. 
One is defined representation theoretically via explicit operators intrinsic
to the symmetric quotient stacks. We then show all three to coincide, 
and thus give alternative descriptions of one natural action. 
For ordinary commutative varities, we give a fourth, geometrical
description by correspondences similar to those used 
by Grojnowski and Nakajima for the Hilbert schemes of points 
on surfaces. 
\end{abstract}
\maketitle
\setcounter{tocdepth}{1}
\tableofcontents
\section{Introduction}
\label{sec:intro}
The Heisenberg algebra originated in quantum mechanics to describe
the commutation relations between position and momentum
operators. The $\infty$-dimensional Heisenberg algebra $\chalg{\mathbb{C}}$
has the generators $\left\{ a(i) \right\}_{ i \in \ZZ \setminus \{0\}}$ 
and the relations $[a(i),\, a(j)] = i\delta_{i,-j}$. It is 
important in many areas of mathematics and physics such as 
conformal field theory, string theory, and representation theory. 

In algebraic geometry, much of its relevance 
is due to the following celebrated result obtained independently by 
Grojnowski and Nakajima in the 1990s:
\begin{Theorem*}[{see
\cite[Theorem~3.1]{nakajima1997heisenberg}, \cite[Theorem~7]{grojnowski1995instantons} and \cite[Theorem~8.13]{nakajima1999lectures}}]
Let $X$ be a smooth projective surface over $\mathbb{C}$. 
Let $X^{[n]}$ be the Hilbert scheme of $n$ points on $X$. Let $\chi$ 
be the pairing on $H^\bullet(X,\mathbb{Q})$ given 
by the cup product and then the direct image along 
$X \rightarrow \text{pt}$. 

For each $\alpha \in H^*(X,\mathbb{Q})$ and $i > 0$, there are
operators 
$A_\alpha(-i)$ and $A_\alpha(i)$ on 
$\bigoplus_{n\geq{0}} H^\bullet(X^{[n]}, \mathbb{Q})$
defined by certain correspondences on $X^{[n]} \times X^{[n-i]}$ and 
$X^{[n]} \times X^{[n+i]}$ for $n \geq 0$. These satisfy 
\begin{small}
\begin{equation}
\label{eqn-intro-symmetric-heisenberg-relation}
A_{\alpha}(i) A_{\beta}(j) = (-1)^{\deg(\alpha)\deg(\beta)}
A_{\beta}(j)  A_{\alpha}(i) + 
\delta_{i,-j}i\langle \alpha,\, \beta\rangle_\chi
\end{equation}
\end{small}
and define an action of the 
Heisenberg algebra $\chalg{H^\bullet(X,\mathbb{Q}), \chi}$ on the
total cohomology $\bigoplus_{n=0}^{\infty} H^\bullet(X^{[n]},
\mathbb{Q})$.  This action 
identifies $\bigoplus_{n=0}^{\infty} H^\bullet(X^{[n]}, \mathbb{Q})$ 
with the Fock space of $\chalg{H^\bullet(X,\mathbb{Q}), \chi}$. 
\end{Theorem*}

Here, the \em Heisenberg algebra $\chalg{V, \chi}$ of a graded 
vector space $V$ with a symmetric bilinear form $\chi$ \rm
is a generalisation introduced in 
\cite{grojnowski1995instantons, nakajima1999lectures}. It has
the generators 
$\left\{ a_{v}(n) \right\}_{v \in V,\; n \in \ZZ \setminus \{0\}}$, 
the relations of linearity in $v$ and the Heisenberg relation
\eqref{eqn-intro-symmetric-heisenberg-relation}. 
%Roughly, in the theorem above the elements $a_{v}(n)$ and $a_{v}(-n)$ act by correspondences which add or remove, respectively, $n$ points belonging to $v$.
 
If $\dim X\geq{3}$ and $n\geq{4}$, $X^{[n]}$ is singular and badly behaved.
Grojnowski conjectured in \cite{grojnowski1995instantons}
that the result should hold in any dimension if one replaces $X^{[n]}$
by the symmetric quotient orbifold $X^n/S_n$ and uses equivariant
$K$-theory. It was later proved by Segal and Wang
\cite{Segal-EquivariantKTheoryAndSymmetricProducts,
Wang-EquivariantKTheoryWreathProductsAndHeisenbergAlgebra}. 

In this paper, for any noncommutative smooth and proper variety
$\basecat$ we construct three actions of the Heisenberg algebra 
$\chalg{\hochhom_\bullet(\basecat), \chi}$ on the total Hochschild 
homology $\bigoplus_{n \geq 0} \hochhom_\bullet(\symbc{n})$ of 
the symmetric quotient stacks $\symbc{n}$. The \em orbifold
decomposition action \rm uses the orbifold decomposition
theorem for Hochschild homology \cite{annobaranovskylogvinenko2023orbifold}
to decompose $\bigoplus_{n \geq 0} \hochhom_\bullet(\symbc{n})$
into the symmetrical algebra $S^*(\hochhom_\bullet(\basecat) \otimes t
\kk[t])$ which carries a canonical Fock space action of 
$\chalg{\hochhom_\bullet(\basecat), \chi}$ (see 
\S\ref{section-fock-space-and-its-symmetric-algebra-structure}). 
For the \em decategorified action\rm, we decategorify
the Heisenberg $2$-category and its categorial action on 
the symmetric quotient stacks constructed 
in \cite{gyenge2021heisenberg}. The 
\em induction-restriction action \rm is defined by explicit 
representation theoretic operators intrinsic to
$\hochhom_\bullet(\symbc{n})$
which resemble the $K$-theoretic operators of Segal and Wang \cite{Segal-EquivariantKTheoryAndSymmetricProducts,
Wang-EquivariantKTheoryWreathProductsAndHeisenbergAlgebra}. Our main 
theorem is that these three actions coincide and are thus 
different descriptions of the same natural action we construct of
$\chalg{\hochhom_\bullet(\basecat), \chi}$ on 
$\bigoplus_{n \geq 0} \hochhom_\bullet(\symbc{n})$. Finally, 
when $\basecat$ is an ordinary (commutative) variety 
we give the fourth description of this action as a \em geometrical action
\rm by correspondences similar to those of 
Grojnowski and Nakajima
\cite{grojnowski1995instantons,nakajima1997heisenberg}.

We now describe the results of this paper in more detail, beginning
with its general context. 
Noncommutative geometry comes in many flavors. Here we follow 
\cite{KontsevichSoibelman-NotesOnAInftyAlgebrasAInftyCategoriesAndNoncommutativeGeometry,KatzarkovKontsevichPantev-HodgeTheoreticAspectsOfMirrorSymmetry, 
Orlov-SmoothAndProperNoncommutativeSchemesAndGluingofDGcategories,
Efimov-HomotopyFinitenessOfSomeDGCategoriesFromAlgebraicGeometry}
where a noncommutative scheme is a small DG category
$\basecat$ over a base field $\mathbb{k}$ considered up to Morita 
equivalence. Such $\basecat$ can be viewed as a DG enhanced 
triangulated category
\cite{BondalKapranov-EnhancedTriangulatedCategories,
Toen-TheHomotopyTheoryOfDGCategoriesAndDerivedMoritaTheory,
Tabuada-InvariantsAdditifsDeDGCategories,
LuntsOrlov-UniquenessOfEnhancementForTriangulatedCategories}.
The triangulated category enhanced by $\basecat$ is $D_c(\basecat)$, 
the compact derived category of $\basecat$-modules. If $\basecat$ 
is a commutative algebra $A$, then this is 
the compact derived category $D_c(\Spec A)$ of quasi-coherent 
sheaves on the scheme $\Spec A$. On the other hand,
the compact derived category of 
any quasi-compact quasi-separated scheme $X$ can be enhanced by 
a noncommutative DG
algebra \cite{BondalVanDenBergh-GeneratorsAndRepresentabilityOfFunctorsInCommutativeAndNoncommutativeGeometry}.

The point of this approach is that we take any enhanced triangulated 
category $\basecat$ and treat it as if it were the compact derived category of a
``noncommutative'' scheme. A number of geometrical features
can be read off at this abstract level:
from smoothness and properness to Hodge and de-Rham cohomologies, 
see \cite{KontsevichSoibelman-NotesOnAInftyAlgebrasAInftyCategoriesAndNoncommutativeGeometry} and a beautiful summary in  
\cite{Kaledin-HomologicalMethodsInNoncommutativeGeometry,
Kaledin-NonCommutativeGeometryFromTheHomologicalPointOfView}. 
In particular, when $\basecat$ is (the derived category of) a smooth
projective variety $X$ over $\mathbb{C}$, by the global version 
of the Hochschild--Kostant--Rosenberg isomorphism 
\cite{Kontsevich-DeformationQuantizationOfPoissonManifolds,
Swan-HochschildCohomologyOfQuasiprojectiveSchemes}
the Hochschild homology $\hochhom_\bullet(\basecat)$ is isomorphic to the
complex cohomology $H^\bullet(X,\mathbb{C})$. Moreover, the symmetric
power $\sym^n \basecat$ is the stack $[X^n/S_n]$
and $\hochhom_\bullet(\symbc{n})$ is isomorphic to the Chen--Ruan
orbifold cohomology $H_{orb}^\bullet\left(X^n/S_n,\mathbb{C}\right)$
\cite{Baranovsky-OrbifoldCohomologyAsPeriodicCyclicHomology}. 

Thus the paradigm of this paper is: 
\begin{align*}
\text{ smooth projective variety $X$ over $\mathbb{C}$ }
&\longrightarrow
\text{ smooth and proper DG category $\basecat$ over $\kk$ }
\\
\text{ symmetric quotient stack $[X^n/S_n]$}
&\longrightarrow
\text{ symmetric power $\symbc{n}$}
\\
\text{ complex cohomology $H^\bullet(X,\mathbb{C})$}
&\longrightarrow
\text{ Hochschild homology $\hochhom_\bullet(\basecat)$}
\\
\text{ Chen--Ruan orbifold cohomology
$H_{orb}^\bullet\left(X^n/S_n,\mathbb{C}\right)$ }
&\longrightarrow
\text{ Hochschild homology $\hochhom_\bullet(\symbc{n})$}
\end{align*}
All the results we prove for the right column specialise
to the left one by setting $\kk := \mathbb{C}, \basecat := D^b_{coh}(X)$. 

A key technical tool we use is the noncommutative orbifold
decomposition theorem established recently by Anno, Baranovsky, 
and the second author \cite{annobaranovskylogvinenko2023orbifold}
motivated by the ideas of  
\cite{GetzlerJones-TheCyclicHomologyOfCrossedProductAlgebras}\cite{Baranovsky-OrbifoldCohomologyAsPeriodicCyclicHomology}. 
For any DG category $\basecat$ and any $n \geq 0$, 
\cite[Theorem 4.1]{annobaranovskylogvinenko2023orbifold}
constructs an explicit isomorphism
\begin{equation}
\label{eqn-intro-noncommutative-orbifold-decomposition-for-sn-v}
\hochhom_\bullet(\symbc{n})
\simeq 
\bigoplus_{\underline{n} \vdash n}
\Sym^{r_1(\underline{n})} \hochhom_\bullet(\basecat)
\otimes 
\Sym^{r_2(\underline{n})} \hochhom_\bullet(\basecat)
\otimes 
\dots 
\otimes 
\Sym^{r_n(\underline{n})} \hochhom_\bullet(\basecat),
\end{equation}
where $\underline{n}$ is an unordered partition of $n$
and $r_i(\underline{n})$ is its number of parts of size $i$. 

For the partition $(n)$ consisting of a single part of size $n$, 
the corresponding summand in 
\eqref{eqn-intro-noncommutative-orbifold-decomposition-for-sn-v} is
just $\hochhom_\bullet(\basecat)$. 
The inclusion of and the projection onto this summand defines 
the linear maps 
\begin{equation}
\label{eqn-linear-map-psi-n-and-kappa-n}
\begin{tikzcd}
\hochhom_\bullet(\basecat) 
\ar[shift left = 1ex]{r}{\psi_n}
&
\hochhom_\bullet(\sym^n \basecat),
\ar[shift left = 1ex]{l}{\kappa_n}
\end{tikzcd}
\end{equation}
which are used below to construct the decategorified and 
the induction-restriction actions. 

Let $\chi$ be the Euler pairing on $\hochhom_\bullet(\basecat)$. 
Decompositions \eqref{eqn-intro-noncommutative-orbifold-decomposition-for-sn-v}
combine into an isomorphism 
\begin{equation}
\label{eqn-intro-symmetric-algebra-homotopy-equivalence-hh}
\bigoplus_{n \geq 0} \hochhom_\bullet\left(\sym^n\basecat\right)
\;\simeq\;
S^*\left(\hochhom_\bullet(\basecat) \otimes t \kk[t]\right). 
\end{equation}
For any $V$ and $\chi$, 
the symmetric algebra $S^*\left(V \otimes t \kk[t]\right)$
is known to be canonically isomorphic to the Fock space
representation of $\chalg{V, \chi}$, see \S\ref{section-fock-space-and-its-symmetric-algebra-structure}. 
Define the \em orbifold decomposition action \rm 
of $\chalg{\hochhom_\bullet(\basecat),\chi}$ on 
$\bigoplus_{n \geq 0} \hochhom_\bullet\left(\sym^n \basecat\right)$ to be 
the action induced 
by \eqref{eqn-intro-symmetric-algebra-homotopy-equivalence-hh}. 
In terms of decompositions
\eqref{eqn-intro-noncommutative-orbifold-decomposition-for-sn-v}, it
has the following simple algebraic form. For any $n \geq 0$, 
any partition $\underline{n} \vdash n$, and any element
$$
\underline{\gamma}_1 \otimes \dots \otimes \underline{\gamma}_i \otimes \dots 
\otimes \underline{\gamma}_n 
\quad \in \quad  
\Sym^{r_1(\underline{n})} \hochhom_\bullet(\basecat)
\otimes 
\dots 
\otimes \Sym^{r_n(\underline{n})} \hochhom_\bullet(\basecat)
\quad \quad \quad $$
define the operators 
$A_\alpha(-i)$ and $A_\alpha(i)$ for 
$\alpha \in \hochhom_\bullet(\basecat)$ and $i > 0$ by
\begin{small}
\begin{align*}
A_\alpha(-i)\left(\underline{\gamma}_1 \otimes \dots \otimes \underline{\gamma}_i \otimes \dots 
\otimes \underline{\gamma}_n\right)
& \; := \;
(-1)^{|\alpha|\left(|\underline{\gamma}_1| + \dots + |\underline{\gamma}_{i-1}|\right)} \;
\underline{\gamma}_1 \otimes \dots \otimes \left<\alpha,
\underline{\gamma}_i\right>_\chi \otimes \dots \otimes \underline{\gamma}_n 
\\
A_\alpha(i)\left(\underline{\gamma}_1 \otimes \dots \otimes \underline{\gamma}_i \otimes \dots 
\otimes \underline{\gamma}_n\right)
& \; := \; 
(-1)^{|\alpha|\left(|\underline{\gamma}_1| + \dots + |\underline{\gamma}_{i-1}|\right)} \;
\underline{\gamma}_1 \otimes \dots \otimes (\alpha \vee
\underline{\gamma}_i) \otimes \dots \otimes \underline{\gamma}_n, 
\end{align*}
\end{small}
where for any $m > 0$ we define the bracket
$\left<\alpha,-\right>_\chi\colon \Sym^m \hochhom_\bullet(\basecat) \rightarrow 
\Sym^{m-1} \hochhom_\bullet(\basecat)$ to send $\beta_1 \vee \dots \vee
\beta_m$ to $\sum_{j = 1}^{m} (-1)^{|\alpha|\left(|\beta_1| + \dots +
|\beta_{j-1}|\right)}
\left<\alpha,\beta_j\right>_\chi
\beta_1 \vee \dots \vee \beta_{j-1} \vee \beta_{j+1} \vee \dots \vee
\beta_m$. 

Next, we decategorify DG 2-categorical constructions of 
\cite{gyenge2021heisenberg} via Hochschild homology. 
In \cite{gyenge2021heisenberg}, we constructed the \em Heisenberg DG 
$2$-category $\hcat\basecat$ \rm of $\basecat$. In the language above, 
$\hcat\basecat$ is a noncommutative scheme version of 
the Heisenberg algebra of $\basecat$. We also constructed its action 
$\Phi_\basecat$ on the $2$-category of the symmetric powers 
$\sym^n \basecat$, the \em categorical Fock space \rm of
$\basecat$. Here, in \S\ref{section-decategorification-map} we define 
the functor $\hochhom_{alg}$ which flattens the Hochschild homology of 
a DG $2$-category into an algebra, see 
Definition~\ref{defn-hhalg-flattening-functor}. We then apply it 
to $\Phi_\basecat$ to obtain an algebra homomorphism 
\begin{equation*}
\hochhom_{alg}(\hcat\basecat) 
\xrightarrow{\eqref{eqn-decategorigication-of-phi-basecat}}
\End \bigl( \bigoplus_{n \geq 0} \hochhom_{\bullet}(\symbc{n}) \bigr). 
\end{equation*}
In the rest of \S\ref{section-hochschild-homology-and-heisenberg-2-category}
we construct the \em decategorification map\rm, that is -- an injective algebra
homomorphism 
\begin{equation*}
\chalg{\hochhom_\bullet(\basecat), \chi}
\xrightarrow{\eqref{eqn-injective-decategorification-map-into-the-flattening}}
\hochhom_{alg}(\hcat\basecat). 
\end{equation*}
We define the \em decategorified action \rm of
$\chalg{\hochhom_\bullet(\basecat),\chi}$ on $\bigoplus_{n \geq 0}
\hochhom_\bullet\left(\sym^n \basecat\right)$ by 
composing
\eqref{eqn-decategorigication-of-phi-basecat} and 
\eqref{eqn-injective-decategorification-map-into-the-flattening}. 

Finally, we define the \em induction-restriction action \rm 
of $\chalg{\hochhom_\bullet(\basecat),\chi}$ on $\bigoplus_{n \geq 0}
\hochhom_\bullet\left(\sym^n \basecat\right)$ by explicit representation 
theoretic formulas below. These are intrinsic to the symmetric
quotient stacks $\symbc{n}$ and use the functors of
induction/restriction with respect to the subgroup 
$S_i \times S_n < S_{n+i}$ which permutes the first $i$ and the last $n$
elements separately. They also use the maps 
$\psi_i$ and $\kappa_i$ of \eqref{eqn-linear-map-psi-n-and-kappa-n}.

Our main theorem is: 
\begin{Theorem}[see Theorem 
\ref{theorem-three-heisenberg-actions-on-hh-of-noncomm-sym-quot-stacks}]
\label{theorem-intro-noncommutative-grojnowski-nakajima-action}
Let $\basecat$ be a smooth and proper DG category over an
algebraically closed field $\kk$ of characteristic $0$.  Let $\chi$ be
the Euler pairing on the Hochschild homology $\hochhom_\bullet(\basecat)$.

For all $\alpha \in \hochhom_\bullet(\basecat)$ and $i > 0$,
the following three definitions of the operators 
$A_\alpha(-i)$ and $A_\alpha(i)$ on 
$\bigoplus_{n\geq{0}} HH_\bullet(\symbc{n})$ are equivalent:
\begin{enumerate}
\item \underline{Orbifold decomposition operators}: define $A_\alpha(-i), A_\alpha(i)$  
in terms of isomorphisms 
\eqref{eqn-intro-noncommutative-orbifold-decomposition-for-sn-v}, 
where for any $n \geq 0$, $\underline{n} \vdash n$, and 
$\underline{\gamma}_1 \otimes \dots \otimes \underline{\gamma}_n 
\in \Sym^{r_1(\underline{n})} \hochhom_\bullet(\basecat)
\otimes 
\dots 
\otimes \Sym^{r_n(\underline{n})} \hochhom_\bullet(\basecat)$ 
\begin{small}
\begin{align*}
A_\alpha(-i)\left(\underline{\gamma}_1 \otimes \dots \otimes \underline{\gamma}_i \otimes \dots 
\otimes \underline{\gamma}_n\right)
& \; := \;
(-1)^{|\alpha|\left(|\underline{\gamma}_1| + \dots + |\underline{\gamma}_{i-1}|\right)} \;
\underline{\gamma}_1 \otimes \dots \otimes \left<\alpha,
\underline{\gamma}_i\right>_\chi \otimes \dots \otimes \underline{\gamma}_n 
\\
A_\alpha(i)\left(\underline{\gamma}_1 \otimes \dots \otimes \underline{\gamma}_i \otimes \dots 
\otimes \underline{\gamma}_n\right)
& \; := \; 
(-1)^{|\alpha|\left(|\underline{\gamma}_1| + \dots + |\underline{\gamma}_{i-1}|\right)} \;
\underline{\gamma}_1 \otimes \dots \otimes (\alpha \vee
\underline{\gamma}_i) \otimes \dots \otimes \underline{\gamma}_n. 
\end{align*}
\end{small}

\item \underline{Decategorified operators}: define 
$A_\alpha(-i), A_\alpha(i)$ by applying to $a_\alpha(-i),
a_\alpha(i) \in \chalg{\hochhom_\bullet(\basecat), \chi}$ 
the decategorification map $\chalg{\hochhom_\bullet(\basecat), \chi}
\xrightarrow{\eqref{eqn-injective-decategorification-map-into-the-flattening}}
\hochhom_{alg}(\hcat\basecat)$ followed by the Hochschild homology
flattening $\hochhom_{alg}(\hcat\basecat) 
\xrightarrow{\eqref{eqn-decategorigication-of-phi-basecat}}
\End \bigl( \bigoplus_{n \geq 0} \hochhom_{\bullet}(\symbc{n}) \bigr)$
of the categorical action of \cite{gyenge2021heisenberg}. 
\item \underline{Induction-restriction operators}: define 
$A_\alpha(-i), A_\alpha(i)$ by setting 
\begin{scriptsize}
\begin{equation*}
A_\alpha(-i) \colon 
\hochhom_\bullet \left(\sym^{n+i} \basecat\right) 
\xrightarrow{\Res^{S_{n+i}}_{S_i \times S_n}}
\hochhom_\bullet \left(\sym^{i} \basecat \otimes \sym^{n} \basecat\right) 
\simeq 
\hochhom_\bullet \left(\sym^{i} \basecat \right) \otimes 
\hochhom_\bullet\left(\sym^{n} \basecat\right) 
\xrightarrow{\left<\alpha, \kappa_i(-)\right>}
\hochhom_\bullet \left(\sym^{n} \basecat \right),
\end{equation*}
\begin{equation*}
A_\alpha(i) \colon 
\hochhom_\bullet \left(\sym^{n} \basecat \right)
\xrightarrow{\psi_i(\alpha)\otimes(-)}
\hochhom_\bullet \left(\sym^{i} \basecat \right)
\otimes 
\hochhom_\bullet \left(\sym^{n} \basecat \right)
\simeq 
\hochhom_\bullet \left(\sym^{i} \basecat 
\otimes \sym^{n} \basecat \right)
\xrightarrow{\Ind^{S_{n + i}}_{S_{i} \times S_n}}
\hochhom_\bullet \left(\sym^{n + i} \basecat \right), 
\end{equation*}
\end{scriptsize}
where $\psi_i$ and $\kappa_i$ are the maps
\eqref{eqn-linear-map-psi-n-and-kappa-n} defined explicitly in 
\S\ref{section-noncommutative-orbifld-decomposition-for-hh}. 
\end{enumerate}

These operators satisfy 
\begin{small}
\begin{equation}
\label{eqn-intro-heisenberg-a-commutation-relation}
A_{\alpha}(i) A_{\beta}(j) - (-1)^{\deg(\alpha)\deg(\beta)}
A_{\beta}(j)  A_{\alpha}(i) 
= 0  \quad \quad \quad \quad \quad \quad \quad\;\;\; i,j > 0 \; \text{ or } \; i,j < 0, 
\end{equation}
\begin{equation}
\label{eqn-intro-heisenberg-a-heisenberg-relation}
A_{\alpha}(-i) A_{\beta}(j) - (-1)^{\deg(\alpha)\deg(\beta)}
A_{\beta}(j)  A_{\alpha}(-i)  = 
 \delta_{i,j} i \langle \alpha,\, \beta\rangle_\chi,
\quad \quad \quad \quad \quad \quad i,j > 0 
\end{equation}
\end{small}
and thus define an action of the 
Heisenberg algebra $H_{HH_\bullet(\basecat), \chi}$ 
on $\bigoplus_{n\geq{0}} HH_\bullet(\symbc{n})$.  This action
identifies $\bigoplus_{n\geq{0}} HH_\bullet(\symbc{n})$ with the
Fock space of $H_{HH_\bullet(\basecat), \chi}$. 
\end{Theorem}

In the commutative case, when $\basecat$ is a smooth
projective variety $X$ of \em any \rm dimension, 
we give a geometrical construction 
of this action. It uses correspondences on fiber products 
similar to those used for Hilbert schemes by 
Grojnowski and Nakajima \cite{grojnowski1995instantons,nakajima1997heisenberg}.
To state our result geometrically, we switch from 
Hochschild homology to Chen--Ruan orbifold cohomology 
\cite{ChenRuan-ANewCohomologyTheoryofOrbifold,
FantechiGottsche-OrbifoldCohomologyForGlobalQuotients}:
\begin{equation*}
 H_{orb}^\bullet\left(X^n/S_n,\mathbb{C}\right)
: = \left(\bigoplus_{\sigma \in S_n}
H^\bullet\left((X^n)^\sigma,\mathbb{C}\right) \right)_{S_n} 
\simeq 
\bigoplus_{\underline{n} \vdash n}
H^\bullet\left(X_{\underline{n}},\mathbb{C}\right)_{S_{\underline{r}(\underline{n})}}, 
\end{equation*}
where $(X^n)^\sigma \subset X^n$ is the fixed point locus of 
$\sigma \in S_n$ and $(-)_{S_n}$ denotes taking the coinvariants. Moreover,
$\underline{n} \vdash n$ is a partition,
$r_i(\underline{n})$ is its number of parts of size $i$, and
$X_{\underline{n}} := X^{r_1(\underline{n})} \times \dots \times 
X^{r_n(\underline{n})}$ is isomorphic to $(X^n)^\sigma$ for any
$\sigma$ of cycle type $\underline{n}$. Finally, 
$S_{\underline{r}(\underline{n})} := S_{r_1(\underline{n})} 
\times \dots \times S_{r_n(\underline{n})}$.

Using the K{\"u}nneth formula to further decompose each
$H^\bullet\left(X_{\underline{n}},\mathbb{C}\right)$ yields an 
isomorphism $\bigoplus_{n\geq{0}}
H_{orb}^\bullet\left(X^n/S_n,\mathbb{C}\right) \simeq 
S^*(H^\bullet(X,\mathbb{C}) \otimes t\mathbb{C}[t])$, see Lemma 
\ref{lemma-symmetrical-algebra-structure-on-total-chen-ruan-cohomology}. 
In view of \eqref{eqn-intro-symmetric-algebra-homotopy-equivalence-hh}, 
we see that the HKR isomorphism $\hochhom_\bullet(\basecat) 
\simeq H^\bullet(X,\mathbb{C})$ identifies the total Hochschild
homology $\bigoplus_{n \geq 0} \hochhom_\bullet\left(\sym^n \basecat\right)$
with the total Chen--Ruan cohomology 
$\bigoplus_{n\geq{0}}
H_{orb}^\bullet\left(X^n/S_n,\mathbb{C}\right)$. 

\begin{Definition*}[see
Definition~\ref{defn-generalised-grojnowski-nakajima-correspondence}]
For any $i > 0$, $n \geq 0$, and a partition $\underline{n} \vdash n$, 
define the subvariety
\begin{equation}
\label{eqn-intro-generalised-grojnowski-nakajima-correspondence}
Z(i) := \left\{ (\underline{x}, \underline{y}, z) 
\;\middle|\;
\begin{matrix}
\forget(\underline{x}_j) = \forget(\underline{y}_j) + z & j = i,
\\
\forget(\underline{x}_j) = \forget(\underline{y}_j) & j \neq i.
\end{matrix}
\right\} \subset 
X_{\underline{n} + (i)} \times X_{\underline{n}}  \times X.
\end{equation}
Here $(i)$ is the single part partition of $i$ 
and $\underline{n} + (i)$ is the union 
of partitions $\underline{n}$ and $(i)$. We write
$\underline{x} \in X_{\underline{n}+(i)}$ as 
$(\underline{x}_1, \dots, \underline{x}_n)$ with 
$\underline{x}_j \in X^{r_j(\underline{n}+(i))}$ and analogously
for $\underline{y} \in X_{\underline{n}}$. For any $m > 0$,
$\forget\colon X^m \rightarrow X^m/S_m$ is the $S_m$-quotient map 
which forgets the ordering on an $m$-tuple of points.
\end{Definition*}
\begin{Theorem}[see Theorem
\ref{theorem-generalised-nakajima-grojnowski-heisenberg-action}]
\label{theorem-intro-generalised-nakajima-grojnowski-heisenberg-action}
Let $X$ be a smooth projective variety over $\mathbb{C}$ and $\chi$
be the pairing on $H^\bullet(X,\mathbb{C})$ given 
by taking the cup product and then the direct image along 
$X \rightarrow \text{pt}$. 

For each $\alpha \in H^\bullet(X, \mathbb{C})$ and $i > 0$,
define operators $A_\alpha(\pm i)$ on 
$\bigoplus_{n\geq{0}}
H_{orb}^\bullet\left(X^n/S_n,\mathbb{C}\right)$ by
\begin{equation*}
A_\alpha(-i)\colon 
H^\bullet(X_{\underline{n} + (i)}, \mathbb{C})_{S_{\underline{r}(\underline{n} + (i))}}
\xrightarrow{\frac{1}{|S_{\underline{r}(\underline{n})}|} \pi_{12 *}
\left([Z(i)] \cap \pi_3^* \alpha\right)
}
H^\bullet(X_{\underline{n}}, \mathbb{C})_{S_{\underline{r}(\underline{n})}}
\end{equation*}
\begin{equation*}
A_\alpha(i)\colon 
H^\bullet(X_{\underline{n}}, \mathbb{C})_{S_{\underline{r}(\underline{n})}}
\xrightarrow{\frac{1}{|S_{\underline{r}(\underline{n} + (i))}|}
\pi_{21 *} \left([Z(i)] \cap \pi_3^* \alpha\right)
}
H^\bullet(X_{\underline{n} + (i)}, \mathbb{C})_{S_{\underline{r}(\underline{n} + (i))}},
\end{equation*}
where $\pi_{12}$, $\pi_{21}$, and $\pi_{3}$ are the 
projections from 
\eqref{eqn-intro-generalised-grojnowski-nakajima-correspondence}
to
$X_{\underline{n} + (i)} \times X_{\underline{n}}$, 
to
$X_{\underline{n}} \times X_{\underline{n} + (i)}$, 
and to $X$. 

These define an action of the Heisenberg 
algebra $\chalg{H^\bullet(X,\mathbb{C}), \chi}$ on  
$\bigoplus_{n\geq{0}} H_{orb}^\bullet\left(X^n/S_n,\mathbb{C}\right)$
which equals its canonical action on 
$S^*\left(H^\bullet(X,\mathbb{C}) \otimes_\mathbb{C}
t\mathbb{C}[t]\right)$ via the isomorphism of Lemma 
\ref{lemma-symmetrical-algebra-structure-on-total-chen-ruan-cohomology}.
\end{Theorem}

The last assertion of Theorem 
\ref{theorem-intro-generalised-nakajima-grojnowski-heisenberg-action}
implies that the HKR isomorphism $\hochhom_\bullet(\basecat) 
\simeq H^\bullet(X,\mathbb{C})$ identifies its geometrical action
with the noncommutative action of Theorem
\ref{theorem-intro-noncommutative-grojnowski-nakajima-action}. 
Indeed, we have 
\begin{small}
$$
\bigoplus_{n \geq 0} \hochhom_\bullet\left(\sym^n \basecat\right)
\overset{\eqref{eqn-intro-symmetric-algebra-homotopy-equivalence-hh}}{\simeq}
S^*\left(\hochhom_\bullet(\basecat) \otimes_\mathbb{C} t\mathbb{C}[t]\right)
\overset{\text{HKR}}{\simeq}
S^*\left(H^\bullet(X,\mathbb{C}) \otimes_\mathbb{C} t\mathbb{C}[t]\right)
\overset{\text{Lemma }\ref{lemma-symmetrical-algebra-structure-on-total-chen-ruan-cohomology}}{\simeq}
\bigoplus_{n\geq{0}}
H_{orb}^\bullet\left(X^n/S_n,\mathbb{C}\right).$$
\end{small}
The action of Theorem
\ref{theorem-intro-noncommutative-grojnowski-nakajima-action} was 
defined by the canonical action (see
\S\ref{section-fock-space-and-its-symmetric-algebra-structure})
on $S^*\left(\hochhom_\bullet(\basecat) \otimes_\mathbb{C} t\mathbb{C}[t]\right)$, while the action of Theorem 
\ref{theorem-intro-generalised-nakajima-grojnowski-heisenberg-action}
was shown to equal the canonical action
on $S^*\left(H^\bullet(X,\mathbb{C}) \otimes_\mathbb{C}
t\mathbb{C}[t]\right)$. Let $\chi_{eul}$ and $\chi_{int}$
be the Euler pairing on $\hochhom_\bullet(\basecat)$
and the intersection product on $H^\bullet(X,\mathbb{C})$. 
We now need the HKR isomorphism to induce an algebra
isomorphism
$\chalg{\hochhom_\bullet(\basecat), \chi_{eul}}
\simeq
\chalg{H^\bullet(X,\mathbb{C}), \chi_{int}} $ identifying 
their canonical actions on the symmetric algebras. This seems obvious, 
until one realises that the HKR isomorphism does not identify 
$\chi_{eul}$ and $\chi_{int}$. However, both pairings are
non-degenerate and in Lemma
\ref{lemma-canonical-actions-of-heisenberg-on-symmetric-are-iso}
for any two non-degenerate pairings 
we give a natural isomorphism between the Heisenberg algebras 
which identifies the canonical actions. This settles the issue. 

To sum up, Theorems 
\ref{theorem-intro-noncommutative-grojnowski-nakajima-action} and \ref{theorem-intro-generalised-nakajima-grojnowski-heisenberg-action}
say that for any smooth and proper noncommutative variety $\basecat$
\\
there is a natural action of $\chalg{\hochhom_\bullet(\basecat), \chi}$ 
on $\bigoplus_{n \geq 0} \hochhom_\bullet\left(\sym^n \basecat\right)$ which 
admits several equivalent descriptions: an algebraic, a
representation-theoretic, and, in commutative case, also a geometric. 
Moreover, this action is the projection to the Hochschild homology of
the categorical action of \cite{gyenge2021heisenberg}.

A few more words about the categorical picture.
The main idea of \cite{gyenge2021heisenberg} 
was that instead of constructing the action of the Heisenberg
algebra separately for each additive invariant 
of the stacks $[X^n/S_n]$ 
($K$-theory, orbifold cohomology, Hochschild homology, etc) 
we would construct the noncommutative Heisenberg scheme of $X$ and 
have it act directly on $[X^n/S_n]$. This universal 
action could then be decategorified via any additive invariant 
to give an analogue of 
Theorem \ref{theorem-intro-noncommutative-grojnowski-nakajima-action}.   

The decategorification process is far from automatic. 
In \cite{gyenge2021heisenberg}
we decategorified our constructions using the numerical Grothendieck 
group $K_0^{\text{num}}$. That meant defining an injective algebra map 
\begin{equation*}
\pi\colon \chalg{\numGgp{\basecat}} \hookrightarrow \numGgp{\hcat\basecat}.
\end{equation*}
Ideally, we want $\pi$ to be an isomorphism, but injectivity is
enough for categorical actions of $\hcat\basecat$
to induce actions of $\chalg{\numGgp{\basecat}}$ on
$\numGgp{-}$. Our universal action $\Phi_\basecat$
induced an action of $\chalg{\numGgp{\basecat}}$ on $\bigoplus_{n
\geq 0} \numGgp{\sym^n \basecat}$. This induced in turn an embedding
of the Fock space of $\chalg{\numGgp{\basecat}}$
$$ \phi\colon \falg{\numGgp{\basecat}} \hookrightarrow 
\bigoplus_{n \geq 0} \numGgp{\sym^n \basecat}. $$
However, $K_0^{\text{num}}$ was not the best invariant to use: 
it fails the K{\"u}nneth formula. This leads to the cases
where $\phi$ is not surjective. By \cite[Theorem 8.13]{gyenge2021heisenberg} 
$\pi$ couldn't then be surjective either. 

The referees of \cite{gyenge2021heisenberg} pointed out that 
to claim our categorical constructions to be a universal version 
of Theorem
\ref{theorem-intro-noncommutative-grojnowski-nakajima-action}
we best show that they can be decategorified with other additive
invariants. We agree and construct in this 
paper the Hochschild homology decategorification map:
\begin{Theorem}[Theorem
\ref{theorem-construction-of-the-decategorification-map} and
Proposition~\ref{prps-decategorification-map-pi-is-injective}]
\label{theorem-intro-decategorification-map-pi}
Let $\basecat$ be a smooth and proper DG category over $\kk$. 
There exists an injective algebra homomorphism 
\begin{equation}
\label{eqn-intro-decategorification-map-pi}
\pi\colon \chalg{\hochhom_\bullet(\basecat)} \hookrightarrow 
\hochhom_{alg}\left(\hcat\basecat\right). 
\end{equation}
\end{Theorem}
We use $\pi$ to construct 
the decategorified action of $\chalg{\hochhom_\bullet(\basecat)}$ on  
$\bigoplus \hochhom_\bullet(\symbc{n})$. This again induces 
an embedding $\phi\colon \falg{\hochhom_\bullet(\basecat)} \hookrightarrow 
\bigoplus_{n \geq 0} \hochhom_\bullet(\sym^n \basecat)$. But now, 
by the dimension count provided by the orbifold
decomposition \eqref{eqn-intro-symmetric-algebra-homotopy-equivalence-hh},
$\phi$ is always an isomorphism. It encourages us to conjecture that:
\begin{Conjecture}
The injective algebra homomorphism
\eqref{eqn-intro-decategorification-map-pi} is always an isomorphism. 
\end{Conjecture}

In proving Theorem 
\ref{theorem-intro-noncommutative-grojnowski-nakajima-action}, 
we were encouraged by recent results of Belmans, Fu, and Krug 
\cite{BelmansFuKrug-HochschildCohomologyOfHilbertSchemesOfPointsOnSurfaces}.
In the commutative case, they computed 
$\dim \bigoplus_{n=0}^{\infty} \hochhom_\bullet\left(\sym^n\basecat\right)$
and showed that it matches the dimension of the Fock space
of $\chalg{\hochhom_\bullet\left(\basecat\right)}$. They conjectured
\cite[Conj.~3.24]{BelmansFuKrug-HochschildCohomologyOfHilbertSchemesOfPointsOnSurfaces} that the same holds in the noncommutative case. 
This matched our expectations 
\cite[Cor.~8.6]{gyenge2021heisenberg} and led 
Anno, Baranovsky and the second author to prove
the conjecture in \cite{annobaranovskylogvinenko2023orbifold}. 
It was also independently proved by Nordstrom 
\cite{Nordstrom-HochschildHomologyOfSymmetricPowersOfDGcategories}
via general considerations which do not yield explicit maps of 
Hochschild complexes. 

It remains to sum up our proof of Theorem 
\ref{theorem-intro-decategorification-map-pi}. 
For all $\alpha \in \hochhom_\bullet(\basecat)$ and $i > 0$ we define
the classes 
$\hAA_\alpha(-i), \hAA_\alpha(i) \in \hochhom_{alg}(\hcat\basecat)$
as the images of $\psi_i(\alpha)$ under  
functors $\Xi_\QQ,\Xi_\PP$ of
\cite[\S6.1]{gyenge2021heisenberg}. 
For these classes, we prove the commutation relation
\eqref{eqn-intro-heisenberg-a-commutation-relation}
and the Heisenberg relation
\eqref{eqn-intro-heisenberg-a-heisenberg-relation}. 
The commutation relation is easy
because it holds for the classes $\psi_i(\alpha)$ which 
live in the \em graded commutative \rm algebra $\bigoplus_{n\geq{0}} \hochhom_\bullet(\sym^n\basecat)$. Proving the Heisenberg relation
was the hardest step. Some $1$-morphism identities 
we need only hold up to homotopy in $\hcat\basecat$. 
This makes explicit computations with Hochschild chains difficult. 
To sidestep this, we construct two functors $\Xi_{\QQ\PP}$ and $\bigoplus_k
\Xi_{\QQ\PP}(\hat{k})$ from $\sym^i \basecat^{\opp} \otimes \sym^j
\basecat$ to $\hcat\basecat$ and show that 
for any $\alpha, \beta \in \hochhom_\bullet(\basecat)$ 
the images of $\psi_i(\alpha) \otimes \psi_j(\beta)$ under these two
functors are the LHS and the RHS of the Heisenberg relation
for $\hAA_\alpha(-i)$ and $\hAA_\beta(j)$. 
Homotopy equivalent functors induce the same map on Hochshchild homology 
\cite[Lemma 3.4]{Keller-OnTheCyclicHomologyOfExactCategories}, so
we complete the argument by constructing in Theorem 
\ref{theorem-functorial-categorification-of-the-PQ-Heisenberg-relation}
a functorial homotopy equivalence 
\begin{equation}
\label{eqn-intro-functorial-categorification-of-the-PQ-Heisenberg-relation}
\bigoplus_k \Xi_{\PP\QQ}(\hat{k})
\longrightarrow 
\Xi_{\QQ\PP}. 
\end{equation}

The functorial homotopy equivalence
\eqref{eqn-intro-functorial-categorification-of-the-PQ-Heisenberg-relation}
can be viewed as a functorial categorification of both
the A-generator Heisenberg relation 
\eqref{eqn-intro-heisenberg-a-heisenberg-relation}
and the PQ-generator Heisenberg relation
of \cite[\S2.2]{gyenge2021heisenberg}. 
A non-functorial categorification 
of the latter appeared in 
\cite[Theorem 6.3]{gyenge2021heisenberg}. It could not be 
made functorial directly as it related
the symmetrised elements $\PP_a^{(i)}$ and $\QQ_b^{(j)}$ which are 
not functorial in $a,b \in \basecat$ for $i,j > 1$. 
Instead, we use its case $i=j=1$ to iteratively construct
\eqref{eqn-intro-functorial-categorification-of-the-PQ-Heisenberg-relation}. 
Applying \eqref{eqn-intro-functorial-categorification-of-the-PQ-Heisenberg-relation} to $\psi_i(\alpha) \otimes \psi_j(\beta)$ yields 
the A-generator Heisenberg relation 
\eqref{eqn-intro-heisenberg-a-heisenberg-relation}
we prove in Theorem 
\ref{theorem-heisenberg-relation-for-the-assignments-of-pi},
while applying it to the product of the symmetrised powers 
$a^{(i)} \in \sym^{i}\basecat^{\opp} $ and $b^{(j)} \in
\sym^{j}\basecat$ recovers the
categorified PQ-generator Heisenberg relation in 
\cite[Theorem 6.3]{gyenge2021heisenberg}. 

\subsection*{Acknowledgments}
We would like to thank Rina Anno, Alexey Bondal, Ian Grojnowski,
Bernhard Keller, Clemens Koppensteiner and Hiraku Nakajima for useful discussions in the 
course of writing this paper. We would like to thank the anonymous
referees of \cite{gyenge2021heisenberg} for their many invaluable
comments and for pushing us to investigate 
other additive invariants and, in particular, Hochschild homology. 
Without them, this paper would have never been written. 
We would like to thank Hiraku Nakajima for his invaluable feedback 
on our first draft. 
We would also like to thank Alexey Bondal for pointing out a silly mistake 
with $\hochhom_0$ in an earlier version of this paper.  
Á.Gy.~was supported by the János Bolyai Research Scholarship and by the ``Lendület (Momentum)'' Program of the Hungarian Academy of Sciences (grant no. LP2026-5/2026.). T.L. would like to thank the organisers 
of the conferences ``Representations, Moduli and Duality'' in 
Lausanne and ``Corti60: A Tour through Algebraic
Geometry'' in Cortona for a stimulating environment in which some of the 
final technical breakthroughs in this project were made. T.L. 
would also like to thank the Max-Planck Institute for Mathematics
in Bonn where a significant revision of the first draft was carried
out.

\section{Preliminaries}
\label{sec:prelim}

Throughout the paper $\kk$ denotes an algebraically closed field of
characteristic $0$.

\subsection{Partitions}
\label{section-partitions}

\begin{Definition}
Let $n \in \mathbb{Z}_{\geq 0}$. An \em (unordered) partition of $n$ \rm 
is an unordered collection 
$\underline{n}\coloneqq \left\{ n_1, \dots, n_m \right\}$ of strictly
positive integers $n_i$ with $\sum_{i=1}^m n_i = n$. We denote this by 
$\underline{n} \vdash n$. 
\end{Definition}

The integers $n_i$ are the \em parts \rm of $\underline{n}$. 
We write $r(\underline{n})$ for the \em length \rm of $\underline{n}$. 
By this we mean the total number of parts in $\underline{n}$, i.e. $m$. 
We write $r_k(\underline{n})$ for
the number of parts of size $k$ in $\underline{n}$, i.e. the number of $i$ such 
that $n_i = k$. Finally, for any $n \in \mathbb{Z}_{\geq 0}$
we write $\partn_n$ to be the set of all partitions of $n$ and
we set $\partn \coloneqq \coprod_{n \in \mathbb{Z}_{\geq 0}} \partn_n$. 

In this paper, we also need the following notion: 
\begin{Definition}
Let $n,m \in \mathbb{Z}_{\geq 0}$. An \em ordered $m$-partition 
$(n_1, n_2, \dots, n_m)$ of $n$ \rm is an ordered $m$-tuple 
of non-negative integers $n_i$ with $\sum_{i=1}^m n_i = n$. 
We denote this by $(n_1, \dots, n_m) \vdash n$. 
\end{Definition}
Note that in these ordered partitions of fixed length we  
allow some parts to be of zero size. 

We write $m\text{-}\ordpartn_{n}$ for the set of all ordered $m$-partitions
of $n$. This set can be viewed as a subset of $\mathbb{Z}^m$ but 
note that there it is neither closed under addition nor multiplication. 
We further set $\ordpartn_{n} \coloneqq \coprod_{m \in \mathbb{Z}_{\geq 0}}
m\text{-}\ordpartn_{n}$ and 
$\ordpartn \coloneqq \coprod_{n \in \mathbb{Z}_{\geq 0} \ordpartn_{n}}$. 

For each $n,m \in \mathbb{Z}_{\geq 0}$ we have a natural forgetful map 
of sets
\begin{equation}
\label{eqn-map-from-m-ordered-to-unordered-partitions-of-n} 
F\colon m\text{-}\ordpartn_n \rightarrow \partn_n 
\end{equation}
which takes an ordered $m$-partition, discards its parts of zero
size, and forgets the ordering on the remaining parts. Its image
in $\partn_n$ are the length $\leq m$ partitions of $n$. Maps 
\eqref{eqn-map-from-m-ordered-to-unordered-partitions-of-n} 
combine into 
\begin{equation}
\label{eqn-map-from-ordered-to-unordered-partitions} 
F\colon \ordpartn \rightarrow \partn.  
\end{equation}

Let $\underline{n} \vdash n$. Then $\bigl(r_1(\underline{n}), \dots, r_n(\underline{n})\bigr)$ is an ordered $n$-partition of $r(\underline{n})$, the length of $\underline{n}$. We
denote by $\underline{r}(\underline{n})$ the resulting unordered partition
$F\bigl(r_1(\underline{n}), \dots, r_n(\underline{n})\bigr)$.  

\begin{Proposition}
Let $n,m \in \mathbb{Z}_{\geq 0}$ and let $\underline{n} \vdash n$ be
a partition of $n$. Let $F^{-1}(\underline{n})$ be the pre-image 
of $\underline{n}$ in $m\text{-}\ordpartn_n$. Then 
$$ |F^{-1}(\underline{n})| = \binom{m}{\underline{r}(\underline{n})}. $$
\end{Proposition}
Note that, by definition, the generalised binomial coefficient
$\binom{z}{k_1, \dots, k_m}$ only depends on the image of $(k_1,
\dots, k_m)$ in $\partn$, and thus we can evaluate it on unordered 
partitions. 
\begin{proof}
A pre-image of $\underline{n}$ is obtained by distributing the parts
of $\underline{n}$ between $m$ ordered positions and filling
the rest with zeroes. As the parts of same size in 
$\underline{n}$ are indistinguishable, the number of such
distributions is the number of ways to choose $r_1(\underline{n})$, 
\dots, $r_n(\underline{n})$ positions out of $m$ available. 
\end{proof}
\begin{Corollary}
\label{cor-ordered-to-undordered-reindexing-in-the-summation}
Let $A$ be an abelian group and let $f\colon \partn_n \rightarrow A$ 
be any maps of sets. Then
$$
\sum_{(n_1, \dots, n_m) \vdash n} f(F(n_1, \dots, n_m))
= 
\sum_{\underline{n} \vdash n} \binom{m}{\underline{r}(\underline{n})}
f(\underline{n}). 
$$
\end{Corollary}

\subsection{Heisenberg algebra}
\label{sec:heisenberg_existing}

\subsubsection{A-generator definition}

This is the original definition in \cite[\S8.1]{nakajima1999lectures}
\cite[\S3]{EllingsrudGottsche-HilbertSchemesOfPointsAndHeisenbergAlgebras}
 systematising the results of \cite{nakajima1997heisenberg}. It takes
the notion of the $\infty$-dimensional Heisenberg Lie algebra
\cite[\S9.13]{kac1990infinite} and extends it to have the generators
parametrised by elements of a vector space with a symmetric bilinear
form. 

We prefer to work with usual algebras, rather than Lie algebras.
Similar to  \cite{cautis2012heisenberg} and \cite{krug2018symmetric},
by \em the Heiseinberg algebra \rm of a lattice or a vector space we
mean the universal enveloping algebra of the corresponding Lie algebra
where we identified the central charge with $1$. This yields:
 
 \begin{Definition}
\label{defn-the-heisenberg-algebra-of-a-vector-space-a-gen-sym}
Let $(V,\chi)$ be a vector space with a symmetric bilinear form. 
The \em Heisenberg algebra \rm
$\chalga{V,\chi}$ is
the unital $\kk$-algebra with generators $a_v(n)$ for $v \in V$ and
integers $n \in \mathbb{Z} \setminus \{0\}$ modulo the following relations 
for all $v,w \in V$, $z \in \kk$ and $n,m \in \mathbb{Z} \setminus \{0\}$:
\begin{equation}\label{eq:heisvectrel1-add-a-gen-sym}
	a_{v+w}(n) = a_v(n) + a_w(n),
\end{equation}
\begin{equation}\label{eq:heisvectrel2-mult-a-gen-sym}
	a_{zv}(n) = z a_v(n), 
\end{equation}
\begin{equation}\label{eq:heisvectrel3-a-gen-sym}
	[a_v(n), a_w(m)] = \delta_{m,-n} m \langle v,w \rangle_{\chi},
\end{equation}
where $[-,-]$ denotes the commutator.
\end{Definition}

The relations \eqref{eq:heisvectrel1-add-a-gen-sym}
and \eqref{eq:heisvectrel2-mult-a-gen-sym} are the relations of
linearity in $v \in V$. The relation \eqref{eq:heisvectrel3-a-gen-sym}
is the \em Heisenberg relation\rm. The 
$\infty$-dimensional Heisenberg algebra (the case $V = \kk$ with
$\langle1,1\rangle_\chi = 1$) is isomorphic to the algebra of
differential operators of the polynomial ring $\kk[x_1, x_2, \dots]$
via
\begin{equation}
	a(n) \mapsto 
	\begin{cases}
		n x_n, \quad\quad n > 0, \\
		\frac{\partial}{\partial x_n}, \quad\quad n < 0,
	\end{cases}
\end{equation}  
The Heisenberg relation then corresponds to the identity 
\begin{equation}
	\frac{\partial}{\partial x_n}\left(nx_n f \right) = 
	nx_n \frac{\partial}{\partial x_n}\left(f \right) + 
	nf.
\end{equation}

For a lattice, we use the same definition but without the scalar multiplication relation \eqref{eq:heisvectrel2-mult-a-gen-sym}:
\begin{Definition}
\label{defn-the-heisenberg-algebra-of-a-lattice-gen-sym}
Let $(M,\chi)$ be a lattice with a symmetric bilinear form. 
The \em Heisenberg algebra \rm
$\chalga{M,\chi}$ 
is the unital $\kk$-algebra with generators $a_v(n)$ for $v \in M$ and 
$n \in \mathbb{Z} \setminus \{0\}$ modulo the following relations 
for all $v,w \in M$, $z \in \kk$ and $n,m \in \mathbb{Z} \setminus \{0\}$:
\begin{equation}\label{eq:heisrel2-add-a-gen-sym}
	a_{v+w}(n) = a_v(n) + a_w(n),
\end{equation}
\begin{equation}\label{eq:heisrel3-a-gen-sym}
	[a_v(n), a_w(m)] = \delta_{m,-n} m \langle v,w \rangle_{\chi}.
\end{equation}
\end{Definition}

As the Heisenberg relation is bilinear in $v,w \in V$ we immediately have 
the \em basis reduction \rm result: 

\begin{Proposition}
\label{theorem-basis-reduction-for-heisenberg-algebra-vector-space-a-gen-sym}
Let $(V,\chi)$ be a vector space (or a lattice) with a
symmetric bilinear form and let $e_1, \dots, e_l$ be a basis of $V$.
The Heisenberg algebra $\chalga{V,\chi}$ is isomorphic to the unital
$\kk$-algebra with generators $a_v(n)$ for all $v \in \left\{ e_1,
\dots, e_l \right\}$ and $n \in \mathbb{Z}\setminus\{0\}$ modulo 
the relations \eqref{eq:heisvectrel1-add-a-gen-sym}
and \eqref{eq:heisvectrel3-a-gen-sym} for all $v,w\in \left\{ e_1, \dots, e_l \right\}$ and $n,m\in \mathbb{Z}$.
\end{Proposition}
\begin{proof}
We give the proof for $V$ being a vector space, the proof for the lattice is similar. The relations  \eqref{eq:heisvectrel1-add-a-gen-sym} and
\eqref{eq:heisvectrel2-mult-a-gen-sym} 
are equivalent to the following basis reduction relation 
\begin{equation}
\label{eq:heisvectrel2-basis-red-a-gen-sym}
a_v(n)  = v_1 a_{e_1}(n) + \dots + v_l a_{e_l}(n) 
\quad \quad \quad \quad 
\forall\; v \in V, n \in \mathbb{Z} \setminus \{0\},
\end{equation}  where $v_i \in \kk$ are the unique coefficients such that $v = v_1 e_1 + \dots + v_l e_l$.

We can thus replace relations \eqref{eq:heisvectrel1-add-a-gen-sym}
and \eqref{eq:heisvectrel2-mult-a-gen-sym} by the relation
\eqref{eq:heisvectrel2-basis-red-a-gen-sym}. Next, as the Heisenberg
relation \eqref{eq:heisvectrel3-a-gen-sym} is bilinear in $v,w \in V$,
it is sufficient to only impose it for all $v,w\in \left\{ e_1, \dots,
e_l \right\}$.  Now, for each $v \in V$ and $n \in
\mathbb{Z}\setminus\{0\}$, the element $a_v(n)$ occurs in precisely 
one relation \eqref{eq:heisvectrel2-basis-red-a-gen-sym}. Hence we can
only take the generators $a_v(n)$ for $v \in \left\{ e_1, \dots, e_l
\right\}$ and the relations 
\eqref{eq:heisvectrel2-basis-red-a-gen-sym}, and 
\eqref{eq:heisvectrel3-a-gen-sym}. Finally, 
the relation \eqref{eq:heisvectrel2-basis-red-a-gen-sym} is tautological for $v \in \left\{ e_1, \dots, e_l \right\}$, so we can get rid of it entirely. 
\end{proof}

\section{Heisenberg algebra of a graded vector space}

In this section, for a graded vector space $V$ with a bilinear form $\chi$
we give a minor generalisation of the standard ($A$-generator) 
definition of the Heisenberg algebra $\chalg{V,\chi}$ which works with 
for all bilinear forms $\chi$. In this generality, we spell out some
well-known structures and results: the bigrading, basis reduction, 
the canonical action on the symmetric algebra $S^*(V \otimes t\kk[t])$. 
We also prove that for non-degenerate $\chi$ our definition is independent 
of the choice of $\chi$. 

\subsection{Definition}
\label{section-heisenberg-algebra-definition}

Throughout this section, a \em graded vector space \rm means a
vector space $V$ with $\mathbb{Z}$- or $\mathbb{Z}_2$-grading. Our
Heisenberg algebra definition only uses the parity of the degree, 
so as an algebra the result only depends on the $\mathbb{Z}_2$-grading. 
However, the grading on $V$, whether by $\mathbb{Z}$ or by
$\mathbb{Z}_2$, is incorporated into the bigrading we consider on the resulting Heisenberg algebra,  cf.~\S\ref{section-heisenberg-algebra-bigrading}. 

By a \em bilinear form \rm $\chi$ on a graded vector space $V$ we mean 
a map $V \otimes_\kk V \rightarrow \kk$ of \em graded \rm vector
spaces. This implies that for homogeneous elements 
$\left<v,w\right>_\chi \neq 0$ if and only $\deg(v) + \deg(w) = 0$. 
Example \ref{example-why-pairing-has-to-be-Z2-graded} below shows that 
if an odd degree $v$ were to pair non-trivially with some even degree
$w$, the Heisenberg relation would imply that $a_v(n) = 0$ for 
all $n \neq 0$. 

\begin{Definition}
\label{defn-the-heisenberg-algebra-of-a-graded-vector-space-a-gen}
Let $(V,\chi)$ be a graded vector space $V$ with a bilinear form
$\chi$. The \em Heisenberg algebra \rm $\chalg{V,\chi}$
is the unital graded $\kk$-algebra with generators $a_v(n)$ 
for all $v \in V$ and $n \in \mathbb{Z} \setminus \{0\}$ modulo the relations:
\begin{equation}\label{eq:vectheisrel1-a-gen-graded}
	a_v(n)a_w(m) = (-1)^{\deg(w)\deg(v)} a_w(m) a_v(n)
\quad \quad v,w \in V \text{ and }  m,n \in \mathbb{Z}_{> 0} \text{ or
} m,n \in \mathbb{Z}_{< 0}
\end{equation}
\begin{equation}\label{eq:vectheisrel2add-a-gen-graded}
	a_{v+w}(n) = a_v(n) + a_w(n)
\quad \quad \quad v,w \in V \text{ and } n \in \mathbb{Z} \setminus \{0\}, 
\end{equation}
\begin{equation}\label{eq:vectheisrel2scalmult-a-gen-graded}
a_{zv}(n) = z a_v(n)
\quad \quad \quad v \in V, z \in \kk, n \in \mathbb{Z} \setminus \{0\}, 
\end{equation}
\begin{equation}\label{eq:vectheisrel3-a-gen-graded}
	a_v(-n) a_w(m)  = (-1)^{\deg(w)\deg(v)} a_w(m) a_v(-n) + \delta_{n,m} 
	m \langle v, w \rangle_{\chi}
\quad v,w \in V \text{ and } n,m \in \mathbb{Z}_{> 0}.
\end{equation}
\end{Definition}

Denote by $V^{odd}$ the sum of all odd degree graded parts of $V$ and
by $V^{even}$ the sum of all even degree ones. We use similar notation
for any graded vector space.
 
The even part $\chalg{V^{even}}$ is the Heisenberg algebra 
in classical (non-graded) sense, while the odd part $\chalg{V^{odd}}$ 
is a Clifford algebra \cite[\S8.1]{nakajima1999lectures}. 
Since $\chi$ is a graded pairing, 
$\chalg{V^{even}}$ and $\chalg{V^{odd}}$ commute. 
The following example shows why $\chi$ has to be at least
$\mathbb{Z}_2$-graded:
\begin{Example}
\label{example-why-pairing-has-to-be-Z2-graded}
Let $V$ be a graded vector space and $\chi$ be a bilinear form on
$V$, not necessarily graded. 
Let $v \in  \chalg{V^{even}}$, $w \in \chalg{V^{odd}}$ and $n \in
\mathbb{Z} \setminus \{0\}$. We have 
\begin{align*}
0 &= a_v(-n) a_w(n) a_w(n) = 
\left( a_w(n) a_v(-n) + n \left<v,w\right>\right) a_w(n)  
= a_w(n) a_v(-n) a_w(n) + n \left<v,w\right> a_w(n) = \\
&= a_w(n)\left( a_w(n) a_v(-n) + n\left<v,w\right> \right) + 
n \left<v,w\right> a_w(n) = 
a_w(n) a_w(n) a_v(-n) + 2n\left<v,w\right> a_w(n) = \\
&= 2n\left<v,w\right> a_w(n). 
\end{align*}
A similar computation with $a_w(-n) a_w(-n) a_v(n)$ shows that 
$2n \left<w,v\right> a_w(n) = 0$. If $\chi$ is graded, 
$\left<w,v\right> = \left<v,w\right> = 0$ 
for any even $v$ and odd $w$. Were $\chi$ not to be graded, 
we must have $a_w(n) = 0$ 
for any odd $w$ which pairs non-trivially, on either side, 
with an even $v$. Thus 
for non-graded $\chi$ the Heisenberg algebra $\chalg{V,\chi}$ 
coincides with $\chalg{V^{even} \oplus K,\chi}$ where $K \subset
V^{odd}$ consists of all elements which pair trivially on both sides
with $V^{even}$. Note that on $V^{even} \oplus K$ the pairing
$\chi$ is $\mathbb{Z}_2$-graded. 
\end{Example}

Since the relations in 
Definition \ref{defn-the-heisenberg-algebra-of-a-graded-vector-space-a-gen}
are still linear, the basis reduction is still immmediate: 

\begin{Proposition}
\label{theorem-basis-reduction-for-heisenberg-algebra-graded-vector-space-pq-gen}
Let $(V,\chi)$ be a graded vector space with bilinear form 
and $\left\{ e_1, \dots, e_l \right\}$ be a homogeneous basis. 
$\chalg{V,\chi}$ is isomorphic to the unital $\kk$-algebra with generators 
$a_{e_i}(n)$ for $n \in \mathbb{Z} \setminus \{0\}$ and the relations 
\eqref{eq:vectheisrel1-a-gen-graded} and  
\eqref{eq:vectheisrel3-a-gen-graded} for $v,w\in \left\{ e_1, \dots,
e_l \right\}$ and $n,m\in \mathbb{Z}\setminus\{0\}$.
\end{Proposition}
\begin{proof}
Same as the proof of Proposition~\ref{theorem-basis-reduction-for-heisenberg-algebra-vector-space-a-gen-sym}. 
\end{proof}

\subsection{Bigrading}
\label{section-heisenberg-algebra-bigrading}

We equip the Heisenberg algebra $\chalg{V,\chi}$ of $(V,\chi)$ with
the following bigrading: 

\begin{Definition}
Let $(V,\chi)$ be a graded vector space with a bilinear form. Define a bigrading on the Heisenberg algebra $\chalg{V,\chi}$ as follows: 
\begin{enumerate}
	\item The \em base \rm grading: for every homogeneous $v \in V$ and 
$n \neq 0$ the degree of the generator $a_{v}(n)$ is $\deg(v)$. 
This grading is of the same type ($\mathbb{Z}$ or $\mathbb{Z}/2$) as 
that on the vector space $V$.
	\item The \em weight \rm grading: the $\mathbb{Z}$-grading 
where for every $v \in V$ and $n \neq 0$ the degree of the generator 
$a_{v}(n)$ is $n$.
\end{enumerate}
\end{Definition}

The relations 
\eqref{eq:vectheisrel1-a-gen-graded}-\eqref{eq:vectheisrel3-a-gen-graded}
in 
Definition \ref{defn-the-heisenberg-algebra-of-a-graded-vector-space-a-gen}
respect both of these gradings, thus we get a well-defined bigrading on the
whole of $\chalg{V,\chi}$. The weight grading agrees with the grading in 
the definition of the
Heisenberg super Lie algebra in \cite[\S8.1]{nakajima1999lectures}. 
Our terminology also agrees with
\cite{EllingsrudGottsche-HilbertSchemesOfPointsAndHeisenbergAlgebras} in 
the following sense: for the Heisenberg algebra of $V$, 
Ellingsrud and G{\"o}ttsche also used the term \em weight \rm 
for the second grading, but used the term 
\em parity \rm for the first grading since they only considered
$\mathbb{Z/2}$-grading on $V$. Our point is that the first grading is
induced from the base vector space $V$, hence we call it 
the \em base grading\rm. 

\subsection{Fock space representation and its symmetric algebra structure}
\label{section-fock-space-and-its-symmetric-algebra-structure}

Heisenberg algebra has the following well-known irreducible representation:

\begin{Definition}
Let $(V,\chi)$ be a graded vector space with a bilinear form. Let
$I_{-} \subset \chalg{V,\chi}$ be the left ideal generated by
$a_v(-n)$ for all $v \in V$ and $n > 0$. The \em Fock space
representation \rm $\falg{V,\chi}$ is the bigraded (left)
$\chalg{V,\chi}$-module $\chalg{V,\chi}/I_{-}$. 
\end{Definition}

By the irreducibility of $\falg{V,\chi}$, for any
representation $W$ of $\chalg{V,\chi}$ and for any $w \in W$ which is
annihilated by $I_{-}$ we obtain a unique injective map 
$\falg{V,\chi} \hookrightarrow W$ sending $1 \in \falg{V,\chi}$ 
to $w \in W$. 

Let $H_+ \subset \chalg{V,\chi}$ be the subalgebra generated by
$a_v(n)$ with $n \geq 0$. It follows from the Heisenberg relation
\eqref{eq:vectheisrel3-a-gen-graded}
that $H_+$ splits the projection to the Fock space. In other words, 
the composition 
\begin{equation}
\label{eqn-identification-of-hplus-and-fock-space}
H_+ \hookrightarrow \chalg{V,\chi} \twoheadrightarrow \falg{V,\chi}
\end{equation}
is an isomorphism. Since Heisenberg relation is empty in absence
of $a_v(n)$ with $n < 0$, $H_+$ is the $\kk$-algebra generated by 
the set $\left\{ a_v(n) \right\}_{v \in V, n > 0}$ subject to
the linearity relations
\eqref{eq:vectheisrel2add-a-gen-graded},
\eqref{eq:vectheisrel2scalmult-a-gen-graded} and the commutation
relation \eqref{eq:vectheisrel3-a-gen-graded}. The linearity relations
give the set $\left\{ a_v(n) \right\}_{v \in V, n > 0}$ the structure
of a bigraded vector space. We identify it with 
the bigraded vector space $V \otimes_\kk t\kk[t]$ via an isomorphism 
\begin{equation}
\label{eqn-isomorphism-sym-algebra-hplus-linear-part}
V \otimes_\kk t\kk[t] 
\xrightarrow{\sim}
\left\{ a_v(n) \right\}_{v \in V, n > 0}
\quad \quad 
v \otimes t^n \mapsto a_v(n). 
\end{equation}
The bigrading on $V \otimes_\kk t\kk[t]$ is defined by setting
$\deg(v \otimes t^n) = (\deg(v), n)$, and as before we use 
the terms \em base grading \rm and \em weight grading \rm 
to refer to these two gradings. 

By the commutation relation \eqref{eq:vectheisrel3-a-gen-graded}
and $\kk$-linearity, it follows that the isomorphism 
\eqref{eqn-isomorphism-sym-algebra-hplus-linear-part}
induces an isomorphism of bigraded algebras
\begin{equation}
\label{eqn-isomorphism-sym-algebra-hplus}
S^* \left(V \otimes_\kk t\kk[t]\right)
\xrightarrow{\sim}
H_+. 
\end{equation}
Note that the skew-commutation sign twist 
in relation \eqref{eq:vectheisrel3-a-gen-graded} only involves the
base grading. Thus the symmetric algebra 
$S^* \left(V \otimes_\kk t\kk[t]\right)$ is the symmetric algebra
of $V \otimes_\kk t\kk[t]$ considered as a single-graded vector space
with respect to the base grading. 

Composing \eqref{eqn-isomorphism-sym-algebra-hplus} with 
\eqref{eqn-identification-of-hplus-and-fock-space}
we obtain the isomorphism of bigraded vector spaces
\begin{equation}
\label{eqn-isomorphism-sym-algebra-fock-space}
S^* \left(V \otimes_\kk t\kk[t]\right)
\xrightarrow{\sim}
\falg{V,\chi}. 
\end{equation}
It gives the Fock space $\falg{V,\chi}$ the structure of a symmetric
algebra and defines the action of $\chalg{V,\chi}$ on 
$S^* \left(V \otimes_\kk t\kk[t]\right)$. This construction is
well-known, see e.g.
\cite[\S3]{EllingsrudGottsche-HilbertSchemesOfPointsAndHeisenbergAlgebras}, 
\cite[\S4.2]{Lehn-LecturesOnHilbertSchemes}, 
\cite[\S3.7]{BelmansFuKrug-HochschildCohomologyOfHilbertSchemesOfPointsOnSurfaces}.

We now spell out the resulting action of 
$\chalg{V,\chi}$ on $S^* \left(V \otimes_\kk t\kk[t]\right)$. 
First, we describe  the structure of $S^* \left(V \otimes_\kk t\kk[t]\right)$
and introduce notation for its elements. 
For any $n > 0$ and any $\underline{n} \vdash n$, write
$$ t^{\underline{n}} = t^{n_1}t^{n_2}\dots{t}^{n_{r(\underline{n})}}. $$
Note that $t^i t^j \neq t^{i+j}$ in $S^* \left(V \otimes_\kk t\kk[t]\right)$ 
as $t\kk[t]$ is only a graded vector space in this context.
For any ordered partition $(n_1, n_2, \dots, n_m) \vdash n$ let 
$$ S_{n_1, n_2, \dots, n_m} := S_{n_1} \times S_{n_2} \times \dots
\times S_{n_m} < S_n $$
$$ \Sym^{n_1, n_2, \dots, n_m} V := V^{n} / S_{(n_1, n_2, \dots, n_m)}, $$
For an unordered partition $\underline{n} \vdash n$ write
$S_{\underline{n}}$ and $\Sym^{\underline{n}} V$ for the (isomorphism
classes of) $S_{(n_1, n_2, \dots, n_m)}$ and $\Sym^{n_1, n_2, \dots,
n_m} V$ where $(n_1, n_2, \dots, n_m)$ is any pre-image of
$\underline{n}$ under the forgetful map. 

Recall that for any $\underline{n} \vdash n$ we have 
ordered partition $\underline{r}(\underline{n}) := (r_1(\underline{n}),
\dots, r_n(\underline{n})) \vdash r(\underline{n})$ where
$r_i(\underline{n})$ is the number of parts of size $i$ in
$\underline{n}$ and $r(\underline{n}) = \sum r_i(\underline{n})$. 
Thus for any $\underline{n} \vdash n$ we have 
$$ \Sym^{\underline{r}(\underline{n})} V
\simeq 
V^{r(\underline{n})} / 
S_{\underline{r}(\underline{n})}
\simeq 
\Sym^{r_1(\underline{n})} V \otimes \dots 
\otimes \Sym^{r_n(\underline{n})} V, $$
where we view each factor $V$ of $V^{r(\underline{n})}$ as
corresponding to a part of $\underline{n}$
and $S_{\underline{r}(\underline{n})}$ acts 
by permuting the factors corresponding
to the parts of the same size.
For any $\underline{v} \in \Sym^{\underline{r}(\underline{n})} V$ we write 
$$ \underline{v} := \underline{v}_1 \otimes \dots \otimes \underline{v}_n
\quad \quad \underline{v}_i \in \Sym^{r_i(\underline{n})} V,
$$
and we further write
$$ \underline{v}_i = v_{i1} \vee \dots \vee v_{ir_i(\underline{n})} \quad \quad v_{ij} \in V. $$
Here we distinguish between the skew-commutative product $\vee$ within
each $\Sym^{r_i(n)} V$ and ordinary tensor product $\otimes$ in 
$\Sym^{r_1(\underline{n})} V \otimes \dots \otimes 
\Sym^{r_n(\underline{n})} V$. 

For example, when $n = 12$ and $\underline{n} = (1 1 1 1 2 3 3)$, we
have
$$ t^{\underline{n}} = t^1 t^1 t^1 t^1 t^2 t^3 t^3, $$
$$ \Sym^{\underline{r}(\underline{n})} V = \Sym^4 V \otimes V \otimes \Sym^2 V, $$
and any $\underline{v} \in \Sym^{\underline{r}(\underline{n})} V$ 
can be written as
$$ \underline{v} = (v_{11} \vee v_{12} \vee v_{13} \vee v_{14})
\otimes (v_{21}) \otimes (v_{31} \vee v_{32}). $$

With this notation in mind:
\begin{Lemma}
\label{lemma-direct-sum-description-of-tkt-sym-algebra-of-v}
We have an isomorphism of bigraded vector spaces
\begin{equation}
\label{eqn-direct-sum-description-of-tkt-sym-algebra-of-v}
S^* \left(V \otimes_\kk t\kk[t]\right) \simeq 
\bigoplus_{n \geq 0} \bigoplus_{\underline{n} \vdash n}
\Sym^{\underline{r}(\underline{n})}V\; t^{\underline{n}}
\end{equation}
which sends any 
$$ v_1 t^{n_1} \vee v_2 t^{n_2} \vee \dots \vee v_{m} t^{n_m} \in
S^*(E \otimes t \kk[t]) $$
with $n_1 \leq n_2 \leq \dots \leq n_m$ to 
$$ (v_1 \vee \dots \vee v_{m_1}) \otimes
(v_{m_1+1} \vee \dots \vee v_{m_1+m_2}) \otimes \dots \otimes
(v_{m_1 + \dots + m_{m-1} + 1} \vee \dots \vee v_m) t^{\underline{\lambda}} $$
where $n = \sum n_i$, $\underline{\lambda} = \left\{n_1, \dots, n_m\right\}$
and $m_i = r_i(\underline{\lambda})$. 
Its inverse sends any 
$$ \underline{v}_1 \otimes \dots \otimes \underline{v}_n
t^{\underline{\lambda}} \in \left(\Sym^{\underline{r}(\underline{\lambda})}
E\right) t^{\underline{\lambda}},
$$
to 
$$ (v_{11} t^1) \vee (v_{12} t^1) \vee \dots \vee (v_{n r_n(\underline{\lambda})} t^n). $$
\end{Lemma}
\begin{proof}
It is straightforward to check that the two maps given above are
bigraded and mutually inverse. For a more conceptual explanation, 
it is well known \cite[\S{XVI}, Proposition~8.2]{Lang-Algebra} that given 
any two graded vector spaces $V$ and $W$ we have a bigraded isomorphism 
$$ S^*(V \oplus W) \simeq S^*(V) \otimes S^*(W) $$
given by 
$$ v_1 \vee \dots \vee v_n \vee w_1 \vee \dots \vee w_m 
\mapsto 
(v_1 \vee \dots \vee v_n) \otimes (w_1 \vee \dots \vee w_m). $$
Applying this iteratively we get 
$$ S^* \left(V \otimes_\kk t\kk[t]\right) \simeq 
S^*(V t^1) \otimes S^*(V t^2) \otimes \dots. $$
For each $n > 0$, the $n$-th graded part of the RHS with respect to 
the weight grading is
$$
\bigoplus_{\underline{n} \vdash n}
S^{r_1(\underline{n})}(V t^1) \otimes  
S^{r_2(\underline{n})}(V t^2) \otimes \dots 
\otimes S^{r_n(\underline{n})}(V t^n)
\simeq  
\bigoplus_{\underline{n} \vdash n}
\Sym^{\underline{r}(\underline{n})}V\; t^{\underline{n}}. 
$$
Composing and simplifying the isomorphisms above gives
the two mutually inverse maps in the assertion of this Lemma. 
\end{proof}

Note that multiplying elements in this notation introduces sign 
twist due to skew-commutativity:
$$ 
\bigl(x\; t^{2}\bigr)
\bigl((v \vee w) \otimes u \; t^{1 1 2}\bigr) 
= (-1)^{x|v,w} (v \vee w) \otimes (x \vee u) \; t^{1 1 2 2}. $$
Here $(-1)^{x|v,w} := (-1)^{\deg(x)(\deg(v) + \deg(w))}$
is the sign-twist for commuting $x$ through $v$ and
$w$. The difference between $\vee$ and $\otimes$ remains important:
$v \vee w\; t^{1 1} = (-1)^{v|w} w \vee v\; t^{1 1}$, 
but $v \otimes w\; t^{2 1}$ and $w \otimes v\; t^{2 1}$ are linearly
independent. 

Finally, for any $w \in V$ and any $n \geq 0$ define the map 
$\left<w,-\right>\colon \Sym^{n+1} V  \;\rightarrow\; \Sym^{n} V$
by
\begin{align}
\label{eqn-pairing-bracket-on-symmetric-powers}
\left<w,-\right>\colon v_1 \vee \dots \vee v_{n+1}  \;\mapsto\;
\sum_{i = 1}^{i+1} (-1)^{w|v_1, \dots, v_{i-1}} \left<w,v_i\right>
v_1 \vee \dots \vee v_{i-1} \vee v_{i+1} \vee \dots v_{n+1}.
\end{align}
\begin{Lemma}
\label{lemma-canonical-action-of-chalg-v-chi-on-symmetric-algebra}
Let $(V,\chi)$ be a graded vector space with a bilinear form. 
The isomorphism \eqref{eqn-isomorphism-sym-algebra-fock-space}
induces the bigraded action of $\chalg{V,\chi}$ on 
$S^* \left(V \otimes_\kk t\kk[t]\right)$ where for every $w \in V$
and every $i > 0$:
\begin{itemize}
\item $a_w(i)$ acts by the sum of the maps
$\Sym^{\underline{r}(\underline{n})}V\; t^{\underline{n}} 
\longrightarrow 
\Sym^{\underline{r}(\underline{n} + (i))}V\; t^{\underline{n} + (i)}$ sending
\begin{equation}
\label{eqn-action-of-a_w(n)-on-symmetric-algebra}
\underline{v}_1 \otimes \dots \otimes \underline{v}_i \otimes \dots 
\otimes \underline{v}_n \; t^{\underline{n}}
\mapsto 
(-1)^{w \;|\; \underline{v}_{1}, \dots, \underline{v}_{i-1}}\;
\underline{v}_1 \otimes \dots \otimes (w \vee \underline{v}_i ) \otimes \dots 
\otimes \underline{v}_n \; t^{\underline{n}+(i)}
\end{equation}
\item $a_w(-i)$ acts by the sum of the maps
$\Sym^{\underline{r}\left(\underline{n} + (i)\right)}V\; t^{\underline{n} + (i)},
 \longrightarrow 
\Sym^{\underline{r}\left(\underline{n}\right)}V\; t^{\underline{n}}$ sending
\begin{equation}	
\label{eqn-action-of-a_w(-n)-on-symmetric-algebra}
\underline{v}_1 \otimes \dots \otimes \underline{v}_i \otimes \dots
\otimes \underline{v}_n\; t^{\underline{n}+(i)}
 \mapsto 
(-1)^{w \;|\; \underline{v}_{1}, \dots, \underline{v}_{i-1}}\;
\underline{v}_1 \otimes \dots \otimes 
\left<w, \underline{v}_i\right>\;
\otimes \dots \otimes \underline{v}_n\; t^{\underline{n}},
\end{equation}
\end{itemize}
\end{Lemma}
\begin{proof}
Let $w \in V$ and $n > 0$. We identified $S^* \left(V \otimes_\kk t\kk[t]\right)$ with 
the lift of the quotient $\falg{V,\chi}$ to a subalgebra $H_+$ of $\chalg{V,\chi}$.  To compute the action 
of $a_w(i)$ and $a_w(-i)$ on  any $x \in H_+$, we need to compute the lifts of the classes $[a_w(i) x], [a_w(-i) x] \in \falg{V,\chi}$ back to $H_+$. 

Since $a_w(i) \in H_+$, the lift of $[a_w(i) x]$ is $a_w(i) x$ itself. In other words, $a_w(i) x$ acts
on $H_+$ by left multiplication within $H_+$. Hence $a_w(i)$  acts on $S^* \left(V \otimes_\kk t\kk[t]\right)$
by the left multiplication by $w t^i$, the pre-image of $a_w(i)$ under isomorphism \eqref{eqn-isomorphism-sym-algebra-fock-space}. We thus have \eqref{eqn-action-of-a_w(n)-on-symmetric-algebra}. 

For the action of $a_w(-i)$,  we note that 
$$ \underline{v}_1 \otimes \dots \otimes \underline{v}_i \otimes \dots \otimes \underline{v}_n t^{\underline{n}+(i)} \in S^* \left(V \otimes_\kk t\kk[t]\right)$$
corresponds to 
\begin{small}
\begin{equation}
\label{eqn-the-term-of-H_+-corresponding-to-t^-part(m)+n}
\left( a_{v_{11}}(1) \dots a_{v_{1r_1(\underline{n})}}(1) \right) \dots 
\left( a_{v_{i1}}(i) \dots a_{v_{i(r_i(\underline{n})+1)}}(i) \right)
\dots
\left( a_{v_{n1}}(n) \dots a_{v_{nr_n(\underline{n})}}(n) \right)
\in H_+.
\end{equation}
\end{small}
Multiplying \eqref{eqn-the-term-of-H_+-corresponding-to-t^-part(m)+n}
on the left by $a_w(-i)$ lands us outside $H_+$. We then use the
Heisenberg relation repeatedly to commute $a_w(i)$ from left to write
past all the terms
in  \eqref{eqn-the-term-of-H_+-corresponding-to-t^-part(m)+n}. Commuting $a_w(-i)$ past any $a_{v_{jk}}(j)$ generates sign twist $(-1)^{w|v_{jk}}$. 
Additionally, when $j = i$ it also generates the additional term where
$a_w(-i) a_{v_{ij}}(i)$ is removed and replaced by $\left<w,v_{ij}\right> \in \kk$. These additional terms all lie 
in $H_+$ and the summands on the RHS of
\eqref{eqn-pairing-bracket-on-symmetric-powers} are the
corresponding elements of $S^* \left(V \otimes_\kk t\kk[t]\right)$.
Finally, after commuting $a_w(-i)$ all the way to the right of
\eqref{eqn-the-term-of-H_+-corresponding-to-t^-part(m)+n}, we obtain
something that lies in $I_{-}$, and thus vanishes in the quotient
$\falg{V,\chi}$. We thus have
\eqref{eqn-action-of-a_w(-n)-on-symmetric-algebra}.
\end{proof}

\subsection{Independence of $\chi$}

For non-degenerate $\chi$, we show that $\chalg{V, \chi}$ and 
its action on the symmetrical  
algebra $S^* \left(V \otimes_\kk t\kk[t]\right)$ are both 
independent of $\chi$.

\begin{Proposition}
\label{prop-heisenberg-algebra-twisting-the-pairing}
Let $(V,\chi)$ be a graded vector space with a graded bilinear form. 
Let $S$ and $T$ be any two graded automorphisms of $V$. Define
the form $S \chi T$ by $\left<v,w\right>_{S \chi T} := \left<Sv,
Tw\right>_{\chi}$. There is a bigraded isomorphism
$$ \chalg{V,\chi} \xrightarrow{\sim}  \chalg{V, S \chi T}. $$
\end{Proposition}
\begin{proof}
Define a map $\chalg{V,\chi} \rightarrow \chalg{V, S \chi T}$ 
on the generators by setting for all $v \in V$ and $n > 0$
\begin{equation}
\label{eqn-heisenberg-algebra-twisting-the-pairing}
a_v(-n) \mapsto a_{S^{-1}v}(-n)
\quad \quad \text{ and } \quad \quad 
a_v(n)\mapsto a_{T^{-1}v(n)}
\end{equation}
This respects the relations \eqref{eq:vectheisrel1-a-gen-graded}-
\eqref{eq:vectheisrel3-a-gen-graded}, 
and hence gives a well-defined map 
$\chalg{V,\chi} \rightarrow \chalg{V, S \chi T}$. As $S$ and $T$
are invertible, this map is a bijection on the sets of generators 
and hence an isomorphism. 
\end{proof}

\begin{Theorem}
\label{theorem-independence-of-chi-for-nondegenerate-chi}
Let $V$ be a graded vector space and $\chi$ and $\xi$ be 
two non-degenerate graded bilinear forms on $V$. Then 
there is a canonical bigraded isomorphism 
\begin{equation}
\rho\colon \chalg{V,\chi} \xrightarrow{\sim} \chalg{V,\xi}. 
\end{equation}
It is canonical in the following sense: it is the unique
bigraded isomorphism $\chalg{V,\chi} \rightarrow \chalg{V,\xi}$
which comes from a morphism of Lie algebras and 
which sends $a_v(n)$ to itself for any $v \in V$ and $n > 0$.  
\end{Theorem}
\begin{proof}
If $\chi$ and $\xi$ are both non-degenerate, there is a unique
graded automorphism $S$ of $V$ such that $\xi = S{\chi}$, i.e. 
$\left<-,-\right>_{\xi} = \left<S-,-\right>_{\chi}$. 
Define $\rho$ to be the bigraded isomorphism 
\eqref{eqn-heisenberg-algebra-twisting-the-pairing} constructed by 
Proposition~\ref{prop-heisenberg-algebra-twisting-the-pairing} with $S$ as
above and $T = \id$.  Since $T = \id$, the isomorphism $\rho$ sends 
$a_v(n)$ to itself
for any $v \in V$ and $n >0$. 

For the uniqueness, it suffices to show that when $\xi = \chi$ 
any map $\rho$ satisfying those assumptions must be the identity map. 
By the linearity relations
\eqref{eq:vectheisrel2add-a-gen-graded}-\eqref{eq:vectheisrel2scalmult-a-gen-graded} and since $\rho$ preserves
the bigrading and comes from a Lie algebra morphism, 
we have 
$$ \rho(a_{v}(n)) = a_{\rho_n(v)}(n)  \quad \quad \quad \forall\; v \in V $$
for some collection of graded linear maps $\rho_n\colon V \rightarrow V$
for all $n \neq 0$. 

Let $n > 0$. Our assumptions on $\rho$ imply that $\rho_n = \id_V$. 
By the Heisenberg relation \eqref{eq:vectheisrel3-a-gen-graded} and 
since $\rho(1) = 1$, we have for all $v,w \in V$
$$ n\left<v,w\right>_\chi = \rho\left(n \left<v,w\right>_\chi\right) = 
  \rho\left([a_{v}(-n), a_w(n)]\right) = 
  [a_{\rho_{-n} v}(-n),a_w(n)] = n \left<\rho_{-n}(v),w\right>_{\chi}. $$ 
By the nondegeneracy of $\chi$, it follows that $\rho_{-n} = \id_V$. 

We conclude that $\rho_n = \id_V$ for all $n \neq 0$, whence 
$\rho = \id_{\chalg{V,\chi}}$. 
\end{proof}

Let $V$ be a graded vector space and $\chi$ and $\xi$ be two 
non-degenerate graded bilinear forms on $V$.  The canonical
isomorphism $\rho\colon \chalg{V,\chi} \xrightarrow{\sim} \chalg{V,\chi}$ 
of Theorem \ref{theorem-independence-of-chi-for-nondegenerate-chi} sends 
the ideal $I_{-} \subset \chalg{V,\chi}$ to the ideal  $I_{-} \subset
\chalg{V,\xi}$. It also sends the subalgebra 
$H_+ \subset \chalg{V,\chi}$ to the subalgebra $H_+ \subset
\chalg{V,\xi}$. We thus have a commutative diagram of isomorphisms:
\begin{equation}
\label{eqn-fock-space-and-hplus-diagram-for-two-forms}
\begin{tikzcd}
H_+
\ar[hook]{r}
\ar{d}{\rho}[']{\sim}
&
\halg{V,\chi} 
\ar[two heads]{r}
\ar{d}{\rho}[']{\sim}
&
\falg{V,\chi} 
\ar{d}{\rho}[']{\sim}
\\
H_+
\ar[hook]{r}
&
\halg{V,\xi} 
\ar[two heads]{r}
&
\falg{V,\xi}.
\end{tikzcd}
\end{equation}

The isomorphism $\rho$ also identifies the canonical 
actions of $\chalg{V,\chi}$ and $\chalg{V,\xi}$ on 
the symmetric algebra $S^* \left(V \otimes_\kk t\kk[t]\right)$ 
described in
\S\ref{section-fock-space-and-its-symmetric-algebra-structure}:
\begin{Lemma}
\label{lemma-canonical-actions-of-heisenberg-on-symmetric-are-iso}
Let $V$ be a graded vector space and $\chi$ and $\xi$ be two 
non-degenerate graded bilinear forms on $V$. The following diagram
commutes:
\begin{equation}
\begin{tikzcd}[column sep = 2cm]
\chalg{V,\chi}
\ar{d}{\rho}[']{\sim}
\ar{rr}{\eqref{eqn-isomorphism-sym-algebra-fock-space}}
& &
\eend_\kk \bigl( S^* \left(V \otimes_\kk t\kk[t]\right) \bigr)
\\
\chalg{V,\xi}
\ar{urr}[',pos=0.25]{\eqref{eqn-isomorphism-sym-algebra-fock-space}}
& &
\end{tikzcd}
\end{equation}
\end{Lemma}
\begin{proof}
By definition, the isomorphism \eqref{eqn-isomorphism-sym-algebra-fock-space}
is a composition of the isomorphism 
\eqref{eqn-isomorphism-sym-algebra-hplus} from 
$S^* \left(V \otimes_\kk t\kk[t]\right)$ to 
$H_+$ with the isomorphism 
\eqref{eqn-identification-of-hplus-and-fock-space}
from $H_+$ to $\falg{V,\chi}$. The latter comprises the horizontal 
rows of \eqref{eqn-fock-space-and-hplus-diagram-for-two-forms}, thus 
suffices to show that the following diagram commutes:
\begin{equation}
\begin{tikzcd}
S^* \left(V \otimes_\kk t\kk[t]\right)
\ar{rr}{\eqref{eqn-isomorphism-sym-algebra-hplus}}
\ar{drr}[']{\eqref{eqn-isomorphism-sym-algebra-hplus}}
& &
H_+ \subset \chalg{V,\chi}
\ar{d}{\rho}[']{\sim}
\\
& &
H_+ \subset \chalg{V,\chi}.
\end{tikzcd}
\end{equation}
Indeed, the isomorphism \eqref{eqn-isomorphism-sym-algebra-hplus}
sends $v \otimes t^n$ to $a_v(n)$, and $\rho|_{H_+}$ 
sends $a_v(n)$ to itself. 
\end{proof}

\section{Generalised Grojnowski-Nakajima action}
\label{section-generalised-nakajima-grojnowski-action}

Geometrical relevance of the Heisenberg algebras came to prominence with
the following famous result by Grojnowski and Nakajima:
 
\begin{Theorem}[see \cite{nakajima1997heisenberg}, Theorem 3.1, \cite{grojnowski1995instantons},
Theorem 7, and \cite{nakajima1999lectures}, Theorem 8.13]
\label{theorem-original-nakajima-grojnowski-heisenberg-action}
Let $X$ be a smooth projective surface over $\mathbb{C}$. 
Let $X^{[n]}$ be the Hilbert scheme of $n$ points on $X$. Let $\chi$ 
be the pairing on $H^\bullet(X,\mathbb{Q})$ given 
by taking the cup product and then the direct image along 
$X \rightarrow \text{pt}$. 

For any $n \geq 0$ and any $i > 0$ define a subvariety 
\begin{equation}
P(i) := \left\{ \mathcal{I}, \mathcal{J}, x) 
\;\middle|\;
\mathcal{I} \subset \mathcal{J},  \supp(\mathcal{J}/\mathcal{I}) =
\left\{ x \right\} \right\}
\subset 
X^{[n]} \times X^{[n + i]} \times X,
\end{equation}
Let $\pi_{12}$, $\pi_{21}$, $\pi_{3}$ be the 
projections from the triple product to 
$X^{[n]} \times X^{[n + i]}$, 
$X^{[n + i]} \times X^{[n]}$, and $X$, respectively. 
For each $\alpha \in H^\bullet(X, \mathbb{Q})$ and $i > 0$,
define operators $A_\alpha(\pm i)$ on 
$\bigoplus_{n=0}^{\infty}
H^\bullet\left(X^{[n]},\mathbb{Q}\right)$ by
\begin{equation}
A_\alpha(-i)\colon 
H^\bullet(X^{[n + i]}, \mathbb{Q})
\xrightarrow{(-1)^m \pi_{21 *} \left([P(i)] \cup \pi_3^* \alpha\right)
}
H^\bullet(X^{[n]}, \mathbb{Q})
\end{equation}
\begin{equation}
A_\alpha(i)\colon 
H^\bullet(X^{[n]}, \mathbb{Q})
\xrightarrow{\pi_{12 *} \left([P(i)] \cup \pi_3^* \alpha\right)
}
H^\bullet(X^{[n+i]}, \mathbb{Q}). 
\end{equation}
These satisfy 
\begin{equation}
\label{eqn-symmetric-heisenberg-relation}
A_{\alpha}(i) A_{\beta}(j) = (-1)^{\deg(\alpha)\deg(\beta)}
A_{\beta}(j)  A_{\alpha}(i) + 
\delta_{i,-j} i\langle \alpha,\, \beta\rangle_\chi
\end{equation}
and thus define an action of the Heisenberg algebra 
$\chalg{H^\bullet(X,\mathbb{Q}), \chi}$ on the
total cohomology $\bigoplus_{n=0}^{\infty} H^\bullet(X^{[n]},
\mathbb{Q})$.
This action
identifies $\bigoplus_{n=0}^{\infty} H^\bullet(X^{[n]}, \mathbb{Q})$ 
with the Fock space $\falg{H^\bullet(X,\mathbb{Q}), \chi}$. 
\end{Theorem}
The Hilbert schemes $X^{[n]}$ are smooth and the above says 
that their cohomology is determined by that of $X$. 
Indeed, the Betti numbers of 
$X^{[n]}$ were computed earlier by G{\"o}ttsche \cite{Gottsche-TheBettiNumbersOfTheHilbertSchemeOfPointsOnASmoothProjectiveSurface}
indicating that the dimension of  $\bigoplus_{n=0}^{\infty}
H^\bullet(X^{[n]}, \mathbb{Q})$ coincided with that of
$\falg{H^\bullet(X,\mathbb{Q}), \chi}$. The results of Grojnowski and
Nakajima gave a natural isomorphism between the two defined by 
geometrical correspondences. 

When $\dim(X) \geq 3$, $X^{[n]}$ are badly singular and their
cohomology is no longer well behaved. Grojnowski conjectured  
\cite[Footnote 3]{grojnowski1995instantons} that 
a similar result 
should hold for a smooth projective variety $X$ of any dimension 
if we replace $X^{[n]}$ by the symmetric quotient orbifold 
$X^n/S_n$ and $H^\bullet(-,\mathbb{Q})$ by equivariant topological
$K$-theory. This was later proved by 
Segal \cite{Segal-EquivariantKTheoryAndSymmetricProducts}
and Wang
\cite{Wang-EquivariantKTheoryWreathProductsAndHeisenbergAlgebra}. 

Baranovsky decomposition
\cite{Baranovsky-OrbifoldCohomologyAsPeriodicCyclicHomology} allows
to translate our results in 
\S\ref{section-action-on-hochschild-homology-of-symmetric-quotient-stacks}
from Hochschild homology to 
the orbifold cohomology introduced by Chen and Ruan 
\cite{ChenRuan-ANewCohomologyTheoryofOrbifold}. They constructed
a rationally graded ring with a rather intricately defined
product and rational grading structures, cf.
\cite[\S4]{AdemLeidaRuan-OrbifoldsAndStringyTopology}
\cite[\S2]{FantechiGottsche-OrbifoldCohomologyForGlobalQuotients}. 
We need it as a module for the Heisenberg algebra, so the product 
structure is irrelevant. Moreover, its grading as such module
is not the rational grading above. 
The following simple definition suffices:

\begin{Definition}[see
\cite{AdemLeidaRuan-OrbifoldsAndStringyTopology}, Remark 4.18, \cite{FantechiGottsche-OrbifoldCohomologyForGlobalQuotients}, Definition~1.1]
Let $Y$ be a smooth complex variety and $G$ a finite group acting on $Y$. 
The \em orbifold cohomology \rm of $Y/G$ is the vector space
$$ H_{orb}^\bullet(Y/G,\mathbb{C}) := \left( \bigoplus_{g \in G} H^\bullet(Y_g, \mathbb{C})\right)_G, $$
with its natural $\mathbb{Z}$-grading. Here 
$Y_g$ is the fixed point locus of $g \in G$ and $\left(-\right)_G$ denotes
taking coinvariants under the action of $G$
induced by each $h \in G$ acting as $X_g \xrightarrow{h.(-)} X_{hgh^{-1}}$. 
\end{Definition}

Recall the notation of 
\S\ref{section-fock-space-and-its-symmetric-algebra-structure}, 
e.g. $S_{\underline{r}(\underline{n})} := 
S_{r_1(n)} \times \dots \times S_{r_n(n)} < S_{r(\underline{n})}$,
and write $X_{\underline{n}} := X^{r(\underline{n})}$. 
We have an isomorphism: 
\begin{equation}
\label{eqn-description-of-cr-cohomology-of-xn-sn-via-partitions}
\left(\bigoplus_{\sigma \in S_n}
H^\bullet\left((X^n)_\sigma,\mathbb{C}\right) \right)_{S_n}
\simeq 
\bigoplus_{\underline{n} \vdash n}
H^\bullet\left(X_{\underline{n}},\mathbb{C}\right)_{S_{\underline{r}(\underline{n})}}.
\end{equation}
For any $\sigma \in S_n$ the fixed point locus 
$(X^n)_\sigma$ is non-canonically
isomorphic to $X_{\underline{n}}$ where $\underline{n} \vdash n$
is the conjugacy class of $\sigma$. The isomorphism 
\eqref{eqn-description-of-cr-cohomology-of-xn-sn-via-partitions}
is obtained by choosing a representative $\sigma_{\underline{n}} \in S_n$ 
of each conjugacy class $\underline{n} \vdash n$
and choosing an isomorphism 
$X_{\underline{n}} \simeq (X^n)_{\sigma_{\underline{n}}}$. This identifies
each $S_{\underline{r}(\underline{n})}$ with the quotient of the centraliser 
$C(\sigma_{\underline{n}})$ by the kernel of its action on 
$X_{\underline{n}}$, whence we obtain 
$\eqref{eqn-description-of-cr-cohomology-of-xn-sn-via-partitions}$. 
The choices do not matter due to taking the coinvariants, 
so the isomorphism 
\eqref{eqn-description-of-cr-cohomology-of-xn-sn-via-partitions}
is canonical. 

The K{\"u}nneth formula then yields:
\begin{Lemma}
\label{lemma-symmetrical-algebra-structure-on-total-chen-ruan-cohomology}
Let $Y$ be a smooth complex variety. We have an isomorphism
\begin{equation}
\bigoplus_{n=0}^{\infty}  H_{orb}^\bullet(X^n/S_n,\mathbb{C}) 
\; \simeq \;
S^*\left( H^\bullet(X,\mathbb{C}) \otimes_\mathbb{C} t\mathbb{C}[t]\right). 
\end{equation}
\end{Lemma}
\begin{proof}
By K{\"u}nneth formula we have  
$$
H^\bullet\left(X_{\underline{n}},\mathbb{C}\right)_{S_{\underline{r}(\underline{n})}}
\simeq
\left(H(X,\mathbb{C})^{r(\underline{n})}\right)_{S_{\underline{r}(\underline{n})}}
= \Sym^{r(\underline{n})} H(X,\mathbb{C}). $$
We thus have a chain of isomorphisms:
\begin{align*}
\bigoplus_{n=0}^{\infty}  H_{orb}^\bullet(X^n/S_n,\mathbb{C})
:= 
& \bigoplus_{n=0}^{\infty} \left(\bigoplus_{\sigma \in S_n}
H^\bullet\left((X^n)_\sigma,\mathbb{C}\right) \right)_{S_n}
\xrightarrow{\eqref{eqn-description-of-cr-cohomology-of-xn-sn-via-partitions}} 
\bigoplus_{n=0}^{\infty} \bigoplus_{\underline{n} \vdash n}
H^\bullet\left(X_{\underline{n}},\mathbb{C}\right)_{S_{\underline{r}(\underline{n})}}
\xrightarrow{\text{K{\"u}nneth}} \\
\rightarrow 
&
\bigoplus_{n=0}^{\infty} \bigoplus_{\underline{n} \vdash n}
\Sym^{r(\underline{n})} H^\bullet(X,\mathbb{C})
\xrightarrow{\eqref{eqn-direct-sum-description-of-tkt-sym-algebra-of-v}}
S^*\left( H(X,\mathbb{C}) \otimes_\mathbb{C} t\mathbb{C}[t]\right).
\end{align*}
\end{proof}

The isomorphism of Lemma
\ref{lemma-symmetrical-algebra-structure-on-total-chen-ruan-cohomology}
identifies the total orbifold cohomology 
$\bigoplus_{n=0}^{\infty}  H_{orb}^\bullet(X^n/S_n,\mathbb{C})$
with the symmetric algebra $S^*\left(H(X,\mathbb{C})
\otimes_\mathbb{C} t\mathbb{C}[t]\right)$. The latter, as explained in 
\S\ref{section-fock-space-and-its-symmetric-algebra-structure}, is
canonically isomorphic to the Fock space representation of the
Heisenberg algebra $\chalg{H(X,\mathbb{C})}$. We thus obtain 
a natural Heisenberg algebra action of $\chalg{H(X,\mathbb{C})}$
on $\bigoplus_{n=0}^{\infty}  H_{orb}^\bullet(X^n/S_n,\mathbb{C})$
which we computed explicitly in
Lemma \ref{lemma-direct-sum-description-of-tkt-sym-algebra-of-v} 
in terms of the decomposition 
$\bigoplus_{n=0}^{\infty} \bigoplus_{\underline{n} \vdash n}
\Sym^{r(\underline{n})} H(X,\mathbb{C}) t^{\underline{n}}$. 

On the other hand, for any DG category $\basecat$ over $\mathbb{C}$
(or, more generally, over any base field $\kk$), the noncommutative 
orbifold decomposition \cite[Theorem
4.1]{annobaranovskylogvinenko2023orbifold} gives an isomorphism 
between the total Hochschild homology 
$\bigoplus_{n \geq 0} \hochhom_\bullet(\sym^n \basecat)$ of the
symmetric powers of $\basecat$ and the symmetric algebra 
$S^*\left(\hochhom_\bullet(\basecat) 
\otimes_\mathbb{C} t\mathbb{C}[t]\right)$. Again, as per 
\S\ref{section-fock-space-and-its-symmetric-algebra-structure}, 
we obtain a natural action of the Heisenberg algebra 
$\chalg{\hochhom_\bullet(\basecat)}$ on 
$\bigoplus_{n \geq 0} \hochhom_\bullet(\sym^n \basecat)$
which we call the \em orbifold decomposition action\rm. 
In \S\ref{section-action-on-hochschild-homology-of-symmetric-quotient-stacks}
we show that this action coincides with the \em decategorified action
\rm obtained in \S\ref{section-hochschild-homology-and-heisenberg-2-category}
from the categorical action of the Heisenberg $2$-category 
of $\basecat$ constructed in \cite[\S4]{gyenge2021heisenberg}. We also 
describe it directly on the symmetric powers $\sym^n \basecat$ as 
the \em induction-restriction action \rm by the operators analogous 
to those in the equvariant topological $K$-theory constructions of 
Segal \cite{Segal-EquivariantKTheoryAndSymmetricProducts} and 
Wang \cite{Wang-EquivariantKTheoryWreathProductsAndHeisenbergAlgebra}. 

When $\basecat$ is the standard DG-enhancement of a smooth projective
variety $X$, we have, as explained in \S\ref{section-hkr-isomorphism}, 
the Hochschild-Kostant-Rosenberg isomorphism 
$\hochhom_\bullet(\basecat) \simeq H^\bullet(X,\mathbb{C})$. 
It doesn't respect the $\mathbb{Z}$-grading, but it respects the
$\mathbb{Z}/2$-grading. It induces bigraded isomorphisms 
$\chalg{\hochhom_\bullet(\basecat)} \simeq
\chalg{H^\bullet(X,\mathbb{C})}$ of the Heisenberg algebras 
and $S^*\left(\hochhom_\bullet(\basecat) 
\otimes_\mathbb{C} t\mathbb{C}[t]\right) 
\simeq S^*\left(H(X,\mathbb{C}) \otimes_\mathbb{C} t\mathbb{C}[t]\right)$
of the symmetric algebras, where we coarsen the base grading to
$\mathbb{Z}/2$-grading. By 
Lemma \ref{lemma-canonical-actions-of-heisenberg-on-symmetric-are-iso}
these isomorphisms identify the natural action of 
of $\chalg{H(X,\mathbb{C})}$
on $\bigoplus_{n=0}^{\infty}  H_{orb}^\bullet(X^n/S_n,\mathbb{C})$
given by Lemma
\ref{lemma-symmetrical-algebra-structure-on-total-chen-ruan-cohomology}
with the orbifold decomposition action in the noncommutative setting, 
and hence also with the decategorified action of 
\cite{gyenge2021heisenberg} and the induction-restriction action. 

\begin{center}
\fbox{
\begin{minipage}{0.9\textwidth}
The main result of this section completes this picture by 
giving a direct description of this action of $\chalg{H(X,\mathbb{C})}$
on $\bigoplus_{n=0}^{\infty}  H_{orb}^\bullet(X^n/S_n,\mathbb{C})$
in terms of geometrical correspondences analogous to those of
Grojnowski \cite{grojnowski1995instantons} and 
Nakajima \cite{nakajima1997heisenberg}. 
\end{minipage}
}
\end{center}

We spell out the action by a correspondence on complex
cohomology to set up the notation: 

\begin{Definition}
Let $X$ and $Y$ be smooth projective varieties. Let $\alpha 
\in H_*(X \times Y)$ be a homology class. The corresponding action 
\begin{equation}
H^*(X) \xrightarrow{\alpha} H^*(Y)
\end{equation}
is given by
$$\beta \mapsto \left( \pi_{X *}\left(\alpha \cap \pi_Y^*
\beta\right)\right)^\vee$$
where $\cap$ is the cap product and $(-)^\vee$ denotes the Poincare duality. 
\end{Definition}

For any $n \geq 1$ and any $\underline{n} \vdash n$, we have 
$X_{\underline{n}} = X^{r_1(\underline{n})} \times 
\dots \times X^{r_n(\underline{n})}$.  
For any $\underline{x} \in X_{\underline{n}}$ we write 
$$ \underline{x} = (\underline{x}_1, \dots, \underline{x}_n)
\in X^{r_1(\underline{n})} \times \dots \times X^{r_n(\underline{n})} $$
where $\underline{x}_j = (x_{j1}, \dots, x_{j r_i(\underline{n})}) \in 
X^{r_j(\underline{n})}$ is the component of $\underline{x}$ corresponding 
to parts of size $j$ in $\underline{n}$. If $r_j(\underline{n}) = 0$, 
i.e. there are not parts of size $j$ in  $\underline{n}$, we have
$\underline{x}_j = \emptyset$. Thus
$$ \underline{x} = (\underline{x}_1, \dots, \underline{x}_n) = 
(x_{jk})_{1 \leq j \leq n, 1 \leq k \leq r_j(\underline{n})} \in
X_{\underline{n}}. $$ 
Finally, for any $m \geq 1$, define the map of sets  
$\forget\colon X^m \rightarrow X^m/S_m$
to be the map which forgets the ordering on an $m$-tuple of points. 

\begin{Definition}
\label{defn-generalised-grojnowski-nakajima-correspondence}
Let $X$ be a smooth projective variety. Let $i > 0$ and write 
$(i)$ for the single part partition of $i$.  
For any $n \geq 0$ and any $\underline{n} \vdash n$, 
define the subvariety
\begin{equation}
\label{eqn-generalised-grojnowski-nakajima-correspondence}
Z(i) := \left\{ (\underline{x}, \underline{y}, z) 
\;\middle|\;
\begin{matrix}
\forget(\underline{x}_j) = \forget(\underline{y}_j) \cup z & j = i,
\\
\forget(\underline{x}_j) = \forget(\underline{y}_j) & j \neq i.
\end{matrix}
\right\} \subset 
X_{\underline{n} + (i)} \times X_{\underline{n}}  \times X,
\end{equation}
where $\underline{n} + (i)$ is the partition of $n+i$
given by the union of $\underline{n}$ and $(i)$.
Let $\pi_{12}$, $\pi_{21}$ and $\pi_{3}$ be the 
projections from the triple product to 
$X_{\underline{n} + (i)} \times X_{\underline{n}}$, 
to 
$X_{\underline{n}} \times X_{\underline{n} + (i)}$, 
and to $X$, respectively. 
\end{Definition}

The following generalises Grojnowski-Nakajima action 
to all smooth projective varieties:

\begin{Theorem}
\label{theorem-generalised-nakajima-grojnowski-heisenberg-action}
Let $X$ be a smooth projective variety over $\mathbb{C}$ and $\chi$
be the pairing on $H^\bullet(X,\mathbb{C})$ given 
by taking the cup product and then the direct image along 
$X \rightarrow \text{pt}$. 

For each $\alpha \in H^\bullet(X, \mathbb{C})$ and $i > 0$,
define operators $A_\alpha(\pm i)$ on 
$\bigoplus_{n=0}^{\infty}
H_{orb}^\bullet\left(X^n/S_n,\mathbb{C}\right)$ by
\begin{equation}
\label{eqn-geometric-heisenberg-action-a-minus-i-operator}
A_\alpha(-i)\colon 
H^\bullet(X_{\underline{n} + (i)}, \mathbb{C})_{S_{\underline{r}(\underline{n} + (i))}}
\xrightarrow{\frac{1}{|S_{\underline{r}(\underline{n})}|} \pi_{12 *}
\left([Z(i)] \cap \pi_3^* \alpha\right)
}
H^\bullet(X_{\underline{n}}, \mathbb{C})_{S_{\underline{r}(\underline{n})}}
\end{equation}
\begin{equation}
\label{eqn-geometric-heisenberg-action-a-plus-i-operator}
A_\alpha(i)\colon 
H^\bullet(X_{\underline{n}}, \mathbb{C})_{S_{\underline{r}(\underline{n})}}
\xrightarrow{\frac{1}{|S_{\underline{r}(\underline{n} + (i))}|}
\pi_{21 *} \left([Z(i)] \cap \pi_3^* \alpha\right)
}
H^\bullet(X_{\underline{n} + (i)}, \mathbb{C})_{S_{\underline{r}(\underline{n} + (i))}}.
\end{equation}
These define an action of the Heisenberg 
algebra $\chalg{H^\bullet(X,\mathbb{C}), \chi}$ on the total cohomology 
$\bigoplus_{n=0}^{\infty} H_{orb}^\bullet\left(X^n/S_n,\mathbb{C}\right)$
which equals the one given by identifying the latter 
with $S^*\left(H^\bullet(X,\mathbb{C}) \otimes_\mathbb{C}
t\mathbb{C}[t]\right)$ via Lemma 
\ref{lemma-symmetrical-algebra-structure-on-total-chen-ruan-cohomology}.
\end{Theorem}

\begin{proof}	
It suffices to show that applying the K{\"u}nneth isomorphisms of 
Lemma \ref{lemma-symmetrical-algebra-structure-on-total-chen-ruan-cohomology}
identifies the operators $A_\alpha(-n)$ and $A_\alpha(n)$ defined in 
\eqref{eqn-geometric-heisenberg-action-a-minus-i-operator}-\eqref{eqn-geometric-heisenberg-action-a-plus-i-operator} with the operators spelled out in 
Lemma 
\ref{lemma-canonical-action-of-chalg-v-chi-on-symmetric-algebra}
which give the natural action of $\chalg{H^\bullet(X,\mathbb{C}),
\chi}$ on $S^*\left(H^\bullet(X,\mathbb{C}) \otimes_\mathbb{C}
t\mathbb{C}[t]\right)$. 

Let $n \geq 0$, $\underline{n} \vdash n$ and $i > 0$. 
By definition, the varieties $X_{\underline{n}}$ and $X_{\underline{n} + (i)}$ differ by a single extra factor of  $X$ in the latter corresponding to its extra part of size $i$. We have $r_i(\underline{n}) + 1$ isomorphisms 
\begin{equation}
\mu_j\colon 
X_{\underline{n} + (i)} 
\xrightarrow{\quad\sim\quad}
X_{\underline{n}}  \times X
\quad\quad 1 \leq j \leq r_i(\underline{n})+1
\end{equation}
identifying $X$ on the RHS with one of the $r_i(\underline{n}) + 1$ size $i$ part factors on the LHS. Explicitly, 
$$ 
\mu_j(\underline{x}_1,  \dots, \underline{x}_n)
:=
\Bigl( \bigl(\underline{x}_1,  \dots ,
\underline{x}_i(\hat{j}), 
\dots \underline{x}_n)\bigr), x_{ij} \Bigr),
$$
where $\underline{x}_i(\hat{j})$ denote $\underline{x}_i$ with $j$-th coordinate skipped:
$$
\underline{x}_i(\hat{j}) := 
(x_{i1}, \dots, x_{i(j-1)}, x_{i(j+1)}, \dots,
x_{ir_i(\underline{n}+(i))}).
$$
Denote by $p_j$ and $q_j$ the compositions of $\mu_j$ with projections to $X_{\underline{n}}$ and $X$, respectively: 
\begin{equation}
\begin{tikzcd}
&
X_{\underline{n} + (i)}
\ar{ld}[']{p_j}
\ar{rd}{q_j}
&
\\
X_{\underline{n}}
&
&
X
\end{tikzcd}
\end{equation}

Let $\sigma \in S_{\underline{r}(\underline{n}+(i))}$ and consider the diagram
\begin{equation}
\label{eqn-diagram-for-computing-correspondence-z(i)}
\begin{tikzcd}
X_{\underline{n} + (i)} 
\ar{r}{\Delta}
&
X_{\underline{n} + (i)} \times X_{\underline{n} + (i)} 
\ar{r}{\sigma \times \mu_1}
&
X_{\underline{n} + (i) } \times X_{\underline{n}} \times X
\ar{ld}{\pi_{12}}
\ar{rd}{\pi_{3}}
&
\\
&
X_{\underline{n} + (i)} \times  X_{\underline{n}}
\ar{ld}{\pi_{1}}
\ar{rd}{\pi_{2}}
&
&
X
\\
X_{\underline{n} + (i)}
&
&
X_{\underline{n}}
&
\end{tikzcd}
\end{equation}

The
closed subvariety 
defined in \eqref{eqn-generalised-grojnowski-nakajima-correspondence}
\begin{equation*}
Z(i) \subset 
X_{\underline{n} + (i)} \times X_{\underline{n}} \times X,
\end{equation*}
is reducible. Its 
irreducible components are 
\begin{equation}
Z(i) = \bigcup_{\sigma \in S_{\underline{r}(\underline{n}+(i))}}
\img((\sigma \times \mu_1) \circ \Delta).
\end{equation}
Therefore in 
$H_*(X_{\underline{n} + (i)} \times X_{\underline{n}} \times X)$ 
we have
$$
[Z(i)] = \sum_{\sigma \in S_{\underline{r}(\underline{n}+(i))}}
(\sigma \times \mu_{1})_* \Delta_* [X_{\underline{n} + (i)}]. 
$$

Let $\alpha \in H^\bullet(X,\mathbb{C})$. 
By the projection formula for cap product
\cite[\S3.3]{Hatcher-AlgebraicTopology} we have 
\begin{align*}
\pi_{12 *} \left( [Z(i)] \cap \pi_3^* \alpha \right) &= 
\sum_{\sigma \in S_{\underline{r}(\underline{n}+(i))}}
\pi_{12 *} \left((\sigma \times \mu_{1})_* \Delta_* [X_{\underline{n} + (i)}]
\cap \pi_3^* \alpha \right) = 
\\
& = 
\sum_{\sigma \in S_{\underline{r}(\underline{n}+(i))}}
\pi_{12 *} (\sigma \times \mu_{1})_* \Delta_* \left( [X_{\underline{n} + (i)}]
\cap \Delta^* (\sigma \times \mu_{1})^* \pi_3^* \alpha \right) = 
\\
& = 
\sum_{\sigma \in S_{\underline{r}(\underline{n}+(i))}}
(\sigma \times p_1)_* \Delta_* \left( [X_{\underline{n} + (i)}]
\cap q_1^* \alpha \right) 
\in H^\bullet(X_{\underline{n} + (i)} \times X_{\underline{n}}, \mathbb{C}). 
\end{align*}

Let now 
$$ \underline{\beta} = \underline{\beta}_1 \times \dots \times
\underline{\beta}_n \in H^\bullet(X_{\underline{n}}, \mathbb{C}) $$
with each $\underline{\beta}_j \in H^\bullet(X^{r_j(\underline{n})})$
further equal to $\beta_{j1} \times \dots \times \beta_{jr_j(\underline{n})}$
with $\beta_{jk} \in H^\bullet(X,\mathbb{C})$. In the notation 
of diagram \eqref{eqn-diagram-for-computing-correspondence-z(i)}
the image of $\underline{\beta}$ under the correspondence defining 
the operator $A_{\alpha}(i)$ in 
\eqref{eqn-geometric-heisenberg-action-a-plus-i-operator} is:
\begin{align}
\nonumber
& \frac{1}{|S_{\underline{r}(\underline{n} + (i))}|}
\sum_{\sigma \in S_{\underline{r}(\underline{n} + (i))}}
\left(\pi_{1*} \Bigl( 
(\sigma \times p_1)_* \Delta_* \left( [X_{\underline{n} + (i)}]
\cap q_1^* \alpha \right)
\cap \pi_2^* \underline{\beta} \Bigr)\right)^\vee = 
\\
\nonumber
=
& \frac{1}{|S_{\underline{r}(\underline{n} + (i))}|}
\sum_{\sigma \in S_{\underline{r}(\underline{n} + (i))}}
\left(\pi_{1*} (\sigma \times p_1)_* \Delta_* \Bigl( 
[X_{\underline{n} + (i)}] \cap 
\left( q_1^* \alpha \cup \Delta^* (\sigma \times p_1)^* \pi_2^* 
\underline{\beta} \right) \Bigr)\right)^\vee = 
\\
\nonumber
=
& \frac{1}{|S_{\underline{r}(\underline{n} + (i))}|}
\sum_{\sigma \in S_{\underline{r}(\underline{n} + (i))}}
\left(\sigma_* \Bigl( 
[X_{\underline{n} + (i)}] \cap 
\left( q_1^* \alpha \cup p_1^* \underline{\beta} \right) \Bigr)\right)^\vee =
\\
\nonumber
=
& \frac{1}{|S_{\underline{r}(\underline{n} + (i))}|}
\sum_{\sigma \in S_{\underline{r}(\underline{n} + (i))}}
\Bigl(\sigma_*[X_{\underline{n} + (i)}] \cap
(\sigma^{-1})^* 
\left( q_1^* \alpha \cup p_1^* \underline{\beta} \right) \Bigr)^\vee =
\\
\label{eqn-a-alpha-i-before-coinvariants-simplified}
= 
& \frac{1}{|S_{\underline{r}(\underline{n} + (i))}|}
\sum_{\sigma \in S_{\underline{r}(\underline{n} + (i))}}
\sigma^* \left(q_1^* \alpha \cup p_1^* \underline{\beta}\right),
\end{align}
where in the last equality we use $\sigma_* [X_{\underline{n} + (i)}]
= [X_{\underline{n} + (i)}]$ and reindex the summation 
from $\sigma$ to $\sigma^{-1}$. 

To see that the resulting map descends to the coinvariants, 
consider the subgroup $S_{\underline{r}(\underline{n})} < 
S_{\underline{r}(\underline{n} + (i))}$
of the elements acting trivially on the leftmost copy of
$X$ in the size $i$ part factor of 
$X_{\underline{n} + (i)}$. Recall, that $X_{\underline{n} + (i)}
:= X^{r_1(\underline{n} + (i))} \times \dots \times 
X^{r_n(\underline{n} + (i))}$ is a product of one copy of $X$
for each part in $\underline{n} + (i)$. The size $i$ part factor
is $X^{r_i(\underline{n} + (i))}$, the copies of $X$ which 
correspond to the parts of size $i$.
This subgroup $S_{\underline{r}(\underline{n})}$ 
has $r_i(\underline{n}+(i))$ left cosets 
whose representatives are $\id, i(12), i(123), \dots, i(12 \dots
r_i(\underline{n}+(i)))$, where $i(1{\dots}j)$ is the cycle 
which acts by $(1{\dots}j)$ on $X^{r_i(\underline{n}+(i))}$ and 
trivially on the rest of  $X_{\underline{n} + (i)}$. We can thus
rewrite  \eqref{eqn-a-alpha-i-before-coinvariants-simplified} as 
\begin{equation*}
\frac{1}{|S_{\underline{r}(\underline{n} + (i))}|}
\sum_{j=1}^{r_i(\underline{n}+(i))}
\sum_{\tau \in S_{\underline{r}(\underline{n})}}
(i(1{\dots}j))^* \tau^* \left(q_1^* \alpha \cup p_1^* \underline{\beta}\right) 
= 
\frac{1}{|S_{\underline{r}(\underline{n} + (i))}|}
\sum_{j=1}^{r_i(\underline{n}+(i))}
\sum_{\tau \in S_{\underline{r}(\underline{n})}}
q_j^* \alpha \otimes p_j^* \tau^* \underline{\beta}, 
\end{equation*}
which clearly kills all the elements of the form $\underline{\beta} - \tau'.
\underline{\beta}$ for $\tau' \in  S_{\underline{r}(\underline{n})}$.  

Thus \eqref{eqn-a-alpha-i-before-coinvariants-simplified} descends to 
the coinvariants. There it becomes the map 
$$ [\beta] \mapsto 
\left[\frac{1}{|S_{\underline{r}(\underline{n} + (i))}|}
\sum_{\sigma \in S_{\underline{r}(\underline{n} + (i))}}
\sigma^* \left(q_1^* \alpha \cup p_1^* \underline{\beta}\right)\right]
= [q_1^* \alpha \cup p_1^* \underline{\beta}]. $$
In components, this is the map 
$$ [(\underline{\beta}_1 \times \dots \times \underline{\beta}_n)]
\mapsto 
= 
(-1)^{\alpha | \underline{\beta}_1, \dots, \underline{\beta}_{i-1}}
[(\underline{\beta}_1 \times \dots \times (\alpha \times
\underline{\beta}_i) \times \dots \times \underline{\beta}_n)], $$
which the K{\"u}nneth isomorphisms of 
Lemma \ref{lemma-symmetrical-algebra-structure-on-total-chen-ruan-cohomology}
identify with the operator 
\eqref{eqn-action-of-a_w(n)-on-symmetric-algebra}, as
required. 

Similarly, let 
$$ \underline{\gamma} = \underline{\gamma}_1 \times \dots \times
\underline{\gamma}_n \in H^\bullet(X_{\underline{n}+(i)}, \mathbb{C}) $$
with each $\underline{\gamma}_j \in H^\bullet(X^{r_j(\underline{n}+(i))})$
further equal to $\gamma_{j1} \times \dots \times
\gamma_{jr_j(\underline{n}+(i))}$
with $\gamma_{jk} \in H^\bullet(X,\mathbb{C})$. 
The image of $\underline{\gamma}$ under the correspondence defining 
the operator $A_{\alpha}(-i)$ in 
\eqref{eqn-geometric-heisenberg-action-a-minus-i-operator} is:
\begin{align}
\nonumber
& \frac{1}{|S_{\underline{r}(\underline{n})}|}
\sum_{\sigma \in S_{\underline{r}(\underline{n}+(i))}}
\left(\pi_{2*} \Bigl( 
(\sigma \times p_1)_* \Delta_* \left( [X_{\underline{n} + (i)}]
\cap q_1^* \alpha \right)
\cap \pi_1^* \underline{\gamma} \Bigr)\right)^\vee = 
\\
\nonumber
=
& \frac{1}{|S_{\underline{r}(\underline{n})}|}
\sum_{\sigma \in S_{\underline{r}(\underline{n} + (i))}}
\left( \pi_{2*} (\sigma \times p_1)_* \Delta_* \Bigl( 
[X_{\underline{n} + (i)}] \cap 
\left( q_1^* \alpha \cup \Delta^* (\sigma \times p_1)^* \pi_1^* 
\underline{\gamma} \right) \Bigr)\right)^\vee = 
\\
\nonumber
=
& \frac{1}{|S_{\underline{r}(\underline{n})}|}
\sum_{\sigma \in S_{\underline{r}(\underline{n} + (i))}}
\left(p_{1 *}
\Bigl( 
[X_{\underline{n} + (i)}] \cap 
\left( q_1^* \alpha \cup \sigma^* \underline{\gamma} \right)
\Bigr)\right)^\vee, 
\end{align}
which is a sum over $\sigma \in S_{\underline{r}(\underline{n} + (i))}$
of expressions in $\sigma^*  \underline{\gamma}$ and
therefore descends to coinvariants in $\underline{\gamma}$. 
It further equals to
\begin{align}
\nonumber
& \frac{1}{|S_{\underline{r}(\underline{n})}|}
\sum_{\sigma \in S_{\underline{r}(\underline{n} + (i))}}
\left(p_{1 *} (\sigma^{-1})_* \sigma_* 
\Bigl([X_{\underline{n} + (i)}] \cap
\sigma^* 
\left((\sigma^{-1})^* q_1^* \alpha \cup \underline{\gamma} \right)
\Bigr)\right)^\vee =
\\
\nonumber
=
& \frac{1}{|S_{\underline{r}(\underline{n})}|}
\sum_{\sigma \in S_{\underline{r}(\underline{n} + (i))}}
\left(p_{1 *} (\sigma^{-1})_*  
\Bigl(\sigma_* [X_{\underline{n} + (i)}] \cap
\left((\sigma^{-1})^* q_1^* \alpha \cup \underline{\gamma} \right)
\Bigr)\right)^\vee =
\\
\nonumber
= 
& \frac{1}{|S_{\underline{r}(\underline{n})}|}
\sum_{\sigma \in S_{\underline{r}(\underline{n} + (i))}}
\left(p_{1 *} \sigma_* 
\Bigl([X_{\underline{n} + (i)}] \cap
\left(\sigma^* q_1^* \alpha \cup \underline{\gamma} \right)
\Bigr)\right)^\vee. 
\end{align}
Using the same left coset decomposition of 
$S_{\underline{r}(\underline{n} + (i))}$ as above, we further rewrite this as:
\begin{align}
\nonumber 
& \frac{1}{|S_{\underline{r}(\underline{n})}|}
\sum_{j=1}^{r_i(\underline{n}+(i))}
\sum_{\tau \in S_{\underline{r}(\underline{n})}}
\left(p_{1 *} \tau_* (i(1{\dots}j)_*  
\Bigl([X_{\underline{n} + (i)}] \cap
\left((i(1{\dots}j)^* \tau^* q_1^* \alpha \cup \underline{\gamma} \right)
\Bigr)\right)^\vee
=
\\
\label{eqn-a-alpha-minus-i-before-coinvariants-simplified}
= & 
\frac{1}{|S_{\underline{r}(\underline{n})}|}
\sum_{j=1}^{r_i(\underline{n}+(i))}
\sum_{\tau \in S_{\underline{r}(\underline{n})}}
\left(\tau_* p_{j *}   
\Bigl([X_{\underline{n} + (i)}] \cap
\left( q_j^* \alpha \cup \underline{\gamma} \right)
\Bigr)\right)^\vee
\end{align}
In components, we have
\begin{small}
$$ q_j^* \alpha \cup \underline{\gamma} 
=
(-1)^{\alpha | \underline{\gamma}_1, \dots, 
\underline{\gamma}_{i-1}}
(-1)^{\alpha | \gamma_{i1}, \dots, 
\gamma_{i(j-1)}}
\underline{\gamma}_1 \times 
\dots 
\times
(\gamma_{i1} \times \dots \times (\alpha \cup \gamma_{ij})
\times \dots \times \gamma_{ir_i(\underline{n}+(i))})
\times 
\dots 
\underline{\gamma}_n
$$
\end{small}
The projection $p_{j *}$ contracts the copy of $X$
which carries $\alpha \cup \gamma_{ij}$ to a point and 
leaves all the other factors untouched. Therefore 
\eqref{eqn-a-alpha-minus-i-before-coinvariants-simplified}
equals to
\begin{align}
\nonumber
&
\frac{1}{|S_{\underline{r}(\underline{n})}|}
\sum_{j=1}^{r_i(\underline{n}+(i))}
\sum_{\tau \in S_{\underline{r}(\underline{n})}}
(-1)^{\alpha | \underline{\gamma}_1, \dots, 
\underline{\gamma}_{i-1}}
(-1)^{\alpha | \gamma_{i1}, \dots, 
\gamma_{i(j-1)}}
\left<\alpha,\gamma_{ij}\right>
\tau_* \gamma(\hat{ij}) 
\end{align}
where
$$
\gamma(\hat{ij}) := 
\underline{\gamma}_1 \times 
\dots 
\times
(\gamma_{i1} \times \dots \times
\gamma_{i(j-1)} \times \gamma_{i(j+1)}  
\times \dots \times \gamma_{ir_i(\underline{n}+(i))})
\times 
\dots 
\underline{\gamma}_n.
$$

We conclude that 
\eqref{eqn-a-alpha-minus-i-before-coinvariants-simplified}
descends to coinvariants as the map 
$$
[\underline{\gamma}] \mapsto 
\sum_{j=1}^{r_i(\underline{n}+(i))}
(-1)^{\alpha | \underline{\gamma}_1, \dots, 
\underline{\gamma}_{i-1}}
(-1)^{\alpha | \gamma_{i1}, \dots, 
\gamma_{i(j-1)}}
\left<\alpha,\gamma_{ij}\right> [\gamma(\hat{ij})], 
$$
which the K{\"u}nneth isomorphisms of 
Lemma \ref{lemma-symmetrical-algebra-structure-on-total-chen-ruan-cohomology}
identify with the operator 
\eqref{eqn-action-of-a_w(-n)-on-symmetric-algebra}, as
required.
\end{proof}

\section{Preliminaries on DG categories and Hochschild homology}

\subsection{DG categories}

A general introduction can be found in 
\cite[Section~2]{AnnoLogvinenko-SphericalDGFunctors}\cite{Toen-LecturesOnDGCategories} and the technicalities relevant to this paper 
in \cite[\S4]{gyenge2021heisenberg}. We 
use freely the concepts and notation introduced therein. 
 
For the convenience of the reader, we summarise some notation below. 
Let $\A$ be a small DG category over $\kk$.  We write $\modA$, $\Amod$, 
and $\AmodA$ to denote 
the DG categories of right and left $\A$-modules and of $\A$-$\A$-bimodules. 
Similarly, we write $\modk$ for the DG category of DG $\kk$-modules. 
We view DG algebras as DG categories with a single object and 
ordinary associative algebras as DG algebras concentrated in degree $0$. 

We view a small DG category as a Morita enhanced triangulated
category. Its underlying triangulated category $D_c(\A)$, 
the derived category of compact DG $\A$-modules. 
When $\A$ is a commutative associative algebra $A$ we 
can equivalently view it as an affine scheme $X = \Spec A$.
The triangulated category $D_c(A)$ is then equivalent to $D_c(X)$, 
the compact derived category of perfect complexes of quasi-coherent 
sheaves on $X$. Following
\cite{KontsevichSoibelman-NotesOnAInftyAlgebrasAInftyCategoriesAndNoncommutativeGeometry}, we also view any small DG category $\A$ as 
a noncommutative scheme whose derived category of perfect complexes 
is $D_c(\A)$. We say that such a noncommutative scheme is
commutative when $D_c(\A) \simeq D_c(X)$ for some scheme $X$ over $\kk$. 

This notion is most meaningful when $\A$ possesses
the following two properties
\cite[\S8.1-8.2]{KontsevichSoibelman-NotesOnAInftyAlgebrasAInftyCategoriesAndNoncommutativeGeometry}: we say that $\A$ is \em smooth \rm if $\A$ is a perfect
object in the category of $\AbimA$-bimodules and that $\A$ is \em proper
\rm when $\A$ is a perfect object in the category of $\kk$-modules, i.e. 
when the total cohomology of any $\homm$-complex in $\A$ is finite dimensional. 
If $\A$ is a commutative scheme $X$, then 
$\A$ is smooth (resp. proper) if and only if $X$ is smooth (resp. proper). 
Thus a smooth and proper noncommutative scheme for us is a smooth 
and proper DG category. These are our main objects of interest 
for which we construct an analogue of
the Grojnowski-Nakajima Heisenberg algera action described  
in \S\ref{section-generalised-nakajima-grojnowski-action}. 

\subsection{Strong group actions and DG categories}
\label{section-strong-group-actions-and-DG-cats}

We need to work with symmetric powers of DG categories. We realise
this technically as follows. 
Let $\A$ be a small \dg category. A \emph{strong} action 
of a finite group $G$ on $\A$ is an embedding of $G$ 
into the group of \dg automorphisms of $\A$.  

\begin{Definition}[\cite{gyenge2021heisenberg}, Definition~4.45]
\label{defn-equivariant-dg-category}
	The \emph{semi-direct product} $\A \rtimes G$ is the following \dg
	category:
	\begin{itemize}
		\item $\obj \A \rtimes G = \obj \A$,
		\item For any $a,b \in \obj(\A \rtimes G)$ their morphism complex is
		\begin{equation*}
			\homm^i_{\A \rtimes G}(a,b) :=
			\left\{ (\alpha, g) \; \middle| \;
			\alpha \in \homm^i_{\A}(g.a,b), g \in G
			\right\}
		\end{equation*}
		with $\deg_{\A \rtimes G} (\alpha,g) = \deg_{\A} \alpha$ and 
		$d_{\A \rtimes G}(\alpha,g) = (d_\A \alpha,g)$,
		\item The composition in $\A \rtimes G$ is given by 
		\begin{equation*}
			(\alpha_1, g_1) \circ (\alpha_2, g_2) =
			(\alpha_1 \circ g_1.\alpha_2,\, g_1 g_2). 
		\end{equation*}
		\item For any $a \in \obj (\A \rtimes G)$ the identity morphism 
		of $a$ is $(\id_a, 1_G)$. 
	\end{itemize}
\end{Definition}

The point of this definition is that modules over $\A \rtimes G$ 
are $G$-equivariant modules over $\A$:
\begin{Lemma}[\cite{gyenge2021heisenberg}, Lemma 4.49]
	\label{lemma-A-rtimes-G-modules-are-G-equiv-A-modules}
	There are mutually inverse isomorphisms of categories
	\[ \rightmod{(\A \rtimes G)} \leftrightarrows \modd^G\mkern-4mu\mhyphen\A, \]
        \[ \hperf{(\A \rtimes G)} \leftrightarrows \hperf^G\mkern-4mu\mhyphen\A.\]
\end{Lemma}

Each autoequivalence $g\colon \A \xrightarrow{\sim} \A$ with which $G$ acts
on $\A$ extends to an autoequivalence $\A \rtimes G$:
\begin{Definition}
\label{defn-group-element-as-autoequivalence-of-the-equivariant-category}
Let $\A$ be a small DG category and $G$ be a group acting strongly on
$\A$. For any $g \in G$ define the autoequivalence 
$$ g\colon \A \rtimes G \xrightarrow{\sim} \A \rtimes G $$
to have the same action on objects as
$g\colon \A \xrightarrow{\sim} \A$ and to act on morphisms 
by $g(\alpha, f) = (g(\alpha), gfg^{-1})$. 
\end{Definition}

In $\A \rtimes G$ the autoequivalence $g$ is
isomorphic to the identity functor. The natural isomorphism 
$\id_{\A \times G} \rightarrow g$ is given for $a \in \A$ by 
the isomorphism $g \colon a \rightarrow g.a$.  
Further technical aspects of semi-direct products $\A \times G$
can be found in \cite[\S4.8]{gyenge2021heisenberg}.
We also need one more technical result:

\begin{Definition}
\label{defn-the-set-of-left-cosets-of-a-subgroup}
Let $G$ be a group and $H \leq G$ be a subgroup. 
Write $Q$ for the set of left cosets of $H$ in $G$. 
The group $G$ acts on $Q$ by left multiplication. For any $g \in G$
and $q \in Q$ write $g.q \in Q$ to denote this action.
For any $g \in G$ write $\fix(g)$ for the fixed point set 
of the action of $g$ on $Q$. 

For every $q \in Q$ choose a representative $r_q \in G$ of 
the corresponding coset. Then for every $g \in G$ and $q \in Q$
write $h_{g,q}$ for the unique element of $H$ such that $g r_q =
r_{g.q} h_{g,q}$. 

\end{Definition}
\begin{Lemma}
\label{lemma-descripton-of-res-G-to-H-functor}
Let $\A$ be a small DG category, $G$ be a group acting strongly on
$\A$, and $H \leq G$ be a subgroup. 
Let $Q$, $\fix(g)$, $r_q$, and $h_{g,q}$ be 
as in Definition \ref{defn-the-set-of-left-cosets-of-a-subgroup}. Let 
$$\Res^G_H\colon \modd(\A \rtimes G) \rightarrow \modd(\A \rtimes H) $$
to be the restriction of scalars functor. Then 
\begin{equation}
\Res^G_H(a) \simeq \bigoplus_{q \in Q}
r_q^{-1}.a
\quad \quad \quad  
\forall\; a \in \A. 
\end{equation}
\begin{equation}
\Res^G_H(\alpha) \colon \bigoplus_{q \in G / H}
r_q^{-1}.a \rightarrow \bigoplus_{q \in G / H}
r_q^{-1}.b
\quad \quad \quad  
\forall\; \alpha \in \homm_\A(a,b),
\end{equation}
is the sum of morphisms $r_q^{-1}(\beta) \colon r_q^{-1}.a \rightarrow
r_q^{-1}.b$ over all $q \in Q$. 
\begin{equation}
\Res^G_H(g) \colon \bigoplus_{q \in Q}
r_q^{-1}.a \rightarrow \bigoplus_{q \in Q} r_q^{-1}.(g.a)
\quad \quad \quad  
\forall\; g \in G
\end{equation}
is the sum of morphisms 
$h_{g,q} \colon 
r_q^{-1}.a \rightarrow r_{g.q}^{-1}.(g.a)$
over all $q \in Q$. 
\end{Lemma}
\begin{proof}
Let $a \in \A$. We have
$ \Res^G_H(a) = \leftidx{_a}{\left(\A \rtimes G\right)} \in
\modd\text{-}\left(\A \rtimes H\right)$.
By \cite[Eq.(4.48)]{gyenge2021heisenberg}, this decomposes as 
$\leftidx{_a}{\left(\A \rtimes G\right)} \simeq \bigoplus_{g \in G}
\leftidx{_a}{\A}{_g}$.
Right action of $h \in H$ sends each $\leftidx{_a}{\A}{_g}$ to $\leftidx{_a}{\A}{_{gh}}$ by precomposing with $h$, so we can regroup this direct sum into a direct sum of $\left(\A \rtimes H\right)$-modules as
$$
\bigoplus_{g \in G} \leftidx{_a}{\A}{_g} \simeq 
\bigoplus_{q \in Q} \left(\bigoplus_{h \in H}
\leftidx{_a}{\A}{_{r_qh}}\right). $$
Since $\homm_{\A}(r_qh.(-),a) \simeq \homm_{\A}(h.(-), r_q^{-1}.a)$, we further have 
$$
\bigoplus_{q \in Q} \left(\bigoplus_{h \in H}
\leftidx{_a}{\A}{_{r_qh}}\right)
\simeq  
\bigoplus_{q \in Q} \left(\bigoplus_{h \in H}
\leftidx{_{r_q^{-1}.a}}{\A}{_h} \right)
 \simeq 
\bigoplus_{q \in Q} 
\leftidx{_{r_q^{-1}.a}}{\left(\A \rtimes H\right)}, 
$$
where the second isomorphism is by \cite[Eq.(4.48)]{gyenge2021heisenberg} again. 

This establishes the first assertion. The remaining assertions are established by chasing 
$$ \id_{\bigoplus r_q^{-1}.a}  = \sum 
\id_{r_q^{-1}.a} \in 
\bigoplus_{q \in Q} 
\leftidx{_{r_q^{-1}.a}}{\left(\A \rtimes H\right)} $$
through the isomorphisms above and the maps 
$\leftidx{_a}{\left(\A \rtimes G\right)}
\rightarrow  \leftidx{_b}{\left(\A \rtimes G\right)}$
and 
$\leftidx{_a}{\left(\A \rtimes G\right)}
\rightarrow  \leftidx{_{g.a}}{\left(\A \rtimes G\right)}$ given by postcompositions with $\alpha$ and 
$g$, respectively. 
\end{proof}

Let $\A$ be a small DG category. For any $n \geq 1$ the permutation group $S_n$ strongly acts on the DG category $\A^{\otimes n}$ by permuting factors of objects and of morphisms:  for any 
$a_1, \dots, a_n \in \A$ 
\begin{equation}
\sigma(a_1 \otimes \dots \otimes a_n) = a_{\sigma^{-1}(1)} \otimes 
\dots a_{\sigma^{-1}(1)}. 
\end{equation}
It is necessary to invert $\sigma$ to ensure that we get a left action 
of $S_n$ and thus an embedding of $S_n$ into the group of 
automorphisms of $\A$.  Similarly, 
given $\alpha_i \in \homm_\A(a_i, a'_i)$ for $1 \leq i \leq n$
we have 
\begin{equation}
\sigma(\alpha_1 \otimes \dots \otimes \alpha_n) 
= 
\alpha_{\sigma^{-1}(1)} 
\otimes \dots \otimes 
\alpha_{\sigma^{-1}(n)}. 
\end{equation}

\begin{Definition}
\label{defn-symmetric-power-of-a-dg-category}
Let $\A$ be a small DG category. For any $n > 0$ define 
the \em $n$-th symmetric power \rm $\sym^{n}\A$ of the enhanced triangulated category $\A$ to be the semi-direct product $\A \rtimes S_n$. 
\end{Definition}

By Lemma \ref{lemma-A-rtimes-G-modules-are-G-equiv-A-modules}, the underlying triangulated category
$D_c(\sym^{n} \A) = 
H^0 \hperf(\sym^{n} \A)$ of $\sym^{n}\A$ 
coincides with other definitions of $n$-th
symmetrical powers for $DG$ categories \cite[Section~2.2.7]{SymCat}.   

\subsection{Hochschild homology}
\label{section-hochschild-homology}

\begin{Definition}
Let $\A$ be a small DG category. Its \em Hochschild homology \rm is 
\begin{equation}
\hochhom_\bullet(\A) := H^\bullet(\A \ldertimes_{\AbimA} \A),
\end{equation}
\end{Definition}
where $\A \ldertimes_{\AbimA} \A \in D(\kk)$ and 
$\ldertimes_{\A\text{-}\A}$ is the derived functor of the DG functor
\begin{equation}
\label{eqn-bimodule-left-right-tensor}
\otimes_{\A\text{-}\A}\colon \AmodA \otimes_k \AmodA \rightarrow \modk 
\end{equation}
where we tensor the left $\A$-action with the right $\A$-action and vice
versa. More precisely, $\A$-$\A$-bimodules are, equivalently, right
$\Aopp \otimes_{\kk} \A$-modules or left $\A \otimes_\kk \Aopp$-modules . 
We can view them as both left and right $\Aopp \otimes_{\kk} \A$-modules 
via the canonical isomorphism
\begin{align*}
\A \otimes_\kk \Aopp & \xrightarrow{\sim} \Aopp \otimes_\kk \A,
\\
a \otimes b \quad & \mapsto \quad  b \otimes a
\end{align*}
and \eqref{eqn-bimodule-left-right-tensor} is the functor of tensoring 
over this module structure. Explicitly, it sends any pair $E, F \in
\AmodA$ to the complex of $\kk$-modules
\begin{equation}
E \otimes_{\AbimA} F \coloneqq E \otimes_\kk F /
\left\{ e \otimes a.f.b - b.e.a \otimes f \;\middle|\; \forall\; e \in E, f \in
F, \text{ and } a,b \in \A \right\}. 
\end{equation}

We compute $\ldertimes_{\A\text{-}\A}$ by taking an h-projective
resolution in either variable. Using the bar-resolution $\barA$ of the diagonal bimodule $\A$, see e.g. 
\cite[Section~2.11]{AnnoLogvinenko-BarCategoryOfModulesAndHomotopyAdjunctionForTensorFunctors},
we see that $\hochhom_\bullet(\A)$ are isomorphic to the cohomologies
of the convolution of the \em Hochschild complex \rm $\hochcx_\bullet(\A)$ 
over $\modk$:
\begin{scriptsize}
\begin{equation}
\label{eqn-hochschild-complex-of-a-category}
\dots
\rightarrow 
\bigoplus_{a,b,c \in \A} 
\homm^\bullet_\A(c,a)
\otimes_k
\homm^\bullet_\A(b,c) 
\otimes_k
\homm^\bullet_\A(a,b)
\rightarrow 
\bigoplus_{a,b \in \A} 
\homm^\bullet_\A(b,a) \otimes_k
\homm^\bullet_\A(a,b)
\rightarrow 
\bigoplus_{a \in \A} \homm^\bullet_\A(a,a), 
\end{equation}
\end{scriptsize}
with the differentials defined by 
\begin{align*}
\alpha_0 \otimes \alpha_1 \otimes \dots \otimes \alpha_n 
& \quad \mapsto \quad 
\sum_{i = 0}^{n-1} (-1)^i \alpha_0 \otimes \dots \otimes \alpha_{i} \alpha_{i+1}
\otimes \dots \otimes \alpha_n + \\ 
& \quad \quad \quad \quad \quad + (-1)^{n + |\alpha_n|(|\alpha_0| + \dots +
|\alpha_{n-1}|)} \alpha_n \alpha_0 \otimes \alpha_1 \otimes \dots
\otimes \alpha_{n-1}
\end{align*}
For example, $\alpha_0 \otimes \alpha_1 \otimes \alpha_2 \mapsto 
\alpha_0 \alpha_1 \otimes \alpha_2 - \alpha_0 \otimes \alpha_1 \alpha_2
+ (-1)^{|\alpha_3|(|\alpha_1||\alpha_2|)}\alpha_2 \alpha_0 \otimes \alpha_1$.

\begin{Example}
	Let $A$ be an associative algebra. Each term of 
	\eqref{eqn-hochschild-complex-of-a-category} is concentrated in degree 0, so the complex is naturally one of $\kk$-modules. 
	Since it is concentrated in non-positive degrees, the same holds for \(\hochhom_\bullet(A)\). In particular,
	\[
	\hochhom_0(A) = A / [A,A],
	\] 
	where $[A,A]$ is the subspace of commutators, not the ideal they generate.
	
	This recovers the classical Hochschild complex 
	\cite[\S{IX.4}]{CartanEilenberg-HomologicalAlgebra} 
	\cite{HochschildKostantRosenberg-DifferentialFormsOnRegularAffineAlgebras} 
	up to a minor adjustment. Originally, the $i$-th Hochschild homology $HH_i(A)$ was defined as $\tor^i_{\AbimA}(A,A)$, and 
	\eqref{eqn-hochschild-complex-of-a-category} was  
	a \em chain \rm complex with $A^{\otimes i}$ in degree $i$. In the DG setting, this grading is unnatural, so we follow Shklyarov 
	\cite{Shklyarov-HirzebruchRiemannRochTypeFormulaForDGAlgebras} 
	and others: $\hochhom_i(A)$ in the classical sense corresponds to $\hochhom_{-i}(A)$ here, consistent with $\tor^i_{\AbimA}(A,A) = H^{-i}(A \ldertimes_{\AbimA} A)$.

\begin{comment}
Let $A$ be an associative algebra. Then each term of 
\eqref{eqn-hochschild-complex-of-a-category} is also concentrated in degree $0$. 
Thus \eqref{eqn-hochschild-complex-of-a-category} is a complex of 
$\kk$-modules.
In particular, since \eqref{eqn-hochschild-complex-of-a-category} is
concentrated in non-positive degress, so is $\hochhom_\bullet(A)$. We
also have $\hochhom_0(A) = A / [A,A]$. This is a vector
space quotient and $[A,A]$ is the subspace of commutators, and not 
the ideal generated by them.  

This recovers the original definition of the Hochschild complex 
\cite[\S{IX.4}]{CartanEilenberg-HomologicalAlgebra}
\cite{HochschildKostantRosenberg-DifferentialFormsOnRegularAffineAlgebras} 
with a minor difference. Originally the $i$-th Hochschild homology 
$HH_i(A)$ was defined as $\tor^i_{\AbimA}(A,A)$ and 
\eqref{eqn-hochschild-complex-of-a-category} was  
a \em chain \rm complex with $A^{\otimes i}$ in degree $i$. In DG setup
this is unnatural, so we follow Shklyarov
\cite{Shklyarov-HirzebruchRiemannRochTypeFormulaForDGAlgebras} and
others in our present conventions. Thus, $\hochhom_{i}(A)$ in 
the original definition for associative algebras is $\hochhom_{-i}(A)$
in our conventions, matching the fact that $\tor^i_{\AbimA}(A,A)$
is $H^{-i}(A \ldertimes_{\AbimA} A)$. 
\end{comment}
\end{Example}

The key properties of Hochschild homology are:
\begin{itemize}
\item \em Self-opposite: \rm 
$\hochhom_\bullet(\A) \simeq \hochhom_\bullet(\Aopp)$.
\item \em Functoriality: \rm 
a DG functor $F\colon \A \rightarrow \B$ induces a map 
$F^{\hochhom}\colon \hochhom_\bullet(\A) \rightarrow
\hochhom_\bullet(\B)$. 
\item \em Homotopy invariance: \rm 
If $F\colon \A \rightarrow \B$ is a quasi-equivalence, 
then $F^{\hochhom}$ is an isomorphism
and if functors $F,G \colon \A \rightarrow \B$ are homotopy
equivalent, then $F^{\hochhom} = G^{\hochhom}$ 
\cite[Lemma 3.4]{Keller-OnTheCyclicHomologyOfExactCategories}. 
\item \em K{\"u}nneth formula: \rm 
\begin{equation}
\label{eqn-kunneth-formula-isomorhism}
\hochhom_{\bullet}(A \otimes_\kk B) \simeq \hochhom_{\bullet}(A) \otimes_\kk
\hochhom_{\bullet}(B). 
\end{equation}
\item \em Morita invariance: \rm 
the Yoneda embedding $\A \hookrightarrow \perfA$ into 
the DG category of perfect $\A$-modules induces an isomorphism 
$\hochhom_\bullet(\A) \simeq \hochhom_\bullet(\perfA)$. 
\end{itemize}

\subsection{Opposite category}
\label{section-hochschild-homology-of-the-opposite-category}

Let $\A$ be a DG category. The isomorphism $\hochhom_\bullet(\A)
\simeq \hochhom_\bullet(\Aopp)$ is induced by the isomorphism of Hochschild
complexes 
\begin{equation}
\label{eqn-category-to-opposite-category-hh-map}
\hochcx_\bullet(\A) 
\xrightarrow{\sim}
\hochcx_\bullet(\Aopp)
\end{equation}
defined for any 
$$ \alpha_0 \otimes \alpha_1 \otimes \dots \otimes \alpha_n 
\in \homm_\A(a_1, a_0) \otimes \homm_\A(a_2, a_1) \otimes \dots
\otimes \homm_\A(a_0, a_n) \in \hochcx_n(\A) $$
by
$\alpha_0 \otimes \alpha_1 \otimes \dots \otimes \alpha_n
\mapsto 
(-1)^{\sum_{i < j} \deg(\alpha_i)\deg(\alpha_j)}
\alpha_n \otimes \alpha_{n-1} \otimes \dots \otimes \alpha_0$.
Note that the image lies in 
\begin{align*}
\homm_\Aopp(a_n, a_0) \otimes \homm_\Aopp(a_{n-1}, a_{n}) \otimes 
\dots \otimes \homm_\Aopp(a_0, a_1) = \hochcx_n(\Aopp). 
\end{align*}

\subsection{Functoriality} 
\label{section-hochschild-homology-functoriality}

Let $F\colon \A \rightarrow \B$ be a DG functor between two DG categories. 
The map $F^\hochhom\colon \hochhom_\bullet(\A) \rightarrow
\hochhom_\bullet(\B)$ is induced by the closed degree zero map 
of Hochschild complexes 
\begin{equation}
\label{eqn-hochschild-homology-functoriality-map}
F\colon \hochcx_\bullet(\A) 
\rightarrow 
\hochcx_\bullet(\B)
\end{equation}
defined for any $\alpha_0 \otimes \alpha_1 \otimes \dots \otimes \alpha_n 
\in \hochcx_n(\A)$ 
by
$\alpha_0 \otimes \alpha_1 \otimes \dots \otimes \alpha_n
\mapsto 
F(\alpha_0) \otimes F(\alpha_1) \otimes \dots \otimes F(\alpha_n)$. 

\subsection{K{\"u}nneth isomorphism}
\label{section-kunneth-isomorphism}

Let $\A$ and $\B$ be two DG categories. The K{\"u}nneth formula isomorphism 
\eqref{eqn-kunneth-formula-isomorhism} is induced by the closed degree
zero map of Hochschild complexes 
\begin{equation}
\label{eqn-kunneth-map}
K\colon \hochcx_\bullet(\A) \otimes_{\kk}  \hochcx_\bullet(\B)
\xrightarrow{\sim}
\hochcx_\bullet(\A \otimes_{\kk} \B)
\end{equation}
defined using the shuffle product as follows. For any Hochschild chains
$$\alpha = \alpha_0 \otimes \dots \otimes \alpha_n \in \hochcx_n(\A)
\quad \text{ and } \beta = \beta_0 \otimes \dots \otimes \beta_m \in \hochcx_m(\A)$$
we have
\begin{small}
\begin{equation}
\label{eqn-explicit-formula-for-kunneth-map}
K\left(\alpha \otimes \beta\right)
= 
\sum_{\sigma \in S_{n,m}} (-1)^{\sigma}
(-1)^{\deg_\sigma(\alpha,\beta)}
(\alpha_0 \otimes \beta_0) \otimes \dots \otimes (\alpha_i \otimes \id)
\otimes \dots (\id \otimes \beta_j) \otimes \dots
\end{equation}
\end{small}
where in the summand indexed by the shuffle $\sigma \in S_{n,m}$ the
factors $(\alpha_i \otimes \id)$  and $(\id \otimes \beta_j)$ occur
in the positions $\sigma(i)$ and $\sigma(n+j)$, respectively. The 
second sign is computed by setting 
$d_\sigma(\alpha,\beta)$ to be the sum of 
$\deg(\alpha_i) \deg(\beta_j)$ for all $0 \leq i \leq n$ and 
$0 \leq j \leq m$ such that the factor containing $\alpha_i$ occurs to 
the right of the factor containing $\beta_j$. 

The K{\"u}nneth map \eqref{eqn-kunneth-map}
is a homotopy equivalence of twisted complexes
over $\modk$
\cite[\S2.4]{Shklyarov-HirzebruchRiemannRochTypeFormulaForDGAlgebras}.
In particular, its induced map on cohomologies is 
an isomorphism, establishing \eqref{eqn-kunneth-formula-isomorhism}.

\subsection{Hochschild homology and direct sums} 
\label{section-hochschild-homology-idempotents-direct-sums}

In this paper, we need to work with Hochschild chains of morphisms 
whose objects decompose as direct sums:

\begin{Definition}
\label{definition-single-component-and-continuous-chains}
Let $\A$ be a small DG category. For every $a \in \A$ fix a direct sum 
decomposition
$$ a \simeq s(a)_{1} \oplus \dots \oplus s(a)_{m(a)}, $$
which may be trivial, i.e. $m(a) = 1$ and $s(a)_1 = a$. 
For each $1 \leq i \leq m(a)$, let $\iota_i\colon s(a)_i \rightarrow a$
and $\pi_i\colon a \rightarrow s(a)_i$ be the inclusion and projection
morphisms. Let $p_i \colon a \rightarrow a$ 
be the idempotent $\iota_i \circ \pi_i$. 

Every $\homm$-complex of $\A$ decomposes as:
\begin{align*}
\homm_\A(a,b) & = \bigoplus_{i,j} p_i \homm_\A(a,b) p_j. 
\end{align*}
The elements of each $p_i \homm_\A(a,b) p_j$ are said to be
\em single component \rm morphisms as they each go from 
a single summand of $a$ to a single summand of $b$. 
Correspondingly, each term of the Hochschild complex $\hochcx_\bullet(\A)$ 
decomposes further into a direct sum of tensor products of 
$p_i \homm_\A(a,b) p_j$'s:
\begin{small}
\begin{equation}
\label{eqn-decomposition-into-single-component-chains}
\hochcx_{n}(\A) =  
\bigoplus_{
\begin{smallmatrix}
a_0, a_1, \dots, a_n \in \A \\
1 \leq s_j,t_j \leq m(a_j)
\end{smallmatrix}
}
p_{t_0} \homm^\bullet_\A(a_1,a_0) p_{s_1}
\otimes
p_{t_1} \homm^\bullet_\A(a_2,a_1) p_{s_2}
\otimes
\dots
\otimes
p_{t_n} \homm^\bullet_\A(a_0,a_n) p_{s_0}. 
\end{equation}
\end{small}
This a \em single component chain \rm decomposition of the terms of
$\hochcx_\bullet(\A)$ and the elements of the summands in
\eqref{eqn-decomposition-into-single-component-chains} are \em single
component chains \rm consisting of single component morphisms. 

For each summand in  
\eqref{eqn-decomposition-into-single-component-chains} we say 
that it is \em continuous \rm if $s_0 = t_0$, \dots, $s_n = t_n$
and \em discontinous \rm otherwise. Its elements are,
correspondingly, \em continuous \rm and \em discontinuous \rm single 
component chains. The Hochschild differential preserves continuity or
discontinuity of a chain and therefore we obtain a decomposition of
the Hochschild twisted complex $\hochcx_\bullet(\A)$ 
into a direct sum of twisted complexes
\begin{equation}
\label{eqn-decomposition-into-continuous-and-discontinuous-chains}
\hochcx_\bullet(\A) = \hochcx_\bullet(\A)_{\text{cts}} \oplus
\hochcx_\bullet(\A)_{\text{dis}}, 
\end{equation}
its continuous and discontinuous components.  

Finally, define a closed, degree zero twisted complex map 
\begin{equation}
\redcts\colon \hochcx_\bullet(\A) \rightarrow \hochcx_\bullet(\A) 
\end{equation}
by setting its action on each summand in 
\eqref{eqn-decomposition-into-single-component-chains} to send each 
$\alpha_0 \otimes \dots \otimes \alpha_n$ with 
$\alpha_i = p_{t_i} \alpha_i p_{s_{i+1}}$ to 
$$ \quad \pi_{t_0} \alpha_0 \iota_{t_1} \otimes \dots \otimes 
   \pi_{t_i} \alpha_i \iota_{t_{i+1}}
\otimes
\dots
\otimes 
\pi_{t_n} \alpha_n \iota_{t_0}. $$
\end{Definition}

In other words, the map $\redcts$ kills the discontinuous chains
and discards the redundant summands in the continuous chains. 
A continuous single component Hochshild chain is a chain 
of morphisms between objects going 
from a single direct summand to a single direct summand, and 
with each morphism starting at the direct summand where the
previous morphism ended. 
The map $\redcts$ discards all the summands not involved 
and reduces a chain to the chain of morphisms between the single
direct summands which are involved. 

\begin{Theorem}
\label{theorem-reduction-to-continuous-subchains-and-direct-summands}
Let $\A$ be a small DG category. For every $a \in \A$ fix its decomposition 
into direct summands as in Definition 
\ref{definition-single-component-and-continuous-chains}.
 
The map 
$\redcts\colon \hochcx_\bullet(\A) \rightarrow \hochcx_\bullet(\A)$
is homotopic to the identity map of twisted complexes. In particular, 
for any chain $\alpha \in \hochcx_\bullet(\A)$ which defines
a class $[\alpha] \in \hochhom_\bullet(\A)$ we have 
$[\alpha] = [\redcts(\alpha)]$. 
\end{Theorem}
\begin{proof}
The desired homotopy between $\id_{\hochcx_{\bullet}}$ and $\redcts$
is the degree $-1$ twisted complex map
$$ h: \hochcx_\bullet(\A) \rightarrow \hochcx_\bullet(\A) $$
whose action on each summand in 
\eqref{eqn-decomposition-into-single-component-chains} sends each 
$ \alpha_0 \otimes \dots \otimes \alpha_n$
with $\alpha_i = p_{t_i} \alpha_i p_{s_{i+1}}$ to
$$ 
\sum_{i=0}^n
(-1)^i
\pi_{t_0} \alpha_0 \otimes \alpha_1 \otimes
\dots
\otimes
\alpha_i \otimes 
\iota_{t_{i+1}}
\otimes 
\pi_{t_{i+1}} \alpha_{i+1} \iota_{t_{i+2}}
\otimes 
\dots
\otimes
\dots
\otimes \pi_{t_n} \alpha_n \iota_{s_0}.
$$
\end{proof}

\subsection{Hochschild homology, DG bicategories and enhanced
triangulated categories}
\label{section-hochschild-homology-bicategories}

Functoriality of Hochschild homology implies that we have a $1$-functor:
\begin{equation}
\label{eqn-1-functor-HH-from-dgcat-to-grvect}
\hochhom\colon \DGCatone \rightarrow \catgrVect,
\end{equation}
sending any small DG category $\A$ to $\hochhom_\bullet(\A)$ 
and any DG functor $\A \rightarrow \B$ to 
the induced map $\hochhom_\bullet(\A) \rightarrow \hochhom_\bullet(\B)$. 
As quasi-equivalences induce isomorphisms on
$\hochhom_\bullet$ we further have a $1$-functor 
\begin{equation}
\label{eqn-1-functor-HH-from-hodgcat-to-grvect}
\hochhom\colon \HoDGCatone \rightarrow \catgrVect.
\end{equation}
This and the K{\"u}nneth formula allow us to define:

\begin{Definition}
\label{defn-hochschild-homology-1-cat-of-a-dg-bicategory}
Let $\HoDGCat\text{-}{\bicat{Bicat}}^1_{hst}$ be 
the $1$-category of $\HoDGCat$-enriched bicategories and homotopy 
strong $2$-functors. Let $\grCatonek$ be the $1$-category of 
$\kk$-linear $1$-categories. Define the following functor
\begin{equation}
\hochhom_\bullet\colon 
\HoDGCat\text{-}{\bicat{Bicat}}^1_{hst} 
\rightarrow 
\grCatonek.
\end{equation}

For any $\HoDGCat$-enriched bicategory $\bicat M$ 
let $\hochhom_\bullet(\bicat M)$ be a graded $\kk$-linear $1$-category 
with the same objects as $\bicat M$ and morphism spaces 
$$ \homm_{\hochhom_\bullet(\bicat M)}(N,N') = 
\hochhom_\bullet\left(\homm_{\bicat M}(N,N')\right). $$
The composition in $\hochhom_\bullet(\bicat M)$ is induced by the
composition in $\bicat M$ as follows
\begin{align*}
& \hochhom_\bullet\left(\homm_{\bicat M}(N',N'')\right) \otimes
\hochhom_\bullet\left(\homm_{\bicat M}(N,N')\right) \xrightarrow{\sim} 
\\
\xrightarrow{\sim} \quad
&
\hochhom_\bullet\left(\homm_{\bicat M}(N',N'') \otimes \homm_{\bicat M}(N,N')\right) \xrightarrow{\text{1-composition in } \bicat M}
\\
\xrightarrow{\text{1-composition in } \bicat M} \quad
&
\hochhom_\bullet\left(\homm_{\bicat M}(N,N'')\right).
\end{align*}
Identity morphisms are given by the classes of identity $1$-morphisms in 
$\bicat M$. 

For any homotopy strong $2$-functor $F: \bicat M \rightarrow \bicat N$
let $\hochhom_\bullet(F)$ to be the functor $\hochhom_\bullet(M) 
\rightarrow \hochhom_\bullet(N)$ which on objects has the same
action as $F$ and on morphism spaces is  
\eqref{eqn-1-functor-HH-from-hodgcat-to-grvect} applied to $F$.
As $F$ is homotopy strong and as homotopy equivalent DG functors
induce the same map on Hochschild homology  
\cite[Lemma 3.4]{Keller-OnTheCyclicHomologyOfExactCategories}, 
$\hochhom_\bullet(F)$ does commute with the compositions in 
$\hochhom_\bullet(\bicat{M})$ and $\hochhom_\bullet(\bicat{N})$.
\end{Definition}

There is a natural $1$-functor 
$$ \DGCatone \rightarrow \hochhom_\bullet(\DGCatdg) $$
which is identity on objects and sends each DG functor 
$\A \rightarrow \B$ to its Euler character in $HH_0(\DGFun(\A,\B))$. The $1$-functor \eqref{eqn-1-functor-HH-from-dgcat-to-grvect} lifts to a $1$-functor
\begin{equation}
\label{eqn-1-functor-from-HH-dgcat-to-grvect}
\hochhom\colon \hochhom_\bullet (\DGCatdg) \rightarrow \catgrVect,
\end{equation}
which sends any small DG category $\A$ to $\hochhom_\bullet(\A)$ and whose
action on the morphism spaces
$$ \hochhom_{\bullet}(\DGFun(\A,\B)) \rightarrow
\homm_\kk(\hochhom_\bullet(\A), \hochhom_\bullet(\B)) $$
is adjoint to the composite map 
$$ \hochhom_{\bullet}\left(\DGFun(\A,\B)\right) 
\otimes_\kk \hochhom_\bullet(\A) 
\xrightarrow{K}
\hochhom_{\bullet}\left(\DGFun(\A,\B) \otimes_\kk \A\right) 
\xrightarrow{\hochhom_{\bullet}(\eval)}
\hochhom_\bullet(\B). $$

To construct a similar lift of \eqref{eqn-1-functor-HH-from-hodgcat-to-grvect}
it is best to restrict to the full subcategory $\HoDGCatone$ 
consisting of all DG categories of the form $\hperf(\A)$.
This subcategory is equivalent to the $1$-categorical truncation 
$\MoDGCatone$ of the strict $2$-category 
$\EnhCatKC$ of Morita enhanced triangulated categories, and 
$\EnhCatKC$ lifts to the DG bicategory $\EnhCatKCdg$ of 
Morita enhanced triangulated categories, 
cf.~\cite[\S4.4]{gyenge2021heisenberg}. The objects of $\EnhCatKCdg$
are small DG categories viewed as enhanced triangulated categories
and its $1$-morphisms categories are
$$ \homm_{\EnhCatKCdg}(\A,\B) := (\AmodbarB)_{\Bperf}, $$ 
where $\AmodbarB$ is the bar-category of $\A$-$\B$-bimodules, cf. 
\cite[\S3.2]{AnnoLogvinenko-BarCategoryOfModulesAndHomotopyAdjunctionForTensorFunctors}. 
We similarly have a $1$-functor
$$ \MoDGCatone \rightarrow \hochhom_\bullet(\EnhCatKCdg) $$
which is identity on objects and sends any Morita quasifunctor
$F\colon \A \rightarrow \B$ to 
the Euler class of any bimodule $M \in \AmodbarB$ representing it. 
This lifts to a $1$-functor 
\begin{equation}
\label{eqn-1-functor-from-HH-enhcatkc-to-grvect}
\hochhom\colon \hochhom_\bullet (\EnhCatKCdg) \rightarrow \catgrVect
\end{equation}
which sends any small DG category $\A$ to 
$$\hochhom_\bullet(\A) \simeq \hochhom_\bullet(\hperfA) \simeq 
\hochhom_\bullet(\barperfA)$$
and whose action on the morphism spaces 
$$ \hochhom_{\bullet}((\AmodbarB)_{\Bperf}) \rightarrow
\homm_\kk(\hochhom_\bullet(\A), \hochhom_\bullet(\B)) $$
is adjoint to the composite map
$$ \hochhom_{\bullet}(\barperfA) 
\otimes_\kk \hochhom_{\bullet}(\AmodbarB)_{\Bperf}  
\xrightarrow{K}
\hochhom_{\bullet}(\barperfA \otimes_\kk 
(\AmodbarB)_{\Bperf})
\xrightarrow{\hochhom_{\bullet}(\bartimes)}
\hochhom_\bullet(\barperfB). $$

\subsection{Hochschild homology with coefficients in a bimodule}
\label{section-twisted-hochschild-homology}

We need the following:
\begin{Definition}
Let $\A$ be a small DG-category and let $M \in \AmodA$. 
The \em Hochschild homology of $\A$ with coefficients in $M$ \rm is 
\begin{equation}
\hochhom_\bullet(\A;M) := H^\bullet(M \ldertimes_{\A\text{-}\A} \A).  
\end{equation}
\end{Definition}
See \cite[\S1.1.3]{Loday-CyclicHomology}
\cite[\S3, Step
4]{Baranovsky-OrbifoldCohomologyAsPeriodicCyclicHomology}\cite[\S3.2]{Nordstrom-FiniteGroupActionsOnDGCategoriesAndHochschildHomology}
for other variations of this notion. 

As before, using the bar-resolution $\barA$ of $\A$ 
we see that $\hochhom_\bullet(\A;M)$ are isomorphic to the cohomologies
of the convolution of the \em Hochschild complex 
$\hochcx_\bullet(\A;M)$ with coefficients in $M$\rm:
\begin{small}
\begin{equation}
\label{eqn-hochschild-complex-of-a-category-with-coefficients-in-bimod}
\dots
\rightarrow 
\bigoplus_{a,b,c \in \A} 
\leftidx{_a}{M}{_c}
\otimes_k
\homm^\bullet_\A(b,c) 
\otimes_k
\homm^\bullet_\A(a,b)
\rightarrow 
\bigoplus_{a,b \in \A} 
\leftidx{_a}{M}{_b}
\otimes_k
\homm^\bullet_\A(a,b)
\rightarrow 
\bigoplus_{a \in \A} 
\leftidx{_a}{M}{_a}
\end{equation}
\end{small}
with the differentials defined by 
\begin{align*}
m_0 \otimes \alpha_1 \otimes \dots \otimes \alpha_n 
& \quad \mapsto \quad 
\quad \quad \quad  m_0.\alpha_1 \otimes \dots
\otimes \alpha_{n-1}\; +
\\
& \quad \quad \quad \quad \quad + \sum_{i = 1}^{n-1} (-1)^i \; m_0 \otimes \dots \otimes \alpha_{i} \alpha_{i+1}
\otimes \dots \otimes \alpha_n\; + \\ 
& \quad \quad \quad \quad \quad + (-1)^{n + |\alpha_n|(|m_0| + \dots +
|\alpha_{n-1}|)} \alpha_n.m_0 \otimes \alpha_1 \otimes \dots
\otimes \alpha_{n-1}
\end{align*}

Let $F\colon \A \rightarrow \A$ be a DG functor. Define the bimodule
$\leftidx{_F}{\A} \in \AmodA$ by setting 
$$ \leftidx{_a}{\left(\leftidx{_F}{\A}\right)}{_b} := 
\homm_{\A}(b,Fa), $$
and letting $\A$ act naturally on the right and via $F$ on the left: 
$$ \alpha.(\beta).\gamma = F(\alpha) \beta \gamma \quad \quad \quad 
\forall \beta \in \homm_{\A}(b,Fa), \alpha \in \homm_{\A}(a,c), \gamma
\in \homm_{\A}(d,b). $$

\begin{Definition}
Let $\A$ be a small DG category and $F\colon \A \rightarrow \A$ be
a DG functor. The 
\em $F$-twisted Hochschild complex $\hochcx_\bullet(\A;F)$ \rm 
and the \em $F$-twisted Hochschild homology $\hochhom_\bullet(\A;F)$ \rm
are the Hochschild complex $\hochhom_\bullet(\A;\leftidx{_F}{\A})$
and the Hochschild homology $\hochhom_\bullet(\A;\leftidx{_F}{\A})$. 
\end{Definition}
\begin{Definition}
\label{defn-map-induced-by-functor-between-twisted-hochschild-homologies}
Let $\A, \B$ be small DG categories and
$F\colon \A \rightarrow \A$, $G\colon \B \rightarrow \B$ be
DG functors. Let $H\colon \A \rightarrow \B$ be a DG functor and 
$\eta\colon HF \rightarrow GH$ be a natural transformation. 
Define the map 
\begin{equation}
(H,\eta)\colon \hochhom_\bullet(\A;F) \rightarrow \hochhom(\B;G) 
\end{equation}
by setting
$$
\alpha_0 \otimes \alpha_1 \otimes \dots \otimes \alpha_n 
\mapsto 
(\eta \circ H\alpha_0) \otimes H\alpha_1 \otimes \dots \otimes
H\alpha_n.  
$$ 
We similarly denote by $(H,\eta)$ the induced map on Hochschild
homologies. If $HF = GH$ and $\eta = \id$, we write $H$ for the map
$(H,\eta)$. It is simply termwise application of $H$ to 
the chains of $\hochhom_\bullet(\A;F)$. 
\end{Definition}

\subsection{Hochschild homology and strong group actions}
\label{section-hochschild-homology-and-strong-group-actions}

We need the following observation: 
\begin{Lemma}
\label{lemma-group-acts-on-equivarant-category-via-hochschild-quasi-isos}
Let $\A$ be a DG category with a strong action of a group $G$. 
For any $g \in G$ the map 
$$ g\colon \hochhom_\bullet(\A \rtimes G) \rightarrow 
\hochhom_\bullet(\A \rtimes G) $$
induced by the autoequivalence $g$ of $\A \rtimes G$ of
Definition 
\ref{defn-group-element-as-autoequivalence-of-the-equivariant-category}
is the identity map. 
\end{Lemma}
\begin{proof}
This follows from  
\cite[Lemma 3.4]{Keller-OnTheCyclicHomologyOfExactCategories}
since the autoequivalence $g$ is isomorphic to the identity functor,
cf.~\S\ref{section-strong-group-actions-and-DG-cats}. 
\end{proof}

The notion of a functor twisted Hochschild homology detailed
in \S\ref{section-twisted-hochschild-homology} is
useful in the context of strong categorical group actions. 
Let $\A$ be a DG category and $G$ a finite group acting strongly on $\A$. 
For any $g,h \in G$ we have $hg = (hgh^{-1}) h$, and so by
Definition 
\ref{defn-map-induced-by-functor-between-twisted-hochschild-homologies}
we get an isomorphism 
\begin{equation}
\label{eqn-group-element-induced-map-on-group-element-twisted-hh}
h\colon \hochcx_\bullet(\A;g) \rightarrow
\hochcx_\bullet(\A;hgh^{-1}). 
\end{equation}

\begin{Definition}
\label{defn-g-twisted-hcx-to-the-regular-hcx-hh}
Let $\A$ be a DG category with a strong action of a group $G$
For any $g \in G$ define 
\begin{equation}
\label{eqn-g-twisted-hcx-to-the-regular-hcx-hh}
\xi_g\colon \hochcx_\bullet(\A;g) \rightarrow \hochcx_\bullet(\A \rtimes G),
\end{equation}
by setting
\begin{equation*}
\alpha_0 \otimes \dots \otimes \alpha_m
\mapsto 
(g^{-1}(\alpha_0), g^{-1}) \otimes (\alpha_1, \id) \otimes \dots \otimes
(\alpha_{m-1},\id) \otimes  (\alpha_m,\id). 
\end{equation*}
\end{Definition}

The maps $\xi_\bullet$ are compatible 
with the shuffle product $K$ detailed in Section
\ref{section-kunneth-isomorphism}:
\begin{Lemma}
\label{lemma-compatibility-of-the-map-xi-with-shuffle-product}
Let $\A$ and $\B$ be DG categories and $G$ and $H$ be  
groups acting strongly on them. Then $G \times H$ acts strongly
on $\A \times \B$, and for any $g \in G$ and $h \in H$
the following square commutes:
\begin{equation}
\begin{tikzcd}
\hochcx^\bullet(\A;g) 
\otimes 
\hochcx^\bullet(\B;h) 
\ar{r}{K}
\ar{d}{\xi_g \otimes \xi_h}
&
\hochcx^\bullet(\A \otimes \B;g \times h) 
\ar{d}{\xi_{g \times h}}
\\
\hochcx^\bullet(\A \rtimes G) 
\otimes 
\hochcx^\bullet(\B \rtimes H) 
\ar{r}{K}
&
\hochcx^\bullet\left((\A \otimes \B) \rtimes (G \times H)\right). 
\end{tikzcd}
\end{equation}
\end{Lemma}
\begin{proof}
A straightforward verification. The key point is that  
$\xi_g$ and $\xi_h$ insert $g$ and $h$ into the first element
of each basic chain, while $K$ tensors the first
elements of the two basic chains and makes this the first element
of each summand of their shuffle product. 
\end{proof}

The following unassuming technical fact lies at the heart of the proof
of Theorem \ref{theorem-heisenberg-relation-for-the-assignments-of-pi}:

\begin{Lemma}
\label{lemma-restriction-functor-composed-with-the-twisted-part}
Let $\A$ be a DG category with a strong action of a group $G$. 
Let $H \leq G$ be a subgroup and let $Q$, $\fix(g)$, $r_q$, and $h_{g,q}$ 
be as in Definition~\ref{defn-the-set-of-left-cosets-of-a-subgroup}. 
For every $g \in G$ the following diagram commutes:
\begin{equation}
\label{eqn-res-G-H-circ-xi_g-square-for-HH}
\begin{tikzcd}[column sep = 2cm]
\hochhom_\bullet(\A;g)
\ar{r}{\sum_{q \in \fix(g)} r^{-1}_q}
\ar{d}{\xi_g}
&   
\displaystyle
\bigoplus_{q \in \fix(g)} \hochhom_\bullet(\A;h_{g,q})
\ar{d}{\sum_{q \in \fix(g)} \xi_{h_{g,q}}}
\\
\hochhom_\bullet(\A \rtimes G)
\ar{r}{\Res^G_H}
&
\hochhom_\bullet(\A \rtimes H). 
\end{tikzcd}
\end{equation}
where $\xi_g$ and $\xi_{h_{g,q}}$ are morphisms
\eqref{eqn-g-twisted-hcx-to-the-regular-hcx-hh}. 
\end{Lemma}

\begin{proof}
Let $\underline{\alpha} \in \hochhom_\bullet(\A;g)$ and lift it to a chain 
$\underline{\alpha} \in \hochcx_\bullet(\A;g)$. Write
$\underline{\alpha}$ as a sum $\sum_i \underline{\alpha}_i$ of 
\begin{equation*}
\underline{\alpha}_i := \alpha_{i0} \otimes \dots \otimes
\alpha_{in_i},  
\quad \quad \quad 
\text{ where }
\alpha_{ij} \text{ is a morphism }
\begin{cases}
a_{i1} \rightarrow g.a_{i0}     & j = 0, \\
a_{i(j+1)} \rightarrow a_{ij}     & j = 1,\dots, n_i-1, \\
a_{i0} \rightarrow a_{in_i}     & j = n_i,
\end{cases}
\end{equation*}
for some objects $a_{ij} \in \A$. 
By Definition \ref{defn-g-twisted-hcx-to-the-regular-hcx-hh}, we have
\begin{equation}
\label{eqn-xi-g-of-alpha}
\xi_g\left(\underline{\alpha}_i\right)
= \left((\id,g^{-1}) \circ (\alpha_{i0}, \id)\right) \otimes (\alpha_{i1}, \id) 
\otimes \dots \otimes (\alpha_{in_i}, \id).  
\end{equation}

We now apply $Res^G_H$ to the chain on the RHS of \eqref{eqn-xi-g-of-alpha}. 
By Lemma \ref{lemma-descripton-of-res-G-to-H-functor}, we have
$$ \Res^G_H(a_{ij}) = \bigoplus_{q \in Q} r_{q}^{-1}.a_{ij},
\quad \quad \quad \text{ and } \quad \quad \quad  
\Res^G_H(\alpha_{ij}, \id) = \sum_{q \in Q} r_{q}^{-1}(\alpha_{ij}). $$
Finally, $\Res^G_H(\id, g^{-1})$ sends each $r_q^{-1}.(g.a_{i0})$, the 
summand indexed by $q \in Q$ in $\Res^G_H(g.a_{i0})$, to 
$r^{-1}_{g^{-1}.q}.(a_{i0})$, the summand indexed by $g^{-1}.q$ in 
$\Res^G_H(a_{i0})$, via the morphism $h^{-1}_{g,g^{-1}.q}$. 
By definition $g r_{g^{-1}.q} = r_{q} h_{g,g^{-1}.q}$, so 
$h^{-1}_{g,g^{-1}.q} = r^{-1}_{g^{-1}.q} g^{-1} r_q$. 
Writing $Q = \left\{q_1, \dots, q_N \right\}$, we can depict 
the way the morphisms of the chain $\Res^G_H(\xi_g(\underline{\alpha}_i))$ act
on the direct summands of its objects as follows 
\begin{equation*}
\begin{tikzcd}[column sep = 0.5em]
\hspace{3cm} & (q_1) & \oplus & (q_2) & \oplus & \dots & \oplus &
(q_{N}) & \hspace{3cm} \\
\ar[phantom]{u}{\sum_{q \in Q} h^{-1}_{g,g^{-1}.q}}
&   \ar{urrrr} & & \ar{urr} & & 
\dots \ar{ullll} \ar{ull} \ar{urr} & & \ar{ull} &
\ar[phantom]{u}{(q_i) \mapsto (g^{-1}.q_i)}
\\
\ar[phantom]{u}{\sum_{q \in Q} r^{-1}_q (\alpha_{i1})} 
& (q_1) \ar{u} & \oplus  & (q_2) \ar{u} & \oplus  &  \dots & \oplus & (q_{N}) \ar{u} & 
\\
\ar[phantom]{u}{\sum_{q \in Q} r^{-1}_q (\alpha_{i2})} \
& \dots \ar{u} & \dots  & \dots \ar{u} & \dots & \dots & \dots & \dots
\ar{u} & \\
\ar[phantom]{u}{\sum_{q \in Q} r^{-1}_q (\alpha_{i(n_i-2)})}
& (q_1) \ar{u} & \oplus & (q_2) \ar{u} & \oplus & \dots & \oplus &
(q_{N}) \ar{u} & \\
\ar[phantom]{u}{\sum_{q \in Q} r^{-1}_q (\alpha_{i(n_i-1)})}
& (q_1) \ar{u} & \oplus & (q_2) \ar{u} & \oplus & \dots & \oplus &
(q_{N}) \ar{u} &\\
\ar[phantom]{u}{\sum_{q \in Q} r^{-1}_q (\alpha_{in_i})}
& (q_1) \ar{u} & \oplus & (q_2) \ar{u} & \oplus & \dots & \oplus &
(q_{N}), \ar{u} &
\end{tikzcd}
\end{equation*}
Each $(q_i)$ represents the $q_i$-indexed summand of the corresponding 
direct sum and each arrow represents a non-zero component of the  
morphism between the corresponding direct sums.  

Decompose, as described in 
\S\ref{section-hochschild-homology-idempotents-direct-sums}, 
the chain $\Res^G_H(\xi_g(\underline{\alpha}_i))$ into a sum of single 
component chains of Definition~\ref{definition-single-component-and-continuous-chains}. By above, 
the continuous chains correspond 
bijectively to $\fix(g) \subseteq Q$, the elements fixed by $g$. 
For each $q \in \fix(g)$, the corresponding chain goes
through all $q$-indexed summands: 
\begin{equation}
\bigl((h^{-1}_{g,q},\id) \circ (r_q^{-1}(\alpha_{i1}), \id)\bigr) \otimes  
(r_q^{-1}(\alpha_{i2}), \id) \otimes \dots \otimes 
(r_q^{-1}(\alpha_{in_i}), \id). 
\end{equation}
This equals the image of $r^{-1}_{q}(\underline{\alpha}_i)$ under 
the map $\xi_{h_{g,q}}$. 

Repeating this for each
$\Res^G_H(\xi_g(\underline{\alpha}_i))$, 
we see that the continous part of $\Res^G_H(\xi_g(\underline{\alpha}))$ is
$$ \Res^G_H(\xi_g(\underline{\alpha}))_{\cts} = 
\sum_{i} \sum_{q \in \fix{g}} \xi_{h_{g,q}}r^{-1}_{q}(\underline{\alpha}_i)
= \sum_{q \in \fix{g}} \xi_{h_{g,q}}r^{-1}_{q}(\underline{\alpha}). $$
This is the image of $\underline{\alpha}$ under  
$(\sum_{q \in \fix(g)} \xi_{h_{g,q}}) \circ 
(\sum_{q \in \fix(g)} r^{-1}_q)$ going around the upper right half 
of \eqref{eqn-res-G-H-circ-xi_g-square-for-HH} on the level of
Hochschild complexes. 
By Theorem  
\ref{theorem-reduction-to-continuous-subchains-and-direct-summands}, 
the classes of $\Res^G_H(\xi_g(\underline{\alpha}))$
and $\Res^G_H(\xi_g(\underline{\alpha}))_{\cts}$ 
are equal in the Hochschild homology. It follows
that the square \eqref{eqn-res-G-H-circ-xi_g-square-for-HH} commutes. 
\end{proof}

\subsection{Hochschild-Kostant-Rosenberg isomorphism}
\label{section-hkr-isomorphism}

The famous HKR (Hochschild-Kostant-Rosenberg) theorem interprets the Hochschild
homology as a non-commutative analogue of the Hodge cohomology of a smooth 
algebraic variety. Its original, local version
\cite[Theorem~5.2]{HochschildKostantRosenberg-DifferentialFormsOnRegularAffineAlgebras} 
states that for a smooth and commutative associative algebra $A$ we have 
$$ HH_{-n}(A) \simeq \Omega^{n}_A, $$
where $\Omega^i_A$ is the $A$-module of $\kk$-linear K{\"a}hler
differential $i$-forms. Viewing $A$ as a smooth affine variety $X = \Spec A$, 
this means that 
$$ HH_{-n}(A) \simeq H^0_X(\Omega^{n}_X), $$
where $\Omega^i_X$ is the sheaf of differential $i$-forms on $X$. 
The global version of the HKR theorem 
\cite[Theorem 4.6.2]{Kontsevich-DeformationQuantizationOfPoissonManifolds} \cite[Corollary~2.6]{Swan-HochschildCohomologyOfQuasiprojectiveSchemes} \cite[Theorem 4.9]{Caldararu-TheMukaiPairingIITheHochschildKostantRosenbergIsomorphism}
states that for any 
DG algebra or category $A$ which Morita enhances a smooth, 
quasi-projective variety $X$ we have  
$$ HH_{-n}(A) \simeq \bigoplus_{i-j = n} H^i_X(\Omega^{j}_X). $$
By Morita enhancement we mean $A$ such that $D_c(A) \simeq D_c(X)$. 

Thus, the total Hochschild homology of $A$ is 
isomorphic as $\mathbb{Z}_2$-graded vector space to the total algebraic 
Hodge cohomology of $X$. Since the Hodge-to-de-Rham spectral sequence
degenerates in characteristic $0$, this is further isomorphic
to the total algebraic de Rham cohomology of $X$:
$$ HH_\bullet(A) \simeq H^{\bullet, \bullet}_{Hodge}(X) 
\simeq H^{\bullet}_{dR}(X). $$
When $\kk = \mathbb{C}$ by Poincaré Lemma this is further 
isomorphic to the topological cohomology $H^\bullet(X,\mathbb{C})$. 

\subsection{Euler pairing on Hochschild homology}
\label{section-euler-pairing-on-hochschild-homology}

If $A$ is proper, then its Hochschild homology carries a graded bilinear form 
$$ \chi: HH_{\bullet}(\A) \otimes_{\kk} HH_{\bullet}(\A) \to \kk, $$
called the \em Euler pairing\rm, see
\cite{Shklyarov-HirzebruchRiemannRochTypeFormulaForDGAlgebras}.
We briefly summarise its cosntruction. 

For a smooth and proper \dg category $\A$, the standard $\homm$-pairing on $\A$
\begin{align*}
        \A \otimes_{\kk} \Aopp \quad & \rightarrow \quad \modk \\
        (a,b) \quad & \mapsto \quad \homm_{\A}(a,b)
\end{align*}
restricts to 
$\A \otimes_\kk \Aopp \longrightarrow \hperf(k)$
and therefore induces a map 
$$ \hochhom_{\bullet}(\A \otimes_\kk \Aopp) \rightarrow
\hochhom_{\bullet}(\hperf(k))$$
which the canonical isomorphisms 
$$ \hochhom_{\bullet}(\A) \otimes_\kk \hochhom_{\bullet}(\A) \simeq  
\hochhom_{\bullet}(\A) \otimes_\kk \hochhom_{\bullet}(\Aopp) \simeq 
\hochhom_{\bullet}(\A \otimes_\kk \Aopp) $$
$$ \hochhom_{\bullet}(\hperf(\kk)) \simeq \hochhom_{\bullet}(\kk) \simeq \kk $$
turn into a $\kk$-bilinear pairing on $\hochhom_{\bullet}(\A)$
known as the \em Euler pairing\rm:
\begin{equation}
        \label{eqn-euler-pairing-on-hh_0-A}
        \chi^H: \hochhom_{\bullet}(\A) \otimes_\kk
\hochhom_{\bullet}(\A) \rightarrow \kk 
\end{equation}
For a smooth and proper $\A$
it is shown in \cite[Theorem 4]{Shklyarov-HirzebruchRiemannRochTypeFormulaForDGAlgebras} that this pairing is non-degenerate. 

\subsection{Euler character}

\begin{Definition}
Let $\A$ be a small DG category. For any $a \in \A$ the \em Euler
character \rm $\euler(a) \in \hochhom_0(\A)$ is the class of 
$\id_a \in \homm^\bullet_\A(a,a)$. 
\end{Definition}

$\hochhom_\bullet$ is stable under Morita equivalences
\cite{Keller-OnTheCyclicHomologyOfExactCategories}. 
Hence $\hochhom_\bullet(\A) \simeq \hochhom_\bullet(\hperfA)$, 
so any element of $\hperfA$ has a class in $\hochhom_0(\A)$. 
An explicit formula for this class was given 
in \cite[Theorem 1]{Shklyarov-HirzebruchRiemannRochTypeFormulaForDGAlgebras} 
and is as follows. Any $F \in \hperfA$ is homotopy equivalent
to a homotopy direct summand of an element of $\pretriag(\A)$. It
thus suffices, for any twisted complex $(a_i, \alpha_{ij}) \in
\pretriag(\A)$ and any homotopy idempotent 
$\pi\colon (a_i, \alpha_{ij}) \rightarrow (a_i, \alpha_{ij})$ 
to give a formula for the class of $\pi$ in $\hochhom_0(\A)$:
\begin{equation}
\euler(\pi) = \sum_{l = 0}^{\infty} (-1)^l \strace(\pi \otimes 
(\alpha_{\bullet \bullet})^{\otimes l}),
\end{equation}
where 
$\strace\colon \hochhom(\pretriag \A) \rightarrow \hochhom(\A)$
is the supertrace map
$$ \beta^1_{\bullet\bullet} \otimes \beta^2_{\bullet\bullet} \otimes \dots \otimes \beta^n_{\bullet\bullet} 
\quad \mapsto \quad \sum_{j \in \mathbb{Z}} (-1)^? \beta^1_{ji_1} \otimes \beta^2_{i_1i_2}
\otimes \dots \otimes \beta^n_{i_n j}, 
$$
where the sign is as described in
\cite[\S3.2]{Shklyarov-HirzebruchRiemannRochTypeFormulaForDGAlgebras}. 

Assignment of the Euler character gives a functorial map  
$$ \euler\colon\quad \hperfA \rightarrow \hochhom_0(\A) $$
which is a non-commutative analogue of the Chern character map
in the following sense. For any $\alpha\colon E \rightarrow F$ in $\hperfA$, 
$\cone(\alpha)$ is the convolution of the twisted complex 
$$ E \xrightarrow{\alpha} \underset{\degzero}{F}. $$
Applying the formula above to this twisted complex we see that
$$ \euler(\cone(\alpha)) = \euler(F) - \euler(E). $$
It follows that the Euler character map induces $\kk$-module morphism 
\begin{equation}
\label{eqn-euler-character-map-from-k(a)-to-hh_0(a)}
\euler\colon \quad \mathrm{K}_0(\A)  \rightarrow \hochhom_0(\A)
\end{equation}
where $\mathrm{K}_0(\A) : = \mathrm{K}_0(H^0(\hperf(\A))$. 
For a scheme $X$, set $\A$ to be the DG category of perfect
complexes of injective sheaves on $X$, then $\mathrm{K}_0(X) \simeq \mathrm{K}_0(\A)$. 
On a smooth variety $X$, the Hochschild-Kostant-Rosenberg theorem identifies 
$HH_0(X)$ with its Hodge cohomology $\bigoplus_{i} H^{i,i}(X)$. The map \eqref{eqn-euler-character-map-from-k(a)-to-hh_0(a)}
is then identified with the Chern character map 
$$ \mathrm{K}_0(X) \rightarrow \bigoplus_{i} H^{i,i}(X). $$ 

The Euler pairing \eqref{eqn-euler-pairing-on-hh_0-A} on $\hochhom_0(\A)$
was defined via the $\homm$-pairing on $\A$, so its composition with 
the Euler character map \eqref{eqn-euler-character-map-from-k(a)-to-hh_0(a)}
gives the usual Euler pairing on $\mathrm{K}_0(\A)$.
As the pairing on $\hochhom_0(\A)$ is non-degenerate, 
the map \eqref{eqn-euler-character-map-from-k(a)-to-hh_0(a)}
kills the kernel of the Euler pairing on $\mathrm{K}_0(\A)$. It 
therefore factors through the projection 
to the numerical Grothendieck group:
\begin{equation}
\label{eqn-euler-character-map-from-knum(a)-to-hh_0(a)}
\euler\colon \numGgp{\A} \rightarrow \hochhom_0(\A). 
\end{equation}

\subsection{Noncommutative orbifold decomposition for Hochschild
homology}
\label{section-noncommutative-orbifld-decomposition-for-hh}

Following the ideas of
\cite{FeiginTsygan-CyclicHomologyOfAlgebrasWithQuadraticRelationsUniversalEnvelopingAlgebrasAndGroupAlgebras}\cite{Brylinski-CyclicHomologyAndEquivariantTheories}\cite{GetzlerJones-TheCyclicHomologyOfCrossedProductAlgebras}, 
Baranovsky gave in \cite{Baranovsky-OrbifoldCohomologyAsPeriodicCyclicHomology}
a decomposition of the Hochschild homology of 
the orbifold stack $[X/G]$ where a finite
group $G$ acts on a smooth, quasi-projective variety $X$ with 
generically trivial stabilisers. In 
\cite{annobaranovskylogvinenko2023orbifold}, Anno, Baranovsky, and 
the second author gave a noncommutative version of this result
for symmetric quotient stacks. It serves as an important 
technical tool for this paper. 

Let $\A$ be any DG category. In 
\cite[Theorem 4.1]{annobaranovskylogvinenko2023orbifold}, there
were constructed two mutually inverse isomorphisms decomposing 
the Hochschild homology of $\sym^n\A$ into 
symmetric powers of $\hochhom_\bullet(\A)$:
\begin{equation}
\label{eqn-noncommutative-orbifold-decomposition-for-sn-v}
\begin{tikzcd}
\hochhom_\bullet(\sym^n\A)
\ar[shift left = 1.24ex]{r}{\zeta_n}[']{\simeq}
&
\bigoplus_{\underline{n} \vdash n}
\Sym^{\underline{r}(\underline{n})}
\hochhom_\bullet(\A)
\ar[shift left = 1.25ex]{l}{\eta_n},
\end{tikzcd}
\end{equation}
where we use the notation of 
\S\ref{section-fock-space-and-its-symmetric-algebra-structure}, 
e.g.~$\Sym^{\underline{r}(\underline{n})} \hochhom_\bullet(\A)
 := 
\Sym^{r_1(n)} \hochhom_\bullet(\A) \times \dots \times \Sym^{r_n(n)} 
\hochhom_\bullet(\A)$.
In \cite[Theorem 4.2]{annobaranovskylogvinenko2023orbifold},
it was further shown that $\zeta_n$ and $\eta_n$ combine
into isomorphisms identifying the total Hochschild homology 
of the symmetric powers of $\A$ with the symmetric algebra of 
\S\ref{section-fock-space-and-its-symmetric-algebra-structure}
which is canonically identified with the Fock space of
$\chalg{\hochhom_\bullet(\A)}$:
\begin{equation}
\label{eqn-symmetric-algebra-homotopy-equivalence-hh}
\begin{tikzcd}
\bigoplus_{n \geq 0} \hochhom_\bullet(\sym^n A)
\ar[shift left = 1.24ex]{r}{\zeta}[']{\simeq}
&
S^*(\hochhom_\bullet(\A) \otimes t \kk[t]). 
\ar[shift left = 1.25ex]{l}{\eta}
\end{tikzcd}
\end{equation}

The decomposition 
\eqref{eqn-noncommutative-orbifold-decomposition-for-sn-v}
has two steps. The first is a pair of quasi-inverse 
quasi-isomorphisms
\begin{equation}
\label{eqn-decomposition-of-orbifold-hh-into-the-sum-of-twisted-hhs}
\begin{tikzcd}
\hochcx_\bullet(\sym^{n}\A)
\ar[shift left = 1ex]{r}{\nu}
&
\left(\bigoplus_{\sigma \in S_n}
\hochcx_{\bullet}(\A^n;\sigma)\right)_{S_n} 
\ar[shift left = 1ex]{l}{\xi}
\end{tikzcd}
\end{equation}
where the subscript denotes the coinvariants under the 
action of $S_n$ by the maps
$$ \tau\colon 
\hochcx_{\bullet}(\A^n;\sigma) 
\xrightarrow{
\eqref{eqn-group-element-induced-map-on-group-element-twisted-hh}
}
\hochcx_{\bullet}\left(\A^n;\tau \sigma \tau^{-1}\right)
\quad \quad \quad \quad \tau \in S_n. $$
The map $\nu$ sends 
\begin{small}
\begin{equation*}
(\alpha_0, \sigma_0) \otimes \dots \otimes (\alpha_m, \sigma_m)
\mapsto 
(\sigma_0\dots\sigma_{m})^{-1}(\alpha_0) \otimes 
(\sigma_1\dots\sigma_{m})^{-1}(\alpha_1) \otimes
\dots \otimes \sigma_{m}^{-1} (\alpha_m)
\in \hochcx_{m}(\A^n;(\sigma_0 \dots \sigma_m)^{-1}). 
\end{equation*}
\end{small}
It is a straightforward generalisation of 
\cite[Proposition 3.5]{Baranovsky-OrbifoldCohomologyAsPeriodicCyclicHomology}
and it works identically for the strong action of any 
finite group $G$ on any DG category $\A$.
In the language of $G$-equivariant objects in $\hperf\A$, it appeared also in  
\cite[Theorem 4.3]{Nordstrom-FiniteGroupActionsOnDGCategoriesAndHochschildHomology}. 
The map $\xi$ is the composition 
$$
\left(\bigoplus_{\sigma \in S_n}
\hochcx_{\bullet}(\A^n;\sigma)\right)_{S_n} 
\xrightarrow{\frac{1}{|S_n|}\sum_{\tau \in S_n} \tau.(-)}
\bigoplus_{\sigma \in S_n}
\hochcx_{\bullet}(\A^n;\sigma)
\xrightarrow{\sum_{\sigma \in S_n} \xi_\sigma}
\hochcx_\bullet(\sym^{n}\A)
$$
where the maps $\xi_\sigma$ are as in Definition 
\ref{defn-g-twisted-hcx-to-the-regular-hcx-hh}. 

For each conjugacy class $\underline{n} \vdash n$ of $S_n$,
choose a representative $\sigma_{\underline{n}}$. 
Taking the coinvariants of the action of $S_n$ on $\bigoplus_{\sigma \in S_n}
\hochcx_{\bullet}(\A^n;\sigma)$ reduces it to
$\bigoplus_{\underline{n} \vdash n} \hochcx_{\bullet}(\A^n;\sigma_{\underline{n}})_{C(\sigma_{\underline{n}})}$, 
where $C(\sigma_{\underline{n}})$ is the centraliser of
${\underline{n}}$. The second step of the
\eqref{eqn-noncommutative-orbifold-decomposition-for-sn-v} is 
a pair of quasi-inverse quasi-isomorphisms
\begin{equation}
\label{eqn-reduction-of-the-twisted-hh-to-the-inertia-stack} 
\begin{tikzcd}
\hochcx_{\bullet}(\A^n;\sigma_{\underline{n}})_{C(\sigma_{\underline{n}})}
\ar[shift left = 1ex]{r}{f}
&
\hochcx_{\bullet}(\A^{r(\underline{n})})_
{S_{r_1(\underline{n})} \times \dots \times S_{r_n(\underline{n})}}. 
\ar[shift left = 1ex]{l}{g}
\end{tikzcd}
\end{equation}
These quasi-isomorphisms are the key result of 
\cite{annobaranovskylogvinenko2023orbifold}. Below, we give 
the explicit formulas for them in the specific case 
of $\underline{n} = (n)$ which we need in this paper. 
General formulas in 
\cite{annobaranovskylogvinenko2023orbifold}
are similar, but involve more cumbersome
notation. We introduce the notation we need.  Let 
$t = (1 \dots n) \in S_n$ and assume without loss 
of generality that $\sigma_{(n)} = t$. 
The long cycle $t$ acts on $\basecat^n$ by sending
$$ a_1 \otimes a_2 \otimes \dots \otimes a_n \rightarrow a_n \otimes a_1 \otimes
\dots a_{n-1} $$
on objects and correspondingly on the morphisms. 
The centraliser $C(t)$ is the cyclic subgroup generated by $t$, 
so for the partition $(n)$ we can rewrite  
\eqref{eqn-reduction-of-the-twisted-hh-to-the-inertia-stack} 
as
\begin{equation}
\label{eqn-reduction-of-the-twisted-hh-to-the-inertia-stack-long-cycle} 
\begin{tikzcd}
\hochcx_{\bullet}(\A^n;t)_{t}
\ar[shift left = 1ex]{r}{f}
&
\hochcx_{\bullet}(\A) 
\ar[shift left = 1ex]{l}{g}
\end{tikzcd}
\end{equation}

Let $m \geq 1$ and let 
$$ \underline{\alpha}_0 \otimes \dots \otimes
\underline{\alpha}_{m-1} \in \hochcx_{m-1}(\A^n;t). $$
By definition of $\hochcx_{m-1}(\A^n;t)$, for $0 \leq i \leq m-1$ 
the morphism $\underline{\alpha}_i$ can be postcomposed with 
$\underline{\alpha}_{i+1}$ in $\A^n$ and $t(\underline{\alpha}_m)$
can be postcomposed with $\underline{\alpha}_0$. 

Every chain in $\hochcx_{m-1}(\A^n;t)$ can be decomposed as
a sum of basic chains where
$$ \underline{\alpha}_i = \alpha_{(i+1)1} \otimes \dots \otimes
\alpha_{(i+1)n} \quad \quad \alpha_{(i+1)j} \in \A. $$
The index shift is because we want to use the matrix notation: 
the chain above can be visualised as 
\begin{small}
\begin{equation}
\begin{tikzcd}[row sep = 0.15cm, column sep = 0.15cm]
(\alpha_{11} 
& \otimes &
\alpha_{12}
& \otimes &
...
& \otimes &
\alpha_{1n})
\\
& & & \bigotimes & & & 
\\
(\alpha_{21} 
& \otimes &
\alpha_{22}
& \otimes &
...
& \otimes &
\alpha_{2n})
\\
& & & \bigotimes & & & 
\\
& & & \dots & & & 
\\
& & & \bigotimes & & & 
\\
(\alpha_{m1} 
& \otimes &
\alpha_{m2}
& \otimes &
...
& \otimes &
\alpha_{mn}). 
\end{tikzcd}
\end{equation}
\end{small}
and we denote it simply by $m \times n$ matrix
\begin{small}
\begin{equation}
\begin{pmatrix}
\alpha_{11} & \alpha_{12} & \dots & \alpha_{1n}  
\\
\alpha_{21} & \alpha_{22} & \dots & \alpha_{2n}  
\\
& & \dots & 
\\
\alpha_{m1} & \alpha_{m2} & \dots & \alpha_{mn}. 
\end{pmatrix}
\end{equation}
\end{small}
For any $1 \leq i \leq m$ and 
$1 \leq j \leq n$ the morphism $\alpha_{ij}$ can be postcomposed
with $\alpha_{(i+1)j}$. The morphism $\alpha_{mj}$ can be
postcomposed with $\alpha_{1(j+1)}$ for $1 \leq j \leq n-1$ and with 
$\alpha_{11}$ for $j = n$.  

Note that for $n = 1$ we simply have 
$\hochcx_{\bullet}(\A^n;t) = \hochcx_{\bullet}(\A)$,
so the chain
$\alpha_1 \otimes \dots \otimes \alpha_{m} \in \hochcx_{m-1}(\A)$
is represented in the above notation by the column vector
\begin{small}
\begin{equation}
\begin{pmatrix}
\alpha_{1} \\
\alpha_{2} \\
\dots
\\
\alpha_{m}. 
\end{pmatrix}
\end{equation}
\end{small}
We thus switch from denoting chains by 
$\alpha_0 \otimes \dots \otimes \alpha_{m-1}$ to denoting
them by $\alpha_1 \otimes \dots \otimes \alpha_{m}$. 
\begin{Definition}
The quasi-isomorphism $f$ in 
\eqref{eqn-reduction-of-the-twisted-hh-to-the-inertia-stack-long-cycle} 
is the map 
\begin{equation}
\label{eqn-quasi-iso-f-long-cycle}
\begin{pmatrix}
\alpha_{11} & \alpha_{12} & \dots & \alpha_{1n}  
\\
\alpha_{21} & \alpha_{22} & \dots & \alpha_{2n}  
\\
& & \dots & 
\\
\alpha_{m1} & \alpha_{m2} & \dots & \alpha_{mn}. 
\end{pmatrix} 
\quad \mapsto \quad 
\frac{1}{n}
\sum_{i = 1}^n
\begin{pmatrix}
\alpha_{1(i+1)} \dots \alpha_{m{i-1}}\alpha_{1i} \\
\alpha_{2i} \\
\dots
\\
\alpha_{mi} 
\end{pmatrix}
\end{equation}
\end{Definition}
\begin{Definition}
\label{defn-quasi-iso-g-in-baranovsky-decomposition-long-cycle} 
The quasi-isomorphism $g$ in 
\eqref{eqn-reduction-of-the-twisted-hh-to-the-inertia-stack-long-cycle} 
is the quotient of the map 
\begin{equation}
\label{eqn-quasi-iso-g-long-cycle}
g\colon \hochcx_\bullet(\A) \rightarrow \hochcx_\bullet(\A^n;t),
\end{equation}
given by
\begin{equation}
\label{eqn-quasi-iso-g-long-cycle-formula}
\begin{pmatrix}
\alpha_{1} \\
\alpha_{2} \\
\dots
\\
\alpha_{m} 
\end{pmatrix}
\quad \mapsto \quad 
\sum_{c \in \left\{1,\dots,n\right\}^{n}
\\
\text{ with } c_1 = 1
}
(-1)^{\sigma_{c}}
\begin{pmatrix}
\beta_{11} & \beta_{12} & \dots & \beta_{1n}  
\\
\beta_{21} & \beta_{22} & \dots & \beta_{2n}  
\\
& & \dots & 
\\
\beta_{m1} & \beta_{m2} & \dots & \beta_{mn}.
\end{pmatrix}
\end{equation}
Here: 
\begin{itemize}
\item The summation is taken over all $c = (c_1, \dots, c_n)$ with 
$c_i \in \left\{1, \dots, n\right\}$ and $c_1 = 1$.
We think of each $c_i$ as a choice of a column position in the $i$-th row 
of the matrix $(\beta_{ij})$ and we always choose 
the first position in the first row.
\item Define $\sigma_c \in S_m$ to send each $i$ to the number of $c_k$ with 
$c_k < c_i$ or with $c_k = c_i$ and $k \leq i$. 
\item Define $\beta_{ij}$ to be $\alpha_{\sigma_c(i)}$ if $j = c_i$
and the identity map otherwise. 
\end{itemize}
In other words, we take the matrix
$(\beta_{ij})$, start at $\beta_{11}$, and follow the composable
order of its entries (down each column and then to the top of 
the column to its right), placing $\alpha_1, \dots, \alpha_n$, 
in that order, in the chosen column position in each row. The remaining
entries are filled with the identity maps. We obtain a matrix 
each whose row contains exactly one $\alpha_i$. 
The permutation $\sigma_{c}$ corresponds to the 
new ordering on $\alpha_i$'s given by counting from the top to the bottom row. 
\end{Definition}
\begin{Example}
For $n = 1$ and $m = 3$ the quasi-isomorphisms $f$ and $g$ are both 
the identity map 
\begin{small}
\begin{equation}
\begin{pmatrix}
\alpha_1
\\
\alpha_2
\\
\alpha_3 
\end{pmatrix}
\mapsto 
\begin{pmatrix}
\alpha_1
\\
\alpha_2
\\
\alpha_3 
\end{pmatrix}.
\end{equation}
\end{small}

For $n = 2$ and $m = 3$ the quasi-isomorphism $g$ is the map 
\begin{small}
\begin{equation}
\begin{pmatrix}
\alpha_1
\\
\alpha_2
\\
\alpha_3 
\end{pmatrix}
\mapsto 
\begin{pmatrix}
\alpha_1 & \id
\\
\alpha_2 & \id
\\
\alpha_3 & \id
\end{pmatrix}
+
\begin{pmatrix}
\alpha_1 & \id
\\
\id & \alpha_2 
\\
\id & \alpha_3 
\end{pmatrix}
+
\begin{pmatrix}
\alpha_1 & \id
\\
\alpha_2 & \id
\\
\id & \alpha_3 
\end{pmatrix}
-
\begin{pmatrix}
\alpha_1 & \id
\\
\id & \alpha_3 
\\
\alpha_2 & \id
\end{pmatrix}, 
\end{equation}
\end{small}
while the quasi-isomorphism $f$ is the map
\begin{small}
\[
\begin{pmatrix}
	\alpha_{11} & \alpha_{12}
	\\
	\alpha_{21} & \alpha_{22}
	\\
	\alpha_{31} & \alpha_{32}
\end{pmatrix} \mapsto
\frac{1}{2}
\left(
\begin{pmatrix}
	\alpha_{12}\alpha_{22}\alpha_{32}\alpha_{11}
	\\
	\alpha_{21}
	\\
	\alpha_{31}
\end{pmatrix}
+
\begin{pmatrix}
	\alpha_{11}\alpha_{21}\alpha_{31}\alpha_{12}
	\\
	\alpha_{22}
	\\
	\alpha_{32}
\end{pmatrix}
\right).
\]
\end{small}

For $n = 3$ and $m = 3$ the quasi-isomorphism $g$ is the map 
\begin{small}
\begin{align*}
\begin{pmatrix}
\alpha_1
\\
\alpha_2
\\
\alpha_3 
\end{pmatrix}
\mapsto 
\quad &
\begin{pmatrix}
\alpha_1 & \id & \id
\\
\alpha_2 & \id & \id
\\
\alpha_3 & \id & \id
\end{pmatrix}
+
\begin{pmatrix}
\alpha_1 & \id & \id
\\
\alpha_2 & \id & \id
\\
\id  & \alpha_3 & \id
\end{pmatrix}
+
\begin{pmatrix}
\alpha_1 & \id & \id
\\
\alpha_2 & \id & \id
\\
\id  & \id & \alpha_3
\end{pmatrix}
-
\begin{pmatrix}
\alpha_1 & \id & \id
\\
\id & \alpha_3 & \id
\\
\alpha_2 & \id & \id 
\end{pmatrix}
+
\\
+ 
& 
\begin{pmatrix}
\alpha_1 & \id & \id
\\
\id & \alpha_2 & \id
\\
\id & \alpha_3 & \id 
\end{pmatrix}
+
\begin{pmatrix}
\alpha_1 & \id & \id
\\
\id & \alpha_2 & \id
\\
\id & \id & \alpha_3  
\end{pmatrix}
-
\begin{pmatrix}
\alpha_1 & \id & \id
\\
\id & \id & \alpha_3
\\
\alpha_2 & \id & \id
\end{pmatrix}
-
\begin{pmatrix}
\alpha_1 & \id & \id
\\
\id & \id & \alpha_3
\\
\id & \alpha_2 & \id
\end{pmatrix}
+
\begin{pmatrix}
\alpha_1 & \id & \id
\\
\id & \id & \alpha_2
\\
\id & \id & \alpha_3 
\end{pmatrix},
\end{align*}
\end{small}
whitle the quasi-isomorphism $f$ is the map 
\begin{small}
\[
\begin{pmatrix}
	\alpha_{11} & \alpha_{12} & \alpha_{13}
	\\
	\alpha_{21} & \alpha_{22} & \alpha_{23}
	\\
	\alpha_{31} & \alpha_{32} & \alpha_{33}
\end{pmatrix} \mapsto
\frac{1}{3}
\left(
\begin{pmatrix}
	\alpha_{12}\alpha_{22}\alpha_{32}\alpha_{13}\alpha_{23}\alpha_{33}\alpha_{11}
	\\
	\alpha_{21}
	\\
	\alpha_{31}
\end{pmatrix}
+
\begin{pmatrix}
	\alpha_{13}\alpha_{23}\alpha_{33}\alpha_{11}\alpha_{21}\alpha_{31}\alpha_{12}
	\\
	\alpha_{22}
	\\
	\alpha_{32}
\end{pmatrix}
+
\begin{pmatrix}
	\alpha_{11}\alpha_{21}\alpha_{31}\alpha_{12}\alpha_{22}\alpha_{32}\alpha_{13}
	\\
	\alpha_{23}
	\\
	\alpha_{33}
\end{pmatrix}
\right).
\]
\end{small}
\end{Example}

Map $g$ in \eqref{eqn-quasi-iso-g-long-cycle} is 
compatible with the shuffle product $K$ detailed in Section
\ref{section-kunneth-isomorphism}:
\begin{Lemma}
\label{lemma-compatibility-of-the-map-g-with-shuffle-product}
Let $\A$ and $\B$ be DG categories and let $n \geq 1$. The following
square commutes:
\begin{equation}
\label{eqn-compatibility-of-the-map-g-with-shuffle-product}
\begin{tikzcd}
\hochcx_\bullet(\A) \otimes \hochcx_\bullet(\B) 
\ar{r}{K}
\ar{d}{g \otimes g}
&
\hochcx_\bullet(\A \otimes \B)
\ar{d}{g}
\\
\hochcx_\bullet(\A^{n}; t_n) \otimes \hochcx_\bullet(\B^{n}; t_n) 
\ar{r}{K}
&
\hochcx_\bullet(\A^n \otimes \B^n; t_n \times t_n) 
\simeq 
\hochcx_\bullet((\A \otimes\B)^{n}; t_n),
\end{tikzcd}
\end{equation}
where $t_n \in S_n$ is the long cycle $(1 \dots n)$. 
\end{Lemma}
\begin{proof}
Let $\alpha \in \hochcx_\bullet(\A)$ and $\beta \in \hochcx_\bullet(\B)$ 
be two basic chains of lengths $p$ and $q$. Going around the lower left half of
\eqref{eqn-compatibility-of-the-map-g-with-shuffle-product} sends 
$\alpha \otimes \beta$ to a sum with the following summands. 
We choose one position in each non-first row of a $p \times n$ matrix 
and of a $q \times n$ matrix, place $\alpha_i$ and $\beta_i$ into 
these chosen positions as described in the definition of $g$, and set 
the remaining entries to be $\id$.  Finally, we choose a shuffle 
$\sigma \in S_{p-1,q-1}$ to shuffle the two matrices together
into a $(p+q-1) \times n$ matrix. 

The upper right half of
\eqref{eqn-compatibility-of-the-map-g-with-shuffle-product} 
yields a sum with the following summands. 
Choose a shuffle $\sigma \in S_{p-1, q-1}$ to shuffle $\alpha$ and
$\beta$ together. Choose one position in each non-first
row of a $(p+q-1) \times n$ matrix, place the elements of 
the shuffled chain into them, and set the remaining entries to $\id$. 

These two sets of choices are naturally bijective, and the summands
produced by corresponding choices are equal. We conclude that 
\eqref{eqn-compatibility-of-the-map-g-with-shuffle-product}
commutes on $\alpha \otimes \beta$. 
\end{proof}

\section{Hochschild homology decategorification of the Heisenberg 2-category}
\label{section-hochschild-homology-and-heisenberg-2-category}

In this section, we decategorify the Heisenberg categorification 
of \cite{gyenge2021heisenberg} via Hochschild homology. 

\subsection{Decategorification map for the Heisenberg $2$-category $\hcat\basecat$}
\label{section-decategorification-map}

Let $\basecat$ be a smooth and proper \dg category.
In \cite{gyenge2021heisenberg}, we constructed for any such $\basecat$:
\begin{enumerate}
\item the $\HoDGCat$-enriched bicategory $\hcat\basecat$, called 
the \em Heisenberg $2$-category \rm of $\basecat$
\cite[\S5]{gyenge2021heisenberg},
\item the $\HoDGCat$-enriched bicategory $\fcat\basecat$, called 
the \em categorical Fock space \rm of $\basecat$
\cite[\S7.2]{gyenge2021heisenberg},
\item the homotopy strong $2$-functor 
$\Phi_{\basecat}\colon \hcat\basecat \rightarrow \fcat\basecat$
\cite[\S7.3-7.5, Theorem 7.38]{gyenge2021heisenberg}. 
\end{enumerate}
As the first step in these constructions, we replace
$\basecat$ by $\hperf(\basecat)$ -- 
this does not change $D_c(\basecat)$. 

In \cite[\S8]{gyenge2021heisenberg} we decategorified these constructions using the numerical Grothendieck group $\numGgp{-}$. 
Here we decategorify them using the Hochschild homology $\hochhom_\bullet(-)$. 
This offers significant simplifications compared to working with 
the numerical Grothendieck groups. There, lack of $\homm$-finiteness of
$K_0(\hcat\basecat, \kk)$ prevented us from passing to 
$\numGgp{\hcat\basecat, \kk}$ by factoring out the kernel of the Euler pairing
\cite[\S6.4]{gyenge2021heisenberg}. 
Instead we had to artificially reproduce that by taking the
further factorisation of $K_0(\hcat\basecat, \kk)$
by the ad-hoc ideal described in \cite[Definition~6.26]{gyenge2021heisenberg}. 

Here, we decategorify by first applying $\hochhom_\bullet(-)$ to a DG bicategory
to get a graded $\kk$-linear $1$-category  as per Definition \ref{defn-hochschild-homology-1-cat-of-a-dg-bicategory}. 
We then flatten this $1$-category into an algebra as follows:
\begin{Definition}
\label{defn-unital-algebra-of-a-k-linear-1-category}
Let $\grCatonek$ be the $1$-category of graded $\kk$-linear $1$-categories
and let $\Algonek$ be the $1$-category of unital $\kk$-algebras. 
Define the flattening functor
\begin{equation}
\algcat\colon \grCatonek \rightarrow \Algonek. 
\end{equation}
as follows: for any $\kk$-linear $1$-category $\C$
$$ \algcat(\C) := 
\bigl\{ (f_{ab}) \in \prod_{a,b \in \C} \homm_{\C}(a,b) \;\bigm|\; 
\forall\; a\in \C \text{ we have } f_{ab} \neq 0 \text { for only finite number of } b \in \C
\bigr\}
$$
with the addition given by the $\kk$-linearity of $\C$, the
multiplication by the composition in $\C$, and the unit 
$1 = (\id_{aa})_{a \in \C}$. For any functor $F\colon \C \rightarrow \D$
betwen $\kk$-linear $1$-categories and any $(f_{ab}) \in \algcat(\C)$ 
$$ \algcat(F)((f_{ab})) := (F(f_{ab})). $$
\end{Definition}

The analogue of \cite[Cor.~8.8]{gyenge2021heisenberg} is now automatic. 
We define:
\begin{Definition}
\label{defn-hhalg-flattening-functor}
Let $\HoDGCat\text{-}{\bicat{Bicat}}^1_{hst}$ be 
the $1$-category of $\HoDGCat$-enriched bicategories and homotopy 
strong $2$-functors. Let $\grCatonek$ be the $1$-category of 
graded $\kk$-linear $1$-categories. Let $\Algonek$ be the $1$-category of
unital $\kk$-algebras. Define the functor
\begin{equation}
\label{eqn-hhalg-flattening-functor}
\hochhom_{alg}\colon \HoDGCat\text{-}{\bicat{Bicat}}^1_{hst}
\rightarrow \Algonek 
\end{equation}
to be the composition 
$$ 
\HoDGCat\text{-}{\bicat{Bicat}}^1_{hst} 
\xrightarrow{\hochhom_\bullet}
\grCatonek
\xrightarrow{\algcat}
\Algonek
$$
of the functor $\hochhom_\bullet$ of taking the Hochschild homology of a 
$\HoDGCat$-enriched bicategory as in
Definition~\ref{defn-hochschild-homology-1-cat-of-a-dg-bicategory}
and the functor $\algcat$ of flattening a $1$-category
into a unital algebra as in 
Definition~\ref{defn-unital-algebra-of-a-k-linear-1-category}. 
\end{Definition}

The $2$-categorical Heisenberg action 
$\Phi_{\basecat}\colon \hcat\basecat \rightarrow \fcat\basecat$ 
of \cite{gyenge2021heisenberg} induces a $1$-functor
$$ \hochhom_\bullet(\hcat\basecat) \rightarrow 
\hochhom_\bullet(\fcat\basecat).$$
By definition, $\fcat\basecat$ is the $1$-full subcategory of $\EnhCatKCdg$ supported at $\symbcn$ for $N \geq 0$.  We can therefore compose 
the above with the $1$-functor \eqref{eqn-1-functor-from-HH-enhcatkc-to-grvect} to get a $1$-functor
$$  
\hochhom_\bullet(\hcat\basecat) 
\rightarrow
\catgrVect
$$
whose image lies in the full subcategory of $\catgrVect$ supported at $\hochhom_\bullet(\symbcn)$ for $N \geq 0$. Hence, applying the flattening $\algcat$ of Definition \ref{defn-unital-algebra-of-a-k-linear-1-category} on this $1$-functor yields 
\begin{equation}
\label{eqn-decategorigication-of-phi-basecat}
\hochhom_{alg}(\Phi_{\basecat})\colon 
\hochhom_{alg}(\hcat\basecat) 
\rightarrow
\End \left( \bigoplus_{\ho} \hochhom_\bullet(\symbcn) \right).
\end{equation}

The main goal of this section is to construct an injective
decategorification map  
\begin{equation}
\label{eqn-injective-decategorification-map-into-the-flattening}
\chalg{\hochhom_\bullet(\basecat)} 
\rightarrow
\hochhom_{alg}(\hcat\basecat),  
\end{equation}
where $\chalg{\hochhom_\bullet(\basecat)}$ is
the Heisenberg algebra of $\hochhom_\bullet(\basecat)$
with the Euler pairing 
described in \S\ref{section-euler-pairing-on-hochschild-homology}. 

\subsection{Reduction to $\alghh_{\hcat\basecat}$ and idempotent modification}
\label{section-reduction-to-alghh}

Applying Definition~\ref{defn-unital-algebra-of-a-k-linear-1-category} to a
$1$-category $\C$ whose set of objects is $\mathbb{Z}$ we obtain a
$\kk$-algebra $\algcat(\C)$ which has a natural $\mathbb{Z}$-grading:
$$
\algcat(\C) = \bigoplus_{n \in \mathbb{Z}} \left( \prod_{i \in \mathbb{Z}} \homm_{\C}(i,i+n) \right). 
$$
We get a functor $\algcat$ from the category whose objects are 
the $1$-categories with the object set $\mathbb{Z}$ and whose morphisms 
are the functors which are identity on objects to the category of $\mathbb{Z}$-graded $\kk$-algebras. This admits a left adjoint:
\begin{Definition}
Let $A$ be a $\mathbb{Z}$-graded algebra. Define 
$\catalg(A)$ to be the $1$-category with:
\begin{itemize}
        \item the set of objects $\mathbb{Z}$, 
        \item $\homm_{\catalg(A)}(i,j) = (A)_{j-i}$,
        \item the composition given by multiplication in $A$,
        \item the identity morphisms given by $1 \in A_0$. 
\end{itemize}
\end{Definition}

The adjunction unit 
\begin{equation}
\label{eqn-adjunction-unit-cat-alg}
A \rightarrow \algcat(\catalg(A))
\end{equation}
sends any $a \in A_n$ for $n \in \mathbb{Z}$ to $(\alpha_{ij})
\in \prod_{i,j \in \mathbb{Z}} A_{i-j}$ with $\alpha_{ij} = a$ if
$j - i = n$ and $\alpha_{ij} = 0$ otherwise. The adjunction counit
\begin{equation}
\label{eqn-adjunction-counit-cat-alg}
\catalg(\algcat(\C)) \rightarrow \C  
\end{equation}
is the functor which on objects is the identity map and 
on morphisms is the map 
$$ \prod_{k \in \mathbb{Z}} \homm_{C}(k,k+j-i)
\rightarrow \homm_{\C}(i,j) $$
which is the projection to the $k = i$ factor. 

\begin{Definition}
Let $\basecat$ be any smooth and proper DG category. 
Define the $\mathbb{Z}$-graded $\kk$-algebra
$$ \alghh_{\hcat\basecat} := \bigoplus_{n \in \mathbb{Z}}
\hochhom_\bullet(\homm_{\hcat\basecat}(0,n)) $$
whose multiplication 
$$ 
\hochhom_\bullet(\homm_{\hcat\basecat}(0,n))
\otimes_\kk
\hochhom_\bullet(\homm_{\hcat\basecat}(0,m))
\rightarrow 
\hochhom_\bullet(\homm_{\hcat\basecat}(0,n+m))
$$
is the corresponding composition map in
$\hochhom_\bullet(\hcat\basecat)$
$$ 
\hochhom_\bullet(\homm_{\hcat\basecat}(0,n))
\otimes_\kk
\hochhom_\bullet(\homm_{\hcat\basecat}(n,n+m))
\rightarrow 
\hochhom_\bullet(\homm_{\hcat\basecat}(0,n+m))
$$
and whose identity element is the Euler class of the identity $1$-morphism 
$\mathbb{1}_0 \in \homm_{\hcat\basecat}(0,0)$. 
\end{Definition}

\begin{Lemma}  
$\hochhom_\bullet(\basecat)$ is isomorphic to $\catalg\left(\alghh_{\hcat\basecat}\right)$. 
\end{Lemma}
\begin{proof}
We observe that, by construction, we have
        $$ \homm_{\hcat\basecat}(i,j) = \homm_{\hcat\basecat}(0,j-i)
        \quad \quad \quad \forall\; i,j \in \mathbb{Z} $$
        for $1$-morphisms categories of $\hcat\basecat$.
\end{proof}

Define a $\mathbb{Z}$-grading on $\chalg{\hochhom_\bullet(\basecat)}$ 
by setting $\deg(a_{\alpha}(n)) = n$ for all $n \in \mathbb{Z}
\setminus \{ 0 \}$ and $\alpha \in \hochhom_\bullet(\basecat)$. 
By adjunction of $\catalg$ and $\algcat$, constructing 
\eqref{eqn-injective-decategorification-map-into-the-flattening}
is equivalent to constructing a $1$-functor
\begin{equation}
        \label{eqn-decategorification-1-functor}
        \catalg(\chalg{\hochhom_\bullet(\basecat)})
        \rightarrow
        \hochhom_\bullet(\hcat\basecat). 
\end{equation}

Thus to construct \eqref{eqn-decategorification-1-functor},
and hence 
\eqref{eqn-injective-decategorification-map-into-the-flattening},
it suffices to construct an algebra homomorphism 
\begin{equation}
\label{eqn-decategorification-map-pi}
\pi \colon 
\chalg{\hochhom_\bullet(\basecat)} \rightarrow \alghh_{\hcat\basecat}. 
\end{equation}
which we also call the decategorification map.  
In the rest of this section, we proceed to construct $\pi$. 
Once constructed, the homomorphism 
\eqref{eqn-injective-decategorification-map-into-the-flattening}
is obtained from $\pi$ as the composition 
\begin{equation}
\label{eqn-injective-decategorification-map-into-the-flattening-via-pi}
\chalg{\hochhom_\bullet(\basecat)}
\xrightarrow{\text{adj.~unit}}
\algcat(\catalg(\chalg{\hochhom_\bullet(\basecat)}))
\xrightarrow{\algcat(\catalg(\pi))}
\algcat(\catalg(\alghh_{\hcat\basecat}))
\simeq \hochhom_{alg}(\hcat\basecat).
\end{equation}
Thus \eqref{eqn-injective-decategorification-map-into-the-flattening} 
is injective if and only if $\pi$ is injective. In 
\S\ref{section-injectivity-of-the-decategorification-map-pi} 
we show that $\pi$ is injective, and the main conjecture of this paper
is that $\pi$ is also surjective. 

On the other hand, 
\eqref{eqn-injective-decategorification-map-into-the-flattening} 
sends $a_{\alpha}(n) \in \chalg{\hochhom_\bullet(\basecat)}$ to 
the element of $\prod_{i \in \mathbb{Z}}
\hochhom_\bullet(\homm_{\hcat\basecat}(i,i+n))$ consisting of 
$\pi(a_{\alpha}(n))$ in every factor. This is clearly never surjective –-
in $\hochhom_{alg}(\hcat\basecat)$ we have $1 = \sum_{n \in
\mathbb{Z}} 1_n$ where $1_n$ are the orthogonal idempotents 
given by the Euler classes of the identity $1$-morphisms
$\mathbb{1}_n \in \homm_{\hcat\basecat}(0,n)$. These idempotents 
never lie in the image of  
\eqref{eqn-injective-decategorification-map-into-the-flattening}. 

This led us in \cite{gyenge2021heisenberg}
to consider the idempotent-modified non-unital version
$\halg{\numGgp\basecat}$ of the Heisenberg algebra
$\chalg{\numGgp\basecat}$.  
In present terms, the idempotent modified $\halg{\hochhom_\bullet(\basecat)}$ 
is the subalgebra
$$ 
\halg{\hochhom_\bullet(\basecat)} :=
\bigoplus_{i,j \in \mathbb{Z}} (\chalg{\hochhom_\bullet(\basecat)})_{j-i} 
\subset \algcat(\catalg(\chalg{\hochhom_\bullet(\basecat)})). 
$$
It is a flattening of the 
$1$-category $\catalg(\chalg{\hochhom_\bullet(\basecat)})$ different
from $\algcat$: we take a direct sum of the $\homm$-spaces 
instead of the direct product with a finiteness condition. The
resulting algebra is not unital as the identity element would
have to be an infinite sum of the identity elements of
$\homm_{\catalg(\chalg{\hochhom_\bullet(\basecat)})}(n,n)$ for all $n
\in \mathbb{Z}$ and this lies outside the direct sum flattening
$\halg{\hochhom_\bullet(\basecat)}$. More generally,
$\halg{\hochhom_\bullet(\basecat)}$ doesn't naturally contain the 
original algebra $\chalg{\hochhom_\bullet(\basecat)}$. 

Enlarging the direct sum flattening to the $\algcat$-flattening
solves this. The adjunction unit embeds 
$\chalg{\hochhom_\bullet(\basecat)}$ as a subalgebra into  
$\algcat(\catalg(\chalg{\hochhom_\bullet(\basecat)}))$.
By above, $\algcat(\catalg(\chalg{\hochhom_\bullet(\basecat)}))$ 
contains $\halg{\hochhom_\bullet(\basecat)}$. 
Constructing the decategorification map $\pi$ as in
$\eqref{eqn-decategorification-map-pi}$, and thus the
$\algcat$-flattening homomorphism 
\begin{equation}
\label{eqn-algcat-flattened-decategorification-map}
\algcat(\catalg(\chalg{\hochhom_\bullet(\basecat)}))
\xrightarrow{\algcat(\catalg(\pi))}
\hochhom_{alg}(\hcat\basecat), 
\end{equation}
we construct homomorphisms from both  
$\chalg{\hochhom_\bullet(\basecat)}$ and its
idempotent-modified version $\halg{\hochhom_\bullet(\basecat)}$ into 
$\hochhom_{alg}(\hcat\basecat)$. Neither would ever be surjective, 
however \eqref{eqn-algcat-flattened-decategorification-map} is
injective if and only $\pi$ is. 

\subsection{Construction of the map $\pi$}
\label{section-overview-of-decategorification-map-pi}

Our Hochschild homology decategorification map 
\begin{equation*}
\pi \colon 
\chalg{\hochhom_\bullet(\basecat)} \rightarrow \alghh_{\hcat\basecat}. 
\end{equation*}
must agree with the numerical Grothendieck group decategorification
map $\pi$ \cite[\S6.1-6.4]{gyenge2021heisenberg}. When it is necessary
to differentiate between the two, we write
$\pi_\hochhom$ for the former and $\pi_{K}$ for the latter. 
The two maps must agree as follows: 
for any $a \in \numGgp\basecat$ and $n
\geq 1$ the Euler classes of the images of $p^{(n)}_{a}$ and
$q^{(n)}_{a}$ under $\pi_{K}$ must coincide with the images of
$p^{(n)}_{\euler(a)}$ and $q^{(n)}_{\euler(a)}$ under $\pi_\hochhom$:
\begin{equation}
\label{eqn-the-agreement-of-two-decategorification-maps}
\begin{tikzcd}
\numGgp{\basecat} 
\ar{d}{\euler}
\ar[shift left = 1ex]{r}{p_{\bullet}^{(n)}}
\ar[shift right = 1ex]{r}[']{q_{\bullet}^{(n)}}
&
\chalg{\numGgp{\basecat}}
\ar{r}{\pi_K}
&
\displaystyle
\bigoplus_{n \in \mathbb{Z}}
\numGgp{\homm_{\hcat\basecat}(0,n)}
\ar{d}{\euler}
\\
\hochhom_\bullet(\basecat)
\ar[shift left = 1ex]{r}{p_{\bullet}^{(n)}}
\ar[shift right = 1ex]{r}[']{q_{\bullet}^{(n)}}
&
\chalg{\hochhom_\bullet(\basecat)}
\ar{r}{\pi_{\hochhom}}
&
\displaystyle
\bigoplus_{n \in \mathbb{Z}}
\hochhom_\bullet(\homm_{\hcat\basecat}(0,n)). 
\end{tikzcd}
\end{equation}

Recall the construction of the map $\pi_{K}$ from 
\cite[\S6.1-6.4]{gyenge2021heisenberg}. For each $n \geq 1$ and each 
$a \in \basecat$, we assign to $a$ a class 
$\psi_n(a) \in \numGgp{\hperf(\symbc{n})}$ which is the class
of the twisted complex which homotopy splits the symmetrising idempotent 
$$ e := \frac{1}{n!}\sum_{\sigma \in S_{n}} \sigma
\quad\in\quad \homm_{\symbc{n}}\left(a^{\otimes n}, a^{\otimes
n}\right). 
$$

We then use the homomorphisms
defined by the $2$-functors 
$\Xi^{\PP}, \Xi^{\QQ}\colon \bihperf(\bicat{Sym}_\basecat) \to \hcat\basecat$
constructed in \cite[\S6.1-6.4]{gyenge2021heisenberg} to obtain 
the classes $\PP^{(n)}_{a}$
and $\QQ^{(n)}_{a}$:
\begin{equation}
\begin{tikzcd}[row sep = 0.5cm]
\label{eqn-k0num-assignment-of-pi(pn)-and-pi(qn)}
a \in \numGgp{\basecat} 
\ar[dashed]{r}{\psi_n}
\ar{d}{\sim}
&
\numGgp{\symbc{n}}
\ar{r}{\Xi^{\PP}}
&
\numGgp{\homm_{\hcat\basecat}(0,n)} \ni \PP^{(n)}_{a}
\\
\quad\quad\quad\numGgp{\basecat^{\opp}} 
\ar[dashed]{r}{\psi_n}
& 
\numGgp{\sym^{n} \basecat^{\opp}}
\ar{r}{\Xi^{\QQ}}
& \numGgp{\homm_{\hcat\basecat}(0,-n)} \ni \QQ^{(n)}_{a}.
\end{tikzcd}
\end{equation}
We choose objects $\left\{a_1, \dots, a_n\right\} \subset \basecat$ 
whose classes give a basis of $\numGgp\basecat$. We can do this
as we replaced $\basecat$ by $\hperf\basecat$ at the outset. 
The algebra $\chalg{\numGgp\basecat}$ is generated by $p_{a_i}^{(n)}$ and
$q_{a_i}^{(n)}$
subject to the relations of \cite[\S2.2]{gyenge2021heisenberg}. 
In \cite[Theorem 6.3]{gyenge2021heisenberg} 
we verify that these hold for $\PP^{(n)}_{a}$ and $\QQ^{(n)}_{a}$ 
and thus \eqref{eqn-k0num-assignment-of-pi(pn)-and-pi(qn)} 
extends to an algebra homomorphism 
$\pi_{K}\colon \chalg{\numGgp\basecat} \rightarrow 
\bigoplus_{n \in \mathbb{Z}} \numGgp{\homm_{\hcat\basecat}(0,n)}$. 

There is no hope of extending $\psi_n$ in 
\eqref{eqn-k0num-assignment-of-pi(pn)-and-pi(qn)} from a set-theoretic 
assignment defined on the classes of objects of $\basecat$ to an additive 
map: in the Heisenberg algebra $\chalg{\numGgp\basecat}$
the parametrisation of the generators $p^{(n)}_v$ and $q^{(n)}_v$
is not additive in $v \in \numGgp\basecat$. The parametrisation 
is additive for $a_\bullet(n)$ and $a_\bullet(-n)$ generators, 
but for them the definition of $\psi_n$ would be very complicated. 

This is another aspect where working with the Hochschild homology
offers a simplification: working with $a_\bullet(n)$ and $a_\bullet(-n)$ 
generators is easy enough. To construct the decategorification map $\pi$, 
we construct below in 
\S\ref{section-maps-psi_n-and-kappa_n} for any $n \geq 1$ $\kk$-linear maps
\begin{equation}
\psi_n\colon \hochhom_\bullet(\basecat) \rightarrow
\hochhom_\bullet(\symbc{n})
\quad \text{ and } \quad 
\psi_n\colon \hochhom_\bullet(\basecat^{\opp}) \rightarrow
\hochhom_\bullet(\sym^{n}\basecat^{\opp}),
\end{equation}
We then define the following assignment:
\begin{Definition}
For any $\alpha \in \hochhom_\bullet(\basecat)$ and $n \geq 1$,
we use the compositions 
\begin{equation}
\label{eqn-hochhom-assignment-of-pi(a(n))-pi(a(-n))}
\begin{tikzcd}[row sep = 0.5cm]
\hochhom_\bullet(\basecat) 
\ar{r}{\psi_n}
\ar{d}{\eqref{eqn-category-to-opposite-category-hh-map}}[']{\sim}
&
\hochhom_\bullet(\symbc{n})
\ar{r}{\Xi^{\PP}}
&
\hochhom_\bullet(\homm_{\hcat\basecat}(0,n)). 
\\
\hochhom_\bullet(\basecat^{\opp}) 
\ar{r}{\psi_n}
& 
\hochhom_\bullet(\sym^{n} \basecat^{\opp})
\ar{r}{\Xi^{\QQ}}
& 
\hochhom_\bullet(\homm_{\hcat\basecat}(0,-n)) 
\end{tikzcd}
\end{equation}
and the isomorphism \eqref{eqn-category-to-opposite-category-hh-map} 
to define 
$$ \hAA_\alpha(n) := \Xi^{\PP}\left(\psi_n(\alpha)\right) \in 
\hochhom_\bullet\left(\homm_{\hcat\basecat}(0,n)\right), $$
$$ \hAA_\alpha(-n) := \Xi^{\QQ}\left(\psi_n(\alpha)\right) \in 
\hochhom_\bullet\left(\homm_{\hcat\basecat}(0,-n)\right), $$
\end{Definition}

We then have:
\begin{Theorem}
\label{theorem-construction-of-the-decategorification-map}
There exists the unique algebra homomorphism
$$\pi_{\hochhom} \colon 
\chalg{\hochhom_\bullet(\basecat)} \rightarrow \alghh_{\hcat\basecat}$$
which sends $a_\alpha(n)$ and
$a_\alpha(-n)$ to $\hAA_\alpha(n)$ and $\hAA_\alpha(-n)$ for any 
$\alpha \in \hochhom_\bullet(\basecat)$ and $n \geq 1$. 
\end{Theorem}
As before, we write simply $\pi$ for $\pi_{\hochhom}$ where no
confusion is possible. 

\begin{proof}
By Definition~\ref{defn-the-heisenberg-algebra-of-a-graded-vector-space-a-gen}, 
$\chalg{\hochhom_\bullet(\basecat)}$ is the unital $\kk$-algebra generated by 
the elements $a_\alpha(n)$ and $a_\alpha(-n)$ subject to the linearity
relations 
\eqref{eq:vectheisrel2add-a-gen-graded},
\eqref{eq:vectheisrel2scalmult-a-gen-graded}, the commutation relation
\eqref{eq:vectheisrel1-a-gen-graded}, and the Heisenberg relation
\eqref{eq:vectheisrel3-a-gen-graded}. For the elements, 
$\hAA_\alpha(n)$ and $\hAA_\alpha(-n)$ the linearity relations
hold by the linearity of \eqref{eqn-hochhom-assignment-of-pi(a(n))-pi(a(-n))}
and in \S\ref{section-commutation-and-heisenberg-relations}
we prove that the commutation and the Heisenberg relations hold as well. 
\end{proof}

\subsection{Linear maps $\psi_n$ and $\kappa_n$}
\label{section-maps-psi_n-and-kappa_n}

In this section, we must construct $\kk$-linear maps
\begin{equation*}
\psi_n\colon \hochhom_\bullet(\basecat) \rightarrow \hochhom_\bullet(\symbc{n})
\end{equation*}
\begin{equation*}
\psi_n\colon \hochhom_\bullet(\basecat^{\opp}) \rightarrow 
\hochhom_\bullet(\sym^{n} \basecat^{\opp})
\end{equation*}
such that the resulting map $\pi_{\hochhom}$ makes 
\eqref{eqn-the-agreement-of-two-decategorification-maps}
commute. The maps $\Xi^{\PP}$ and $\Xi^{\QQ}$ 
in \eqref{eqn-k0num-assignment-of-pi(pn)-and-pi(qn)} and 
\eqref{eqn-hochhom-assignment-of-pi(a(n))-pi(a(-n))}
commute with taking the Euler class, so for 
any $a \in \basecat$, the class of the symmetrising idempotent
$e_n$ of $a^{\otimes n}$ in $HH_0(\symbc{n})$ must be related 
in the algebra $\bigoplus_{m \geq 0} \hochhom_\bullet(\symbc{m})$
to $\psi_m(\id_a)$ for $m \leq n$ via the same combinatorial formulas 
which express $p_a^{(n)}$ in terms of $a_a(m)$:
\begin{align}
\label{eqn-formula-for-e_n-in-terms-of-psi_n_i}
[e_n] = \frac{1}{n!} 
\sum_{\underline{n}  \vdash n} |\underline{n}| \psi_{n_1}(\id_a)\dots
\psi_{n_{r(\underline{n})}}(\id_a), 
\end{align}
where $|\underline{n}|$ is the size of the conjugacy class
$\underline{n}$.
 
In $HH_0(\symbc{n})$, the class $[e_n]$ is the sum of the classes
$\frac{1}{n!} [\sigma]$ for each $\sigma \in S_n$. Since $HH_0$ is
commutative, the classes of conjugate $\sigma \in S_n$ are equal and thus
\begin{equation}
\label{eqn-formula-for-e_n-in-terms-of-sigma_n_i}
[e_n] =  \frac{1}{n!} 
\sum_{\underline{n}  \vdash n} |\underline{n}| 
[\sigma_{\underline{n}}],
\end{equation}
where we choose a representative $\sigma_{\underline{n}}$
of every conjugacy class $\underline{n} \in n$.
The first step
\eqref{eqn-decomposition-of-orbifold-hh-into-the-sum-of-twisted-hhs}
of the orbifold decomposition 
implies that $[\sigma_{\underline{n}}]$ 
lies in the summand of the decomposition indexed by $\underline{n}$.  
The algebra structure on $\bigoplus_{m \geq 0} \hochhom_\bullet(\symbc{m})$
is induced by the functors $\symbc{m_1} \otimes \symbc{m_2}
\rightarrow \symbc{m_1 + m_2}$ which are in turn induced by 
the inclusion $S_{m_1} \times S_{m_2} \hookrightarrow S_{m_1 + m_2}$. 
For the orbifold decomposition, this means that the product
of summands indexed by $\underline{m_1} \vdash m_1$ and 
$\underline{m_2} \vdash m_2$ lies in the summand indexed by 
$\underline{m_1} \cup \underline{m_1} \vdash m_1 + m_2$. 
Comparing \eqref{eqn-formula-for-e_n-in-terms-of-psi_n_i}
with \eqref{eqn-formula-for-e_n-in-terms-of-sigma_n_i}
and keeping the algebra structure in mind, we conclude that we must have 
$\psi_n([\id_a]) = [t]$, where $t := (1 \dots n)$ is the long cycle 
representing the conjugacy class $(n)$, the partition with 1 
part of size $n$. 

In particular, $\psi_a(\id_v)$ lies in the summand of the orbifold 
decomposition of $\hochhom_\bullet(\symbc{n})$ indexed by $(n)$. This
summand is just $\hochhom_\bullet(\basecat)$. This motivates
the following definition: 
\begin{Definition}
\label{defn-linear-maps-psi-n-and-kappa-n}
Define the $\kk$-linear maps 
$$\psi_n\colon \hochhom_\bullet(\basecat) \rightarrow
\hochhom_\bullet(\symbc{n})
\quad \quad \psi_n\colon \hochhom_\bullet(\basecat^{\opp}) \rightarrow
\hochhom_\bullet(\sym^{n}\basecat^{\opp})$$
$$\kappa_n\colon \hochhom_\bullet(\symbc{n}) \rightarrow
\hochhom_\bullet(\basecat)
\quad \quad \kappa_n\colon  
\hochhom_\bullet(\sym^{n}\basecat^{\opp})
\rightarrow
\hochhom_\bullet(\basecat^{\opp})
$$
to be the inclusions of and the projections onto the summand indexed by $(n)$ 
in the orbifold decomposition 
\eqref{eqn-noncommutative-orbifold-decomposition-for-sn-v}. 
\end{Definition}

Explicitly, we have 
\begin{equation}
\label{eqn-linear-maps-psi-n}
\begin{tikzcd}
\psi_n := \quad 
\hochcx_{\bullet}(\basecat)
\ar{r}{\eqref{eqn-quasi-iso-g-long-cycle}}
&
\hochcx_{\bullet}(\basecat^n;t)
\ar{r}{\xi_t}
&
\hochcx_{\bullet}(\symbc{n}), 
\end{tikzcd}
\end{equation}
\begin{equation}
\label{eqn-linear-maps-kappa-n}
\begin{tikzcd}
\kappa := \quad 
\hochcx_{\bullet}(\symbc{n}), 
\ar{r}{\nu}
&
\bigoplus_{\underline{n} \vdash n}
\hochcx_\bullet(\basecat^n;\sigma_{\underline{n}})_{C(\sigma_{\underline{n}})}
\ar[two heads]{r}
&
\hochcx_{\bullet}(\basecat^n;t)_t
\ar{r}{\eqref{eqn-quasi-iso-f-long-cycle}}
&
\hochcx_{\bullet}(\basecat)
\end{tikzcd}
\end{equation}
and similarly for $\basecat^{\opp}$ variants. 
In particular, on $\hochcx_0$ the map $\psi_n$ is 
given by
$$ 
\alpha \mapsto (\id \otimes \dots \otimes \id \otimes \alpha, t^{-1})
\quad \quad \forall\; a\in\basecat, \alpha \in \homm_{\basecat}(a,a). 
$$
Thus we have, as desired
$$ \psi_n([\id_a]) = [t^{-1}] = [t]. $$

Note that for any $\sigma \in S_n$, we have $[\sigma] = [\sigma^{-1}]$
in $\hochhom_0(\symbc{n})$. This is because any $\hochhom_0$ is 
a commutative algebra, and if we decompose $\sigma$ as a product of
transpositions, then $\sigma^{-1}$ is the product of the same
transpositions in the reverse order. 
 
\subsection{Commutation and Heisenberg relations}
\label{section-commutation-and-heisenberg-relations}

We now prove the commutation relations
\eqref{eq:vectheisrel1-a-gen-graded} and the Heisenberg relation
\eqref{eq:vectheisrel3-a-gen-graded} for the elements 
$\hAA_{\alpha}(n)$ and $\hAA_{\alpha}(-n)$ defined in 
Definition \ref{eqn-hochhom-assignment-of-pi(a(n))-pi(a(-n))}. 

The commutation relations \eqref{eq:vectheisrel1-a-gen-graded} are
easy to establish because they only involve the images of the elements
$\psi_n(\alpha)$ under one of the functors $\Xi^{\PP}$ and $\Xi^{\QQ}$. 
Since these functors are monoidal, it is enough to verify 
the corresponding relations for the elements $\psi_n(\alpha)$ in
$\bigoplus_{n \geq 0} \hochhom_\bullet(\symbc{n})$:

\begin{Proposition}
For any $n,m \geq 1$ and $\alpha, \beta \in \hochhom_\bullet(\basecat)$
the following relations hold in $\alghh_{\hcat\basecat}$:
\begin{align}
\label{eqn-commutation-relations-for-the-assignments-of-pi}
\hAA_{\alpha}(n) \hAA_{\beta}(m) &= 
(-1)^{\deg(\alpha)\deg(\beta)}
\hAA_{\beta}(m) \hAA_{\alpha}(n), \\
\hAA_{\alpha}(-n) \hAA_{\beta}(-m) &= 
(-1)^{\deg(\alpha)\deg(\beta)}
\hAA_{\beta}(-m) \hAA_{\alpha}(-n). 
\end{align}
\end{Proposition}
\begin{proof}
By the definition of $\hAA_{\alpha}(\pm n)$ and $\hAA_{\beta}(\pm m)$
we need to establish that:
\begin{align*}
\Xi^{\PP}\left(\psi_n(\alpha)\right)
\Xi^{\PP}\left(\psi_m(\beta)\right)
 &= 
(-1)^{\deg(\alpha)\deg(\beta)}
\Xi^{\PP}\left(\psi_m(\beta)\right)
\Xi^{\PP}\left(\psi_n(\alpha)\right),
\\
\Xi^{\QQ}\left(\psi_n(\alpha)\right)
\Xi^{\QQ}\left(\psi_m(\beta)\right)
 &= 
(-1)^{\deg(\alpha)\deg(\beta)}
\Xi^{\QQ}\left(\psi_m(\beta)\right)
\Xi^{\QQ}\left(\psi_n(\alpha)\right),
\end{align*}
where again we use implicitly the isomorphism
$\hochcx_\bullet(\basecat) \simeq \hochcx_\bullet(\basecat^{\opp})$
established in \S\ref{section-hochschild-homology-of-the-opposite-category}

The maps $\Xi^{\PP}$ and $\Xi^{\QQ}$ in 
\eqref{eqn-hochhom-assignment-of-pi(a(n))-pi(a(-n))}
are induced by the DG $2$-functors 
$$ \Xi^{\PP} \colon \bihperf(\bicat{Sym}_\basecat) 
\to \hcat\basecat
\quad \text{ and } \quad 
\Xi^{\QQ} \colon \bihperf(\bicat{Sym}_{\basecat^{\opp}}) 
\to \hcat\basecat
$$
constructed in 
\cite[\S6.1]{gyenge2021heisenberg}. They thus package up
into algebra homomorphisms
$$ \Xi^{\PP} \colon 
\bigoplus_{n \geq 0} \hochhom_{\bullet}(\symbc{n})
\longrightarrow \alghh_{\hcat\basecat}
\quad\text{ and }\quad
\Xi^{\PP} \colon 
\bigoplus_{n \geq 0} \hochhom_{\bullet}(\sym^{n} \basecat^{\opp})
\longrightarrow \alghh_{\hcat\basecat}.
$$
It suffices therefore to establish
in $\bigoplus_{n \geq 0} \hochhom_{\bullet}(\symbc{n})$ 
and 
$\bigoplus_{n \geq 0} \hochhom_{\bullet}(\sym^{n} \basecat^{\opp})$
the relationship 
\begin{equation}
\psi_n(\alpha)
\psi_m(\beta)
= 
(-1)^{\deg(\alpha)\deg(\beta)}
\psi_m(\beta)
\psi_n(\alpha),
\end{equation}

This is trivial, as $\psi_\bullet$ are degree preserving
and the algebra $\bigoplus_{n \geq 0} \hochhom_{\bullet}(\symbc{n})$ 
is super-commutative for any small DG category $\basecat$. 
For super-commutativity, recall that the algebra structure 
$$ \hochhom_{\bullet}(\symbc{n}) \otimes 
\hochhom_{\bullet}(\symbc{m})
\rightarrow 
\hochhom_{\bullet}(\symbc{n+m}) $$
is induced by the functor 
$\symbc{n} \otimes \symbc{n} \rightarrow \symbc{n+m}$
induced by the subgroup inclusion $S_n \times S_m \subset S_{n+m}$. 

Let $\tau \in S_{n+m}$ be the permutation 
$$
\tau(i) := 
\begin{cases}
m + i, \quad \quad 1 \leq i \leq n, \\
i - n, \quad \quad n+1 \leq i \leq n+m, 
\end{cases}
$$
which swaps the first $n$ and the last $m$ elements. 
For any 
$$ a_1 \otimes \dots \otimes a_n \otimes b_1 \otimes \dots \otimes b_m
\in \symbc{n+m} $$
we have in $\symbc{n+m}$ an isomorphism 
$$ a_1 \otimes \dots \otimes a_n \otimes b_1 \otimes \dots \otimes b_m
\xrightarrow{(\id, \tau)}
b_1 \otimes \dots \otimes b_m \otimes a_1 \otimes \dots \otimes a_m. $$
These give an endofunctor 
$T\colon \symbc{n+m} \rightarrow \symbc{n+m}$
naturally isomorphic to $\id_{\symbc{n+m}}$. The induced map 
$$ T\colon \hochhom_\bullet(\symbc{n+m}) \rightarrow
\hochhom_\bullet(\symbc{n+m}) $$
for any $\eta \in \hochhom_{\bullet}(\symbc{n})$ and $\zeta \in 
\hochhom_{\bullet}(\symbc{m})$ sends 
$\eta \zeta \mapsto (-1)^{\deg(\eta)} (-1)^{\deg(\zeta)} \zeta \eta$.
But as the functor $T$ is isomorphic to $\id_{\symbc{n+m}}$ the induced map 
must the identity map.  
\end{proof}

The Heisenberg relation requires more work to establish:
\begin{Theorem}
\label{theorem-heisenberg-relation-for-the-assignments-of-pi}
For any $n,m \geq 1$ and $\alpha, \beta \in \hochhom_\bullet(\basecat)$
the following relation hold in $\alghh_{\hcat\basecat}$:
\begin{align}
\label{eqn-heisenberg-relation-for-the-assignments-of-pi}
\hAA_{\alpha}(-n) \hAA_{\beta}(m) &= 
(-1)^{\deg(\alpha)\deg(\beta)}
\hAA_{\beta}(m) \hAA_{\alpha}(-n)
+ \delta_{n,m} m \left<\alpha, \beta\right>.
\end{align}
\end{Theorem}

Since this relation involves the images of the elements
$\psi_n(\alpha)$ and $\psi_m(\beta)$ under both the functors $\Xi^{\PP}$ and $\Xi^{\QQ}$, 
it can not come from some relation already existing in 
$\hochhom_\bullet(\symbc{n+m})$. Ideally, we would have wanted
to prove it by constructing two homotopy equivalent functors 
$$\basecat^{\opp} \otimes \basecat \rightarrow
\homm_{\hcat\basecat}(0,m-n)$$
whose induced maps on Hochschild homology send $\alpha \otimes \beta
\in \hochhom_\bullet(\basecat^{\opp}) \otimes
\hochhom_\bullet(\basecat)$ to the LHS and the RHS of 
\eqref{eqn-heisenberg-relation-for-the-assignments-of-pi}. 
The Heisenberg relation would 
then follow since homotopy equivalent functors induce the
same map on Hochshchild homology
\cite[Lemma 3.4]{Keller-OnTheCyclicHomologyOfExactCategories}.  
But this is impossible: our construction of
$\pi$, and hence of classes $\hAA_\alpha(-n)$ and $\hAA_\beta(m)$, 
is not functorial. It involves the $\kk$-linear maps 
$\psi_n\colon \hochhom_\bullet(\basecat) \rightarrow
\hochhom_\bullet(\symbc{n})$ which do not come from any functors
$\basecat \rightarrow \symbc{n}$ when $n > 1$. 

However, once we have the classes $\psi_n(\alpha)$ the rest
of $\pi$ is functorial. We get $\hAA_\alpha(-n)$ and 
$\hAA_\alpha(n)$ by applying to $\psi_n(\alpha)$ the maps induced
on the Hochschild homology by the functors $\Xi_{\QQ}$ and
$\Xi_{\PP}$, respectively. To prove the Heisenberg relation
we construct two homotopy equivalent functors 
\begin{equation}
\label{eqn-two-functors-SnSm-H}
\sym^{n}\basecat^{\opp} \otimes \sym^{m}\basecat
\rightarrow 
\homm_{\hcat\basecat}(0,m-n)
\end{equation}
whose induced maps send $\psi_n(\alpha) \otimes \psi_m(\beta)$
for any $\alpha, \beta \in \hochcx_\bullet(\basecat)$ to 
the LHS and the RHS of 
\eqref{eqn-heisenberg-relation-for-the-assignments-of-pi}. 
Here and below we implicitly use the isomorphism 
$\hochcx_\bullet(\basecat) \simeq \hochcx_\bullet(\basecat^{\opp})$
defined in Section
\ref{section-hochschild-homology-of-the-opposite-category}. 

The first functor, $\Xi_{\QQ \PP}$, is straightforward: it is 
a $1$-composition of $2$-functors $\Xi_{\QQ}$ and $\Xi_{\PP}$. 
To be precise, note that our $2$-category $\hcat\basecat$ 
has object set $\mathbb{Z}$, and by construction 
$$ \homm_{\hcat\basecat}(n,m) = \homm_{\hcat\basecat}(0,m-n)
\quad \quad \quad \forall\; n,m \in \mathbb{Z}. $$
It has a natural structure of monoidal $2$-category 
given on the object set $\mathbb{Z}$ by addition 
$n \otimes n' = n + n'$
and on $1$-morphism categories by the $1$-composition functor
\begin{equation*}
\begin{tikzcd}[row sep = 0.25cm]
\homm_{\hcat\basecat}(n,m) \otimes \homm_{\hcat\basecat}(n',m') 
\ar[equals]{d}
\\
\homm_{\hcat\basecat}(n+m',m+m') \otimes
\homm_{\hcat\basecat}(n'+n,m'+n) 
\ar{d}{\quad 1\text{-comp}}
\\
\homm_{\hcat\basecat}(n+n',m+m'). 
\end{tikzcd}
\end{equation*}
We therefore have a $2$-functor
\begin{equation}
\label{eqn-two-functor-QP}
\Xi_{\QQ\PP}\colon 
\bicat{Sym}_{{\basecat}^{\opp}} \otimes \bicat{Sym}_{\basecat}
\xrightarrow{\Xi_{\QQ} \otimes \Xi_{\PP}}
\hcat\basecat \otimes \hcat\basecat
\rightarrow \hcat\basecat.
\end{equation}
\begin{Definition}
For any $n,m \geq 0$ define the DG functor
\begin{equation}
\Xi_{\QQ\PP}\colon 
\sym^{n}\basecat^{\opp} \otimes \sym^{m}\basecat 
\rightarrow 
\homm_{\hcat\basecat}(0,m-n)
\end{equation}
to be the action of the $2$-functor \eqref{eqn-two-functor-QP}
on the $1$-morphism category 
$\homm_{\sym^{n}\basecat^{\opp} \otimes \sym^{m}\basecat}(0 \otimes 0, n
\otimes m)$. 
\end{Definition}
\begin{Example}
Let 
$(a_n \otimes \dots \otimes a_1) \otimes (b_1 \otimes \dots \otimes
b_m) \in \sym^{n}\basecat^{\opp} \otimes \sym^{m}\basecat.$
The functor $\Xi_{\QQ\PP}$ takes it to the $1$-composition of
$\Xi_{\QQ}(a_n \otimes \dots \otimes a_1)$ and 
$\Xi_{\PP}(b_1 \otimes \dots \otimes b_m)$. That is 
$$ 
\Xi_{\QQ\PP}((a_n \otimes \dots \otimes a_1)\otimes 
(b_1 \otimes \dots \otimes b_m)) = 
\QQ_{a_n} \dots \QQ_{a_1} \PP_{b_1} \dots \PP_{b_m} 
\in \homm_{\hcat\basecat}(0,m-n). $$
Similarly, given any morphisms 
$$ \alpha \colon a_n \otimes \dots \otimes a_1 \rightarrow 
a'_n \otimes \dots \otimes a'_1 $$
$$ \beta \colon b_1 \otimes \dots \otimes b_m \rightarrow 
b'_1 \otimes \dots \otimes b'_m $$
in $\sym^{n}\basecat^{\opp}$ and $\sym^{m}\basecat$ we have
$$
\Xi_{\QQ\PP}(\alpha \otimes \beta) = \Xi_{\QQ}(\alpha) \circ_1 
\Xi_{\PP}(\beta). 
$$

Let $\sigma = (1234) \in S_4$, $\tau = (12)(45) \in S_5$.
Let $a_i, a'_i, b_j, b'_j \in \basecat$, 
$\alpha_i \in \homm_\basecat \bigl(a'_i, a_{\sigma^{-1}(i)}\bigr)$
and $\beta_i \in \homm_\basecat  \bigr(b_{\tau^{-1}(i)}, b'_{i}\bigr)$.
Then 
$\Xi_{\QQ\PP}\bigl((\alpha_4 \otimes \alpha_3 \otimes \alpha_2 \otimes
\alpha_1, \sigma) \otimes (\beta_1 \otimes \beta_2 \otimes \beta_3 \otimes
\beta_4 \otimes \beta_5, \tau)\bigr)$
is 
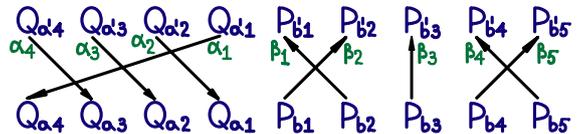
\begin{center}
\begin{tikzpicture}[scale=0.73]
	\draw  (0,0)  node {$\QQ_{a_4}$};
	\draw   (2,0)  node {$\QQ_{a_3}$};
	\draw   (4,0)   node {$\QQ_{a_2}$};
	\draw   (6,0)  node {$\QQ_{a_1}$};
	\draw   (8,0)  node {$\PP_{b_1}$};
	\draw  (10,0) node  {$\PP_{b_2}$};
	\draw  (12,0) node  {$\PP_{b_3}$};
	\draw  (14,0) node  {$\PP_{b_4}$};
	\draw  (16,0) node  {$\PP_{b_5}$};
	\draw   (0,2)  node{$\QQ_{a_4'}$};
	\draw   (2,2)  node{$\QQ_{a_3'}$};
	\draw   (4,2)   node {$\QQ_{a_2'}$};
	\draw   (6,2)  node {$\QQ_{a_1'}$};
	\draw   (8,2)  node {$\PP_{b_1'}$};
	\draw   (10,2)  node {$\PP_{b_2'}$};
	\draw   (12,2)  node {$\PP_{b_3'}$};
	\draw   (14,2)  node {$\PP_{b_4'}$};
	\draw   (16,2)  node {$\PP_{b_5'}$};
	\draw[->] (5.7,1.7) 
	-- node[label=below:{$\alpha_1$}, dot, pos=0.1] {} (0.3,0.3);
	\draw[->] (0.3,1.7) 
	-- node[label=left:{$\alpha_4$}, dot, pos=0.2] {} (1.7,0.3); 
\draw[->] (2.3,1.7) 
-- node[label=left:{$\alpha_3$}, dot, pos=0.2] {} (3.7,0.3); 
\draw[->] (4.3,1.7) 
-- node[label=left:{$\alpha_2$}, dot, pos=0.1] {} (5.7,0.3); 
\draw[->] (8.3,0.3) 
-- node[label=right:{$\beta_2$}, dot, pos=0.8] {} (9.7,1.7); 
\draw[->] (9.7,0.3) 
-- node[label=left:{$\beta_1$}, dot, pos=0.8] {} (8.3,1.7); 
\draw[->] (12,0.3) 
-- node[label=right:{$\beta_3$}, dot, pos=0.8] {} (12,1.7);
\draw[->] (14.3,0.3) 
-- node[label=right:{$\beta_5$}, dot, pos=0.8] {} (15.7,1.7); 
\draw[->] (15.7,0.3) 
-- node[label=left:{$\beta_4$}, dot, pos=0.8] {} (14.3,1.7); 
	%node[below] {$\QQ_{a_4}$};
\end{tikzpicture}
\end{center}
\end{Example}

Before defining the second functor \eqref{eqn-two-functors-SnSm-H} we  
give some intuition for its construction. It corresponds
to the following iterative procedure:

\begin{Definition}[Commutation-annihilation procedure]
\label{defn-commutation-annihilation-procedure}
Let 
$(a_n \otimes \dots \otimes a_1)\otimes (b_1 \otimes \dots \otimes b_m)
\in \sym^{n}\basecat^{\opp} \otimes \sym^{m}\basecat$ and take
the $1$-morphism the $2$-functor $\Xi_{\QQ\PP}$ sends it to:
\begin{equation}
\label{eqn-QP-word-iterative-procedure-step-1}
\QQ_{a_n} \dots \QQ_{a_1} \PP_{b_1} \dots \PP_{b_m} 
\in \homm_{\hcat\basecat}(0,m-n). 
\end{equation}

We now apply to it the following iterative procedure. Locate the
rightmost $\QQ$ which has $\PP$ to its right. At this first step, 
it is $\QQ_{a_1} \PP_{b_1}$.  Take this pair, 
and apply the homotopy equivalence 
\begin{equation}
\label{eq:dg-baby-Heisenberg-morphism}
        \PP_{b}\QQ_a \oplus \bigl(\Hom(a,b) \otimes_\kk \hunit \bigr) 
        \xrightarrow{
                \left[
                \begin{tikzpicture}[baseline={(0,0.15)}, scale=0.5]
                        \draw[->] (0.5,0) -- (1.5,1);
                        \draw[->] (0.5,1) -- (1.5,0);
                \end{tikzpicture}
                \,,\,
                \,\psi_2
                \right]
        }
        \QQ_a \PP_{b},
\end{equation}
defined in \cite[\S5.4]{gyenge2021heisenberg}. Here, it breaks 
\eqref{eqn-QP-word-iterative-procedure-step-1} into two summands. 
In one $\QQ_{a_1}$ is commuted past $\PP_{b_1}$, while in the other
$\QQ_{a_1}$ and $\PP_{b_1}$ annihilate each other and the rest 
is tensored by $\homm_{\basecat}(a_1, b_1)$:
\begin{equation}
\label{eqn-QP-word-iterative-procedure-step-2}
\begin{tikzcd}[row sep = 0.25cm]
\QQ_{a_n} \dots \QQ_{a_1} \PP_{b_1} \dots \PP_{b_m} 
\\
\QQ_{a_n} \dots \QQ_{a_2} \PP_{b_1}\QQ_{a_1} \PP_{b_2}  \dots \PP_{b_m} 
\oplus 
\QQ_{a_n} \dots \QQ_{a_2} \PP_{b_2}  \dots \PP_{b_m} \otimes
\homm_{\basecat}(a_1, b_1)
\ar{u}
\end{tikzcd}
\end{equation}

We now take each of the newly obtained summands and apply the same
procedure to them. We locate the rightmost $\QQ\PP$ subword, 
if one exists, and apply to it the homotopy equivalence 
\eqref{eq:dg-baby-Heisenberg-morphism}. This replaces the old summand
by two new ones, where the $\QQ\PP$ pair in question is commuted and
is annihilated, respectively. 
We repeat this until none of the summands contain a $\QQ\PP$ subword. 

The result is a homotopy equivalence into 
\eqref{eqn-QP-word-iterative-procedure-step-1} from a direct sum 
whose summands have form 
\begin{small}
\begin{equation}
\label{eqn-the-summands-of-sumPQcheck}
\PP_{b_1} \dots 
\widehat{\PP_{b_{j_1}}} \dots \widehat{\PP_{b_{j_k}}}
\dots \PP_{b_m}
\QQ_{a_n} \dots 
\widehat{\QQ_{a_{i_k}}} \dots \widehat{\QQ_{a_{i_1}}} 
\dots \QQ_{a_1}
\otimes \homm_{\basecat}(a_{i_1}, b_{j_{\sigma(1)}}) \otimes \dots \otimes 
\homm_{\basecat}(a_{i_k}, b_{j_{\sigma(k)}}). 
\end{equation}
\end{small}
where a hat indicates that we skip this element and where $\sigma \in S_k$. 
We obtain each summand by starting at 
\eqref{eqn-QP-word-iterative-procedure-step-1}, 
taking each $\QQ$ (from right to left) and starting to move it 
past all the $\PP$s. At each $\PP$ we choose whether to commute 
the $\QQ$ past it or to annihilate the $\QQ$ against it and start on 
the next $\QQ$. In \eqref{eqn-the-summands-of-sumPQcheck} the elements $\QQ_{a_{i_1}}$, \dots, $\QQ_{a_{i_k}}$ 
were annihilated against the elements $\PP_{b_{j_1}}$, \dots,
$\PP_{b_{j_k}}$ in the order specified by $\sigma$. 
The remaining $\QQ$s were succesfully commuted past all the $\PP$s. 
\end{Definition}

\begin{Example}
\label{exmpl-commutation-annihilation-procedure-for-n-equals-three}
Let $a_1, a_2, a_3, b_1, b_2, b_3 \in \basecat$. 
Write $(a_i,b_j)$ for $\homm_\basecat(a_i,b_j)$. For 
$$ \QQ_{a_3}\QQ_{a_2}\QQ_{a_1} \PP_{b_1}\PP_{b_2}\PP_{b_3} \in
\homm_{\hcat\basecat}(0,0) $$
the first few steps of
the commutation-annihilation procedure described above are 
\begin{tiny}
\begin{equation*}
\begin{tikzcd}[row sep = 0.25cm]
\QQ_{a_3}\QQ_{a_2}\QQ_{a_1} \PP_{b_1}\PP_{b_2}\PP_{b_3}
\\
\left(\QQ_{a_3}\QQ_{a_2}\PP_{b_1}\QQ_{a_1} \PP_{b_2}\PP_{b_3} \right)
\oplus 
\left(
\QQ_{a_3}\QQ_{a_2}\PP_{b_2}\PP_{b_3} \otimes (a_1, b_1)
\right)
\ar{u}
\\
\begin{matrix}
\left( 
\QQ_{a_3}\QQ_{a_2}\PP_{b_1}\PP_{b_2}\QQ_{a_1}\PP_{b_3} 
\right)
\oplus 
\left(
\QQ_{a_3}\QQ_{a_2}\PP_{b_1}\PP_{b_3} \otimes (a_1, b_2)
\right)
\oplus 
\left(
\QQ_{a_3}\PP_{b_2}\QQ_{a_2}\PP_{b_3} \otimes (a_1, b_1)
\right)
\oplus
\left( 
\QQ_{a_3}\PP_{b_3} \otimes (a_1, b_1) \otimes
(a_2, b_2)
\right)
\end{matrix}
\ar{u} 
\\
\begin{matrix}
\left( 
\QQ_{a_3}\QQ_{a_2}\PP_{b_1}\PP_{b_2}\PP_{b_3}\QQ_{a_1} 
\right)
\oplus
\left(
\QQ_{a_3}\QQ_{a_2}\PP_{b_1}\PP_{b_2} \otimes (a_1, b_3) 
\right)
\oplus 
\left(
\QQ_{a_3}\PP_{b_1}\QQ_{a_2}\PP_{b_3} \otimes (a_1, b_2)
\right)
\oplus 
\left(
\QQ_{a_3}\PP_{b_3} \otimes (a_1, b_2) \otimes
(a_2, b_1)
\right)
\oplus
\\
\oplus
\left(
\QQ_{a_3}\PP_{b_2}\PP_{b_3}\QQ_{a_2} \otimes (a_1, b_1)
\right)
\oplus
\left(
\QQ_{a_3}\PP_{b_2}\otimes (a_1, b_1) \otimes
(a_2, b_3)
\right)
\oplus
\left( 
\PP_{b_3}\QQ_{a_3} \otimes (a_1, b_1) \otimes
(a_2, b_2)
\right)
\oplus
\bigl( 
(a_1, b_1) \otimes (a_2, b_2)
\otimes (a_3, b_3)
\bigr)
\end{matrix}
\ar{u} 
\\
\begin{matrix}
\left( 
\QQ_{a_3}\PP_{b_1}\QQ_{a_2}\PP_{b_2}\PP_{b_3}\QQ_{a_1} 
\right)
\oplus
\left( 
\QQ_{a_3}\PP_{b_2}\PP_{b_3}\QQ_{a_1} \otimes (a_2, b_1)
\right)
\oplus
\left(
\QQ_{a_3}\PP_{b_1}\QQ_{a_2}\PP_{b_2} \otimes (a_1, b_3) 
\right)
\oplus
\left(
\QQ_{a_3}\PP_{b_2} \otimes (a_1, b_3) \otimes
(a_2, b_1)
\right)
\oplus 
\\
\oplus
\left(
\QQ_{a_3}\PP_{b_1}\PP_{b_3}\QQ_{a_2} \otimes (a_1, b_2)
\right)
\oplus 
\left(
\QQ_{a_3}\PP_{b_1} \otimes (a_1, b_2) \otimes
(a_2, b_3)
\right)
\oplus
\left(
\PP_{b_3}\QQ_{a_3} \otimes (a_1, b_2) \otimes
(a_2, b_1)
\right)
\oplus
\bigl(
(a_1, b_2) \otimes (a_2, b_1) \otimes (a_3, b_3)
\bigr)
\oplus
\\
\oplus
\left(
\PP_{b_2}\QQ_{a_3}\PP_{b_3}\QQ_{a_2} \otimes (a_1, b_1)
\right)
\oplus
\left(
\PP_{b_3}\QQ_{a_2} \otimes (a_1, b_1) \otimes
(a_3, b_2)
\right)
\oplus
\left(
\PP_{b_2}\QQ_{a_3}\otimes (a_1, b_1) \otimes
(a_2, b_3)
\right)
\oplus
\bigl(
(a_1, b_1) \otimes
(a_2, b_3) \otimes (a_3, b_2)
\bigr)
\oplus
\\
\oplus
\left( 
\PP_{b_3}\QQ_{a_3} \otimes (a_1, b_1) \otimes (a_2, b_2)
\right)
\oplus
\bigl( 
(a_1, b_1) \otimes (a_2, b_2) \otimes (a_3, b_3)
\bigr)
\end{matrix}
\ar{u}
\end{tikzcd}
\end{equation*}
\end{tiny}

Its result is the direct sum whose summands, grouped by 
the number of annihilations, are:
\begin{enumerate}
\item $\;\;\;\PP_{b_1}\PP_{b_2}\PP_{b_3}\QQ_{a_3}\QQ_{a_2}\QQ_{a_1}$, 
\item 
$\;\;\;\bigl(
\PP_{b_2}\PP_{b_3}\QQ_{a_3}\QQ_{a_2} \otimes (a_1, b_1)
\bigr)
\oplus 
\bigl(
\PP_{b_1}\PP_{b_3}\QQ_{a_3}\QQ_{a_2} \otimes (a_1, b_2)
\bigr)
\oplus 
\bigl(
\PP_{b_1}\PP_{b_2}\QQ_{a_3}\QQ_{a_2} \otimes (a_1, b_3)
\bigr)
\oplus$ \\
$\oplus 
\bigl(
\PP_{b_2}\PP_{b_3}\QQ_{a_3}\QQ_{a_1} \otimes (a_2, b_1)
\bigr)
\oplus 
\bigl(
\PP_{b_1}\PP_{b_3}\QQ_{a_3}\QQ_{a_1} \otimes (a_2, b_2)
\bigr)
\oplus 
\bigl(
\PP_{b_1}\PP_{b_2}\QQ_{a_3}\QQ_{a_1} \otimes (a_2, b_3)
\bigr)
\oplus
$\\
$
\oplus
\bigl(
\PP_{b_2}\PP_{b_3}\QQ_{a_2}\QQ_{a_1} \otimes (a_3, b_1)
\bigr)
\oplus
\bigl(
\PP_{b_1}\PP_{b_3}\QQ_{a_2}\QQ_{a_1} \otimes (a_3, b_2)
\bigr)
\oplus
\bigl(
\PP_{b_1}\PP_{b_2}\QQ_{a_2}\QQ_{a_1} \otimes (a_3, b_3),
\bigr)
$ 

\item 
$\;\;\;\bigl
(\PP_{b_3}\QQ_{a_3} \otimes (a_1, b_1) \otimes (a_2, b_2)
\bigr)
\oplus 
\bigl(
\PP_{b_2}\QQ_{a_3} \otimes (a_1, b_1) \otimes (a_2, b_3)
\bigr)
\oplus 
\bigl(
\PP_{b_3}\QQ_{a_2} \otimes (a_1, b_1) \otimes (a_3, b_2)
\bigr)
\oplus$ \\
$\oplus 
\bigl(
\PP_{b_2}\QQ_{a_2} \otimes (a_1, b_1) \otimes (a_3, b_3)
\bigr)
\oplus 
\bigl(
\PP_{b_3}\QQ_{a_3} \otimes (a_1, b_2) \otimes (a_2, b_1)
\oplus 
\bigl(
\PP_{b_1}\QQ_{a_3} \otimes (a_1, b_2) \otimes (a_2, b_3)
\bigr)
\oplus$ \\
$\oplus 
\bigl(
\PP_{b_3}\QQ_{a_2} \otimes (a_1, b_2) \otimes (a_3, b_1)
\bigr)
\oplus 
\bigl(
\PP_{b_1}\QQ_{a_2} \otimes (a_1, b_2) \otimes (a_3, b_3)
\bigr)
\oplus 
\bigl(
\PP_{b_2}\QQ_{a_3} \otimes (a_1, b_3) \otimes (a_2, b_1)
\bigr)
\oplus$ \\
$\oplus 
\bigl(
\PP_{b_1}\QQ_{a_3} \otimes (a_1, b_3) \otimes (a_2, b_2)
\bigr)
\oplus 
\bigl(
\PP_{b_2}\QQ_{a_2} \otimes (a_1, b_3) \otimes (a_3, b_1)
\bigr)
\oplus 
\bigl(
\PP_{b_1}\QQ_{a_2} \otimes (a_1, b_3) \otimes (a_3, b_2)
\bigr)
\oplus$ \\
$\oplus 
\bigl(
\PP_{b_3}\QQ_{a_1} \otimes (a_2, b_1) \otimes (a_3, b_2)
\bigr)
\oplus 
\bigl(
\PP_{b_2}\QQ_{a_1} \otimes (a_2, b_1) \otimes (a_3, b_3)
\bigr)
\oplus 
\bigl(
\PP_{b_3}\QQ_{a_1} \otimes (a_2, b_2) \otimes (a_3, b_1)
\bigr)
\oplus$ \\
$\oplus 
\bigl(
\PP_{b_1}\QQ_{a_1} \otimes (a_2, b_2) \otimes (a_3, b_3)
\bigr)
\oplus 
\bigl(
\PP_{b_2}\QQ_{a_1} \otimes (a_2, b_3) \otimes (a_3, b_1)
\bigr)
\oplus 
\bigl(
\PP_{b_1}\QQ_{a_1} \otimes (a_2, b_3) \otimes (a_3, b_2)
\bigr)$
\item 
$\;\;\; 
\bigl(
(a_1, b_1) \otimes (a_2, b_2) \otimes (a_3, b_3)
\bigr)
\oplus 
\bigl(
(a_1, b_1) \otimes (a_2, b_3) \otimes (a_3, b_2)
\bigr)
\oplus 
\bigl(
(a_1, b_2) \otimes (a_2, b_1) \otimes (a_3, b_3)
\bigr)
\oplus$ \\
$\oplus 
\bigl(
(a_1, b_2) \otimes (a_2, b_3) \otimes (a_3, b_2)
\bigr)
\oplus 
\bigl(
(a_1, b_3) \otimes (a_2, b_1) \otimes (a_3, b_2)
\bigr)
\oplus 
\bigl(
(a_1, b_3) \otimes (a_2, b_2) \otimes (a_3, b_1)
\bigr)
$. 
\end{enumerate}
\end{Example}

We now make this procedure functorial. We begin with the annihilation:
\begin{Definition}
\label{defn-annihilation-functor}
Let $n,m > 0$ and let $\min(n,m) \geq k \geq 0$. Define the DG functor
\begin{equation}
(\hat{k})\colon 
\sym^{n}\basecat^{\opp} \otimes \sym^{m}\basecat
\longrightarrow 
\hperf\left(
\sym^{n-k}\basecat^{\opp} \otimes \sym^{m-k}\basecat \right)
\end{equation}
as the composition 
\begin{align*}
\sym^{n}\basecat^{\opp} \otimes \sym^{m}\basecat
\xrightarrow{
\Res_{S_{n-k} \times S_{k} \times S_{m-k}}^{S_n \times S_m}
}
& \hperf \left(
\sym^{n-k}\basecat^{\opp} 
\otimes 
\left(\left((\basecat^{\opp})^{\otimes k} \otimes \basecat^{\otimes k}\right)
\rtimes S_k\right) 
\otimes 
\sym^{m-k}\basecat
\right)
\simeq 
\\
\simeq 
& \hperf \left(
\sym^{n-k}\basecat^{\opp} 
\otimes 
\sym^k(\basecat^{\opp} \otimes \basecat)
\otimes 
\sym^{m-k}\basecat
\right)
\xrightarrow{\homm_\basecat(-,-)}
\\
\rightarrow 
&
\hperf \left(
\sym^{n-k}\basecat^{\opp} 
\otimes 
\sym^k(\hperf \kk)
\otimes 
\sym^{m-k}\basecat
\right)
\xrightarrow{(-) \otimes \dots \otimes (-)}
\\
\rightarrow 
& \hperf \left(
\sym^{n-k}\basecat^{\opp} 
\otimes 
\hperf \kk
\otimes 
\sym^{m-k}\basecat
\right)
\simeq 
\\
\simeq 
&
\hperf \left(
\sym^{n-k}\basecat^{\opp} 
\otimes 
\sym^{m-k}\basecat
\right).
\end{align*}
The first composant
is the restriction of scalars functors induced by the group embedding 
$$ S_{n-k} \times S_{k} \times S_{m-k} \hookrightarrow 
S_{n-k} \times S_{k} \times S_{k} \times S_{m-k}
\hookrightarrow  S_{n} \times S_{m} $$
where $S_k$ embeds into $S_k \times S_k$ as $\sigma \mapsto (\iota
\sigma \iota, \sigma)$ with $\iota$  the order reversing involution 
$1 \dots n \mapsto n \dots 1$. 

The second composant is induced by the factor permuting isomorphism 
$$ \left(\basecat^{\opp}\right)^{\otimes k} \otimes \basecat^{\otimes k} 
\simeq (\basecat^{\opp} \otimes \basecat)^{\otimes k} $$
which sends $a_k \otimes \dots \otimes a_1 \otimes b_1 \otimes \dots
\otimes b_k$ to $a_1 \otimes b_1 \otimes 
\dots \otimes a_k \otimes b_k$. 

The third composant is induced by the functor 
$$ \homm_\basecat(-,-) \colon \basecat^{\opp} \otimes \basecat
\rightarrow \hperf \kk $$
which sends $a \otimes b \mapsto \homm_{\basecat}(a,b)$. 

The fourth composant is induced by the functor 
$$ 
(-) \otimes \dots \otimes (-) \colon 
\sym^k \hperf \kk \rightarrow \hperf \kk $$
which tensors $k$ complexes of $\kk$-modules together. 

The last composant is the isomorphism given by 
the restriction of scalars along the Yoneda embedding 
$\kk \rightarrow \hperf \kk$. 
\end{Definition}

\begin{Lemma}
\label{lemma-explicit-description-of-the-functor-hatk}
Let $n,m > 0$ and $\min(n,m) \geq k \geq 0$. Let
$a_1, \dots, a_n \in \basecat$ and $b_1, \dots, b_m \in \basecat$. 
Then
$$
(\hat{k})\left((a_n \otimes \dots \otimes a_1) \otimes (b_1 \otimes
\dots \otimes b_m)\right) \simeq
$$
\begin{equation}
\label{eqn-direct-sum-decomposition-of-hatk}
\bigoplus_{
\begin{smallmatrix}
1 \leq i_1 < \dots < i_k \leq n, \\
1 \leq j_1 < \dots < j_k \leq m, \\
\sigma \in S_k
\end{smallmatrix}
}
\begin{matrix}
(a_n \otimes \dots \widehat{a_{i_k}} \dots \widehat{a_{i_1}} \dots \otimes a_1) \otimes 
(b_1 \otimes \dots \widehat{b_{j_1}} \dots \widehat{b_{j_k}} \dots \otimes b_m) \otimes 
\\
\otimes \homm_{\basecat}(a_{i_1}, b_{j_{\sigma(1)}}) 
\otimes 
\dots 
\otimes 
\homm_{\basecat}(a_{i_k}, b_{j_{\sigma(k)}}). 
\end{matrix}
\end{equation}

Furthermore, let $\alpha_i \in \homm_\basecat(a'_i, a_i)$ for $1 \leq i \leq n$ and $\beta_j \in \homm_\basecat(b_j, b'_j)$ for $1 \leq j \leq m$. Then 
in terms of the direct sum decomposition  
\eqref{eqn-direct-sum-decomposition-of-hatk}
$$ (\hat{k})\bigl((\alpha_n \otimes \dots \otimes \alpha_1) \otimes (\beta_1 \otimes
\dots \otimes \beta_m)\bigr) $$
maps the summand indexed by 
$1 \leq i_1 < \dots < i_k \leq n $,
$1 \leq j_1 < \dots < j_k  \leq m$, 
and $\sigma \in S_k$ to the summand with the same index via the map 
$$ (\alpha_m \otimes \dots \widehat{\alpha_{i_k}} \dots \widehat{\alpha_{i_1}} \dots \otimes \alpha_1) \otimes 
(\beta_1 \otimes \dots \widehat{\beta_{j_{\sigma(1)}}} \dots \widehat{\beta_{j_{\sigma(k)}}} 
\dots \otimes \beta_m ) 
\otimes 
(\beta_{j_{\sigma(1)}} \circ \alpha_1) \otimes \dots \otimes (\beta_{j_{\sigma(k)}} \circ \alpha_k).   
$$

Similarly, let $\eta \in S_n$ and $\zeta \in S_m$. Then 
$$ (\hat{k})\bigl(\eta \otimes \zeta \bigr) $$
maps the summand indexed by $1 \leq i_1 < \dots < i_k \leq n$,
$1 \leq j_1 < \dots < j_k  \leq m$, and $\sigma \in S_k$ 
to the one indexed by 
$1 \leq \eta(i_{\tau^{-1}(1)}) < \dots < \eta(i_{\tau^{-1}(k)}) \leq n$, 
$1 \leq \zeta(j_{\upsilon^{-1}(1)}) < \dots <
\zeta(j_{\upsilon^{-1}(k)}) \leq m$, 
and $\upsilon \sigma \tau^{-1} \in S_k$ via the isomorphism
$$ \phi \otimes \tau \otimes \chi $$
where
\begin{itemize}
\item $\tau \in S_k$ reorders  
$\eta(i_1)$, \dots, $\eta(i_k)$ in the increasing order,
\item $\upsilon \in S_k$ reorders 
$\zeta(j_1)$, \dots, $\zeta(j_k)$ in the increasing order, 
\item $\phi \in S_{n-k}$ reorders 
$\eta(n) \dots \widehat{\eta(i_k)} \dots  \widehat{\eta(i_1)} \dots
\eta(1)$ 
in the decreasing order, 
\item $\chi \in S_{m-k}$ reorders
$\zeta(1) \dots \widehat{\zeta(j_1)} \dots  \widehat{\zeta(j_k)} \dots
\zeta(m)$ in the increasing order. 
\end{itemize}
\end{Lemma}
\begin{proof}
For any $p > q > 0$ the left cosets of 
$S_{p-q} \otimes S_q \hookrightarrow S_p$
are enumerated by the partitions of $\left\{1,\dots, p\right\}$ 
into $q$ and $p-q$ elements. This can be viewed as choices of
$1 \leq i_1 < \dots < i_q \leq p$ and as a representative of each coset 
we can take either of the permutations
$$ 1\dots p \mapsto 1\dots \widehat{i_1} \dots  \widehat{i_q} \dots p \;
i_1 \dots i_q, $$
$$ 1\dots p \mapsto p\dots \widehat{i_q} \dots  \widehat{i_1} \dots 1 \;
i_q \dots i_1. $$  
For any $p > 0$ the left cosets of 
$S_p \hookrightarrow S_p \times S_p$
with the embedding as in Definition \ref{defn-annihilation-functor} 
are enumerated by permutations $\sigma \in S_p$. As a representative
of each coset we can take $\id \times \sigma$ or $\sigma \times \id$.  

Hence the left cosets of 
$S_{n-k} \times S_{k} \times S_{m-k} \hookrightarrow  S_{n} \times S_{m}$
are given by choices of $1 \leq i_1 < \dots < i_k \leq n$, 
$1 \leq j_1 < \dots < j_k \leq m$, $\sigma \in S_k$. 
The coset corresponding to $(\underline{i}, \underline{j}, \sigma)$
has the representative 
\begin{equation}
\rho_{\underline{i}, \underline{j}, \sigma}
:= 
\text{ the product of }
\begin{cases}
1\dots n \mapsto n\dots \widehat{i_k} \dots  \widehat{i_1} \dots
1\; i_k \dots i_1 & \quad \in S_n, \\
1\dots m \mapsto j_{\sigma(1)}  \dots j_{\sigma(k)} 
\; 1\dots \widehat{j_1} \dots  \widehat{j_k} \dots m & \quad \in S_m. 
\end{cases}
\end{equation}

The group $S_n \times S_m$ acts on the set of left
cosets of $S_{n-k} \times S_{k} \times S_{m-k}$ by left
multiplication. On the coset representatives $\rho_{\underline{i},
\underline{j}, \sigma}$, for any 
$\eta \times \zeta \in S_n \times S_m$ we have 
$$ (\eta \times \zeta) \rho_{(\underline{i}, \underline{j}, \sigma)} = 
\rho_{(\eta(\underline{i}), \zeta(\underline{j}), \upsilon\sigma\tau^{-1})} 
(\phi \times \tau \times \chi)
\quad \quad \quad 
(\phi \times \tau \times \chi) \in S_{n-k} \times S_k \times S_{m-k}
$$
with $\tau$, $\upsilon$, $\phi$, and $\chi$ as above. 
By Lemma  
\ref{lemma-descripton-of-res-G-to-H-functor}
the first composant
$\Res_{S_{n-k} \times S_{k} \times S_{m-k}}^{S_n \times S_m}$ of
$(\hat{k})$ sends 
$$  (a_n \otimes \dots \otimes a_1) \otimes (b_1 \otimes \dots \otimes b_m) $$
to the direct sum of representable modules
$$ \bigoplus_{
\begin{smallmatrix}
1 \leq i_1 < \dots < i_k \leq n, \\
1 \leq j_1 < \dots < j_k \leq m, \\
\sigma \in S_k
\end{smallmatrix}
}
(a_n \otimes \dots \widehat{a_{i_k}} \dots \widehat{a_{i_1}} \dots
\otimes a_1) \otimes a_{i_k} \otimes \dots \otimes a_{i_1} \otimes 
b_{j_{\sigma(1)}} \otimes \dots \otimes b_{j_{\sigma(k)}} \otimes 
(b_1 \times \dots \widehat{b_{j_1}} \dots \widehat{b_{j_k}} \dots
\otimes b_m). 
$$

The second composant sends each summand to  
$$ 
(a_n \otimes \dots \widehat{a_{i_k}} \dots \widehat{a_{i_1}} \dots
\otimes a_1) \otimes (a_{i_1} \otimes b_{j_{\sigma(1)}}) 
\otimes \dots \otimes (a_{i_k} \otimes b_{j_{\sigma(k)}}) \otimes 
(b_1 \times \dots \widehat{b_{j_1}} \dots \widehat{b_{j_k}} \dots
\otimes b_m). 
$$
The third sends it further to  
$$
(a_n \otimes \dots \widehat{a_{i_k}} \dots \widehat{a_{i_1}} \dots
\otimes a_1) \otimes \homm_{\basecat}(a_{i_1}, b_{j_{\sigma(1)}}) 
\otimes \dots \otimes \homm_{\basecat}(a_{i_k}, b_{j_{\sigma(k)}}) 
\otimes (b_1 \otimes \dots \widehat{b_{j_1}} \dots \widehat{b_{j_k}} \dots
\otimes b_m). 
$$
The fourth views the factor
$$ \homm_{\basecat}(a_{i_1}, b_{j_{\sigma(1)}}) 
\otimes \dots \otimes \homm_{\basecat}(a_{i_k}, b_{j_{\sigma(k)}}) $$
as a representable $(\hperf \kk)$-module instead of a representable 
$\sym^k(\hperf \kk)$-module. Finally the fifth views it as a
perfect $\kk$-module instead of a representable $(\hperf \kk)$-module. 

The remaining assertions 
follow similarly by Lemma \ref{lemma-descripton-of-res-G-to-H-functor}. 
\end{proof}

We now define a $2$-functor corresponding to each group of summands 
with the same number of annihilations in the result of the
commutation-annihilation procedure of
Definition \ref{defn-commutation-annihilation-procedure}. 

\begin{Definition}
\label{defn-two-functor-PQhatk}
For any $k \geq 0$ define a $2$-functor 
\begin{equation}
\label{eqn-two-functor-PQhatk}
\Xi_{\PP\QQ}(\hat{k})\colon 
\bicat{Sym}_{{\basecat}^{\opp}} \otimes \bicat{Sym}_{\basecat}
\rightarrow 
\hcat\basecat \otimes \hcat\basecat
\rightarrow \hcat\basecat
\end{equation}
on objects to be the map $(n,m) \mapsto m - n$ for all $n,m \in \mathbb{Z}$
and on $1$-morphisms to be the functor
$$ \homm_{\bicat{Sym}_{{\basecat}^{\opp}} \otimes
\bicat{Sym}_{\basecat}}\bigl((r,s), (r+n, s+m)\bigr) \rightarrow 
\homm_{\hcat\basecat}\bigl(s-r, (s-r) + (m-n)\bigr) 
\quad \quad \forall\; r,s,n,m \in \mathbb{Z} $$
defined for $k > \min(n,m)$ to be zero and 
for $k \leq \min(n,m)$ to be the composition 
\begin{align*}
\sym^{n}\basecat^{\opp} \otimes \sym^{m}\basecat
\xrightarrow{(\hat{k})}
& \hperf\left(
\sym^{n-k}\basecat^{\opp} \otimes \sym^{m-k}\basecat \right)
\simeq 
\hperf\left( \sym^{n-k}\basecat^{\opp} \right) 
\otimes 
\hperf\left( \sym^{m-k}\basecat \right)\simeq
\\
\simeq 
& \hperf\left( \sym^{m-k}\basecat \right)
\otimes 
\hperf\left( \sym^{n-k}\basecat^{\opp} \right)
\xrightarrow{\Xi_\PP \circ_1 \Xi_\QQ}
\homm_{\hcat\basecat}\bigl(s-r,(s-r)+(m-n)\bigr). 
\end{align*}
\end{Definition}

We finally define the second functor \eqref{eqn-two-functors-SnSm-H}
which we show to be homotopy equivalent to $\Xi_{\QQ\PP}$:

\begin{Definition}
\label{defn-dg-functor-PQhatk}
Define the functor
\begin{equation}
\bigoplus \Xi_{\PP\QQ}(\hat{k})\colon
\sym^{n}\basecat^{\opp} \otimes \sym^{m}\basecat
\rightarrow 
\homm_{\hcat\basecat}(0,m-n)
\end{equation}
to be the action of the $2$-functor 
$\bigoplus_{k=0}^{\infty} \Xi_{PQ}(\hat{k})$
on the $1$-morphism category $\homm_{\bicat{Sym}_{{\basecat}^{\opp}}
\otimes \bicat{Sym}_{\basecat}}((0,0), (n,m))$. 
\end{Definition}

\begin{Corollary}
For any object
$$ (a_n \otimes \dots \otimes a_1) \otimes (b_1 \otimes \dots \otimes
b_m) \in \sym^{n}\basecat^{\opp} \otimes \sym^{m}\basecat $$
the $1$-morphism 
$$ 
\bigoplus \Xi_{\PP\QQ}(\hat{k})\bigl(
(a_n \otimes \dots \otimes a_1) \otimes (b_1 \otimes \dots \otimes
\bigr)
$$
is isomorphic to the result of the commutation-annihilation procedure
of Definition \ref{defn-commutation-annihilation-procedure} applied
to 
$$ 
\Xi_{\QQ\PP}
\bigl( (a_n \otimes \dots \otimes a_1) \otimes (b_1 \otimes \dots
\otimes b_n)
\bigr).
$$
\end{Corollary}
\begin{proof}
This follows from Lemma 
\ref{lemma-explicit-description-of-the-functor-hatk}
and the description of the summands of the result of the
commutation-annihilation procedure given in 
\eqref{eqn-the-summands-of-sumPQcheck}. 
\end{proof}

We need to show that the homotopy equivalence
given by the commutation-annihilation procedure 
$$ \bigoplus \Xi_{\PP\QQ}(\hat{k})\bigl(
(a_n \otimes \dots \otimes a_1) \otimes (b_1 \otimes \dots \otimes b_n
\bigr)
\xrightarrow{\sim}
\Xi_{QP}
\bigl( (a_n \otimes \dots \otimes a_1) \otimes (b_1 \otimes \dots
\otimes b_n)
\bigr)
$$ 
define a natural transformation 
$ \bigoplus \Xi_{\PP\QQ}(\hat{k}) \rightarrow \Xi_{\QQ\PP}$. 
First, we give a direct definition:
\begin{Definition}
\label{defn-natural-transformation-phi-termwise} 
Let 
$$ \underline{a} = (a_n \otimes \dots \otimes a_1) \in
\sym^{n}\basecat^{\opp}, $$
$$ \underline{b} = (b_1 \otimes \dots \otimes b_m) \in  
\otimes \sym^{m}\basecat. $$
Define a $2$-morphism in $\hcat\basecat$ 
\begin{equation}
\label{eqn-natural-transformation-phi-termwise} 
\phi\colon \bigoplus \Xi_{\PP\QQ}(\hat{k})\bigl(\underline{a} \otimes \underline{b} \bigr)
\longrightarrow 
\Xi_{\PP\QQ}\bigl(\underline{a} \otimes \underline{b}\bigr)
\end{equation}
by setting for 
$0 \leq k \leq \min(n,m)$, $1 \leq i_1 < \dots < i_k \leq n$, 
$1 \leq j_1 < \dots < j_k \leq m$, and $\sigma \in S_k$
the component of $\phi$ on the corresponding summand of 
$\Xi_{\PP\QQ}(\hat{k})\bigl(\underline{a} \otimes \underline{b} \bigr)$
to be the adjoint of the morphism
\begin{equation}
\label{eqn-the-adjoint-of-the-natural-transformation-phi}
\begin{tikzcd}[row sep=0.25cm]
\homm_{\basecat}(a_{i_1}, b_{j_{\sigma(1)}}) \otimes \dots \otimes 
\homm_{\basecat}(a_{i_k}, b_{j_{\sigma(k)}})
\ar{d}
\\
\homm_{\hcat\basecat}\left(
\PP_{b_1} \dots 
\widehat{\PP_{b_{j_\bullet}}} 
\dots \PP_{b_m}
\QQ_{a_n} \dots 
\widehat{\QQ_{a_{i_\bullet}}} 
\dots \QQ_{a_1}, \quad 
\QQ_{a_n} \dots \QQ_{a_1} \PP_{b_1} \dots \PP_{b_m}\right)
\end{tikzcd}
\end{equation}
which sends each $\gamma_1 \otimes \dots \otimes \gamma_k$
to the $2$-morphism in $\hcat\basecat$ 
defined by the following planar diagram:
\begin{enumerate}
\item First, draw the annihilation strands. Start at the top of
the diagram, and for each $1 \leq l \leq k$ draw a circular strand 
counterclockwise from $\QQ_{i_l}$ to $\PP_{j_{\sigma(l)}}$ and
decorate it by $\gamma_l$. Draw these
so that each annihilation strand dips below the previous one, i.e. 
the height of the $l$-th strand arc is less than that of the $(l+1)$-st
for all $1 \leq l < k$. 
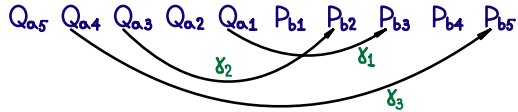
\begin{center}
\begin{tikzpicture}[scale=0.65]
		\draw   (0,2)  node{$\QQ_{a_5}$};
		\draw   (1,2)  node{$\QQ_{a_4}$};
		\draw   (2,2)   node {$\QQ_{a_3}$};
		\draw   (3,2)  node {$\QQ_{a_2}$};
		\draw   (4,2)  node {$\QQ_{a_1}$};
		\draw   (5,2)  node {$\PP_{b_1}$};
		\draw   (6,2)  node {$\PP_{b_2}$};
		\draw   (7,2)  node {$\PP_{b_3}$};
		\draw   (8,2)  node {$\PP_{b_4}$};
		\draw   (9,2)  node {$\PP_{b_5}$};
			\draw[->] (1.3,1.7) 
		to [out=-40, in=180] (5,0.3) to [out=0, in=220] node[label=below:{$\gamma_3$}, dot, pos=0.4] {} (8.7,1.7); 
		\draw[->] (2.3,1.7) 
		to [out=-50, in=180] node[label=above:{$\gamma_2$}, dot, pos=1] {} (4,0.6) to [out=0, in=230] (5.7,1.7); 
		\draw[->] (4.3,1.7) 
		to [out=-50, in=180]  (5.5,1) to [out=0, in=230] node[label=below:{$\gamma_1$}, dot, pos=0.3] {} (6.7,1.7); 
	\end{tikzpicture}
\end{center}
\item Next, begin drawing the commutation strands. From each of
the remaining $\QQ$s and $\PP$s draw a vertical line down to some
horizontal level line below all the annihilation strands. 
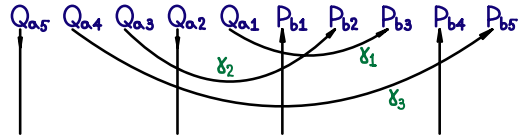
\begin{center}
\begin{tikzpicture}[scale=0.65]
		\draw   (0,2)  node{$\QQ_{a_5}$};
	\draw   (1,2)  node{$\QQ_{a_4}$};
	\draw   (2,2)   node {$\QQ_{a_3}$};
	\draw   (3,2)  node {$\QQ_{a_2}$};
	\draw   (4,2)  node {$\QQ_{a_1}$};
	\draw   (5,2)  node {$\PP_{b_1}$};
	\draw   (6,2)  node {$\PP_{b_2}$};
	\draw   (7,2)  node {$\PP_{b_3}$};
	\draw   (8,2)  node {$\PP_{b_4}$};
	\draw   (9,2)  node {$\PP_{b_5}$};
	\draw[->] (0,1.7) -- (0,0); 
	\draw[->] (3,1.7) -- (3,0); 
	\draw[->] (5,0) -- (5,1.7); 
	\draw[->] (8,0) -- (8,1.7); 
	\draw[->] (1.3,1.7) 
to [out=-40, in=180] (5,0.3) to [out=0, in=220] node[label=below:{$\gamma_3$}, dot, pos=0.4] {} (8.7,1.7); 
\draw[->] (2.3,1.7) 
to [out=-50, in=180] node[label=above:{$\gamma_2$}, dot, pos=1] {} (4,0.6) to [out=0, in=230] (5.7,1.7); 
\draw[->] (4.3,1.7) 
to [out=-50, in=180]  (5.5,1) to [out=0, in=230] node[label=below:{$\gamma_1$}, dot, pos=0.3] {} (6.7,1.7);
\end{tikzpicture}
\end{center}
\item Finish by starting at that level and joing the commutation strand  
of each $\QQ_{a_i}$ or $\PP_{b_j}$ at the top of the diagram
to the same $\QQ_{a_i}$ or $\PP_{b_j}$ at the bottom of the diagram in a
straight line. 
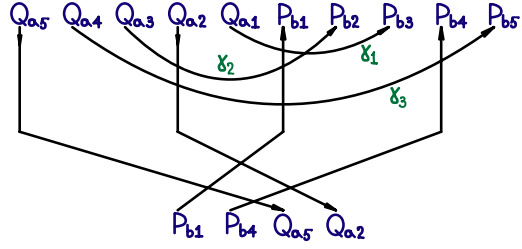
\begin{center}
\begin{tikzpicture}[scale=0.65]
	\draw   (0,2)  node{$\QQ_{a_5}$};
\draw   (1,2)  node{$\QQ_{a_4}$};
\draw   (2,2)   node {$\QQ_{a_3}$};
\draw   (3,2)  node {$\QQ_{a_2}$};
\draw   (4,2)  node {$\QQ_{a_1}$};
\draw   (5,2)  node {$\PP_{b_1}$};
\draw   (6,2)  node {$\PP_{b_2}$};
\draw   (7,2)  node {$\PP_{b_3}$};
\draw   (8,2)  node {$\PP_{b_4}$};
\draw   (9,2)  node {$\PP_{b_5}$};
	\draw (3,-2.5) node {$\PP_{b_1}$};
	\draw (4,-2.5) node {$\PP_{b_4}$};
	\draw (5,-2.5) node {$\QQ_{a_5}$};
	\draw (6,-2.5) node {$\QQ_{a_2}$};
	\draw[->] (0,1.7) -- (0,0) -- (4.7,-2); 
	\draw[->] (3,1.7) -- (3,0) -- (5.7,-2); 
	\draw[->] (3,-2) -- (5,0) -- (5,1.7) ; 
	\draw[->] (4,-2)  -- (8,0) -- (8,1.7) ; 
		\draw[->] (1.3,1.7) 
	to [out=-40, in=180] (5,0.3) to [out=0, in=220] node[label=below:{$\gamma_3$}, dot, pos=0.4] {} (8.7,1.7); 
	\draw[->] (2.3,1.7) 
	to [out=-50, in=180] node[label=above:{$\gamma_2$}, dot, pos=1] {} (4,0.6) to [out=0, in=230] (5.7,1.7); 
	\draw[->] (4.3,1.7) 
	to [out=-50, in=180]  (5.5,1) to [out=0, in=230] node[label=below:{$\gamma_1$}, dot, pos=0.3] {} (6.7,1.7);
\end{tikzpicture}
\end{center}
\end{enumerate}
\end{Definition}

Denote the component of $\phi$ on each summand of 
$\Xi_{\PP\QQ}(\hat{k})\bigl(\underline{a} \otimes \underline{b} \bigr)$
by the same diagram as its adjoint
\eqref{eqn-the-adjoint-of-the-natural-transformation-phi} only
with each annihilation strand decorated 
by $(a_{i_l}, b_{j_{\sigma(l)}})$ instead of $\gamma_l$:
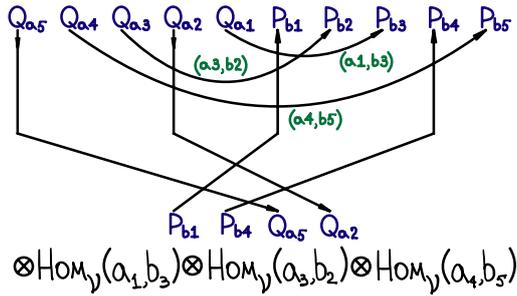
\begin{center}
\begin{tikzpicture}[scale=0.65]
		\draw   (0,2)  node{$\QQ_{a_5}$};
	\draw   (1,2)  node{$\QQ_{a_4}$};
	\draw   (2,2)   node {$\QQ_{a_3}$};
	\draw   (3,2)  node {$\QQ_{a_2}$};
	\draw   (4,2)  node {$\QQ_{a_1}$};
	\draw   (5,2)  node {$\PP_{b_1}$};
	\draw   (6,2)  node {$\PP_{b_2}$};
	\draw   (7,2)  node {$\PP_{b_3}$};
	\draw   (8,2)  node {$\PP_{b_4}$};
	\draw   (9,2)  node {$\PP_{b_5}$};
	\draw (3,-2.5) node {$\PP_{b_1}$};
	\draw (4,-2.5) node {$\PP_{b_4}$};
	\draw (5,-2.5) node {$\QQ_{a_5}$};
	\draw (6,-2.5) node {$\QQ_{a_2}$};
	\draw[->] (0,1.7) -- (0,0) -- (4.7,-2); 
	\draw[->] (3,1.7) -- (3,0) -- (5.7,-2); 
	\draw[->] (3,-2) -- (5,0) -- (5,1.7) ; 
	\draw[->] (4,-2)  -- (8,0) -- (8,1.7) ; 
		\draw[->] (1.3,1.7) 
	to [out=-40, in=180] (5,0.3) to [out=0, in=220] node[label={[label distance=-1mm]270:{${\scriptstyle (a_4,b_5)}$}}, dot, pos=0.4] {} (8.7,1.7); 
	\draw[->] (2.3,1.7) 
	to [out=-50, in=180] node[label={[label distance=-1mm]90:{${\scriptstyle (a_3,b_2)}$}}, dot, pos=1] {} (4,0.6) to [out=0, in=230] (5.7,1.7); 
	\draw[->] (4.3,1.7) 
	to [out=-50, in=180]  (5.5,1) to [out=0, in=230] node[label={[label distance=-1mm]270:{${\scriptstyle (a_1,b_3)}$}}, dot, pos=0.3] {} (6.7,1.7);
	\draw (4.5,-3.5) node {$	\otimes \Hom_\basecat(a_1,b_3)\otimes \Hom_\basecat(a_3,b_2)\otimes \Hom_\basecat(a_4,b_5) $};

\end{tikzpicture}

\end{center}

\begin{Proposition}
\label{prps-natural-transformation-phi-is-homotopy-equivalence}
For any 
$\underline{a} = (a_n \otimes \dots \otimes a_1) \in
\sym^{n}\basecat^{\opp}$ and 
$\underline{b} = (b_1 \otimes \dots \otimes b_m) \in  
\otimes \sym^{m}\basecat$, 
the $2$-morphism 
$\phi\colon \bigoplus \Xi_{\PP\QQ}(\hat{k})\bigl(\underline{a} \otimes \underline{b} \bigr)
\longrightarrow 
\Xi_{\QQ\PP}\bigl(\underline{a} \otimes \underline{b}\bigr)$
of Definition~\ref{defn-natural-transformation-phi-termwise}
equals in $\hcat\basecat$ 
the one constructed by the commutation-annihilation procedure
of Definition~\ref{defn-commutation-annihilation-procedure}. In 
particular, it is a homotopy equivalence. 
\end{Proposition}
\begin{proof}
This follows from the pitchfork and triple move relations 
in $\hcat\basecat$ \cite[Lemma 5.5]{gyenge2021heisenberg}. 

We can prove it separately for each summand of 
$\bigoplus \Xi_{\PP\QQ}(\hat{k})\bigl(\underline{a} \otimes \underline{b}
\bigr)$. For a given summand,
the commutation-annihilation procedure construct the following
$2$-morphism into 
$\Xi_{\QQ\PP}\bigl(\underline{a} \otimes \underline{b}\bigr)$. 
We begin at the top of the diagram at $\QQ_{a_1}$ and 
perform rightward crossings 
$\begin{tikzpicture}[baseline={(0,0.15)}, scale=0.5]
\draw[->] (0.5,0) -- (1.5,1);
\draw[->] (0.5,1) -- (1.5,0);
\end{tikzpicture}$ until: 
\begin{itemize}
\item If the $\QQ_{a_1}$-strand is a commutation strand -- 
until it moves to the right of all $\PP$-strands, 
\item If the $\QQ_{a_1}$-strand is an annihilation strand -- 
until it reaches the $\PP_{b_j}$-strand it is paired with. 
The two strands are then terminated by a cup marked by $(a_1, b_j)$. 
\end{itemize}
We then repeat the same for $\QQ_{a_2}$, etc. 
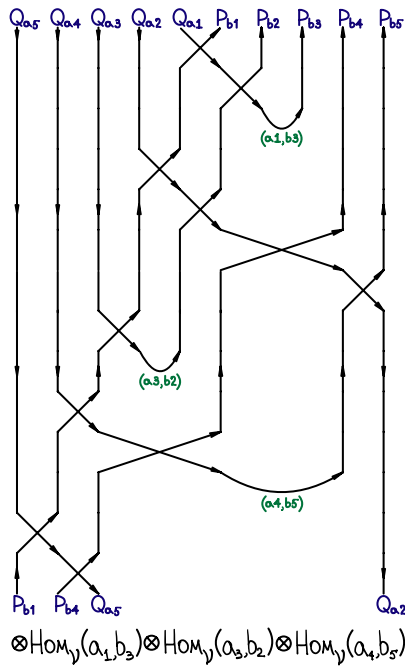
\begin{equation*}
\begin{minipage}{0.45\textwidth}
	\begin{tikzpicture}[scale=0.65]
		\draw   (0,2.25)  node{$\QQ_{a_5}$};
		\draw   (1,2.25)  node{$\QQ_{a_4}$};
		\draw   (2,2.25)   node {$\QQ_{a_3}$};
		\draw   (3,2.25)  node {$\QQ_{a_2}$};
		\draw   (4,2.25)  node {$\QQ_{a_1}$};
		\draw   (5,2.25)  node {$\PP_{b_1}$};
		\draw   (6,2.25)  node {$\PP_{b_2}$};
		\draw   (7,2.25)  node {$\PP_{b_3}$};
		\draw   (8,2.25)  node {$\PP_{b_4}$};
		\draw   (9,2.25)  node {$\PP_{b_5}$};
		\draw (0,-12.4) node {$\PP_{b_1}$};
		\draw (1,-12.4) node {$\PP_{b_4}$};
		\draw (2,-12.4) node {$\QQ_{a_5}$};
		\draw (9,-12.4) node {$\QQ_{a_2}$};
		\draw[->] (0,2) -- (0,-10);
		\draw[->] (0,-10)-- (1,-11);
		\draw[->] (1,-11) -- (2,-12); 
		\draw[->] (3,2) -- (3,-1)-- (4,-2);
		 \draw[->] (4,-2) -- (5,-3);
		 \draw[->] (5,-3) -- (8,-4);
		 \draw[->] (8,-4) -- (9,-5);
		 \draw[->] (9,-5) -- (9,-12);
		\draw[->] (0,-12) -- (0, -11);
		\draw[->] (0,-11) -- (1,-10); 
		\draw[->] (1,-10) -- (1,-8) -- (2,-7) ;
		\draw[->] (2,-7) -- (2,-6) --(3,-5) ; 
		\draw[->] (3,-5) -- (3,-2) --(4,-1) ;
		\draw[->] (4,-1) -- (4,1) -- (5,2); 
		\draw[->] (1,-12) -- (2,-11);
		\draw[->] (2,-11) -- (2,-9) --(5,-8) ; 
		\draw[->] (5,-8) -- (5,-4) -- (8,-3);
		\draw[->] (8,-3)  -- (8,2);
		%\draw[->] (4,-2)  -- (6,0) -- (8,1.5) -- (8,1.7) ; 
		\draw[->] (2,2) -- (2,-5) -- (3,-6);
		\draw[->] (4,-6) -- (4,-3) -- (5,-2);
		\draw[->] (5,-2) -- (5,0) -- (6,1);
		\draw[->] (6,1) -- (6,2);
		\draw[->] (1,2) -- (1,-7) -- (2,-8);
		\draw[->] (2,-8) -- (5,-9);
		\draw[->] (8,-9) -- (8,-5) -- (9,-4);
		\draw[->] (9,-4) -- (9,2);
		\draw[->] (4,2) -- (5,1);
		\draw[->] (5,1) -- (6,0);
		\draw[->] (7,0) -- (7,2);
		\draw[->] (5,-9) 
		to [out=-50, in=180] (6.5,-9.7) to [out=0, in=230] node[label={[label distance=-1mm]270:{${\scriptstyle (a_4,b_5)}$}}, dot, pos=0] {} (8,-9); 
		\draw[->] (3,-6) 
		to [out=-50, in=180] node[label={[label distance=-1mm]270:{${\scriptstyle (a_3,b_2)}$}}, dot, pos=1] {} (3.5,-6.4) to [out=0, in=230] (4,-6); 
		\draw[->] (6,0) 
		to [out=-50, in=180]  (6.5,-0.4) to [out=0, in=230] node[label={[label distance=-1mm]270:{${\scriptstyle (a_1,b_3)}$}}, dot, pos=0] {} (7,0);
		\draw (4.5,-13) node {${\scriptstyle	\otimes \Hom_\basecat(a_1,b_3)\otimes \Hom_\basecat(a_3,b_2)\otimes \Hom_\basecat(a_4,b_5) }$};		
	\end{tikzpicture} 
\end{minipage}
=
\begin{minipage}{0.45\textwidth}
	
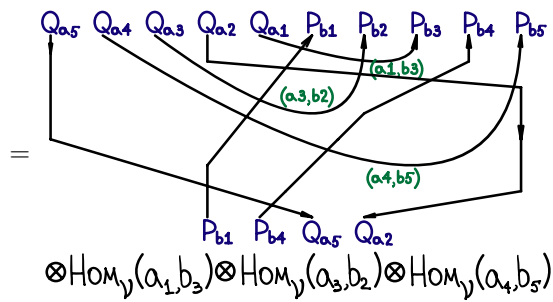
\begin{tikzpicture}[scale=0.65]
		\draw   (0,2.25)  node{$\QQ_{a_5}$};
		\draw   (1,2.25)  node{$\QQ_{a_4}$};
		\draw   (2,2.25)   node {$\QQ_{a_3}$};
		\draw   (3,2.25)  node {$\QQ_{a_2}$};
		\draw   (4,2.25)  node {$\QQ_{a_1}$};
		\draw   (5,2.25)  node {$\PP_{b_1}$};
		\draw   (6,2.25)  node {$\PP_{b_2}$};
		\draw   (7,2.25)  node {$\PP_{b_3}$};
		\draw   (8,2.25)  node {$\PP_{b_4}$};
		\draw   (9,2.25)  node {$\PP_{b_5}$};
		\draw (3,-2.4) node {$\PP_{b_1}$};
		\draw (4,-2.4) node {$\PP_{b_4}$};
		\draw (5,-2.4) node {$\QQ_{a_5}$};
		\draw (6,-2.4) node {$\QQ_{a_2}$};
		\draw[->] (0,2) -- (0,0) -- (4.7,-2); 
		\draw[->] (3,2) -- (3,1)-- (9,0.2) -- (9,-1.3)-- (5.7,-2); 
		\draw[->] (3,-2) -- (3,0) -- (5,2) ; 
		\draw[->] (4,-2)  -- (6,0) -- (8,1.5) -- (8,2) ; 
		\draw[->] (1,2) 
		to [out=-50, in=180] (5,-0.8) to [out=0, in=230] node[label={[label distance=-1mm]270:{${\scriptstyle (a_4,b_5)}$}}, dot, pos=0.4] {} (9,2); 
		\draw[->] (2,2) 
		to [out=-50, in=180] node[label={[label distance=-1mm]90:{${\scriptstyle (a_3,b_2)}$}}, dot, pos=1] {} (4,0.1) to [out=0, in=230] (6,2); 
		\draw[->] (4,2) 
		to [out=-50, in=180]  (5.5,1.2) to [out=0, in=230] node[label={[label distance=-1mm]270:{${\scriptstyle (a_1,b_3)}$}}, dot, pos=0.3] {} (7,2);
		\draw (4.5,-3) node {${\scriptstyle	\otimes \Hom_\basecat(a_1,b_3)\otimes \Hom_\basecat(a_3,b_2)\otimes \Hom_\basecat(a_4,b_5) }$};
		
	\end{tikzpicture}
\end{minipage}
\end{equation*}

If we take this planar diagram and consider all the annihilation strands
separately and all the commutation strands separately, then 
the two configurations match, possibly up to some pitchfork relations, 
the configurations of the annihilation and of the commutation strands 
in the planar diagram defining the corresponding component of 
the $2$-morphism $\phi$ of Definition 
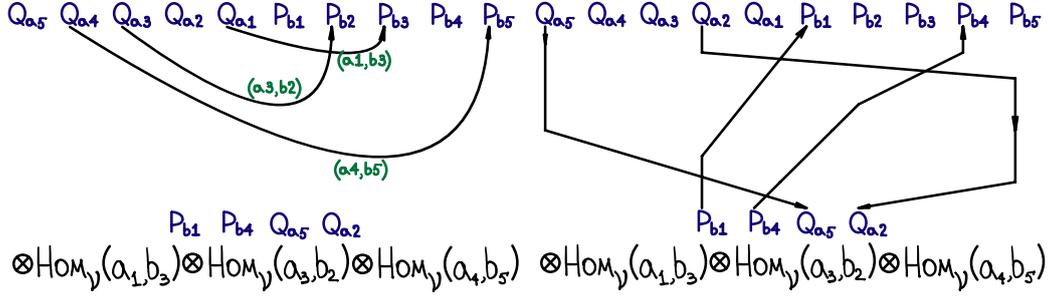
\ref{defn-natural-transformation-phi-termwise}. 
\begin{equation*}
\begin{minipage}{0.45\textwidth}
	\begin{tikzpicture}[scale=0.65]
	\draw   (0,2.25)  node{$\QQ_{a_5}$};
	\draw   (1,2.25)  node{$\QQ_{a_4}$};
	\draw   (2,2.25)   node {$\QQ_{a_3}$};
	\draw   (3,2.25)  node {$\QQ_{a_2}$};
	\draw   (4,2.25)  node {$\QQ_{a_1}$};
	\draw   (5,2.25)  node {$\PP_{b_1}$};
	\draw   (6,2.25)  node {$\PP_{b_2}$};
	\draw   (7,2.25)  node {$\PP_{b_3}$};
	\draw   (8,2.25)  node {$\PP_{b_4}$};
	\draw   (9,2.25)  node {$\PP_{b_5}$};
	\draw (3,-2.4) node {$\PP_{b_1}$};
	\draw (4,-2.4) node {$\PP_{b_4}$};
	\draw (5,-2.4) node {$\QQ_{a_5}$};
	\draw (6,-2.4) node {$\QQ_{a_2}$};
%	\draw[->] (0,2) -- (0,0) -- (4.7,-2); 
%	\draw[->] (3,2) -- (3,1)-- (9,0.2) -- (9,-1.3)-- (5.7,-2); 
%	\draw[->] (3,-2) -- (3,0) -- (5,2) ; 
%	\draw[->] (4,-2)  -- (6,0) -- (8,1.5) -- (8,2) ; 
	\draw[->] (1,2) 
	to [out=-50, in=180] (5,-0.8) to [out=0, in=230] node[label={[label distance=-1mm]270:{${\scriptstyle (a_4,b_5)}$}}, dot, pos=0.4] {} (9,2); 
	\draw[->] (2,2) 
	to [out=-50, in=180] node[label={[label distance=-1mm]90:{${\scriptstyle (a_3,b_2)}$}}, dot, pos=1] {} (4,0.1) to [out=0, in=230] (6,2); 
	\draw[->] (4,2) 
	to [out=-50, in=180]  (5.5,1.2) to [out=0, in=230] node[label={[label distance=-1mm]270:{${\scriptstyle (a_1,b_3)}$}}, dot, pos=0.3] {} (7,2);
	\draw (4.5,-3) node {${\scriptstyle	\otimes \Hom_\basecat(a_1,b_3)\otimes \Hom_\basecat(a_3,b_2)\otimes \Hom_\basecat(a_4,b_5) }$};
	
\end{tikzpicture}\end{minipage}
	\begin{minipage}{0.45\textwidth}
	\begin{tikzpicture}[scale=0.65]
	\draw   (0,2.25)  node{$\QQ_{a_5}$};
	\draw   (1,2.25)  node{$\QQ_{a_4}$};
	\draw   (2,2.25)   node {$\QQ_{a_3}$};
	\draw   (3,2.25)  node {$\QQ_{a_2}$};
	\draw   (4,2.25)  node {$\QQ_{a_1}$};
	\draw   (5,2.25)  node {$\PP_{b_1}$};
	\draw   (6,2.25)  node {$\PP_{b_2}$};
	\draw   (7,2.25)  node {$\PP_{b_3}$};
	\draw   (8,2.25)  node {$\PP_{b_4}$};
	\draw   (9,2.25)  node {$\PP_{b_5}$};
	\draw (3,-2.4) node {$\PP_{b_1}$};
	\draw (4,-2.4) node {$\PP_{b_4}$};
	\draw (5,-2.4) node {$\QQ_{a_5}$};
	\draw (6,-2.4) node {$\QQ_{a_2}$};
	\draw[->] (0,2) -- (0,0) -- (4.7,-2); 
	\draw[->] (3,2) -- (3,1)-- (9,0.2) -- (9,-1.3)-- (5.7,-2); 
	\draw[->] (3,-2) -- (3,0) -- (5,2) ; 
	\draw[->] (4,-2)  -- (6,0) -- (8,1.5) -- (8,2) ; 
%	\draw[->] (1,2) 
%	to [out=-50, in=180] (5,-0.8) to [out=0, in=230] node[label={[label distance=-1mm]270:{${\scriptstyle (a_4,b_5)}$}}, dot, pos=0.4] {} (9,2); 
%	\draw[->] (2,2) 
%	to [out=-50, in=180] node[label={[label distance=-1mm]90:{${\scriptstyle (a_3,b_2)}$}}, dot, pos=1] {} (4,0.1) to [out=0, in=230] (6,2); 
%	\draw[->] (4,2) 
%	to [out=-50, in=180]  (5.5,1.2) to [out=0, in=230] node[label={[label distance=-1mm]270:{${\scriptstyle (a_1,b_3)}$}}, dot, pos=0.3] {} (7,2);
	\draw (4.5,-3) node {${\scriptstyle	\otimes \Hom_\basecat(a_1,b_3)\otimes \Hom_\basecat(a_3,b_2)\otimes \Hom_\basecat(a_4,b_5) }$};
	
\end{tikzpicture}
\end{minipage}
\end{equation*}

Now taking all the commutation strand crossings in the
commutation-annihilation planar diagram and using triple moves to
commute them downwards past all the annihilation strands produces the
planar diagram in the definition of $\phi$.
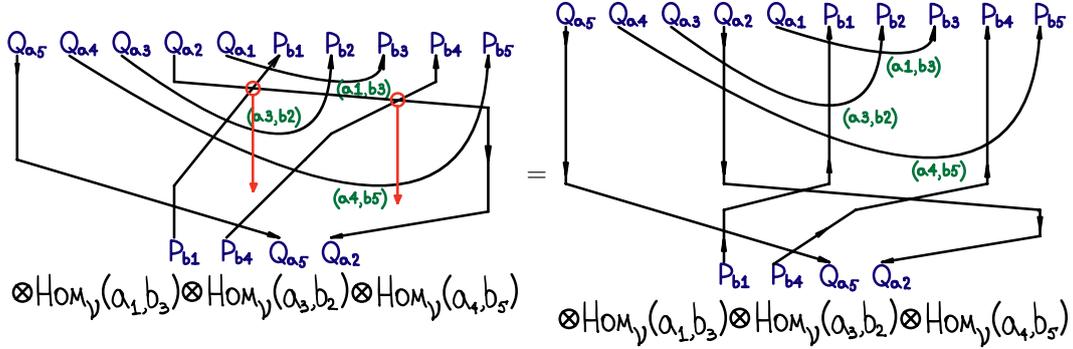
\begin{equation*}
\begin{minipage}{0.45\textwidth}
	\begin{tikzpicture}[scale=0.65]
	\draw   (0,2.25)  node{$\QQ_{a_5}$};
	\draw   (1,2.25)  node{$\QQ_{a_4}$};
	\draw   (2,2.25)   node {$\QQ_{a_3}$};
	\draw   (3,2.25)  node {$\QQ_{a_2}$};
	\draw   (4,2.25)  node {$\QQ_{a_1}$};
	\draw   (5,2.25)  node {$\PP_{b_1}$};
	\draw   (6,2.25)  node {$\PP_{b_2}$};
	\draw   (7,2.25)  node {$\PP_{b_3}$};
	\draw   (8,2.25)  node {$\PP_{b_4}$};
	\draw   (9,2.25)  node {$\PP_{b_5}$};
	\draw (3,-2.4) node {$\PP_{b_1}$};
	\draw (4,-2.4) node {$\PP_{b_4}$};
	\draw (5,-2.4) node {$\QQ_{a_5}$};
	\draw (6,-2.4) node {$\QQ_{a_2}$};
	\draw[->] (0,2) -- (0,0) -- (4.7,-2); 
	\draw[->] (3,2) -- (3,1)-- (9,0.2) -- (9,-1.3)-- (5.7,-2); 
	\draw[->] (3,-2) -- (3,0) -- (5,2) ; 
	\draw[->] (4,-2)  -- (6,0) -- (8,1.5) -- (8,2) ; 
	\draw[->] (1,2) 
	to [out=-50, in=180] (5,-0.8) to [out=0, in=230] node[label={[label distance=-1mm]270:{${\scriptstyle (a_4,b_5)}$}}, dot, pos=0.4] {} (9,2); 
	\draw[->] (2,2) 
	to [out=-50, in=180] node[label={[label distance=-1mm]90:{${\scriptstyle (a_3,b_2)}$}}, dot, pos=1] {} (4,0.1) to [out=0, in=230] (6,2); 
	\draw[->] (4,2) 
	to [out=-50, in=180]  (5.5,1.2) to [out=0, in=230] node[label={[label distance=-1mm]270:{${\scriptstyle (a_1,b_3)}$}}, dot, pos=0.3] {} (7,2);
	\draw (4.5,-3) node {${\scriptstyle	\otimes \Hom_\basecat(a_1,b_3)\otimes \Hom_\basecat(a_3,b_2)\otimes \Hom_\basecat(a_4,b_5) }$};
	\draw[color=red] (3.9,0.9) circle (0.2);
	\draw[->,color=red] (3.9,0.7) -- (3.9,-1);
	\draw[color=red] (6.7,0.5) circle (0.2);
	\draw[->,color=red] (6.7,0.3) -- (6.7,-1.2);
\end{tikzpicture}
\end{minipage}
=
\begin{minipage}{0.45\textwidth}
	\begin{tikzpicture}[scale=0.65]
		\draw   (0,2.25)  node{$\QQ_{a_5}$};
		\draw   (1,2.25)  node{$\QQ_{a_4}$};
		\draw   (2,2.25)   node {$\QQ_{a_3}$};
		\draw   (3,2.25)  node {$\QQ_{a_2}$};
		\draw   (4,2.25)  node {$\QQ_{a_1}$};
		\draw   (5,2.25)  node {$\PP_{b_1}$};
		\draw   (6,2.25)  node {$\PP_{b_2}$};
		\draw   (7,2.25)  node {$\PP_{b_3}$};
		\draw   (8,2.25)  node {$\PP_{b_4}$};
		\draw   (9,2.25)  node {$\PP_{b_5}$};
		\draw (3,-3.4) node {$\PP_{b_1}$};
		\draw (4,-3.4) node {$\PP_{b_4}$};
		\draw (5,-3.4) node {$\QQ_{a_5}$};
		\draw (6,-3.4) node {$\QQ_{a_2}$};
		\draw[->] (0,2) -- (0,-1) -- (4.7,-3); 
		\draw[->] (3,2) -- (3,-1.2)-- (9,-2) -- (9,-2.3)-- (5.7,-3); 
		\draw[->] (3,-3) -- (3,-2) -- (5,-1)-- (5,2) ; 
		\draw[->] (4,-3)  -- (5,-2) -- (8,-1) -- (8,2) ; 
		\draw[->] (1,2) 
		to [out=-50, in=180] (5,-0.8) to [out=0, in=230] node[label={[label distance=-1mm]270:{${\scriptstyle (a_4,b_5)}$}}, dot, pos=0.4] {} (9,2); 
		\draw[->] (2,2) 
		to [out=-50, in=180] node[label={[label distance=-1mm]90:{${\scriptstyle (a_3,b_2)}$}}, dot, pos=1] {} (4,0.1) to [out=0, in=230] (6,2); 
		\draw[->] (4,2) 
		to [out=-50, in=180]  (5.5,1.2) to [out=0, in=230] node[label={[label distance=-1mm]270:{${\scriptstyle (a_1,b_3)}$}}, dot, pos=0.3] {} (7,2);
		\draw (4.5,-4) node {${\scriptstyle	\otimes \Hom_\basecat(a_1,b_3)\otimes \Hom_\basecat(a_3,b_2)\otimes \Hom_\basecat(a_4,b_5) }$};
		
	\end{tikzpicture}
\end{minipage}
\end{equation*}
We conclude that the two planar diagrams define the same $2$-morphism 
in $\hcat\basecat$.  
\end{proof}

The following can be viewed as a functorial 
categorification of the PQ Heisenberg relation
\cite[\S2.2]{gyenge2021heisenberg}:

\begin{Theorem}
\label{theorem-functorial-categorification-of-the-PQ-Heisenberg-relation}
The $2$-morphisms $\phi$ of Definition 
\ref{defn-natural-transformation-phi-termwise}
define a DG natural transformation 
\begin{equation}
\label{eqn-functorial-categorification-of-the-PQ-Heisenberg-relation}
\phi\colon \bigoplus \Xi_{\PP\QQ}(\hat{k})
\longrightarrow 
\Xi_{\QQ\PP}
\end{equation}
of DG functors
$\sym^{n}\basecat^{\opp} \otimes \sym^{m}\basecat \rightarrow 
\homm_{\hcat\basecat}(0,m-n)$. By Proposition~\ref{prps-natural-transformation-phi-is-homotopy-equivalence}, 
$\phi$ is a homotopy equivalence. 
\end{Theorem}

A non-functorial categorification of the Heisenberg relation 
appeared in \cite[Theorem 6.3]{gyenge2021heisenberg}. However, it relates
the symmetrised elements $\PP_a^{(n)}$ and $\QQ_b^{(m)}$ and
as these are not functorial in $a,b \in \basecat$ for $n,m > 1$, 
there was little hope of making it functorial directly. Instead, 
we used its basic case $n = m = 1$ to iteratively construct
the present, functorial categorification
\eqref{eqn-functorial-categorification-of-the-PQ-Heisenberg-relation}. 
It is clear that applying  
\eqref{eqn-functorial-categorification-of-the-PQ-Heisenberg-relation}
to the product of the symmetrised powers 
$a^{(n)} \in \sym^{n}\basecat^{\opp} $ and $b^{(m)} \in
\sym^{m}\basecat$, defined via the twisted complexes analogous
to those defining $\PP_a^{(n)}$ and $\QQ_b^{(m)}$ in 
\cite[Definition 6.2]{gyenge2021heisenberg}, recovers the
non-functorial categorification of \cite[Theorem 6.3]{gyenge2021heisenberg}. 

\begin{proof}[Proof of Theorem
\ref{theorem-functorial-categorification-of-the-PQ-Heisenberg-relation}]

We need to show that for all 
$$ \underline{a} = (a_n \otimes \dots \otimes a_1) \text{ and }
\underline{a}' = (a'_n \otimes \dots \otimes a'_1) \in
\sym^{n}\basecat^{\opp}, $$
$$ \underline{b} = (b_1 \otimes \dots \otimes b_m) \text{ and }
\underline{b'} = (b'_1 \otimes \dots \otimes b'_m)
 \in \otimes \sym^{m}\basecat $$
and each morphism 
\begin{equation}
\label{eqn-morphisms-in-S^n-basecatopp-S^m-basecat-for-functoriality-of-phi}
\underline{\alpha} \otimes \underline{\beta} 
\colon \underline{a} \otimes \underline{b}
\rightarrow \underline{a}' \otimes \underline{b}'
\quad \quad \text{ in } \sym^{n}\basecat^{\opp} \otimes
\sym^{m}\basecat
\end{equation}
the corresponding square commutes in
$\homm_{\hcat\basecat}(0,m-n)$:
\begin{equation}
\label{eqn-commutative-square-for-the-functoriality-of-phi}
\begin{tikzcd}[column sep = 3cm]
\Xi_{\QQ\PP}\left(\underline{a} \otimes \underline{b}\right)
\ar{r}{\Xi_{\QQ\PP}(\underline{\alpha} \otimes \underline{\beta})}
&
\Xi_{\QQ\PP}\left(\underline{a}' \otimes \underline{b}' \right)
\\
\bigoplus \Xi_{\PP\QQ}(\hat{k})\left(\underline{a} \otimes \underline{b}\right)
\ar{r}{\bigoplus \Xi_{\PP\QQ}(\hat{k})(\underline{\alpha} \otimes \underline{\beta})}
\ar{u}{\phi_{\underline{a} \otimes \underline{b}}}
&
\bigoplus \Xi_{\PP\QQ}(\hat{k})\left(\underline{a}' \otimes
\underline{b}'\right).
\ar{u}{\phi_{\underline{a}' \otimes \underline{b}'}}
\end{tikzcd}
\end{equation}
The action of $\Xi_{\QQ\PP}$ on the 
morphism spaces of $\sym^{n}\basecat^{\opp} \otimes
\sym^{m}\basecat$ is clear. The functor $\Xi_{\PP\QQ}(\hat{k})$
is defined in terms of the functor $(\hat{k})$ whose
definition is more involved. However, the action of $(\hat{k})$
on the morphism spaces of $\sym^{n}\basecat^{\opp} \otimes
\sym^{m}\basecat$ is described explicitly by Lemma
\ref{lemma-explicit-description-of-the-functor-hatk}. We thus
prove the functoriality of $\phi$ by explicitly verifying for 
each generating morphism of 
$\sym^{n}\basecat^{\opp} \otimes \sym^{m}\basecat$ that 
\eqref{eqn-commutative-square-for-the-functoriality-of-phi}
commutes
in $\homm_{\hcat\basecat}(0,m-n)$. This reduces, once more,
to applying the pitchfork and triple move relations
\cite[Lemma 5.5]{gyenge2021heisenberg} and the dot sliding 
relations \cite[Lemma 5.4]{gyenge2021heisenberg}. 

Indeed, it suffices to prove that
\eqref{eqn-commutative-square-for-the-functoriality-of-phi}
commutes for each direct summand of $\bigoplus
\Xi_{\PP\QQ}(\hat{k})\left(\underline{a} \otimes \underline{b}\right)$. 
Fix $0 \leq k \leq \min(n,m)$. 
By Lemma \ref{lemma-explicit-description-of-the-functor-hatk}, 
each $\Xi_{\PP\QQ}(\hat{k})\left(\underline{a} \otimes
\underline{b}\right)$ is itself a direct sum 
indexed by choices of 
$1 \leq i_1 < \dots < i_k \leq n$, $1 \leq j_1 < \dots < j_k
\leq m$ and $\sigma \in S_k$. Fix any such 
$(\underline{\iota}, \underline{j}, \sigma)$. 

The morphisms 
\eqref{eqn-morphisms-in-S^n-basecatopp-S^m-basecat-for-functoriality-of-phi}
are generated by composition from the following four basic types:
\begin{enumerate}
\item
\label{item-four-basic-morphisms-in-S^nVoppS^mV-Vopp-morphism}
 $\underline{\alpha} \otimes \id$
with $\underline{\alpha} = \id^{\otimes (n-i)} \otimes \alpha_i
\otimes \id^{\otimes(i-1)}$ for some $1 \leq i \leq n$ and
$\alpha_i \in \homm_{\basecat}(a'_i, a_i)$, 
\item 
\label{item-four-basic-morphisms-in-S^nVoppS^mV-V-morphism}
$\id \otimes \underline{\beta}$
with $\underline{\beta} = \id^{\otimes (j-1)} \otimes \beta_j
\otimes \id^{\otimes(m-j)}$ for some $1 \leq j \leq m$ and
$\beta \in \homm_{\basecat}(b_i, b'_i)$, 
\item 
\label{item-four-basic-morphisms-in-S^nVoppS^mV-S_n-element}
$\eta \otimes \id$ with $\eta = (i(i+1)) \in S_n$ for some 
$1 \leq i \leq n-1$, 
\item 
\label{item-four-basic-morphisms-in-S^nVoppS^mV-S_m-element}
$\id \otimes \zeta$ with $\zeta = (j(j+1)) \in S_m$
for some $1 \leq j \leq m-1$. 
\end{enumerate}
We only give the proofs for the types
\ref{item-four-basic-morphisms-in-S^nVoppS^mV-V-morphism}
and 
\ref{item-four-basic-morphisms-in-S^nVoppS^mV-S_m-element}, the
proofs for the other two types are similar. 

\underline{Morphisms of type
\ref{item-four-basic-morphisms-in-S^nVoppS^mV-V-morphism}}:

Consider $\phi_{\underline{a} \otimes \underline{b}}$ restricted
to the $(\underline{\iota}, \underline{j}, \sigma)$ direct summand 
of $\Xi_{\PP\QQ}(\hat{k})\left(\underline{a} \otimes
\underline{b}\right)$. There are two cases:

\begin{enumerate}
\item {\em The strand entering $\PP_{b_j}$ is a commutation strand ($j
\notin \{j_1, \dots, j_k\})$: }

Then on the $(\underline{\iota}, \underline{j}, \sigma)$ direct
summand
$\Xi_{\QQ\PP}(\id \otimes \underline{\beta}) \circ
\phi_{\underline{a} \otimes \underline{b}}$ and
$\phi_{\underline{a}' \otimes \underline{b}'} \circ \bigoplus
\Xi_{\PP\QQ}(\hat{k})(\id \otimes \underline{\beta})$, 
the two compositions around the upper-left and the lower-right
halves of the square 
\eqref{eqn-commutative-square-for-the-functoriality-of-phi}, 
are defined by two planar diagrams which only differ in the followings parts:
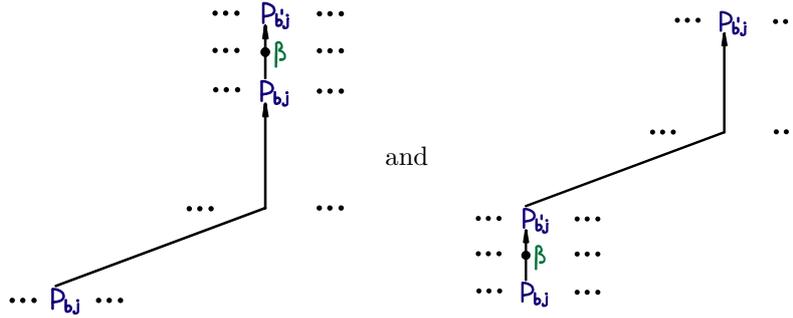
\begin{equation*}
	\begin{minipage}{0.35\textwidth}
		\begin{tikzpicture}[scale=0.7]
		\draw   (0,0)  node {$\PP_{b_j}$};
		\draw (4,4) node {$\PP_{b_j}$};	
		\draw (4,5.5) node {$\PP_{b'_j}$};	
		\draw (3,4) node {$\dots$};
		\draw (5,4) node {$\dots$};
		\draw (3,5.5) node {$\dots$};
		\draw (5,5.5) node {$\dots$};
		\draw (-1,0) node {$\dots$};
		\draw (1,0) node {$\dots$};
			\draw[->] (0.3,0.3) -- (4,2) -- (4,3.7); 
			\draw[->] (4,4.3) -- node[label={0:{$\beta$}}, dot, pos=0.5] {} (4,5.2); 
		\end{tikzpicture}
	\end{minipage}
	\quad \text{ and } \quad 
	\begin{minipage}{0.35\textwidth}
		\begin{tikzpicture}[scale=0.7]
			\draw   (0,0)  node {$\PP_{b'_j}$};
			\draw (4,4) node {$\PP_{b'_j}$};	
			\draw (0,-1.5) node {$\PP_{b_j}$};	
			\draw (-1,0) node {$\dots$};
			\draw (1,0) node {$\dots$};
			\draw (-1,-1.5) node {$\dots$};
			\draw (1,-1.5) node {$\dots$};
			\draw (3,4) node {$\dots$};
			\draw (5,4) node {$\dots$};
			\draw[->] (0.3,0.3) -- (4,2) -- (4,3.7); 
			\draw[->] (0,-1.2) -- node[label={0:{$\beta$}}, dot, pos=0.5] {} (0,-0.3); 
		\end{tikzpicture}
	\end{minipage}
\end{equation*}
Sliding the dot down the strand turns the left diagram 
into the right one. The identical remainder of the two diagrams, 
not depicted above, crosses the depicted parts. 
As it slides, the dot travels through these crossings. By 
the dot sliding relations \cite[Lemma 5.4]{gyenge2021heisenberg}, 
this doesn't change the $2$-morphism defined by the diagram. We
conclude that the corresponding $2$-morphisms are equal and thus 
\eqref{eqn-commutative-square-for-the-functoriality-of-phi}
commutes on this direct summand. 

Here and below, we only draw the relevant parts in which the diagrams 
differ nontrivially. For example, here the rest 
is identical since $\Xi_{\QQ\PP}(\id \otimes \underline{\beta})$ and 
$\Xi_{\PP\QQ}(\hat{k})(\id \otimes \underline{\beta})$
both consist of a row of parallel vertical strands one of which is
adorned by $\beta$, while $\phi_{\underline{a} \otimes \underline{b}}$
and $\phi_{\underline{a}' \otimes \underline{b}'}$ are restricted to 
the same indexed summand and hence are defined by the identical
diagrams. Thus, away from the $\beta$-adorned strand,  
we are post-composing the diagram of 
$\phi_{\underline{a} \otimes \underline{b}}$ and 
pre-composing the same diagram of 
$\phi_{\underline{a}' \otimes \underline{b}'}$ 
with a row of unadorned parallel vertical strands: 
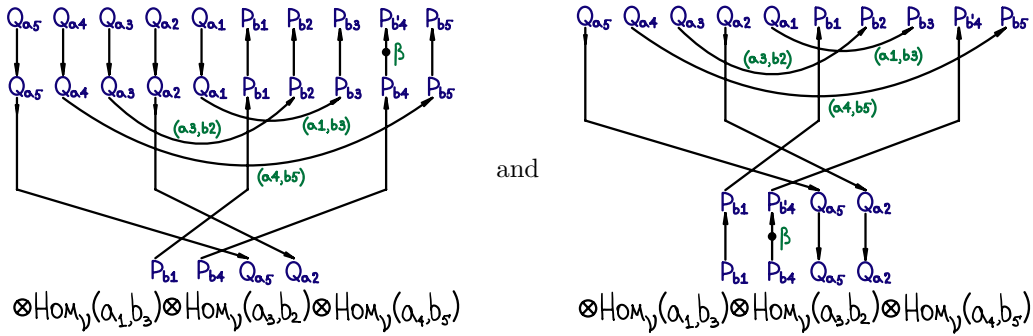
\begin{equation*}
\begin{minipage}{0.45\textwidth}
	\begin{tikzpicture}[scale=0.65]
		\draw   (0,3.5)  node{$\QQ_{a_5}$};
		\draw   (1,3.5)  node{$\QQ_{a_4}$};
		\draw   (2,3.5)   node {$\QQ_{a_3}$};
		\draw   (3,3.5)  node {$\QQ_{a_2}$};
		\draw   (4,3.5)  node {$\QQ_{a_1}$};
		\draw   (5,3.5)  node {$\PP_{b_1}$};
		\draw   (6,3.5)  node {$\PP_{b_2}$};
		\draw   (7,3.5)  node {$\PP_{b_3}$};
		\draw   (8,3.5)  node {$\PP_{b'_4}$};
		\draw   (9,3.5)  node {$\PP_{b_5}$};
		%\draw   (0,2.25)  node{$\QQ_{a_5}$};
		\draw[->] (0,3.2) -- (0,2.3);
		\draw[->] (1,3.2) -- (1,2.3);
		\draw[->] (2,3.2) -- (2,2.3);
		\draw[->] (3,3.2) -- (3,2.3);
		\draw[->] (4,3.2) -- (4,2.3);
		\draw[->] (5,3.2) -- (5,2.3);
		\draw[->] (6,2.3) -- (6,3.2);
		\draw[->] (7,2.3) -- (7,3.2);
		\draw[->] (8,2.3) -- node[label={[label distance=-1mm]0:{$\scriptstyle  \beta$}}, dot, pos=0.5] {} (8,3.2);
		\draw[->] (9,2.3) -- (9,3.2);
		\draw   (0,2)  node{$\QQ_{a_5}$};
\draw   (1,2)  node{$\QQ_{a_4}$};
\draw   (2,2)   node {$\QQ_{a_3}$};
\draw   (3,2)  node {$\QQ_{a_2}$};
\draw   (4,2)  node {$\QQ_{a_1}$};
\draw   (5,2)  node {$\PP_{b_1}$};
\draw   (6,2)  node {$\PP_{b_2}$};
\draw   (7,2)  node {$\PP_{b_3}$};
\draw   (8,2)  node {$\PP_{b_4}$};
\draw   (9,2)  node {$\PP_{b_5}$};
\draw (3,-2.5) node {$\PP_{b_1}$};
\draw (4,-2.5) node {$\PP_{b_4}$};
\draw (5,-2.5) node {$\QQ_{a_5}$};
\draw (6,-2.5) node {$\QQ_{a_2}$};
\draw[->] (0,1.7) -- (0,0) -- (4.7,-2); 
\draw[->] (3,1.7) -- (3,0) -- (5.7,-2); 
\draw[->] (3,-2) -- (5,0) -- (5,1.7) ; 
\draw[->] (4,-2)  -- (8,0) -- (8,1.7) ; 
\draw[->] (1.3,1.7) 
to [out=-40, in=180] (5,0.3) to [out=0, in=220] node[label={[label distance=-1mm]270:{${\scriptstyle (a_4,b_5)}$}}, dot, pos=0.4] {} (8.7,1.7); 
\draw[->] (2.3,1.7) 
to [out=-50, in=180] node[label={[label distance=-1mm]90:{${\scriptstyle (a_3,b_2)}$}}, dot, pos=1] {} (4,0.6) to [out=0, in=230] (5.7,1.7); 
\draw[->] (4.3,1.7) 
to [out=-50, in=180]  (5.5,1) to [out=0, in=230] node[label={[label distance=-1mm]270:{${\scriptstyle (a_1,b_3)}$}}, dot, pos=0.3] {} (6.7,1.7);
\draw (4.5,-3.2) node {${\scriptstyle	\otimes \Hom_\basecat(a_1,b_3)\otimes \Hom_\basecat(a_3,b_2)\otimes \Hom_\basecat(a_4,b_5) }$};
		
	\end{tikzpicture}
\end{minipage}
\quad \text{ and } \quad 
\begin{minipage}{0.45\textwidth}
	\begin{tikzpicture}[scale=0.65]
	\draw (3,-4) node {$\PP_{b_1}$};
	\draw (4,-4) node {$\PP_{b_4}$};
	\draw (5,-4) node {$\QQ_{a_5}$};
	\draw (6,-4) node {$\QQ_{a_2}$};
	\draw[->] (3,-3.7) -- (3,-2.8);
	\draw[->] (4,-3.7) -- node[label={[label distance=-1mm]0:{$\scriptstyle  \beta$}}, dot, pos=0.5] {} (4,-2.8);
	\draw[->] (5,-3.7) -- (5,-2.8);
	\draw[->] (6,-3.7) -- (6,-2.8);
	\draw   (0,2)  node{$\QQ_{a_5}$};
	\draw   (1,2)  node{$\QQ_{a_4}$};
	\draw   (2,2)   node {$\QQ_{a_3}$};
	\draw   (3,2)  node {$\QQ_{a_2}$};
	\draw   (4,2)  node {$\QQ_{a_1}$};
	\draw   (5,2)  node {$\PP_{b_1}$};
	\draw   (6,2)  node {$\PP_{b_2}$};
	\draw   (7,2)  node {$\PP_{b_3}$};
	\draw   (8,2)  node {$\PP_{b'_4}$};
	\draw   (9,2)  node {$\PP_{b_5}$};
	\draw (3,-2.5) node {$\PP_{b_1}$};
	\draw (4,-2.5) node {$\PP_{b'_4}$};
	\draw (5,-2.5) node {$\QQ_{a_5}$};
	\draw (6,-2.5) node {$\QQ_{a_2}$};
	\draw[->] (0,1.7) -- (0,0) -- (4.7,-2); 
	\draw[->] (3,1.7) -- (3,0) -- (5.7,-2); 
	\draw[->] (3,-2) -- (5,0) -- (5,1.7) ; 
	\draw[->] (4,-2)  -- (8,0) -- (8,1.7) ; 
	\draw[->] (1.3,1.7) 
	to [out=-40, in=180] (5,0.3) to [out=0, in=220] node[label={[label distance=-1mm]270:{${\scriptstyle (a_4,b_5)}$}}, dot, pos=0.4] {} (8.7,1.7); 
	\draw[->] (2.3,1.7) 
	to [out=-50, in=180] node[label={[label distance=-1mm]90:{${\scriptstyle (a_3,b_2)}$}}, dot, pos=1] {} (4,0.6) to [out=0, in=230] (5.7,1.7); 
	\draw[->] (4.3,1.7) 
	to [out=-50, in=180]  (5.5,1) to [out=0, in=230] node[label={[label distance=-1mm]270:{${\scriptstyle (a_1,b_3)}$}}, dot, pos=0.3] {} (6.7,1.7);
	\draw (4.5,-4.8) node {${\scriptstyle	\otimes \Hom_\basecat(a_1,b_3)\otimes \Hom_\basecat(a_3,b_2)\otimes \Hom_\basecat(a_4,b_5) }$};
	
\end{tikzpicture}
\end{minipage}
\end{equation*}

\item {\em The strand entering $\PP_{b_j}$ is an annihilation strand ($j
\in \underline{j})$: }

Then on the $(\underline{\iota}, \underline{j}, \sigma)$
direct summand 
the two compositions around 
\eqref{eqn-commutative-square-for-the-functoriality-of-phi}
are defined by two planar diagrams which only differ in:
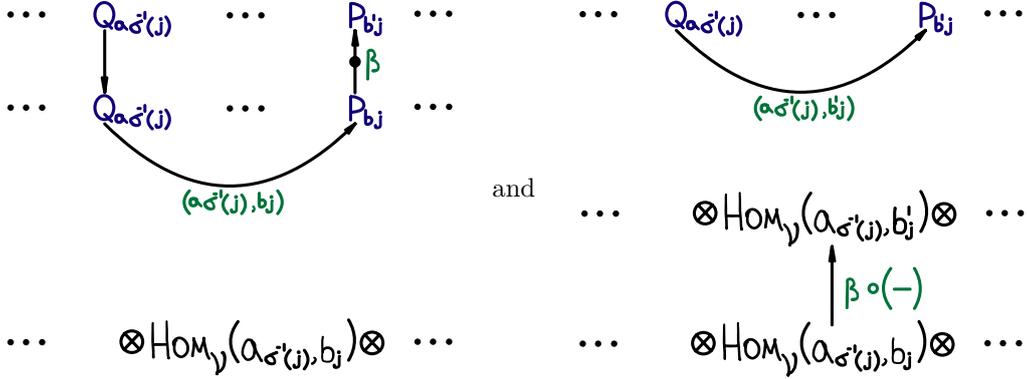
\begin{equation*}
	\begin{minipage}{0.4\textwidth}
			\begin{tikzpicture}[scale=0.65]
			\draw   (0,1.5)  node{$\QQ_{a_{\sigma^{-1}(j)}}$};
			\draw   (0,0)  node{$\QQ_{a_{\sigma^{-1}(j)}}$};
			\draw   (4,0)  node{$\PP_{b_{j}}$};
			\draw   (-2,0)  node{$\cdots$};
			\draw   (2,0)  node{$\cdots$};
			\draw   (6,0)  node{$\cdots$};
			\draw   (-2,1.5)  node{$\cdots$};
			\draw   (2,1.5)  node{$\cdots$};
			\draw   (6,1.5)  node{$\cdots$};
			\draw   (4,1.5)  node{$\PP_{b'_{j}}$};
				\draw[->] (4,0.3) -- node[label={[label distance=-1mm]0:{$\scriptstyle  \beta$}}, dot, pos=0.5] {} (4,1.2);
				\draw[->] (0,1.2) -- (0,0.3);
					\draw[->] (0,-0.3) 
				to [out=-50, in=180] (2,-1.2) to [out=0, in=230] node[label={[label distance=-1mm]270:{${\scriptstyle (a_{\sigma^{-1}(j)},b_j)}$}}, dot, pos=0] {} (4,-0.3); 
			\draw   (2,-3)  node{$\cdots \otimes \Hom_\basecat(a_{\sigma^{-1}(j)},b_j)\otimes \cdots$};
		\end{tikzpicture}
\end{minipage} 
\quad \text{ and } \quad 
\begin{minipage}{0.4\textwidth}
\begin{tikzpicture}[scale=0.65]
	\draw   (0,1.5)  node{$\QQ_{a_{\sigma^{-1}(j)}}$};
%	\draw   (0,0)  node{$\QQ_{a_{\sigma^{-1}(j)}}$};
%	\draw   (4,0)  node{$\PP_{b_{j}}$};
%	\draw   (-2,0)  node{$\cdots$};
%	\draw   (2,0)  node{$\cdots$};
%	\draw   (6,0)  node{$\cdots$};
	\draw   (-2,1.5)  node{$\cdots$};
	\draw   (2,1.5)  node{$\cdots$};
	\draw   (6,1.5)  node{$\cdots$};
	\draw   (4,1.5)  node{$\PP_{b'_{j}}$};
	\draw[->] (2,-2.7) -- node[label={[label distance=-1mm]0:{$\scriptstyle  \beta \circ ( - )$}}, dot, pos=0.5] {} (2,-1.8);
	%\draw[->] (0,1.2) -- (0,0.3);
	\draw[->] (0,1.2) 
	to [out=-50, in=180] (2,0.3) to [out=0, in=230] node[label={[label distance=-1mm]270:{${\scriptstyle (a_{\sigma^{-1}(j)},b_j)}$}}, dot, pos=0] {} (4,1.2); 
	\draw   (2,-3)  node{$\cdots \otimes \Hom_\basecat(a_{\sigma^{-1}(j)},b_j)\otimes \cdots$};
	\draw   (2,-1.5)  node{$\cdots \otimes \Hom_\basecat(a_{\sigma^{-1}(j)},b'_j)\otimes \cdots$};
\end{tikzpicture}
\end{minipage} 
\end{equation*}
Here and below, by abuse of notation, by $\sigma^{-1}(j)$ we 
denote the pre-image of $j$ under the permutation of $\left\{ j_1,
\dots, j_k \right\}$ defined by $\sigma \in S_k$. In other words, 
if $j = j_l$, we mean $j_{\sigma^{-1}(l)}$.  

By definition, these 
$2$-morphisms are the left adjoints of the maps which send
each $\gamma \in \homm_{\basecat}(a_{\sigma^{-1}(j)}, b_j)$ to 
the $2$-morphisms defined by the diagrams only in
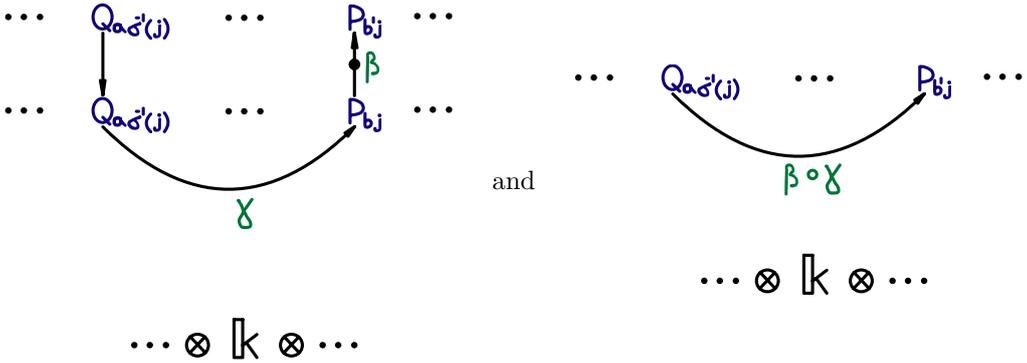
\begin{equation*}
	\begin{minipage}{0.4\textwidth}
	\begin{tikzpicture}[scale=0.65]
	\draw   (0,1.5)  node{$\QQ_{a_{\sigma^{-1}(j)}}$};
	\draw   (0,0)  node{$\QQ_{a_{\sigma^{-1}(j)}}$};
	\draw   (4,0)  node{$\PP_{b_{j}}$};
	\draw   (-2,0)  node{$\cdots$};
	\draw   (2,0)  node{$\cdots$};
	\draw   (6,0)  node{$\cdots$};
	\draw   (-2,1.5)  node{$\cdots$};
	\draw   (2,1.5)  node{$\cdots$};
	\draw   (6,1.5)  node{$\cdots$};
	\draw   (4,1.5)  node{$\PP_{b'_{j}}$};
	\draw[->] (4,0.3) -- node[label={[label distance=-1mm]0:{$\scriptstyle  \beta$}}, dot, pos=0.5] {} (4,1.2);
	\draw[->] (0,1.2) -- (0,0.3);
	\draw[->] (0,-0.3) 
	to [out=-50, in=180] (2,-1.2) to [out=0, in=230] node[label={[label distance=0mm]270:{${ \gamma}$}}, dot, pos=0] {} (4,-0.3); 
	\draw   (2,-3)  node{$\cdots \otimes \kk \otimes \cdots$};
\end{tikzpicture}
	\end{minipage}
	\quad \text{ and } \quad 
	\begin{minipage}{0.4\textwidth}
	\begin{tikzpicture}[scale=0.65]
	%\draw   (0,1.5)  node{$\QQ_{a_{\sigma^{-1}(j)}}$};
	\draw   (0,0)  node{$\QQ_{a_{\sigma^{-1}(j)}}$};
	\draw   (4,0)  node{$\PP_{b'_{j}}$};
	\draw   (-2,0)  node{$\cdots$};
	\draw   (2,0)  node{$\cdots$};
	\draw   (6,0)  node{$\cdots$};
	%\draw   (-2,1.5)  node{$\cdots$};
	%\draw   (2,1.5)  node{$\cdots$};
	%\draw   (6,1.5)  node{$\cdots$};
	%\draw   (4,1.5)  node{$\PP_{b'_{j}}$};
	%\draw[->] (4,0.3) -- node[label={[label distance=-1mm]0:{$\scriptstyle  \beta$}}, dot, pos=0.5] {} (4,1.2);
	%\draw[->] (0,1.2) -- (0,0.3);
	\draw[->] (0,-0.3) 
	to [out=-50, in=180] (2,-1.2) to [out=0, in=230] node[label={[label distance=0mm]270:{${ \beta \circ \gamma}$}}, dot, pos=0] {} (4,-0.3); 
	\draw   (2,-3)  node{$\cdots \otimes \kk \otimes \cdots$};
\end{tikzpicture}
	\end{minipage} 
\end{equation*}
Slide down the $\beta$-dot down the string and merge it with the
$\gamma$-dot turns the left diagram into the right diagram.  
By the dot sliding relations, the corresponding $2$-morphisms are equal. 
Hence their left adjoints are also equal, 
and \eqref{eqn-commutative-square-for-the-functoriality-of-phi}
commutes on this direct summand. 
\end{enumerate}

\underline{Morphisms of type
\ref{item-four-basic-morphisms-in-S^nVoppS^mV-S_m-element}}:

Consider $\phi_{\underline{a} \otimes \underline{b}}$ restricted
to the $(\underline{\iota}, \underline{j}, \sigma)$ direct summand 
of $\Xi_{\PP\QQ}(\hat{k})\left(\underline{a} \otimes
\underline{b}\right)$. There are five cases:

\begin{enumerate}
\item {\em The strand entering $\PP_{b_j}$ is a commutation strand ($j
\notin \underline{j}$), \\
the strand entering $\PP_{b_{j+1}}$ is a commutation strand ($j + 1
\notin \underline{j}$):}

Then on the $(\underline{\iota}, \underline{j}, \sigma)$
direct summand 
the two compositions around 
\eqref{eqn-commutative-square-for-the-functoriality-of-phi}
are defined by two planar diagrams which only differ in:
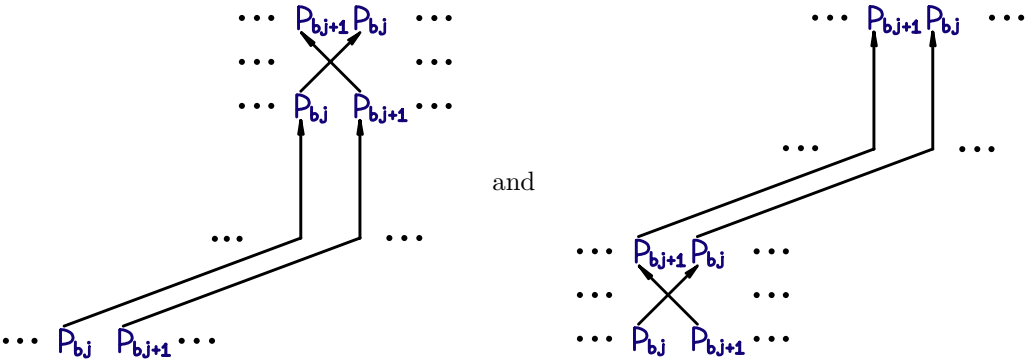
\begin{equation*}
	\begin{minipage}{0.4\textwidth}
		\begin{tikzpicture}[scale=0.7]
			\draw   (0,0)  node {$\PP_{b_j}$};
			\draw (4,4) node {$\PP_{b_j}$};	
			\draw (4,5.5) node {$\PP_{b_{j+1}}$};	
			\draw (3,4) node {$\dots$};
			\draw (6,4) node {$\dots$};
			\draw (3,5.5) node {$\dots$};
			\draw (6,5.5) node {$\dots$};
			\draw[->] (4,4.3) -- (5,5.2); 
			\draw[->] (0.3,0.3) -- (4,2) -- (4,3.7); 
			\draw   (1,0)  node {$\PP_{b_{j+1}}$};
			\draw (5,4) node {$\PP_{b_{j+1}}$};	
			\draw (5,5.5) node {$\PP_{b_{j}}$};	
			\draw (-1,0) node {$\dots$};
			\draw (2,0) node {$\dots$};
			\draw[->] (1.3,0.3) -- (5,2) -- (5,3.7); 
			\draw[->] (5,4.3) -- (4,5.2); 
		\end{tikzpicture}
	\end{minipage}
	\quad \text{ and } \quad 
	\begin{minipage}{0.4\textwidth}
		\begin{tikzpicture}[scale=0.7]
			\draw   (0,0)  node {$\PP_{b_{j+1}}$};
			\draw (4,4) node {$\PP_{b_{j+1}}$};	
			\draw (0,-1.5) node {$\PP_{b_j}$};	
			\draw (-1,0) node {$\dots$};
			\draw (2,0) node {$\dots$};
			\draw (-1,-1.5) node {$\dots$};
			\draw (2,-1.5) node {$\dots$};
			\draw[->] (0.3,0.3) -- (4,2) -- (4,3.7); 
			\draw   (1,0)  node {$\PP_{b_j}$};
			\draw (5,4) node {$\PP_{b_j}$};	
			\draw (1,-1.5) node {$\PP_{b_{j+1}}$};	
			\draw (3,4) node {$\dots$};
			\draw (6,4) node {$\dots$};
			\draw[->] (1.3,0.3) -- (5,2) -- (5,3.7); 
			\draw[->] (0,-1.2) -- (1,-0.3); 
			\draw[->] (1,-1.2) -- (0,-0.3); 
		\end{tikzpicture}
	\end{minipage}
\end{equation*}
Sliding the upward crossing down the parallel strands turns
the left diagram into the right one. On its way down, the crossing 
travels through the identical remainder of the two diagrams.
By the triple move relations \cite[Lemma 5.5]{gyenge2021heisenberg}, 
this doesn't change the corresponding $2$-morphism. Thus 
\eqref{eqn-commutative-square-for-the-functoriality-of-phi}
commutes on this direct summand.

\item {\em The strand entering $\PP_{b_j}$ is a commutation strand ($j
\notin \underline{j})$, \\
the strand entering $\PP_{b_{j+1}}$ is an annihilation 
($j + 1 \in \underline{j})$:}

Then on the $(\underline{\iota}, \underline{j}, \sigma)$
direct summand 
the two compositions around 
\eqref{eqn-commutative-square-for-the-functoriality-of-phi}
are defined by two planar diagrams which only differ in:
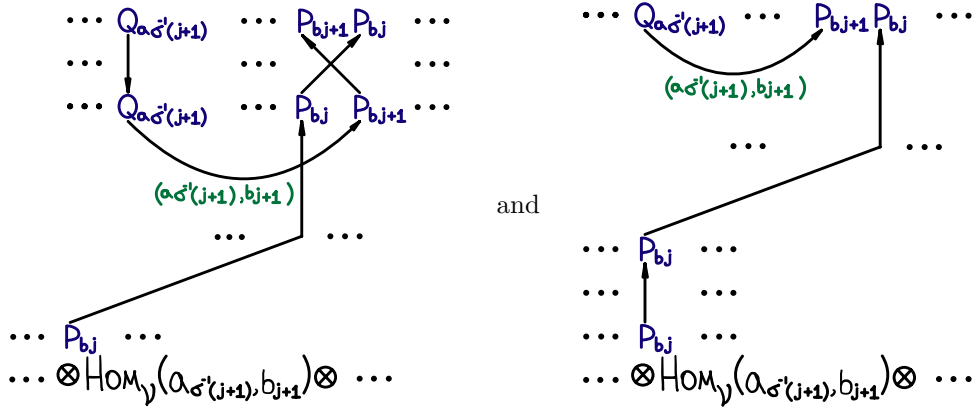
\begin{equation*}
	\begin{minipage}{0.38\textwidth}
		\begin{tikzpicture}[scale=0.7]
				\draw   (1,5.5)  node{$\QQ_{a_{\sigma^{-1}(j+1)}}$};
			\draw   (1,4)  node{$\QQ_{a_{\sigma^{-1}(j+1)}}$};
			\draw   (0,0)  node {$\PP_{b_j}$};
			\draw (4,4) node {$\PP_{b_j}$};	
			\draw (4,5.5) node {$\PP_{b_{j+1}}$};	
			\draw (2.75,4) node {$\dots$};
			\draw (6,4) node {$\dots$};
			\draw (2.75,5.5) node {$\dots$};
			\draw (6,5.5) node {$\dots$};
				\draw[->] (1,5.2) -- (1,4.3); 
			\draw[->] (4,4.3) -- (5,5.2); 
			\draw[->] (0.3,0.3) -- (4,2) -- (4,3.7); 
			%\draw   (1,0)  node {$\PP_{b_{j+1}}$};
			\draw (5,4) node {$\PP_{b_{j+1}}$};	
			\draw (5,5.5) node {$\PP_{b_{j}}$};	
			\draw (-1,0) node {$\dots$};
			\draw (1,0) node {$\dots$};
			\draw (-1,4) node {$\dots$};
			\draw (-1,5.5) node {$\dots$};
			%\draw[->] (1.3,0.3) -- (5,2) -- (5,3.7); 
			\draw[->] (5,4.3) -- (4,5.2); 
				\draw[->] (1,3.7) 
			to [out=-50, in=180] (3,2.8) to [out=0, in=230] node[label={[label distance=-1mm]260:{${\scriptstyle (a_{\sigma^{-1}(j+1)},b_{j+1})}$}}, dot, pos=0] {} (5,3.7); 
			\draw   (2,-1.5)  node{$\cdots \otimes \Hom_\basecat(a_{\sigma^{-1}(j+1)},b_{j+1})\otimes \cdots$};
		\end{tikzpicture}
	\end{minipage}
	\quad \text{ and } \quad 
	\begin{minipage}{0.38\textwidth}
		\begin{tikzpicture}[scale=0.7]
				\draw   (0,4)  node{$\QQ_{a_{\sigma^{-1}(j+1)}}$};
			%\draw   (0,0)  node {$\PP_{b_{j+1}}$};
			\draw (4,4) node {$\PP_{b_{j+1}}$};	
			%\draw (0,-1.5) node {$\PP_{b_j}$};	
			\draw (0,0) node {$\dots$};
			\draw (2,0) node {$\dots$};
			\draw (0,-1.5) node {$\dots$};
			\draw (2,-1.5) node {$\dots$};
			%\draw[->] (0.3,0.3) -- (4,2) -- (4,3.7); 
			\draw   (1,0)  node {$\PP_{b_j}$};
			\draw (5,4) node {$\PP_{b_j}$};	
			\draw (1,-1.5) node {$\PP_{b_{j}}$};	
			\draw (-2,4) node {$\dots$};
			\draw (2,4) node {$\dots$};
			\draw (6,4) node {$\dots$};
			\draw[->] (1.3,0.3) -- (5,2) -- (5,3.7); 
			%\draw[->] (0,-1.2) -- (0,-0.3); 
			\draw[->] (1,-1.2) -- (1,-0.3); 
			\draw[->] (0,3.7) 
			to [out=-50, in=180] (2,2.8) to [out=0, in=230] node[label={[label distance=-1mm]270:{${\scriptstyle (a_{\sigma^{-1}(j+1)},b_{j+1})}$}}, dot, pos=0] {} (4,3.7); 
			\draw   (2,-3)  node{$\cdots \otimes \Hom_\basecat(a_{\sigma^{-1}(j+1)},b_{j+1})\otimes \cdots$};
		\end{tikzpicture}
	\end{minipage}
\end{equation*}
Here for the first time $\id \otimes \zeta$
acts non-trivially on our chosen direct summand. 
We have $\zeta = (j(j+1))$, so by Lemma 
\ref{lemma-explicit-description-of-the-functor-hatk} 
the morphism $\bigoplus \Xi_{\PP\QQ}(\hat{k})(\id \otimes \zeta)$
maps the $(\underline{i}, \underline{j}, \sigma)$ summand
of $\bigoplus \Xi_{\PP\QQ}(\underline{a} \otimes \underline{b})$
to the $(\underline{i}, \zeta(\underline{j}), \zeta(\sigma))$ summand of 
$\bigoplus \Xi_{\PP\QQ}(\underline{a} \otimes \zeta(\underline{b}))$. 
Recall that $\underline{j}$ is a choice of $1 \leq j_1 < \dots < j_k
\leq m$ and $\zeta(\underline{j})$ is $\zeta(j_1)$, \dots,
$\zeta(j_k)$ reordered in the increasing order. As 
$j \notin \underline{j}$ and $j + 1 \in \underline{j}$,   
$\zeta(j_1)$, \dots, $\zeta(j_k)$ is  
$j_1$, \dots, $j_k$ with $j + 1$ replaced by $j$. No reordering
needed, so $\zeta(\sigma) = \sigma$.  Similarly, 
$\zeta(1), \dots, \widehat{\zeta(j_{1})}, \dots, 
\widehat{\zeta(j_k)}, \dots, \zeta(m)$ are  
$1, \dots, \widehat{j_{1}}, \dots, 
\widehat{j_k}, \dots, m$ with $j+1$ replaced by $j$, with 
no reordering needed. Thus, by Lemma
\ref{lemma-explicit-description-of-the-functor-hatk}, 
$\bigoplus \Xi_{\PP\QQ}(\hat{k})(\id \otimes \zeta)$
maps one summand to the other by the identity map.

On the left diagram, slide the lower crossing up until 
it reaches the other crossing. By the triple move relations 
this does not change the corresponding $2$-morphism. 
Once one crossing reaches the other, 
replace their composition by two parallel upward strands to
get the right diagram. By the symmetric group relations for upward 
strands \cite[Lemma 5.5]{gyenge2021heisenberg}, this doesn't change
the $2$-morphism either. Thus 
\eqref{eqn-commutative-square-for-the-functoriality-of-phi}
commutes on this direct summand.

\item {\em The strand entering $\PP_{b_j}$ is an annihilation strand ($j
\in \underline{j})$, \\
the strand entering $\PP_{b_{j+1}}$ is an commutation strand 
($j + 1 \notin \underline{j})$:}

Then on the $(\underline{\iota}, \underline{j}, \sigma)$
direct summand 
the two compositions around 
\eqref{eqn-commutative-square-for-the-functoriality-of-phi}
are defined by two planar diagrams which only differ in:
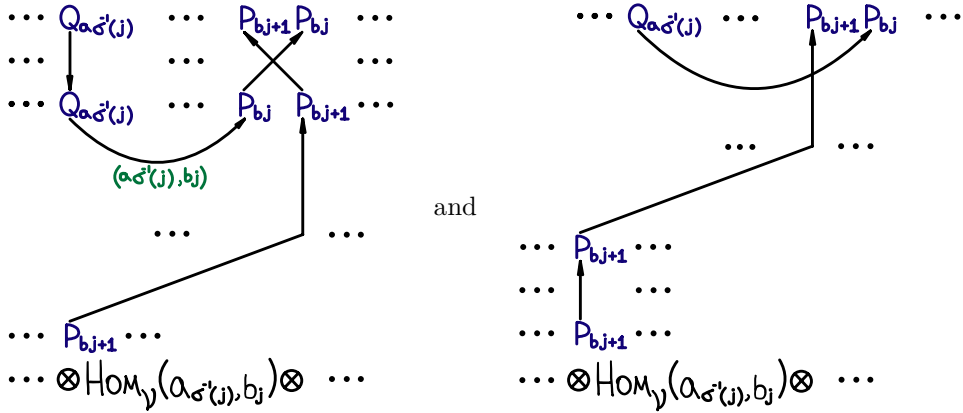
\begin{equation*}
	\begin{minipage}{0.38\textwidth}
		\begin{tikzpicture}[scale=0.7]
			\draw   (1,5.5)  node{$\QQ_{a_{\sigma^{-1}(j)}}$};
			\draw   (1,4)  node{$\QQ_{a_{\sigma^{-1}(j)}}$};
			\draw   (1,0)  node {$\PP_{b_{j+1}}$};
			\draw (4,4) node {$\PP_{b_j}$};	
			\draw (4,5.5) node {$\PP_{b_{j+1}}$};	
			\draw (2.5,4) node {$\dots$};
			\draw (6,4) node {$\dots$};
			\draw (2.5,5.5) node {$\dots$};
			\draw (6,5.5) node {$\dots$};
			\draw[->] (1,5.2) -- (1,4.3); 
			\draw[->] (4,4.3) -- (5,5.2); 
			\draw[->] (1.3,0.3) -- (5,2) -- (5,3.7); 
			%\draw   (1,0)  node {$\PP_{b_{j+1}}$};
			\draw (5,4) node {$\PP_{b_{j+1}}$};	
			\draw (5,5.5) node {$\PP_{b_{j}}$};	
			\draw (0,0) node {$\dots$};
			\draw (2,0) node {$\dots$};
			\draw (-1,4) node {$\dots$};
			\draw (-1,5.5) node {$\dots$};
			%\draw[->] (1.3,0.3) -- (5,2) -- (5,3.7); 
			\draw[->] (5,4.3) -- (4,5.2); 
			\draw[->] (1,3.7) 
			to [out=-50, in=180] (2.5,2.8) to [out=0, in=230] node[label={[label distance=-1mm]270:{${\scriptstyle (a_{\sigma^{-1}(j)},b_{j})}$}}, dot, pos=0] {} (4,3.7); 
			\draw   (2,-1.5)  node{$\cdots \otimes \Hom_\basecat(a_{\sigma^{-1}(j)},b_{j})\otimes \cdots$};
		\end{tikzpicture}
	\end{minipage}
	\quad \text{ and } \quad 
	\begin{minipage}{0.38\textwidth}
		\begin{tikzpicture}[scale=0.7]
			\draw   (0,4)  node{$\QQ_{a_{\sigma^{-1}(j)}}$};
			%\draw   (0,0)  node {$\PP_{b_{j+1}}$};
			\draw (4,4) node {$\PP_{b_{j+1}}$};	
			%\draw (0,-1.5) node {$\PP_{b_j}$};	
			\draw (-1,0) node {$\dots$};
			\draw (1,0) node {$\dots$};
			\draw (-1,-1.5) node {$\dots$};
			\draw (1,-1.5) node {$\dots$};
			%\draw[->] (0.3,0.3) -- (4,2) -- (4,3.7); 
			\draw   (0,0)  node {$\PP_{b_{j+1}}$};
			\draw (5,4) node {$\PP_{b_j}$};	
			\draw (0,-1.5) node {$\PP_{b_{j+1}}$};	
			\draw (-2,4) node {$\dots$};
			\draw (2,4) node {$\dots$};
			\draw (6,4) node {$\dots$};
			\draw[->] (0.3,0.3) -- (4,2) -- (4,3.7); 
			%\draw[->] (0,-1.2) -- (0,-0.3); 
			\draw[->] (0,-1.2) -- (0,-0.3); 
			\draw[->] (0,3.7) 
			to [out=-50, in=180] (2.5,2.8) to [out=0, in=230] node[label={[label distance=-1mm]270:{${\scriptstyle (a_{\sigma^{-1}(j)},b_{j})}$}}, dot, pos=0] {} (5,3.7); 
			\draw   (2,-3)  node{$\cdots \otimes \Hom_\basecat(a_{\sigma^{-1}(j)},b_{j})\otimes \cdots$};
		\end{tikzpicture}
	\end{minipage}
\end{equation*}

Similarly, $\zeta$ sends $(\underline{i}, \underline{j}, \sigma)$
to $(\underline{i}, \zeta(\underline{j}), \sigma)$ where 
$\zeta(\underline{j})$ is $\underline{j}$ with 
$j$ replaced by $j+1$. Again, 
$\bigoplus \Xi_{\PP\QQ}(\hat{k})(\id \otimes \zeta)$  acts by the identity map. 

Take the left diagram and slide the upward crossing down the parallel strands 
to obtain the right diagram. By the triple move relations, this 
doesn't change the corresponding $2$-morphism. Thus 
\eqref{eqn-commutative-square-for-the-functoriality-of-phi}
commutes on this direct summand.

\item {\em The strand entering $\PP_{b_j}$ is an annihilation strand ($j
\in \underline{j})$, \\
the strand entering $\PP_{b_{j+1}}$ is an annihilation strand 
($j + 1 \in \underline{j})$, \\
and these two strands cross ($\sigma^{-1}(j) > \sigma^{-1}(j+1)$):}

Then on the $(\underline{\iota}, \underline{j}, \sigma)$
direct summand the two compositions around 
\eqref{eqn-commutative-square-for-the-functoriality-of-phi}
are defined by two planar diagrams which only differ in:
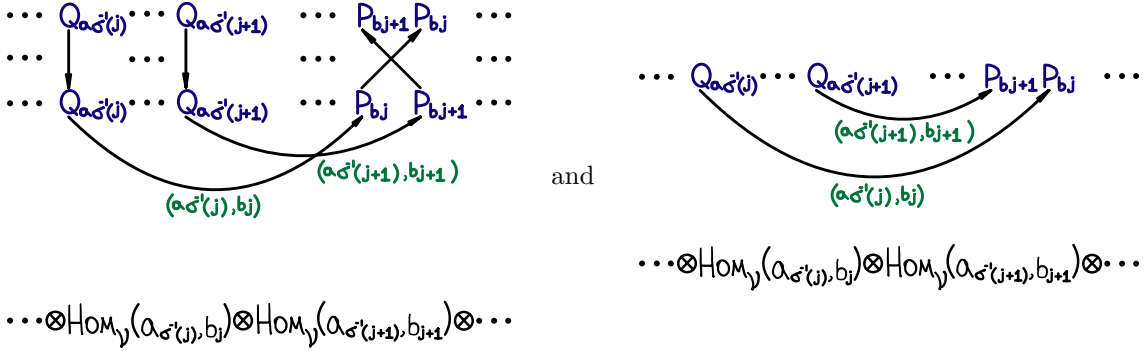
\begin{equation*}
	\begin{minipage}{0.45\textwidth}
		\begin{tikzpicture}[scale=0.64]
			\draw   (1,5.5)  node{$\QQ_{a_{\sigma^{-1}(j+1)}}$};
			\draw   (1,4)  node{$\QQ_{a_{\sigma^{-1}(j+1)}}$};
				\draw   (-2,5.5)  node{$\QQ_{a_{\sigma^{-1}(j)}}$};
			\draw   (-2,4)  node{$\QQ_{a_{\sigma^{-1}(j)}}$};
			%\draw   (0,0)  node {$\PP_{b_j}$};
			\draw (4,4) node {$\PP_{b_j}$};	
			\draw (4,5.5) node {$\PP_{b_{j+1}}$};	
			\draw (2.75,4) node {$\dots$};
			\draw (6,4) node {$\dots$};
			\draw (2.75,5.5) node {$\dots$};
			\draw (6,5.5) node {$\dots$};
			\draw[->] (1,5.2) -- (1,4.3); 
			\draw[->] (-2,5.2) -- (-2,4.3); 
			\draw[->] (4,4.3) -- (5,5.2); 
			%\draw[->] (0.3,0.3) -- (4,2) -- (4,3.7); 
			%\draw   (1,0)  node {$\PP_{b_{j+1}}$};
			\draw (5,4) node {$\PP_{b_{j+1}}$};	
			\draw (5,5.5) node {$\PP_{b_{j}}$};	
			%\draw (-1,0) node {$\dots$};
			%\draw (1,0) node {$\dots$};
			\draw (-3.5,4) node {$\dots$};
			\draw (-3.5,5.5) node {$\dots$};
			\draw (-0.75,4) node {$\dots$};
			\draw (-0.75,5.5) node {$\dots$};
			%\draw[->] (1.3,0.3) -- (5,2) -- (5,3.7); 
			\draw[->] (5,4.3) -- (4,5.2); 
			\draw[->] (1,3.7) 
			to [out=-50, in=180] (3,2.8) to [out=0, in=230] node[label={[label distance=-1mm]300:{${\scriptstyle (a_{\sigma^{-1}(j+1)},b_{j+1})}$}}, dot, pos=0.1] {} (5,3.7); 
			\draw[->] (-2,3.7) 
			to [out=-50, in=180] (1,2.6) to [out=0, in=230] node[label={[label distance=-1mm]270:{${\scriptstyle (a_{\sigma^{-1}(j)},b_{j})}$}}, dot, pos=0] {} (4,3.7); 
			\draw   (1.5,1.5)  node{${\scriptstyle \cdots  \otimes \Hom_\basecat(a_{\sigma^{-1}(j)},b_{j}) \otimes \Hom_\basecat(a_{\sigma^{-1}(j+1)},b_{j+1})\otimes \cdots}$};
		\end{tikzpicture}
	\end{minipage}
	\quad \text{ and } \quad 
	\begin{minipage}{0.45\textwidth}
		\begin{tikzpicture}[scale=0.64]
			\draw   (0,4)  node{$\QQ_{a_{\sigma^{-1}(j+1)}}$};
			\draw   (-3,4)  node{$\QQ_{a_{\sigma^{-1}(j)}}$};
			%\draw   (0,0)  node {$\PP_{b_{j+1}}$};
			\draw (3.5,4) node {$\PP_{b_{j+1}}$};	
			%\draw (0,-1.5) node {$\PP_{b_j}$};	
			%\draw (0,0) node {$\dots$};
			%\draw (2,0) node {$\dots$};
			%\draw (0,-1.5) node {$\dots$};
			%\draw (2,-1.5) node {$\dots$};
			%\draw[->] (0.3,0.3) -- (4,2) -- (4,3.7); 
			%\draw   (1,0)  node {$\PP_{b_j}$};
			\draw (4.5,4) node {$\PP_{b_j}$};	
			%\draw (1,-1.5) node {$\PP_{b_{j}}$};	
			\draw (-4.5,4) node {$\dots$};
				\draw (-1.75,4) node {$\dots$};
			\draw (1.75,4) node {$\dots$};
			\draw (5.5,4) node {$\dots$};
			%\draw[->] (1.3,0.3) -- (5,2) -- (5,3.7); 
			%\draw[->] (0,-1.2) -- (0,-0.3); 
			%\draw[->] (1,-1.2) -- (1,-0.3); 
			\draw[->] (0,3.7) 
			to [out=-40, in=180] (2,3.1) to [out=0, in=220] node[label={[label distance=-1mm]270:{${\scriptstyle (a_{\sigma^{-1}(j+1)},b_{j+1})}$}}, dot, pos=0] {} (3.5,3.7); 
			\draw[->] (-3,3.7) 
			to [out=-50, in=180] (1,1.8) to [out=0, in=230] node[label={[label distance=-1mm]270:{${\scriptstyle (a_{\sigma^{-1}(j)},b_{j})}$}}, dot, pos=0] {} (4.5,3.7); 
			\draw   (0.5,0.5)  node{${\scriptstyle \cdots  \otimes \Hom_\basecat(a_{\sigma^{-1}(j)},b_{j}) \otimes \Hom_\basecat(a_{\sigma^{-1}(j+1)},b_{j+1})\otimes \cdots}$};
		\end{tikzpicture}
	\end{minipage}
\end{equation*}

Here applying $\zeta$ to $\underline{j}$ swaps $j$ and $j+1$, 
and we need to reorder by applying $\zeta$ again. 
Thus $\zeta(\underline{j}) = \underline{j}$, but $\sigma$ becomes 
$\zeta \sigma$. Moreover, applying $\zeta$ to 
$1, \dots, \widehat{j_{1}}, \dots, \widehat{j_k}, \dots, m$
changes nothing. Thus 
$\bigoplus \Xi_{\PP\QQ}(\hat{k})(\id \otimes \zeta)$ maps the $(\underline{i}, \underline{j}, \sigma)$ summand
of $\bigoplus \Xi_{\PP\QQ}(\underline{a} \otimes \underline{b})$
to the $(\underline{i}, \underline{j}, \zeta\sigma)$ summand of 
$\bigoplus \Xi_{\PP\QQ}(\underline{a} \otimes \zeta(\underline{b}))$
by the identity map. 

Take the left diagram, slide the lower crossing up 
until it reaches the other crossing, and then replace the composition
of the two crossings by two upward parallel strands to obtain the
right diagram. 
By the triple move and the symmetric group relations, this doesn't
change the corresponding $2$-morphism.  
Thus \eqref{eqn-commutative-square-for-the-functoriality-of-phi}
commutes on this direct summand.

\item {\em The strand entering $\PP_{b_j}$ is an annihilation strand ($j
\in \underline{j})$, \\
the strand entering $\PP_{b_{j+1}}$ is an annihilation strand 
($j + 1 \in \underline{j})$, \\
and these two strands do not cross ($\sigma^{-1}(j) < \sigma^{-1}(j+1)$):}
\begin{equation*}
	\begin{minipage}{0.45\textwidth}
		\begin{tikzpicture}[scale=0.64]
			\draw   (1,5.5)  node{$\QQ_{a_{\sigma^{-1}(j)}}$};
			\draw   (1,4)  node{$\QQ_{a_{\sigma^{-1}(j)}}$};
			\draw   (-2,5.5)  node{$\QQ_{a_{\sigma^{-1}(j+1)}}$};
			\draw   (-2,4)  node{$\QQ_{a_{\sigma^{-1}(j+1)}}$};
			%\draw   (0,0)  node {$\PP_{b_j}$};
			\draw (4,4) node {$\PP_{b_j}$};	
			\draw (4,5.5) node {$\PP_{b_{j+1}}$};	
			\draw (2.75,4) node {$\dots$};
			\draw (6,4) node {$\dots$};
			\draw (2.75,5.5) node {$\dots$};
			\draw (6,5.5) node {$\dots$};
			\draw[->] (1,5.2) -- (1,4.3); 
			\draw[->] (-2,5.2) -- (-2,4.3); 
			\draw[->] (4,4.3) -- (5,5.2); 
			%\draw[->] (0.3,0.3) -- (4,2) -- (4,3.7); 
			%\draw   (1,0)  node {$\PP_{b_{j+1}}$};
			\draw (5,4) node {$\PP_{b_{j+1}}$};	
			\draw (5,5.5) node {$\PP_{b_{j}}$};	
			%\draw (-1,0) node {$\dots$};
			%\draw (1,0) node {$\dots$};
			\draw (-3.5,4) node {$\dots$};
			\draw (-3.5,5.5) node {$\dots$};
			\draw (-0.5,4) node {$\dots$};
			\draw (-0.5,5.5) node {$\dots$};
			%\draw[->] (1.3,0.3) -- (5,2) -- (5,3.7); 
			\draw[->] (5,4.3) -- (4,5.2); 
			\draw[->] (1,3.7) 
			to [out=-40, in=180] (2.5,3.2) to [out=0, in=220] node[label={[label distance=-1mm]220:{${\scriptstyle (a_{\sigma^{-1}(j)},b_{j})}$}}, dot, pos=0.1] {} (4,3.7); 
			\draw[->] (-2,3.7) 
			to [out=-50, in=180] (1.5,2.2) to [out=0, in=230] node[label={[label distance=-1mm]270:{${\scriptstyle (a_{\sigma^{-1}(j+1)},b_{j+1})}$}}, dot, pos=0] {} (5,3.7); 
			\draw   (1.5,1.2)  node{${\scriptstyle \cdots  \otimes \Hom_\basecat(a_{\sigma^{-1}(j)},b_{j}) \otimes \Hom_\basecat(a_{\sigma^{-1}(j+1)},b_{j+1})\otimes \cdots}$};
		\end{tikzpicture}
	\end{minipage}
	\quad \text{ and } \quad 
	\begin{minipage}{0.45\textwidth}
		\begin{tikzpicture}[scale=0.64]
			\draw   (0,4)  node{$\QQ_{a_{\sigma^{-1}(j)}}$};
			\draw   (-3,4)  node{$\QQ_{a_{\sigma^{-1}(j+1)}}$};
			%\draw   (0,0)  node {$\PP_{b_{j+1}}$};
			\draw (3.5,4) node {$\PP_{b_{j+1}}$};	
			%\draw (0,-1.5) node {$\PP_{b_j}$};	
			%\draw (0,0) node {$\dots$};
			%\draw (2,0) node {$\dots$};
			%\draw (0,-1.5) node {$\dots$};
			%\draw (2,-1.5) node {$\dots$};
			%\draw[->] (0.3,0.3) -- (4,2) -- (4,3.7); 
			%\draw   (1,0)  node {$\PP_{b_j}$};
			\draw (4.5,4) node {$\PP_{b_j}$};	
			%\draw (1,-1.5) node {$\PP_{b_{j}}$};	
			\draw (-4.75,4) node {$\dots$};
			\draw (-1.5,4) node {$\dots$};
			\draw (1.75,4) node {$\dots$};
			\draw (5.5,4) node {$\dots$};
			%\draw[->] (1.3,0.3) -- (5,2) -- (5,3.7); 
			%\draw[->] (0,-1.2) -- (0,-0.3); 
			%\draw[->] (1,-1.2) -- (1,-0.3); 
			\draw[->] (0,3.7) 
			to [out=-30, in=180] (2.5,3.1) to [out=0, in=210] node[label={[label distance=-1mm]250:{${\scriptstyle (a_{\sigma^{-1}(j)},b_{j})}$}}, dot, pos=0] {} (4.5,3.7); 
			\draw[->] (-3,3.7) 
			to [out=-50, in=180] (0.5,2) to [out=0, in=230] node[label={[label distance=-1mm]270:{${\scriptstyle (a_{\sigma^{-1}(j+1)},b_{j+1})}$}}, dot, pos=0] {} (3.5,3.7); 
			\draw   (0.5,0.5)  node{${\scriptstyle \cdots  \otimes \Hom_\basecat(a_{\sigma^{-1}(j)},b_{j}) \otimes \Hom_\basecat(a_{\sigma^{-1}(j+1)},b_{j+1})\otimes \cdots}$};
		\end{tikzpicture}
	\end{minipage}
\end{equation*}

Again, $\zeta$ sends $(\underline{i},
\underline{j}, \sigma)$ to $(\underline{i}, \underline{j}, \zeta\sigma)$
and $\bigoplus \Xi_{\PP\QQ}(\hat{k})(\id \otimes \zeta)$ acts by the
identity map.  

Take the left diagram and slide the upward crossing down the parallel strands 
to obtain the right diagram. By the triple move relations, 
this doesn't change the corresponding $2$-morphism. Thus 
\eqref{eqn-commutative-square-for-the-functoriality-of-phi}
commutes on this direct summand.

\end{enumerate}
\end{proof}

Its functoriality makes the Heisenberg relation categorification
\eqref{eqn-functorial-categorification-of-the-PQ-Heisenberg-relation}
applicable to Hochschild homology. We now complete the proof 
of the desired Heisenberg relation 
\eqref{eqn-heisenberg-relation-for-the-assignments-of-pi}
by applying 
\eqref{eqn-functorial-categorification-of-the-PQ-Heisenberg-relation}
to the product of the Hochschild homology classes $\psi_n(\alpha)$ and 
$\psi_m(\beta)$:

\begin{proof}[Proof of Theorem
\ref{theorem-heisenberg-relation-for-the-assignments-of-pi}]

Let $\alpha, \beta \in \hochcx_\bullet(\basecat)$. 
Implicitly using 
the isomorphism $\hochcx_\bullet(\basecat) \simeq \hochcx_\bullet(\basecat^{\opp})$
defined in Section
\ref{section-hochschild-homology-of-the-opposite-category}, 
we also denote by $\alpha$ the corresponding element in 
$\hochcx_\bullet(\basecat^{\opp})$. 

By \cite[Lemma 3.4]{Keller-OnTheCyclicHomologyOfExactCategories}, 
it follows from
Theorem \ref{theorem-functorial-categorification-of-the-PQ-Heisenberg-relation}
that in $\hochhom \homm_{\hcat\basecat}(0,m-n)$ we have
$$ 
\Xi_{\QQ\PP}\bigl(K(\psi_n(\alpha) \otimes \psi_m(\beta))\bigr)
= \bigoplus_{k = 0}^{k = \min(n,m)}
\Xi_{\PP\QQ}(\hat{k})\bigl(K(\psi_n(\alpha) \otimes
\psi_m(\beta))\bigr), $$
where $K$ is the shuffle product map 
$$ 
\hochcx_\bullet\left(\sym^{n} \basecat^{\opp}\right) 
\otimes 
\hochcx_\bullet\left(\sym^{m} \basecat\right)
\xrightarrow{\eqref{eqn-explicit-formula-for-kunneth-map}}
\hochcx_\bullet\left(\sym^{n} \basecat^{\opp} \otimes \sym^{m}
\basecat\right). 
$$

It follows from the definition \eqref{eqn-two-functor-QP}
of the functor $\Xi_{\QQ\PP}$ that 
$$
\Xi_{\QQ\PP}\bigl(K(\psi_n(\alpha) \otimes \psi_m(\beta)\bigr))
= \Xi_{\QQ}\bigl(\psi_n(\alpha)\bigr)
\Xi_\PP\bigl(\psi_m(\beta)\bigr) = 
\hAA_{\alpha}(-n) \hAA_{\beta}(m). $$ 
Since $\hat{0}$ is the identity functor, we have 
$\Xi_{\PP\QQ}\bigl(\hat{0}\bigr) = \Xi_{\PP\QQ}$ and similarly
$$
\Xi_{\PP\QQ}
\bigl(K(\psi_n(\alpha) \otimes \psi_m(\beta)) \bigr)
= 
\Xi_{\PP}\bigl(\psi_m(\beta)\bigr) 
\Xi_{\QQ}\bigl(\psi_n(\alpha)\bigr)
= 
\hAA_{\beta}(m) \hAA_{\alpha}(-n). $$ 

It now suffices to establish the following claim: in the Hochschild
homology 
\begin{equation}
\label{eqn-key-claim-for-proving-the-heisenberg-relation} 
\forall\; k \geq 1 \quad \quad \quad  
\Xi_{\PP\QQ}\bigl(\hat{k}\bigr)
\bigl(K(\psi_n(\alpha) \otimes \psi_m(\beta)) \bigr)
= 
\begin{cases}
n \left<\alpha, \beta\right> & \text{ if } k = n = m, 
\\
0 & \text{ otherwise. }
\end{cases}
\end{equation}

By definition, $K\left(\psi_n(\alpha) \otimes \psi_m(\beta)\right)$ is 
the image of $\alpha \otimes \beta$ under the composition 
\begin{align*}
\hochcx_\bullet(\basecat^{\opp}) \otimes \hochcx_\bullet(\basecat)
\xrightarrow{g_n \otimes g_m}
& \hochcx_\bullet\left((\basecat^{\opp})^{\otimes n};t_n\right) 
\otimes 
\hochcx_\bullet\left(\basecat^{\otimes m};t_m\right) 
\xrightarrow{\xi_{t_n} \otimes \xi_{t_m}}
\\
\rightarrow\;
& 
\hochcx_\bullet\left(\sym^n \basecat^{\opp} \right) 
\otimes 
\hochcx_\bullet\left(\sym^m \basecat \right) 
\xrightarrow{K}
\hochcx_\bullet\left(\sym^n \basecat^{\opp} \otimes \sym^m \basecat
\right), 
\end{align*}
where $t_n = (1 \dots n) \in S_n$ and $t_m = (1 \dots m) \in S_m$
are the long cycles. 
By Lemma \ref{lemma-compatibility-of-the-map-g-with-shuffle-product}, 
this equals 
\begin{align*}
\hochcx_\bullet(\basecat^{\opp}) \otimes \hochcx_\bullet(\basecat)
\xrightarrow{g_n \otimes g_m}
& \hochcx_\bullet\left((\basecat^{\opp})^{\otimes n};t_n\right) 
\otimes 
\hochcx_\bullet\left(\basecat^{\otimes m};t_m\right) 
\xrightarrow{K} \\
\rightarrow\;
& \hochcx_\bullet\left((\basecat^{\opp})^{\otimes n} \otimes
\basecat^{\otimes m};t_n \times t_m\right) 
\xrightarrow{\xi_{t_n \times t_m}}
\hochcx_\bullet\left(\sym^n \basecat^{\opp} \otimes \sym^m \basecat
\right), 
\end{align*}
where we implicitly identify 
$\sym^n \basecat^{\opp} \otimes \sym^m \basecat$
with $(\basecat^{\opp})^{\otimes n} \otimes \basecat^{\otimes m})
\rtimes (S_n \times S_m)$. 

By definition, $\Xi_{\PP\QQ}\bigl(\hat{k}\bigr) = (\Xi_{\PP} \circ_1
\Xi_{\QQ}) \circ  \bigl(\hat{k}\bigr)$. The annihilation functor
$\bigl(\hat{k}\bigr)$ was itself defined in Definition 
\ref{defn-annihilation-functor} as a composition of functors
which begins with the functor 
$$ \sym^{n}\basecat^{\opp} \otimes \sym^{m}\basecat
\xrightarrow{
\Res_{S_{n-k} \times S_{k} \times S_{m-k}}^{S_n \times S_m}
}
 \hperf \left(
\sym^{n-k}\basecat^{\opp} 
\otimes 
\left(\left((\basecat^{\opp})^{\otimes k} \otimes \basecat^{\otimes k}\right)
\rtimes S_k\right) 
\otimes 
\sym^{m-k}\basecat
\right), $$
where $S_{n-k} \times S_{k} \times S_{m-k}$ embeds into 
$S^n \times S^m$ as described in Definition \ref{defn-annihilation-functor}. 

We conclude that $\Xi_{\PP\QQ}\bigl(\hat{k}\bigr)
\bigl(K(\psi_n(\alpha) \otimes \psi_m(\beta)) \bigr)$ is the image 
of $\alpha \otimes \beta \in \hochcx_\bullet(\basecat^{\opp}) \otimes
\hochcx_\bullet(\basecat)$ under a composition of maps which includes
\begin{equation}
\label{eqn-composition-xi-and-Res-for-t_n-times-t_m}
\begin{tikzcd}
\hochcx_\bullet\left((\basecat^{\opp})^{\otimes n} \otimes
\basecat^{\otimes m};t_n \times t_m\right) 
\ar{d}{\xi_{t_n \times t_m}}
\\
\hochcx_\bullet\left(\sym^n \basecat^{\opp} \otimes \sym^m \basecat
\right) 
\ar{d}{\Res_{S_{n-k} \times S_{k} \times S_{m-k}}^{S_n \times S_m}}
\\
\hochcx_\bullet
\left(
\sym^{n-k}\basecat^{\opp} 
\otimes 
\left(\left((\basecat^{\opp})^{\otimes k} \otimes \basecat^{\otimes k}\right)
\rtimes S_k\right) 
\otimes 
\sym^{m-k}\basecat
\right).
\end{tikzcd}
\end{equation}
By Lemma
\ref{lemma-restriction-functor-composed-with-the-twisted-part}, 
on $\hochhom_\bullet$ the map 
\eqref{eqn-composition-xi-and-Res-for-t_n-times-t_m}
is a sum indexed by the elements of 
$\fix_Q(t_n \times t_m)$, 
the fixed set of the action of $t_n \times t_m$ on the set 
$Q$ of the left cosets of $S_{n-k} \times S_{k} \times S_{m-k}$ in 
$S_n \times S_m$. 
However, $\fix_Q(t_n \times t_m) = \emptyset$ unless $k = 0$ or $k = n = m$. 
This can be seen from the explicit description of the action of 
$S_n \times S_m$ on $Q$ in the proof of Lemma 
\ref{lemma-explicit-description-of-the-functor-hatk}. Alternatively,
since 
$$ S_{n-k} \times S_{k} \times S_{m-k} 
\leq 
S_{n-k} \times S_{k} \times S_{k} \times S_{m-k} 
\leq 
S_n \times S_m, $$
it suffices to consider the action of $t_n \times t_m$ on 
the cosets of $S_{n-k} \times S_{k} \times S_{k} \times S_{m-k}$. 
We have 
$$ \fix_{S_n \times S_m / S_{n-k} \times S_{k} \times S_{k} \times S_{m-k}}
(t_n \times t_m) = 
\fix_{S_n / S_{n-k} \times S_{k}}(t_n) \times 
\fix_{S_m / S_{k} \times S_{m-k}}(t_m), $$
and unless $k = 0$ or $k = n = m$ one of the two sets in the cartesian
product on the RHS is empty. 

Thus the chain 
$\Xi_{\PP\QQ}\bigl(\hat{k}\bigr)\bigl(K(\psi_n(\alpha) \otimes
\psi_m(\beta)) \bigr)$ vanishes in $\hochhom_\bullet$ 
unless $k = 0$ or $k = n = m$. 
This shows most of
\eqref{eqn-key-claim-for-proving-the-heisenberg-relation}. 
It remains to show that when $n = m$ we have 
\begin{equation}
\Xi_{\PP\QQ}\bigl(\hat{n}\bigr)\bigl(K(\psi_n(\alpha) \otimes
\psi_n(\beta)) \bigr) = \left< \alpha, \beta \right> [\hunit] 
\in \hochhom_\bullet(\homm_{\hcat\basecat}(0,0),
\end{equation}
where $\hunit \in \homm_{\hcat\basecat}(0,0)$ is the identity
$1$-morphism. 
By Definition \ref{defn-dg-functor-PQhatk}, the functor 
$\Xi_{\PP\QQ}\bigl(\hat{n}\bigr)$ is  
\begin{align*}
\sym^{n}\basecat^{\opp} \otimes \sym^{n}\basecat
\xrightarrow{(\hat{n})}
\hperf\left( \kk \right)
\xrightarrow{\sim}
\hperf\left( \kk \right) \otimes \hperf\left( \kk \right)
\xrightarrow{\Xi_\PP \circ_1 \Xi_\QQ}
\homm_{\hcat\basecat}\bigl(0,0)\bigr). 
\end{align*}
Each of $\Xi_\PP$ and $\Xi_\QQ$ send $\kk \in \hperf\left(\kk\right)$
to $\hunit \in \homm_{\hcat\basecat}(0,0)$. Thus the latter
two terms in the composition send $\kk \in \hperf\left(\kk\right)$
to $\hunit \circ_1 \hunit = \hunit \in \homm_{\hcat\basecat}(0,0)$. 
Therefore, on the level of Hochschild homology they become 
the linear map 
$\kk  \rightarrow \hochhom_\bullet\left(\homm_{\hcat\basecat}(0,0)\right)$
which sends $1$ to the $[\hunit]$. 

It remains to show that in $\hochhom_\bullet(\hperf \kk) = \kk$ we have
$$ \bigl(\hat{n}\bigr)\bigl(K(\psi_n(\alpha) \otimes
\psi_n(\beta)) \bigr) = n \left< \alpha, \beta \right>. $$
By Definition \ref{defn-annihilation-functor}, the map induced
by $(\hat{n})$ on the Hochschild homology is 
\begin{align*}
\hochhom_\bullet \left(\sym^{n}\basecat^{\opp} \otimes
\sym^{n}\basecat\right)
\xrightarrow{ \Res_{S_{n}}^{S_n \times S_n}}
& \hochhom_\bullet \left( \sym^n(\basecat^{\opp} \otimes \basecat) \right)
\xrightarrow{\homm_\basecat(-,-)}
\\
\xrightarrow{\quad \quad}\;
& 
\hochhom \left(
\sym^n(\hperf \kk)
\right)
\xrightarrow{(-) \otimes \dots \otimes (-)}
\hochhom_\bullet \left(\hperf \kk\right). 
\end{align*}
Consider the diagram
\begin{equation}
\label{eqn-key-commutative-diagram-for-heisen-rel-n-equals-m}
\begin{tikzcd}[row sep = 1.25cm, column sep = 3cm]
\hochhom_\bullet\left(\basecat^{\opp}\right)
\otimes 
\hochhom_\bullet\left(\basecat \right)
\ar{r}{\psi_n \otimes \psi_n}
\ar{dd}{K}
&
\hochhom_\bullet\left(\sym^n \basecat^{\opp}\right) \otimes 
\hochhom_\bullet\left(\sym^n \basecat \right)
\ar{d}{K}
\\
&
\hochhom_\bullet\left(\sym^n\basecat^{\opp}
\otimes \sym^n \basecat \right)
\ar{d}{\Res_{S_{n}}^{S_n \times S_n}}
\\
\hochhom_\bullet\left(\basecat^{\opp} \otimes \basecat \right)
\ar{r}{\psi_n}
\ar{d}{\homm_\basecat(-,-)}
&
\hochhom_\bullet\left(\sym^n\left(\basecat^{\opp} \otimes \basecat\right)\right)
\ar{d}{\sym^n \homm_\basecat(-,-)}
\\
\hochhom_\bullet\left(\hperf \kk \right)
\ar{r}{\psi_n}
\ar[equals]{rd}
&
\hochhom_\bullet\left(\sym^n\left( \hperf \kk \right)\right)
\ar{d}{(-) \otimes \dots \otimes (-)}
\\
&
\hochhom_\bullet\left(\hperf \kk \right). 
\end{tikzcd}
\end{equation}
If we start with $\alpha \otimes \beta \in
\hochhom_\bullet\left(\basecat^{\opp}\right) \otimes
\hochhom_\bullet\left(\basecat \right)$, then going around the
upper right perimeter of the diagram produces 
$\bigl(\hat{n}\bigr)\bigl(K(\psi_n(\alpha) \otimes \psi_n(\beta))
\bigr)$, 
while going around the lower left perimeter produces the Euler pairing
$\left< \alpha, \beta \right>$, as per its definition in 
\S\ref{section-euler-pairing-on-hochschild-homology}. 

The triangle at the bottom of 
\eqref{eqn-key-commutative-diagram-for-heisen-rel-n-equals-m}
commutes because $\hochhom_\bullet(\hperf) \simeq \kk$, so it suffices
to check that it commutes on $1 \in \kk$, i.e. on the class
$[\id_\kk]$, which is trivial. The middle square in 
\eqref{eqn-key-commutative-diagram-for-heisen-rel-n-equals-m}
commutes because our construction of the map $\psi_n\colon
\hochcx_\bullet(\A) \rightarrow \hochcx_\bullet(\sym^n \A)$
is functorial in $\A$. 

It remains to show that the top square in
\eqref{eqn-key-commutative-diagram-for-heisen-rel-n-equals-m} commutes
up to the factor of $n$. For this, consider its refinement to the diagram
 \begin{small}
\begin{equation*}
\begin{tikzcd}[row sep = 1.25cm, column sep = 1.15cm]
\hochhom_\bullet\left(\basecat^{\opp}\right)
\otimes 
\hochhom_\bullet\left(\basecat \right)
\ar{r}{g_n \otimes g_n}
\ar{dd}{K}
&
\hochhom_\bullet\left(\left(\basecat^{\opp}\right)^{\otimes n}; t \right)
\otimes 
\hochhom_\bullet\left(\basecat^{\otimes n}; t \right)
\ar{r}{\xi_{t} \otimes \xi_{t}}
\ar{d}{K}
&
\hochhom_\bullet\left(\sym^n \basecat^{\opp}\right) \otimes 
\hochhom_\bullet\left(\sym^n \basecat \right)
\ar{d}{K}
\\
&
\hochhom_\bullet\left(\left(\basecat^{\opp}\right)^{\otimes n} \otimes 
\basecat^{\otimes n}; t \times t \right)
\ar{r}{\xi_{t \times t}}
\ar[<->]{d}{\sim}
&
\hochhom_\bullet\left(\sym^n\basecat^{\opp}
\otimes \sym^n \basecat \right)
\ar{d}{\Res_{S_{n}}^{S_n \times S_n}}
\\
\hochhom_\bullet\left(\basecat^{\opp} \otimes \basecat \right)
\ar{r}{g_n}
&
\hochhom_\bullet\left(\left(\basecat^{\opp} \otimes
\basecat\right)^{\otimes n}; t \right)
\ar{r}{\xi_{t}}
&
\hochhom_\bullet\left(\sym^n\left(\basecat^{\opp} \otimes
\basecat\right)\right), 
\end{tikzcd}
\end{equation*}
\end{small}
where $t$ is the long cycle $(1\dots{n}) \in S_n$. 

The left square and the top right square commute by Lemmas
\ref{lemma-compatibility-of-the-map-g-with-shuffle-product} 
and \ref{lemma-compatibility-of-the-map-xi-with-shuffle-product},
respectively. It remains to show that the bottom right square
commutes up to the factor of $n$. For this, we apply Lemma 
\ref{lemma-restriction-functor-composed-with-the-twisted-part} to 
the composition $\Res_{S_{n}}^{S_n \times S_n} \circ \xi_{t \times
t}$.

The set $Q$ of the left cosets of $S_n$ in $S_n \times S_n$ can be
identified with $S_n$. The coset corresponding to each 
$\sigma \in Q$ has 
a representative $r_\sigma = \id \times \sigma \in S_n \times S_n$. 
Moreover, for every $\sigma \in Q$ we have 
$$ (t \times t)(\id \times \sigma) = 
(\id \times t \sigma t^{-1}) (t \times t). $$
Thus the action of $t \times t$ on $Q$ is given by 
$(t \times t).\sigma = t \sigma t^{-1}$ and 
we can set $h_{t \times t, \sigma} = t \times t$ in 
the terminology of Definition
\ref{defn-the-set-of-left-cosets-of-a-subgroup}. 
The fixed set $\fix(t \times t)$ of this action is the 
subset $\left\{\id, t, t^2, \dots, t^{n-1}\right\} \subset Q$. 
The representatives of these fixed cosets are the elements 
\[\id, \id \times t, \id \times t^2, \dots, \id \times t^{n-1} \in S_n
\times S_n.\] These give isomorphisms 
\eqref{eqn-group-element-induced-map-on-group-element-twisted-hh}
on twisted Hochschild complexes
\begin{equation*}
(\id \times t^i)\colon
\hochcx_\bullet
\left(\left(\basecat^{\opp}\right)^{\otimes n} \otimes 
\basecat^{\otimes n}; t \times t \right)
\longrightarrow
\hochcx_\bullet
\left(\left(\basecat^{\opp}\right)^{\otimes n} \otimes 
\basecat^{\otimes n}; t \times t \right)
\end{equation*}
whose induced maps on the Hochschild homology are the identity maps. 
 
Therefore, applying
Lemma \ref{lemma-restriction-functor-composed-with-the-twisted-part} we
get that on the Hochschild homology 
$$ 
\Res_{S_{n}}^{S_n \times S_n} \circ \xi_{t \times t}
= \sum_{i = 0}^{n-1} \xi_t \circ (1 \times t^{i})^{-1} = 
\sum_{i = 0}^{n-1} \xi_t = n \xi_t, $$ as desired. 
\end{proof}

\subsection{Injectivity of $\pi$}
\label{section-injectivity-of-the-decategorification-map-pi}

Having proved for $\hAA_\bullet(n) \in \chalg{\hochhom_\bullet(\basecat)}$ 
the commutation and the Heisenberg relations, we 
obtain by Theorem \ref{theorem-construction-of-the-decategorification-map}
the desired Hochschild homology categorification map. 
It is injective for the same tautological reason as its 
numerical Grothendieck group counterpat in \cite{gyenge2021heisenberg}:
\begin{Proposition}
\label{prps-decategorification-map-pi-is-injective}
The decategorification map 
$$
\pi \colon \chalg{\hochhom_\bullet(\basecat)} 
\rightarrow
\alghh_{\hcat\basecat}
$$
constructed in
\S\ref{section-overview-of-decategorification-map-pi}-\ref{section-commutation-and-heisenberg-relations} is injective. 	
 \end{Proposition}
 \begin{proof}
The $2$-functor $\Phi_{\basecat}$ of 
\cite[Theorem 7.30]{gyenge2021heisenberg} sends each 
object $n \in \mathbb{Z}$ of $\hcat\basecat$ to the object 
$\symbc{n}$ of $\fcat\basecat$. Since for $n < 0$ these are zero 
categories, $\Phi_{\basecat}$ sends all
$\homm_{\hcat\basecat}(n,n-m)$ to zero for $m > n$. 
It particular, it kills $\homm_{\hcat\basecat}(0,-n)$ for $n > 0$. 
Let $I_{-}$ be the left ideal in $\chalg{\hochhom_\bullet(\basecat)}$
generated by $a_{\alpha}(-n)$
for $n > 0$ and any $\alpha \in \hochhom_\bullet(\basecat)$.
By construction, $\pi$ maps each $a_{\alpha}(-n)$ to
$\hAA_{\alpha}(-n) \in
\hochhom_\bullet\left(\homm_{\hcat\basecat}(0,-n)\right)$. 
It follows that the image of $I_{-}$ under the $\kk$-algebra homomorphism  
\begin{equation}
\label{eqn-rep-of-H-on-fock-space-hochschild}
\chalg{\hochhom_\bullet(\basecat)}
\xrightarrow{
\eqref{eqn-injective-decategorification-map-into-the-flattening}
}
\hochhom_{alg}(\hcat\basecat) 
\xrightarrow{\hochhom_{alg}\left(\Phi_{\basecat}\right)}
\End \left( \bigoplus_{n \geq 0} \hochhom_{\bullet}(\symbc{n}) \right)
\end{equation}
kills $1 \in \kk \simeq \hochhom_\bullet(\symbc{0})$. Since the quotient  
$\chalg{\hochhom_\bullet(\basecat)}/I_{-}$ is the Fock space representation $\falg{\hochhom_\bullet(\basecat)}$, 
the action of $\chalg{\hochhom_\bullet(\basecat)}$ on $1 \in \kk \simeq \hochhom_\bullet(\symbc{0})$
induces a map of $\chalg{\hochhom_\bullet(\basecat)}$ representations
\begin{equation}
\phi\colon \falg{\hochhom_\bullet(\basecat)} \rightarrow
\bigoplus_{n \geq 0}
\hochhom_{\bullet}(\symbc{n}).
\end{equation}
The map $\phi$ is non-zero since
\eqref{eqn-rep-of-H-on-fock-space-hochschild} is a unital algebra
homomorphism, so $1 \in \chalg{\hochhom_\bullet(\basecat)}$ acts as
the identity map. By irreducibility of the Fock space representation,
$\phi$ is injective. By faithfulness of the Fock space space
representation, 
the morphism \eqref{eqn-rep-of-H-on-fock-space-hochschild} is injective, 
and hence so is the morphism 
\eqref{eqn-injective-decategorification-map-into-the-flattening}.
Since
\eqref{eqn-injective-decategorification-map-into-the-flattening}
is the composition
\eqref{eqn-injective-decategorification-map-into-the-flattening-via-pi}, 
it follows that $\pi$ is also injective. 
\end{proof}

\section{Heisenberg algebra and the Hochschild homology of 
symmetric quotient stacks}
\label{section-action-on-hochschild-homology-of-symmetric-quotient-stacks}

Let $\basecat$ be a smooth and proper DG category over 
$\kk$. Let $\chi$ be the Euler pairing on the Hochschild homology
$\hochhom_\bullet(\basecat)$. 
We view such $\basecat$ as a noncommutative 
smooth and proper variety. We view its symmetric powers
$\symbc{n}$ (see Definition~\ref{defn-symmetric-power-of-a-dg-category})
as noncommutative symmetric quotient stacks of $\basecat$. 
By Lemma \ref{lemma-A-rtimes-G-modules-are-G-equiv-A-modules}, 
when $\basecat$ is a DG enhancement of a smooth and proper variety $X$,
$\symbc{n}$ is a DG enhancement of the symmetric quotient stack $[X^n/S_n]$. 

In this section, we give three different constructions of an action 
of the Heisenberg algebra $\chalg{\hochhom_\bullet(\basecat), \chi}$
on the total Hochschild homology 
$\bigoplus_{n=0}^{\infty} \hochhom_\bullet(\symbc{n})$ of the
symmetric powers of $\basecat$. We then prove that the three resulting
actions coincide and can thus be viewed as different descriptions 
of the same action. 

The first action is obtained using the noncommutative orbifold
decomposition for Hochschild homology established recently in 
\cite{annobaranovskylogvinenko2023orbifold} and summarised in 
\S\ref{section-noncommutative-orbifld-decomposition-for-hh}. In 
\cite[Theorem 4.2]{annobaranovskylogvinenko2023orbifold}
it was shown that the individual decompositions of each
$\hochhom_\bullet(\symbc{n})$ combine into the decomposition 
of $\bigoplus_{n=0}^{\infty} \hochhom_\bullet(\symbc{n})$ into the
symmetric algebra $S^*(\hochhom_\bullet(\A) \otimes t \kk[t])$. 
We use this to translate the canonical action of 
$\chalg{\hochhom_\bullet(\basecat), \chi}$ 
on $S^*(\hochhom_\bullet(\A) \otimes t \kk[t])$ described in 
\S\ref{section-fock-space-and-its-symmetric-algebra-structure}
into an action of $\chalg{\hochhom_\bullet(\basecat), \chi}$ on 
$\bigoplus_{n=0}^{\infty} \hochhom_\bullet(\symbc{n})$:

\begin{Definition}[Orbifold decomposition action]
\label{defn-orbifold-decomposition-action}
Let $\basecat$ be a smooth and proper DG category over $\kk$
and let $\chi$ be the Euler pairing on $\hochhom_\bullet(\basecat)$. 

For each $\alpha \in \hochhom_\bullet(\basecat)$ and $i > 0$,
define operators $A_\alpha(-i)$ and $A_\alpha(i)$ on the 
total Hochschild homology 
$\bigoplus_{n=0}^{\infty} \hochhom_\bullet(\symbc{n})$
to be the image of the operators of the Lemma 
\ref{lemma-canonical-action-of-chalg-v-chi-on-symmetric-algebra}
under mutually inverse isomorphisms of 
\cite[Theorem 4.2]{annobaranovskylogvinenko2023orbifold}:
\begin{equation*}
\begin{tikzcd}
\bigoplus_{n \geq 0} \hochhom_\bullet(\sym^n A)
\ar[shift left = 1.24ex]{r}{\zeta}[']{\simeq}
&
S^*(\hochhom_\bullet(\A) \otimes t \kk[t]). 
\ar[shift left = 1.25ex]{l}{\eta}
\end{tikzcd}
\end{equation*}
\end{Definition}

Since the operators of Lemma 
\ref{lemma-canonical-action-of-chalg-v-chi-on-symmetric-algebra}
do define the action of $\chalg{\hochhom_\bullet(\basecat), \chi}$ 
on $S^*(\hochhom_\bullet(\A) \otimes t \kk[t])$, we know apriori 
that the operators of Definition~\ref{defn-orbifold-decomposition-action}
satisfy the relations 
\eqref{eq:vectheisrel1-a-gen-graded}-\eqref{eq:vectheisrel3-a-gen-graded}
of Definition \ref{defn-the-heisenberg-algebra-of-a-graded-vector-space-a-gen}
and define the action of $\chalg{\hochhom_\bullet(\basecat), \chi}$ on 
$\bigoplus_{n=0}^{\infty} \hochhom_\bullet(\symbc{n})$. 

The second action is a decategorification of the categorical action of 
\cite{gyenge2021heisenberg}. There, we constructed for any smooth and 
proper $\basecat$ its Heisenberg $2$-category $\hcat\basecat$, 
its $2$-categorical Fock space $\fcat\basecat$, and the $2$-categorical action 
$\Phi_\basecat$ of $\hcat\basecat$ on $\fcat\basecat$. 
In \S\ref{section-decategorification-map}, we defined the functor 
$\hochhom_{alg}$ of taking the Hochschild homology of a $2$-category 
and flattening it into an algebra. We also demonstrated that applying it 
to $\Phi_\basecat$ yields an algebra homomorphism 
\begin{equation*}
\hochhom_{alg}(\hcat\basecat) 
\xrightarrow{\eqref{eqn-decategorigication-of-phi-basecat}}
\End \left( \bigoplus_{n \geq 0} \hochhom_{\bullet}(\symbc{n}) \right). 
\end{equation*}
In \S\ref{section-reduction-to-alghh} we constructed from the
decategorification map $\pi$ (see
\S\ref{section-overview-of-decategorification-map-pi}-\ref{section-commutation-and-heisenberg-relations})
an algebra homomorphism 
\begin{equation}
\chalg{\hochhom_\bullet(\basecat), \chi}
\xrightarrow{\eqref{eqn-injective-decategorification-map-into-the-flattening}}
\hochhom_{alg}(\hcat\basecat). 
\end{equation}
The composition of 
\eqref{eqn-injective-decategorification-map-into-the-flattening}
and \eqref{eqn-decategorigication-of-phi-basecat} gives the action 
of $\chalg{\hochhom_\bullet(\basecat), \chi}$  on 
$\bigoplus_{n=0}^{\infty} \hochhom_\bullet(\symbc{n})$. 

\begin{Definition}[Decategorified action]
\label{defn-decategorified-action}
Let $\basecat$ be a smooth and proper DG category over $\kk$
and let $\chi$ be the Euler pairing on $\hochhom_\bullet(\basecat)$. 

For each $\alpha \in \hochhom_\bullet(\basecat)$ and $i > 0$,
define operators $A_\alpha(-i)$ and $A_\alpha(i)$ on the 
total Hochschild homology 
$\bigoplus_{n=0}^{\infty} \hochhom_\bullet(\symbc{n})$
to be the images of $a_\alpha(-i), a_\alpha(i)
\in \chalg{\hochhom_\bullet(\basecat), \chi}$ 
under the composition 
\begin{equation}
\label{eqn-decategorified-action-on-the-total-hochschild-homology}
\chalg{\hochhom_\bullet(\basecat), \chi}
\xrightarrow{\eqref{eqn-injective-decategorification-map-into-the-flattening}}
\hochhom_{alg}(\hcat\basecat)
\xrightarrow{\eqref{eqn-decategorigication-of-phi-basecat}}
\End \left( \bigoplus_{n \geq 0} \hochhom_{\bullet}(\symbc{n}) \right). 
\end{equation}
of the algebra homomorphism 
\eqref{eqn-injective-decategorification-map-into-the-flattening}
constructed in \S\ref{section-hochschild-homology-and-heisenberg-2-category}
and the decategorification
$\eqref{eqn-decategorigication-of-phi-basecat}$
of the categorical action of \cite{gyenge2021heisenberg}.  
\end{Definition}
By construction, since 
$\eqref{eqn-decategorified-action-on-the-total-hochschild-homology}$ 
is an algebra homomorphism, 
the operators of Definition~\ref{defn-decategorified-action}
also satisfy the relations 
\eqref{eq:vectheisrel1-a-gen-graded}-\eqref{eq:vectheisrel3-a-gen-graded}
of Definition \ref{defn-the-heisenberg-algebra-of-a-graded-vector-space-a-gen}
and define the action of $\chalg{\hochhom_\bullet(\basecat), \chi}$ on 
$\bigoplus_{n=0}^{\infty} \hochhom_\bullet(\symbc{n})$.

The third action is by explicitly defined operators which resemble  
the topological $K$-theoretic operators of  
Segal \cite{Segal-EquivariantKTheoryAndSymmetricProducts} and 
Wang \cite{Wang-EquivariantKTheoryWreathProductsAndHeisenbergAlgebra}. 
Unlike the previous two actions, defined via a chain of
identifications, here the operators are defined
intrinsically on $\bigoplus_{n=0}^{\infty} \hochhom_\bullet(\symbc{n})$
using the representation theoretic restriction and induction functors:
\begin{Definition}[Induction-restriction action]
\label{defn-induction-restriction-action}
Let $\basecat$ be a smooth and proper DG category over $\kk$
and let $\chi$ be the Euler pairing on $\hochhom_\bullet(\basecat)$. 

For each $\alpha \in \hochhom_\bullet(\basecat)$ and $i > 0$,
define operators $A_\alpha(-i)$ and $A_\alpha(i)$ on the 
total Hochschild homology 
$\bigoplus_{n=0}^{\infty} \hochhom_\bullet(\symbc{n})$
as follows:
\begin{scriptsize}
\begin{equation}
\label{eqn-noncommutative-operator-A-alpha-(-i)}
A_\alpha(-i) \colon 
\hochhom_\bullet \left(\sym^{n+i} \basecat\right) 
\xrightarrow{\Res^{S_{n+i}}_{S_{i} \times S_n}}
\hochhom_\bullet \left(\sym^{i} \basecat \otimes \sym^{n} \basecat\right) 
\simeq 
\hochhom_\bullet \left(\sym^{i} \basecat \right) \otimes 
\hochhom_\bullet\left( \sym^{n} \basecat\right) 
\xrightarrow{\left<\alpha, \kappa_i(-)\right>}
\hochhom_\bullet \left(\sym^{n} \basecat \right),
\end{equation}
\begin{equation}
\label{eqn-noncommutative-operator-A-alpha-(i)}
A_\alpha(i) \colon 
\hochhom_\bullet \left(\sym^{n} \basecat \right)
\xrightarrow{(-) \otimes \psi_i(\alpha)}
\hochhom_\bullet \left(\sym^{i} \basecat \right)
\otimes 
\hochhom_\bullet \left(\sym^{n} \basecat \right)
\simeq 
\hochhom_\bullet \left(\sym^{i} \basecat 
\otimes \sym^{m} \basecat \right)
\xrightarrow{\Ind^{S_{n+i}}_{S_{i} \times S_n}}
\hochhom_\bullet \left(\sym^{n+i} \basecat \right). 
\end{equation}
\end{scriptsize}
\end{Definition}
Unlike in Definitions \ref{defn-orbifold-decomposition-action} and 
\ref{defn-decategorified-action}, here we do not apriori know that the 
operators $A_\alpha(-i)$ and $A_\alpha(i)$ satisfy the relations 
\eqref{eq:vectheisrel1-a-gen-graded}-\eqref{eq:vectheisrel3-a-gen-graded}
of Definition \ref{defn-the-heisenberg-algebra-of-a-graded-vector-space-a-gen}
and define the action of $\chalg{\hochhom_\bullet(\basecat), \chi}$ on 
$\bigoplus_{n=0}^{\infty} \hochhom_\bullet(\symbc{n})$. The authors
would rather not try and verify the relations 
\eqref{eq:vectheisrel1-a-gen-graded}-\eqref{eq:vectheisrel3-a-gen-graded} 
by means of a direct computation in the Hochschild homology. 

Instead, we show that Definitions \ref{defn-orbifold-decomposition-action}, 
\ref{defn-decategorified-action}, and \ref{defn-induction-restriction-action}
all define the same set of the operators:
\begin{Theorem}
\label{theorem-three-heisenberg-actions-on-hh-of-noncomm-sym-quot-stacks}
Let $\basecat$ be a smooth and proper DG category over an
algebraically closed field $\kk$ of characteristic $0$.  Let $\chi$ be
the Euler pairing on the Hochschild homology $\hochhom_\bullet(\basecat)$.

For all $\alpha \in \hochhom_\bullet(\basecat)$ and $i > 0$
the following three definitions of the operators $A_\alpha(-i)$ and 
$A_\alpha(i)$ on $\bigoplus_{n=0}^{\infty} HH_\bullet(\symbc{n})$
are equivalent:
\begin{enumerate}
\item Orbifold decomposition operators of Definition \ref{defn-orbifold-decomposition-action}, 
\item Decategorified operators of Definition \ref{defn-decategorified-action}, 
\item Induction-restriction operators of Definition
\ref{defn-induction-restriction-action}.  
\end{enumerate}

The operators $A_\alpha(-i)$ and $A_\alpha(i)$ defined by these
equivalent definitions satisfy the relations 
\eqref{eq:vectheisrel1-a-gen-graded}-\eqref{eq:vectheisrel3-a-gen-graded}
and thus define an action of the Heisenberg algebra 
$H_{HH_\bullet(\basecat), \chi}$ on 
$\bigoplus_{n=0}^{\infty} HH_\bullet(\symbc{n})$ which identifies
$\bigoplus_{n=0}^{\infty} HH_\bullet(\symbc{n})$ with the
Fock space of $H_{HH_\bullet(\basecat), \chi}$. 
\end{Theorem}
\begin{proof}
The orbifold decomposition operators satisfy the relations 
\eqref{eq:vectheisrel1-a-gen-graded}-\eqref{eq:vectheisrel3-a-gen-graded}
and the resulting action of $H_{HH_\bullet(\basecat), \chi}$
identifies $\bigoplus_{n=0}^{\infty} HH_\bullet(\symbc{n})$ 
with the Fock space. This is because the analogous claim holds for
the operators of Lemma
\ref{lemma-canonical-action-of-chalg-v-chi-on-symmetric-algebra}
and the orbifold decomposition operators are the images 
of these under the isomorphism 
$\bigoplus_{n \geq 0} \hochhom_\bullet(\symbc{n}) 
\simeq S^*(\hochhom_\bullet(\basecat) \otimes t \kk[t])$ of 
\cite[Theorem 4.2]{annobaranovskylogvinenko2023orbifold}.

It remains to show that the orbifold decomposition operators of
Definition~\ref{defn-orbifold-decomposition-action}, 
the decategorified operators of Definition~\ref{defn-decategorified-action}, 
and the induction-restriction operators of 
Definition~\ref{defn-induction-restriction-action} all coincide. 

\underline{Orbifold decomposition operators $=$ Induction-restriction operators}:

Let $E$ be any graded vector space. The symmetric algebra $S^*(E)$ 
has natural structure of a graded Hopf algebra. cf. 
\cite[\S4.0]{Sweedler-HopfAlgebras}. Its multiplication
$\mu\colon S^*(E) \otimes S^*(E) \rightarrow S^*(E)$ sends 
$$ (s_1 \vee \dots \vee s_n) \otimes (t_1 \vee \dots \vee t_m) 
\rightarrow s_1 \vee \dots \vee s_n \vee t_1 \vee \dots \vee t_m, $$
and its
comultiplication  $\Delta\colon S^*(E) \rightarrow S^*(E) \otimes S^*(E) $
sends 
$$ (s_1 \vee \dots \vee s_n) 
\rightarrow \sum_{
\left\{ i_1, \dots, i_m \right\} 
\subseteq \left\{ 1, \dots, n \right\} }
(-1)^{p}
(s_{i_1} \vee \dots \vee s_{i_m}) \otimes (s_{j_1} \vee \dots \vee
s_{j_{n-m}}), $$
where $\left\{ j_1, \dots, j_{n-m} \right\} = \left\{1,\dots, n\right\}
\setminus \left\{ i_1, \dots, i_m \right\}$ and the sign twist $(-1)^p$ 
is the product of $(-1)^{\deg(s_{i_k})\deg(s_{j_{l}})}$ for 
each $i_k > j_l$. 

In case of the symmetric algebra 
$S^*(\hochhom_\bullet(\basecat) \otimes t \kk[t])$,
for any $\alpha \in \hochhom_\bullet(\basecat)$ and $i > 0$
we can rewrite the operators 
\eqref{eqn-action-of-a_w(n)-on-symmetric-algebra} 
and \eqref{eqn-action-of-a_w(-n)-on-symmetric-algebra}
of Lemma \ref{lemma-canonical-action-of-chalg-v-chi-on-symmetric-algebra}
in terms $\mu$ and $\delta$ as:
\begin{small}
\begin{align*}
a_\alpha(-i)\colon  \Sym^{\underline{r}(\underline{n} + (i))} 
\hochhom_\bullet(\basecat)\; t^{\underline{n} + (i)}
\xrightarrow{\Delta}
\hochhom_\bullet(\basecat)\; t^{(i)}
\otimes 
\Sym^{\underline{r}(\underline{n})} 
\hochhom_\bullet(\basecat)\; t^{\underline{n}}
\xrightarrow{\left<\alpha,-\right> \otimes \id}
\Sym^{\underline{r}(\underline{n})} 
\hochhom_\bullet(\basecat)\; t^{\underline{n}}, 
\end{align*}
\begin{align*}
a_\alpha(i)\colon  
\Sym^{\underline{r}(\underline{n})} 
\hochhom_\bullet(\basecat) \; t^{\underline{n}}
\xrightarrow{\alpha\; t^{(i)} \otimes \id}
\hochhom_\bullet(\basecat) \; t^{(i)}
\otimes 
\Sym^{\underline{r}(\underline{n}} 
\hochhom_\bullet(\basecat) \; t^{\underline{n}}
\xrightarrow{\mu}
\Sym^{\underline{r}(\underline{n} + (i))} 
\hochhom_\bullet(\basecat) \; t^{\underline{n} + (i)}.
\end{align*}
\end{small}

Recall that \eqref{eqn-action-of-a_w(n)-on-symmetric-algebra} 
and \eqref{eqn-action-of-a_w(-n)-on-symmetric-algebra}
are the actions of $a_\alpha(-i)$ and $a_\alpha(i)$ on
the individual summands in the direct sum decomposition
\eqref{eqn-direct-sum-description-of-tkt-sym-algebra-of-v}
of $S^*(\hochhom_\bullet(\basecat) \otimes t \kk[t])$. 
The total operators by which $a_\alpha(-i)$ and $a_\alpha(i)$ act 
are the sums of \eqref{eqn-action-of-a_w(n)-on-symmetric-algebra} 
and \eqref{eqn-action-of-a_w(-n)-on-symmetric-algebra}
over all $n \geq 0$ and all partitions $\underline{n} \vdash n$. 

Fix $n \geq 0$. By above we can rewrite the actions 
\begin{equation}
\label{eqn-action-of-canonical-a_w(-i)-operator-on-weight-n+i-component}
\underset{{\underline{n+i} \vdash n+i}}{\bigoplus}
\Sym^{\underline{r}(\underline{n+i})}
\hochhom_\bullet(\basecat)\; t^{\underline{n + i }}
\xrightarrow{\quad a_\alpha(-i) \quad}
\underset{{\underline{n} \vdash n}}{\bigoplus}
\Sym^{\underline{r}(\underline{n})}
\hochhom_\bullet(\basecat)\; t^{\underline{n}},
\end{equation}
\begin{equation}
\label{eqn-action-of-canonical-a_w(i)-operator-on-weight-n-component}
\underset{{\underline{n} \vdash n}}{\bigoplus}
\Sym^{\underline{r}(\underline{n})}
\hochhom_\bullet(\basecat)\; t^{\underline{n}}
\xrightarrow{\quad a_\alpha(i) \quad}
\underset{{\underline{n+i} \vdash n + i}}{\bigoplus}
\Sym^{\underline{r}(\underline{n+i})}
\hochhom_\bullet(\basecat)\; t^{\underline{n+i}}, 
\end{equation}
as the compositions
\begin{small}
\begin{align*}
& \underset{{\underline{n+i} \vdash n+i}}{\bigoplus}
\Sym^{\underline{r}(\underline{n+i})}
\hochhom_\bullet(\basecat)\; t^{\underline{n + i }}
\; \xrightarrow{\quad \Delta \quad}
\left(\underset{\underline{i} \vdash i}{\bigoplus}
\Sym^{\underline{r}(\underline{i})} 
\hochhom_\bullet(\basecat)\;
t^{\underline{i}}\right)
\otimes 
\left(\underset{\underline{n} \vdash n}{\bigoplus}
\Sym^{\underline{r}(\underline{n})} 
\hochhom_\bullet(\basecat)\;
t^{\underline{n}}\right)
\twoheadrightarrow 
\\
\twoheadrightarrow \;
& 
\hochhom_\bullet(\basecat)\; t^{(i)}
\otimes 
\left(\underset{\underline{n} \vdash n}{\bigoplus}
\Sym^{\underline{r}(\underline{n})} 
\hochhom_\bullet(\basecat)\;
t^{\underline{n}}\right)
\xrightarrow{\left<\alpha,-\right>}
\rightarrow \;
\underset{\underline{n} \vdash n}{\bigoplus}
\Sym^{\underline{r}(\underline{n})} 
\hochhom_\bullet(\basecat)\;
t^{\underline{n}},
\end{align*}

\begin{align*}
& \underset{\underline{n} \vdash n}{\bigoplus}
\Sym^{\underline{r}(\underline{n})} 
\hochhom_\bullet(\basecat)\;
t^{\underline{n}}
\; \xrightarrow{\quad \alpha\; t^{(i)} \otimes (-) \quad}
\hochhom_\bullet(\basecat)\; t^{(i)}
\otimes 
\left(\underset{\underline{n} \vdash n}{\bigoplus}
\Sym^{\underline{r}(\underline{n})} 
\hochhom_\bullet(\basecat)\;
t^{\underline{n}}\right)
\hookrightarrow 
\\
\hookrightarrow \;
& 
\left(\underset{\underline{i} \vdash i}{\bigoplus}
\Sym^{\underline{r}(\underline{i})} 
\hochhom_\bullet(\basecat)\;
t^{\underline{i}}\right)
\otimes 
\left(\underset{\underline{n} \vdash n}{\bigoplus}
\Sym^{\underline{r}(\underline{n})} 
\hochhom_\bullet(\basecat)\;
t^{\underline{n}}\right)
\xrightarrow{\quad \mu \quad } \;
\underset{{\underline{n+i} \vdash n+i}}{\bigoplus}
\Sym^{\underline{r}(\underline{n+i})}
\hochhom_\bullet(\basecat)\; t^{\underline{n + i }},
\end{align*}
\end{small}whose midlle maps are the projection to and the inclusion of
the direct summand corresponding to the long cycle partition $(i)$ of $i$. 

The orbifold decomposition operators $A_\alpha(-i)$ and $A_\alpha(i)$
are the images of the operators $a_\alpha(-i)$ and $a_\alpha(i)$ 
of Lemma \ref{lemma-canonical-action-of-chalg-v-chi-on-symmetric-algebra}
under the mutually inverse isomorphisms $\zeta$ and $\eta$ 
of \eqref{eqn-symmetric-algebra-homotopy-equivalence-hh}
whose individual components are 
mutally inverse isomorphisms $\zeta_n$ and $\eta_n$ of
\eqref{eqn-noncommutative-orbifold-decomposition-for-sn-v}. 

In \cite[Theorem 5.1]{annobaranovskylogvinenko2023orbifold},  
one computed the Hopf algebra structure induced from 
$S^*(\hochhom_\bullet(\basecat) \otimes t \kk[t])$ onto 
$\bigoplus_{n \geq 0} \hochhom_\bullet(\symbc{n})$
via $\zeta$ and $\eta$. By 
\cite[Theorem 5.1]{annobaranovskylogvinenko2023orbifold}
the following diagram commutes 
\begin{small}
\begin{equation*}
\begin{tikzcd}[column sep = 0.25cm]
\underset{{\underline{n+i} \vdash n+i}}{\bigoplus}
\Sym^{\underline{r}(\underline{n+i})}
\hochhom_\bullet(\basecat)\; t^{\underline{n + i }}
\ar{rr}{\Delta}
& &
\left(\underset{\underline{i} \vdash i}{\bigoplus}
\Sym^{\underline{r}(\underline{i})} 
\hochhom_\bullet(\basecat)\;
t^{\underline{i}}\right)
\otimes 
\left(\underset{\underline{n} \vdash n}{\bigoplus}
\Sym^{\underline{r}(\underline{n})} 
\hochhom_\bullet(\basecat)\;
t^{\underline{n}}\right)
\ar{d}{\eta_i \otimes \eta_n}.
\\
\hochhom_\bullet (\symbc{n+i})
\ar{r}{\Res^{S_{n+i}}_{S_{i} \times S_n}}
\ar{u}{\zeta_{n+i}}
&
\hochhom_\bullet (\symbc{i} \otimes \symbc{n})
\ar{r}{\sim}
&
\hochhom_\bullet (\symbc{i})
\otimes 
\hochhom_\bullet (\symbc{n})
\end{tikzcd}
\end{equation*}
\end{small}
It follows that isomorphisms $\eta$ and $\zeta$ identify the action 
\eqref{eqn-action-of-canonical-a_w(-i)-operator-on-weight-n+i-component}
of $a_\alpha(-i)$
with the composition 
\begin{align*}
&\hochhom_\bullet (\symbc{n+i})
\xrightarrow{\Res^{S_{n+i}}_{S_{i} \times S_n}}
\hochhom_\bullet (\symbc{i} \otimes \symbc{n})
\xrightarrow{\sim}
\hochhom_\bullet (\symbc{i})
\otimes 
\hochhom_\bullet (\symbc{n})
\xrightarrow{\eta_i}
\\
\rightarrow 
&
\left(\underset{\underline{i} \vdash i}{\bigoplus}
\Sym^{\underline{r}(\underline{i})} 
\hochhom_\bullet(\basecat)\;
t^{\underline{i}}\right)
\otimes 
\hochhom_\bullet (\symbc{n})
\twoheadrightarrow
\hochhom_\bullet(\basecat)\;
t^{(i)}
\otimes 
\hochhom_\bullet (\symbc{n})
\xrightarrow{\left<\alpha,-\right>}
\hochhom_\bullet (\symbc{n}). 
\end{align*}
In Definition~\ref{defn-linear-maps-psi-n-and-kappa-n} we defined 
map $\kappa_i$ to be the composition 
$$ \hochhom_\bullet (\symbc{i}) 
\xrightarrow{\eta_i} 
\underset{\underline{i} \vdash i}{\bigoplus}
\Sym^{\underline{r}(\underline{i})} 
\hochhom_\bullet(\basecat)\;
t^{\underline{i}}
\twoheadrightarrow 
\hochhom_\bullet(\basecat)\; t^{(i)}.
$$ 
We conclude that isomorphisms $\eta$ and $\zeta$ identify the action 
\eqref{eqn-action-of-canonical-a_w(-i)-operator-on-weight-n+i-component}
of $a_w(-i)$ in 
Lemma \ref{lemma-canonical-action-of-chalg-v-chi-on-symmetric-algebra}
with 
\begin{small}
\begin{align*}
\hochhom_\bullet (\symbc{n+i})
\xrightarrow{\Res^{S_{n+i}}_{S_{i} \times S_n}}
\hochhom_\bullet (\symbc{i} \otimes \symbc{n})
\xrightarrow{\sim}
\hochhom_\bullet (\symbc{i})
\otimes 
\hochhom_\bullet (\symbc{n})
\xrightarrow{\left<\alpha,\kappa_i(-)\right>}
\hochhom_\bullet (\symbc{n}), 
\end{align*}
\end{small}
which is the composition \eqref{eqn-noncommutative-operator-A-alpha-(-i)}
defining the induction-restriction operator $A_{\alpha}(-i)$. 

Similarly, by 
\cite[Theorem 5.1]{annobaranovskylogvinenko2023orbifold}
the following diagram commutes 
\begin{small}
\begin{equation*}
\begin{tikzcd}[column sep = 0.25cm]
\left(\underset{\underline{i} \vdash i}{\bigoplus}
\Sym^{\underline{r}(\underline{i})} 
\hochhom_\bullet(\basecat)\;
t^{\underline{i}}\right)
\otimes 
\left(\underset{\underline{n} \vdash n}{\bigoplus}
\Sym^{\underline{r}(\underline{n})} 
\hochhom_\bullet(\basecat)\;
t^{\underline{n}}\right)
\ar{rr}{\mu}
& &
\underset{{\underline{n+i} \vdash n+i}}{\bigoplus}
\Sym^{\underline{r}(\underline{n+i})}
\hochhom_\bullet(\basecat)\; t^{\underline{n + i }}
\ar{d}{\eta_{n+i}}.
\\
\hochhom_\bullet (\symbc{i})
\otimes 
\hochhom_\bullet (\symbc{n})
\ar{r}{\sim}
\ar{u}{\zeta_{i} \otimes \zeta_{n}}
&
\hochhom_\bullet (\symbc{i} \otimes \symbc{n})
\ar{r}{\Ind^{S_{n+i}}_{S_{i} \times S_n}}
&
\hochhom_\bullet (\symbc{n+i})
\end{tikzcd}
\end{equation*}
\end{small}
and thus $\eta$ and $\zeta$ identify the action
\eqref{eqn-action-of-canonical-a_w(i)-operator-on-weight-n-component}
of $a_w(i)$ in 
Lemma \ref{lemma-canonical-action-of-chalg-v-chi-on-symmetric-algebra}
with the composition
\begin{align*}
&\hochhom_\bullet (\symbc{n}). 
\xrightarrow{\alpha\; t^{(i)} \otimes (-)}
\hochhom_\bullet(\basecat)\;
t^{(i)}
\otimes 
\hochhom_\bullet (\symbc{n})
\hookrightarrow
\left(\underset{\underline{i} \vdash i}{\bigoplus}
\Sym^{\underline{r}(\underline{i})} 
\hochhom_\bullet(\basecat)\;
t^{\underline{i}}\right)
\otimes 
\hochhom_\bullet (\symbc{n})
\xrightarrow{\zeta_i}
\\
\rightarrow \;
&
\hochhom_\bullet (\symbc{i})
\otimes 
\hochhom_\bullet (\symbc{n})
\xrightarrow{\sim}
\hochhom_\bullet (\symbc{i} \otimes \symbc{n})
\xrightarrow{\Ind^{S_{n+i}}_{S_{i} \times S_n}}
\hochhom_\bullet (\symbc{n+i})
\end{align*}
which by definition of $\psi_i$ matches  
the definition \eqref{eqn-noncommutative-operator-A-alpha-(i)}
of the induction-restriction operator $A_{\alpha}(i)$.

\underline{Decategorified operators = Induction-restriction operators}:

We need to show that for any  
$\alpha \in \hochhom_\bullet(\basecat)$ and $i > 0$
the images of the generators $a_\alpha(-i)$ and $a_{\alpha}(i)$
of $\chalg{\hochhom_\bullet(\basecat), \chi}$ under the composition 
of \eqref{eqn-decategorigication-of-phi-basecat} and
\eqref{eqn-injective-decategorification-map-into-the-flattening}
are the induction-restriction operators $A_\alpha(-i)$ and $A_\alpha(i)$ 
defined in \eqref{eqn-noncommutative-operator-A-alpha-(-i)}
and
\eqref{eqn-noncommutative-operator-A-alpha-(i)}. 
By construction, these images are
$\Phi_\basecat(\Xi_\QQ(\psi_i(\alpha))$ and 
$\Phi_\basecat(\Xi_\PP(\psi_i(\alpha))$. 
By \cite[Lemma 8.4]{gyenge2021heisenberg}, the composition
$\Phi_\basecat \circ \Xi_\PP$ is homotopy equivalent to
the functor 
\begin{small}
\begin{align}
\label{eqn-filtering-phi-through-yoneda-embedding}
\begin{tikzcd}[row sep=0.5cm]
\hperf(\symbc{i})
\ar{d}{\text{unit}}
\\
\DGFun\left(\hperf(\symbc{n}),\, \hperf(\symbc{i}) \otimes \hperf(\symbc{n})\right)
\ar{d}{\otimes_\kk \circ (-)}
\\
\DGFun\left(\hperf(\symbc{n}),\, \hperf(\symbc{i} \otimes \symbc{n})\right)
\ar{d}{\Ind_{S_i \times S_n}^{S_{n+i}} \circ (-)}
\\
\DGFun\left(\hperf(\symbc{n}),\, \hperf(\symbc{n+i})\right). 
\end{tikzcd}
\end{align}
\end{small}

By its definition in \S\ref{section-hochschild-homology-bicategories}, 
it follows that the corresponding map 
\begin{equation}
\hochhom_\bullet(\symbc{i})
\rightarrow 
\homm_\kk\left(\hochhom_\bullet(\symbc{n}),
\hochhom_\bullet(\symbc{n+i})\right)
\end{equation}
sends $\psi_i(\alpha)$ to the definition 
\eqref{eqn-noncommutative-operator-A-alpha-(i)}
of the induction-restriction operator $A_\alpha(i)$.
By construction
of $\Phi_\basecat$ in \cite[\S7]{gyenge2021heisenberg}, for 
any $E \in \hperf(\symbc{i})$ the composition 
$\Phi_\basecat \circ \Xi_\QQ(E)$ is the right adjoint of 
$\Phi_\basecat \circ \Xi_\PP(E)$. It follows that 
$\Phi_\basecat \circ \Xi_\QQ$ sends 
$\psi_i(\alpha)$ to the
definition \eqref{eqn-noncommutative-operator-A-alpha-(-i)}
of the induction-restriction operator $A_\alpha(-n)$. 
\end{proof}

\bibliographystyle{amsplain}
\bibliography{references}

\providecommand{\bysame}{\leavevmode\hbox to3em{\hrulefill}\thinspace}
\providecommand{\MR}{\relax\ifhmode\unskip\space\fi MR }
\providecommand{\MRhref}[2]{%
  \href{http://www.ams.org/mathscinet-getitem?mr=#1}{#2}
}
\providecommand{\href}[2]{#2}
\begin{thebibliography}{10}

\bibitem{AdemLeidaRuan-OrbifoldsAndStringyTopology}
A.~Adem, J.~Leida, and Y.~Ruan, \emph{Orbifolds and stringy topology},
  Cambridge Tracts in Mathematics, no. 171, Cambridge University Press, 2007.

\bibitem{annobaranovskylogvinenko2023orbifold}
R.~Anno, V.~Baranovsky, and T.~Logvinenko, \emph{{The Hochschild homology of a
  noncommutative symmetric quotient stack}}, in preparation.

\bibitem{AnnoLogvinenko-SphericalDGFunctors}
R.~Anno and T.~Logvinenko, \emph{{Spherical DG-functors}}, J. Eur. Math. Soc.
  \textbf{49} (2017), no.~9, 2577--2656.

\bibitem{AnnoLogvinenko-BarCategoryOfModulesAndHomotopyAdjunctionForTensorFunctors}
\bysame, \emph{Bar category of modules and homotopy adjunction for tensor
  functors}, Int. Math. Res. Not. \textbf{2021} (2021), no.~2, 1353--1462.

\bibitem{Baranovsky-OrbifoldCohomologyAsPeriodicCyclicHomology}
V.~Baranovsky, \emph{Orbifold cohomology as periodic cyclic homology}, Int. J.
  Math. \textbf{14} (2003), no.~08, 791--812.

\bibitem{BelmansFuKrug-HochschildCohomologyOfHilbertSchemesOfPointsOnSurfaces}
P.~Belmans, L.~Fu, and A.~Krug, \emph{Hochschild cohomology of {H}ilbert
  schemes of points on surfaces}, arXiv:2309.06244, 2023.

\bibitem{BondalKapranov-EnhancedTriangulatedCategories}
A.~Bondal and M.~Kapranov, \emph{Enhanced triangulated categories}, Mat. Sb.
  \textbf{181} (1990), no.~5, 669--683.

\bibitem{BondalVanDenBergh-GeneratorsAndRepresentabilityOfFunctorsInCommutativeAndNoncommutativeGeometry}
A.~Bondal and M.~van~den Bergh, \emph{Generators and representability of
  functors in commutative and noncommutative geometry}, Mosc. Math. J
  \textbf{3} (2003), no.~1, 1--36, arXiv:math/0204218.

\bibitem{Brylinski-CyclicHomologyAndEquivariantTheories}
J.~Brylinski, \emph{Cyclic homology and equivariant theories}, Ann. Inst.
  Fourier \textbf{37} (1987), no.~4, 15--28.

\bibitem{Caldararu-TheMukaiPairingIITheHochschildKostantRosenbergIsomorphism}
A.~Caldararu, \emph{The {M}ukai pairing {II}: the
  {H}ochschild--{K}ostant--{R}osenberg isomorphism}, Advances in Mathematics
  \textbf{194} (2005), no.~1, 34--66.

\bibitem{CartanEilenberg-HomologicalAlgebra}
H.~Cartan and S.~Eilenberg, \emph{Homological algebra}, Princeton Mathematical
  Series, no.~19, Princeton University Press, Princeton, N. J., 1956.

\bibitem{cautis2012heisenberg}
S.~Cautis and A.~Licata, \emph{{Heisenberg categorification and Hilbert
  schemes}}, Duke Mathematical Journal \textbf{161} (2012), no.~13, 2469--2547.

\bibitem{ChenRuan-ANewCohomologyTheoryofOrbifold}
W.~Chen and Y.~Ruan, \emph{A new cohomology theory of orbifold}, Communications
  in Mathematical Physics \textbf{248} (2004), no.~1, 1--31.

\bibitem{Efimov-HomotopyFinitenessOfSomeDGCategoriesFromAlgebraicGeometry}
A.~Efimov, \emph{Homotopy finiteness of some {DG} categories from algebraic
  geometry}, J. Eur. Math. Soc. \textbf{22} (2020), no.~9, 2897--2942,
  pre-print arXiv:1308.0135.

\bibitem{EllingsrudGottsche-HilbertSchemesOfPointsAndHeisenbergAlgebras}
G.~Ellingsrud and L.~G{\"o}ttsche, \emph{Hilbert schemes of points and
  {H}eisenberg algebras}, School on {A}lgebraic {G}eometry ({T}rieste, 1999),
  ICTP Lect. Notes, vol.~1, ICTP, Trieste, 2000, pp.~60--100.

\bibitem{FantechiGottsche-OrbifoldCohomologyForGlobalQuotients}
B.~Fantechi and L.~G{\"o}ttsche, \emph{Orbifold cohomology for global
  quotients}, Duke Math. J. \textbf{117} (2003), no.~2, 197--227.

\bibitem{FeiginTsygan-CyclicHomologyOfAlgebrasWithQuadraticRelationsUniversalEnvelopingAlgebrasAndGroupAlgebras}
B.~Feigin and B.~Tsygan, \emph{Cyclic homology of algebras with quadratic
  relations, universal enveloping algebras and group algebras}, K-Theory,
  Arithmetic and Geometry: Seminar, Moscow University, 1984--1986 (Yuri~I.
  Manin, ed.), Lecture Notes in Math., vol. 1289, Springer-Verlag, Berlin,
  Heidelberg, 1987, pp.~210--239.

\bibitem{SymCat}
N.~Ganter and M.~Kapranov, \emph{Symmetric and exterior powers of categories},
  Transformation Groups \textbf{19} (2014), no.~1, 57--103.

\bibitem{GetzlerJones-TheCyclicHomologyOfCrossedProductAlgebras}
E.~Getzler and J.~Jones, \emph{The cyclic homology of crossed product
  algebras}, J. Reine Angew. Math. \textbf{445} (1993), 161--174.

\bibitem{Gottsche-TheBettiNumbersOfTheHilbertSchemeOfPointsOnASmoothProjectiveSurface}
L.~G{\"o}ttsche, \emph{The {B}etti numbers of the {H}ilbert scheme of points on
  a smooth projective surface.}, Math. Ann. \textbf{286} (1990), no.~1-3,
  193--208.

\bibitem{grojnowski1995instantons}
I.~Grojnowski, \emph{Instantons and affine algebras. {I}. {T}he {H}ilbert
  scheme and vertex operators}, Math. Res. Lett. \textbf{3} (1996), no.~2,
  275--291.

\bibitem{gyenge2021heisenberg}
Á. Gyenge, C.~Koppensteiner, and T.~Logvinenko, \emph{The {H}eisenberg
  category of a category}, Mem. Am. Math. Soc. \textbf{320} (2026), no.~1633.

\bibitem{Hatcher-AlgebraicTopology}
A.~Hatcher, \emph{{Algebraic Topology}}, Cambridge University Press, 2002.

\bibitem{HochschildKostantRosenberg-DifferentialFormsOnRegularAffineAlgebras}
G.~Hochschild, B.~Kostant, and A.~Rosenberg, \emph{Differential forms on
  regular affine algebras}, Trans. Amer. Math. Soc. \textbf{102} (1962),
  383--408.

\bibitem{kac1990infinite}
V.~G. Kac, \emph{{Infinite-dimensional Lie algebras}}, Cambridge University
  Press, 1990.

\bibitem{Kaledin-HomologicalMethodsInNoncommutativeGeometry}
D.~Kaledin, \emph{Homological methods in non-commutative geometry},
  \url{http://imperium.lenin.ru/~kaledin/tokyo/}, 2007.

\bibitem{Kaledin-NonCommutativeGeometryFromTheHomologicalPointOfView}
\bysame, \emph{Non-commutative geometry from the homological point of view},
  \url{http://imperium.lenin.ru/~kaledin/seoul/}, 2009.

\bibitem{KatzarkovKontsevichPantev-HodgeTheoreticAspectsOfMirrorSymmetry}
L.~Katzarkov, M.~Kontsevich, and T.~Pantev, \emph{Hodge theoretic aspects of
  mirror symmetry}, From Hodge Theory to Integrability and TQFT: tt*-geometry
  (Ron~Y. Donagi and Katrin Wendland, eds.), Proceedngs of Symposia in Pure
  Mathematics, vol.~78, Amer. Math. Soc., Providence, RI, 2007, pp.~87--174.

\bibitem{Keller-OnTheCyclicHomologyOfExactCategories}
B.~Keller, \emph{On the cyclic homology of exact categories}, Journal of Pure
  and Applied Algebra \textbf{136} (1999), no.~1, 1--56.

\bibitem{Kontsevich-DeformationQuantizationOfPoissonManifolds}
M.~Kontsevich, \emph{Deformation quantization of {P}oisson manifolds}, Lett.
  Math. Phys. \textbf{66} (2003), no.~3, 157--216.

\bibitem{KontsevichSoibelman-NotesOnAInftyAlgebrasAInftyCategoriesAndNoncommutativeGeometry}
M.~Kontsevich and Y.~Soibelman, \emph{Notes on ${A}_\infty$-algebras,
  ${A}_\infty$-categories and {N}on-{C}ommutative {G}eometry}, Homological
  Mirror Symmetry, Lecture Notes in Physics, vol. 757, Springer Berlin
  Heidelberg, 2009, pp.~1--67.

\bibitem{krug2018symmetric}
A.~Krug, \emph{{Symmetric quotient stacks and Heisenberg actions}},
  Mathematische Zeitschrift \textbf{288} (2018), no.~1-2, 11--22.

\bibitem{Lang-Algebra}
S.~Lang, \emph{Algebra}, {R}ev. 3rd ed., {G}raduate {T}exts in {M}athematics,
  vol. 211, Springer, 2002.

\bibitem{Lehn-LecturesOnHilbertSchemes}
M.~Lehn, \emph{Lectures on {H}ilbert schemes}, Algebraic Structures and Moduli
  Spaces, CRM Proceedings \& Lecture Notes, vol.~38, American Mathematical
  Society, 2004, pp.~1--30.

\bibitem{Loday-CyclicHomology}
J.~Loday, \emph{Cyclic homology}, Grundlehren der mathematischen
  {W}issenschaften, no. 301, Springer, 1997.

\bibitem{LuntsOrlov-UniquenessOfEnhancementForTriangulatedCategories}
V.~Lunts and D.~Orlov, \emph{Uniqueness of enhancement for triangulated
  categories}, J. Amer. Math. Soc. \textbf{23} (2010), 853--908,
  arXiv:0908.4187.

\bibitem{nakajima1997heisenberg}
H.~Nakajima, \emph{{Heisenberg algebra and Hilbert schemes of points on
  projective surfaces}}, Annals of Mathematics \textbf{145} (1997), no.~2,
  379--388.

\bibitem{nakajima1999lectures}
\bysame, \emph{{Lectures on Hilbert schemes of points on surfaces}}, no.~18,
  American Mathematical Soc., 1999.

\bibitem{Nordstrom-HochschildHomologyOfSymmetricPowersOfDGcategories}
V.~Nordstr{\"o}m, \emph{{A decomposition theorem for the Hochschild homology of
  symmetric powers of a dg category}}, arXiv:2511.03269, 2025.

\bibitem{Nordstrom-FiniteGroupActionsOnDGCategoriesAndHochschildHomology}
\bysame, \emph{Finite group actions on {DG} categories and {H}ochschild
  homology}, Can. Math. Bull. (2025), 1--21.

\bibitem{Orlov-SmoothAndProperNoncommutativeSchemesAndGluingofDGcategories}
D.~Orlov, \emph{Smooth and proper noncommutative schemes and gluing of {DG}
  categories}, Adv. Math. \textbf{302} (2016), 59--105, arXiv:1402.7364.

\bibitem{Segal-EquivariantKTheoryAndSymmetricProducts}
G.~Segal, \emph{{Equivariant K-theory and symmetric products}}, pre-print,
  1996.

\bibitem{Shklyarov-HirzebruchRiemannRochTypeFormulaForDGAlgebras}
D.~Shklyarov, \emph{Hirzebruch--{R}iemann--{R}och-type formula for {DG}
  algebras}, Proceedings of the London Mathematical Society \textbf{106}
  (2013), no.~1, 1--32, arXiv:0710.1937.

\bibitem{Swan-HochschildCohomologyOfQuasiprojectiveSchemes}
R.~Swan, \emph{Hochschild cohomology of quasiprojective schemes}, J. Pure Appl.
  Algebra \textbf{110} (1996), 57--80.

\bibitem{Sweedler-HopfAlgebras}
M.~Sweedler, \emph{Hopf algebras}, W.A. Benjamin Inc., 1969.

\bibitem{Tabuada-InvariantsAdditifsDeDGCategories}
G.~Tabuada, \emph{{Invariants additifs de dg-cat{\'e}gories}}, International
  Mathematics Research Notices \textbf{2005} (2005), no.~53, 3309--3339.

\bibitem{Toen-TheHomotopyTheoryOfDGCategoriesAndDerivedMoritaTheory}
B.~To{\"e}n, \emph{The homotopy theory of {\em dg}-categories and derived
  {M}orita theory}, Invent. Math. \textbf{167} (2007), no.~3, 615--667.

\bibitem{Toen-LecturesOnDGCategories}
\bysame, \emph{Lectures on {DG}-categories}, Topics in Algebraic and
  Topological K-Theory, Lecture Notes in Mathematics, vol. 2008, Springer
  Berlin Heidelberg, 2011, pp.~243--302.

\bibitem{Wang-EquivariantKTheoryWreathProductsAndHeisenbergAlgebra}
W.~Wang, \emph{{Equivariant $K$-theory, wreath products, and Heisenberg
  algebra}}, Duke Mathematical Journal \textbf{103} (2000), no.~1, 1--23.

\end{thebibliography}

\end{document}